\documentclass[11pt]{article}

\usepackage{amsthm,amsmath,amssymb,bbm}
\usepackage{natbib}
\usepackage{multirow}
\usepackage[pdftex]{graphicx}
\usepackage{subfigure}
\usepackage{makecell}
\usepackage{booktabs}
\usepackage{array}
\usepackage{tabularx}
\usepackage{tabulary}
\usepackage{caption}
\usepackage{booktabs}
\usepackage{url}
\usepackage{algorithm}
\usepackage{algorithmic}
\usepackage{bm}
\usepackage{wrapfig}
\usepackage{lipsum}
\usepackage{mathrsfs}
\usepackage{dsfont}
\usepackage{titling}
\usepackage{txfonts}
\usepackage{setspace}

\usepackage{xr}
\usepackage[usenames,dvipsnames,svgnames,table]{xcolor}
\usepackage[colorlinks,
linkcolor=red,
anchorcolor=blue,
citecolor=blue
]{hyperref}

\usepackage{smile}

\newcommand*{\var}{\textnormal{var}}

\newcommand{\nn}{\nonumber}

\def\##1\#{\begin{align}#1\end{align}}
\def\$#1\${\begin{align*}#1\end{align*}}

\usepackage{relsize}
\newcommand{\T}{{\mathsmaller {\rm T}}}

\def\sn{\sum_{i=1}^n}
\def\Sb{\mathbf{S}}

\newcommand{\BB}{\mathbb{B}}

\newcommand{\wt}{\widetilde}

\newcommand{\bfsym}[1]{\ensuremath{\boldsymbol{#1}}}
\def \balpha   {\bfsym{\alpha}}       \def \bbeta    {\bfsym{\beta}}
       \def \bdelta   {\bfsym{\delta}}

\newcommand{\Rom}[1]{\text{\uppercase\expandafter{\romannumeral #1\relax}}}

\usepackage{geometry}
\usepackage{enumitem}

\newcommand{\dd}{{\rm d}}

\newcommand{\est}{{\rm est}}

\numberwithin{equation}{section}
\numberwithin{figure}{section}

\begin{document}

\title{Smoothed Quantile Regression with Large-Scale Inference}

\author{Xuming He,
%\thanks{Department of Statistics, University of Michigan, Ann Arbor, MI 48109, USA. E-mail:\href{xmhe@umich.edu}{\textsf{xmhe@umich.edu}}.}
~~Xiaoou Pan,
%\thanks{Department of Mathematics, University of California, San Diego, La Jolla, CA, 92093, USA. E-mail:\href{mailto:xip024@ucsd.edu}{\textsf{xip024@ucsd.edu}}.}
~~Kean Ming Tan,
%\thanks{Department of Statistics, University of Michigan, Ann Arbor, MI 48109, USA. E-mail:\href{mailto:keanming@umich.edu}{\textsf{keanming@umich.edu}}.}
~~and~~Wen-Xin Zhou
%\thanks{Department of Mathematics, University of California, San Diego, La Jolla, CA 92093, USA. E-mail:\href{mailto:wez243@ucsd.edu}{\textsf{wez243@ucsd.edu}}.} 
}

\date{}
\maketitle

\vspace{-0.5in}

\begin{abstract}
Quantile regression is a powerful tool for learning the relationship between a response variable and a multivariate predictor while exploring heterogeneous effects. In this paper, we consider statistical inference for quantile regression with large-scale data in the ``increasing dimension" regime. We provide a comprehensive and in-depth analysis of a convolution-type smoothing approach that achieves adequate approximation to computation and inference for quantile regression. This method, which we refer to as {\it{conquer}}, turns the non-differentiable quantile loss function into a twice-differentiable, convex and locally strongly convex surrogate, which admits a fast and scalable  Barzilai-Borwein gradient-based algorithm to perform optimization, and multiplier bootstrap for statistical inference. Theoretically, we establish explicit non-asymptotic bounds on both estimation and Bahadur-Kiefer linearization errors, from which we show that the asymptotic normality of the conquer estimator holds under a weaker requirement on the number of the regressors than needed for conventional quantile regression. Moreover, we prove the validity of multiplier bootstrap confidence constructions. Our numerical studies confirm the conquer estimator as a practical and reliable approach to large-scale inference for quantile regression. Software implementing the methodology is available in the \texttt{R} package \texttt{conquer}.
\end{abstract}

\noindent
{\bf Keywords}: Bahadur-Kiefer representation;  convolution; quantile regression; gradient descent; multiplier bootstrap; non-asymptotic statistics.

 \let\thefootnote\relax\footnotetext{Xuming He: Department of Statistics, University of Michigan, Ann Arbor, MI.  Email:\href{xmhe@umich.edu}{\textsf{xmhe@umich.edu}}.   Xiaoou Pan: Department of Mathematics, University of California, San Diego, La Jolla, CA.  Email:\href{mailto:xip024@ucsd.edu}{\textsf{xip024@ucsd.edu}}. 
 Kean Ming Tan: Department of Statistics, University of Michigan, Ann Arbor, MI.  Email:\href{mailto:keanming@umich.edu}{\textsf{keanming@umich.edu}}.  
 Wen-Xin Zhou: Department of Mathematics, University of California, San Diego, La Jolla, CA.  Email:\href{mailto:wez243@ucsd.edu}{\textsf{wez243@ucsd.edu}}.
We thank four anonymous referees,  the associate editor,  and editor (Xiaohong Chen) for many valuable  suggestions that help improve the overall quality of the paper. }

%%%%%%%%%%%%%%%%%%%%%%%%%%%%%%%%%%%%%%%%%%%
%%%%%%%%%%%%%%%%%%%%%%%%%%%%%%%%%%%%%%%%%%%
% Introduction
%%%%%%%%%%%%%%%%%%%%%%%%%%%%%%%%%%%%%%%%%%%
%%%%%%%%%%%%%%%%%%%%%%%%%%%%%%%%%%%%%%%%%%%
\section{Introduction}
\label{sec:1}
Quantile regression (QR) is a useful statistical tool for modeling and inferring the  relationship between a scalar response $y$ and a $p$-dimensional predictor $\bx$ \citep{KB1978}. 
Compared to the least squares regression that focuses on modeling the conditional mean of $y$ given $\bx$, QR allows modeling of the entire conditional distribution of $y$ given $\bx$, and thus provides valuable insights into heterogeneity in the relationship between $\bx$ and $y$.
Moreover, quantile regression is robust against outliers and can be performed for skewed or heavy-tailed response distributions without correct specification of the likelihood. 
These advantages make QR an appealing method to explore data features that are invisible to the least squares regression. We refer to \citet{K2005} and \cite{KCHP2017} for an extensive overview of QR  from methods, theory, computation, to various extensions under complex data structures.

Quantile regression involves a convex optimization problem with a piecewise linear loss function, also known as the \emph{check function}.  One can reformulate the QR problem as a linear program (LP), solvable by the Frisch-Newton algorithm with an average-case computational complexity that grows as a cubic function of $p$, i.e., $\cO_{\PP}(n^{1+\alpha} p^3 \log n)$ for some constant $\alpha \in (0, 1/2)$ \citep{PK1997}, where $n$ is the sample size and $p$ is the parametric dimension.  However, when applied to large-scale problems---both $n$ and $p$ are large,  QR computation via LP reformulation tends to be slow or too memory-intensive.
To better appreciate such a challenge,  we take the empirical study of U.S.  equities from \cite{GKX2020} as an example.  The dataset consists of monthly total individual equity returns,  which begins in March 1957 and ends in December 2016,  from CRSP for all firms listed in the NYSE,  AMEX,  and NASDAQ.  Within this span of 60 years,  the average number of stocks considered  is around 6200 per month.  After processing, the number of observations (over 60 years) in the entire panel exceeds 4 million, and the number of  stock-level covariates is 920.  Using half of the data (1957--1986) as the training sample,  the training sample size is still as large as 2 million.  Even with preprocessing, the interior point QR solver in \texttt{R} \citep{K2019} may either  run out of memory or take too much time on a personal computer.  This shortcoming arguably makes QR less attractive among various machine learning tools. We refer to Chapter~5 of \cite{KCHP2017} for an overview of the prevailing computational methods for quantile regression, such as simplex-based algorithms \citep{BR1974,KO1987}, interior point methods \citep{PK1997}, and alternating direction method of multipliers among other first-order proximal methods \citep{PB2014}.
 
We consider conducting large-scale inference for quantile regression under the ``increasing dimension" regime, namely,  the dimension $p=p_n$ is subject to the growth condition $p \asymp n^a$ for some $a\in (0,1)$. Two general principles have been widely used to suit this purpose.  The first uses a nonparametric estimate of the asymptotic variance \citep{GJ1992} that involves the conditional density of the response given the covariates, yet such an estimate can be fairly unstable. Even if the asymptotic variance is well estimated, its approximation accuracy to the finite-sample variance depends on the design matrix and the quantile level. Resampling methods, on the other hand, provide a more reliable approach to inference for QR under a wide variety of settings \citep{PWY1994,HH2002,KHH2005,FHH2011}.
Inevitably, the resampling approach requires repeatedly computing QR estimates up to thousands of times, and therefore is unduly expensive  for large-scale data.

 Theoretically, valid statistical inference is often justified by asymptotic normal approximations to QR estimators. The Bahadur-Kiefer representation of the nonlinear QR estimators are essential to this end, as shown in \cite{A1996} and \cite{HS1996}.  In large-$p$ (non)asymptotic settings in which the parametric dimension $p$ may tend to infinity with the sample size, we refer to \cite{W1989}, \cite{HS2000}, \cite{BCCF2019}, and \cite{PZ2019} for normal approximation results of the QR estimators under fixed and random designs. The question of how large $p$ can be relative to $n$ to ensure asymptotic normality has been addressed by those authors. It is now recognized that we may have to pay a price here as compared to $M$-estimators with smooth loss functions that are at least twice continuously differentiable.

To circumvent the non-differentiability of the QR loss function, \cite{H1998} proposed to smooth the indicator part of the check function via the survival function of a kernel.
 This smoothing method, which we refer to as \emph{Horowitz's smoothing} throughout, has been widely used for various QR-related problems \citep{W2006, KS2017,GK2016,WSZ2012, WMY2015,CGKL2019,CLZ2019}. 
However, Horowitz's smoothing gains smoothness at the cost of convexity, which inevitably raises optimization-related issues.
In general, computing a global minimum of a non-convex function is intractable: finding an $\epsilon$-suboptimal point for a $k$-times continuously differentiable function $f:\mathbb{R}^p \to \RR$ requires at least as many as $(1/\epsilon)^{p/k}$ evaluations of the function and its first $k$ derivatives \citep{NY1983}.
As we shall see from the numerical studies in Section~\ref{sec:numerical}, the convergence of gradient-based algorithms can be relatively slow  for high and low quantile levels. To address the aforementioned issue, \cite{FGH2019} proposed a convolution-type smoothing method that yields a convex and twice differentiable loss function, and studied the asymptotic properties of the smoothed estimator when $p$ is fixed.   To distinguish this approach from Horowitz's smoothing, we adopt the term \emph{conquer} for \underbar{con}volution-type smoothed \underbar{qu}antil\underbar{e} \underbar{r}egression.

In this paper, we first provide an in-depth statistical analysis of conquer under various nonstandard asymptotics settings in which $p$ increases with $n$. Our results reveal a key feature of the smoothing parameter, often referred to as the  bandwidth: the bandwidth adapts to both the sample size $n$ and dimensionality $p$, so as to achieve a tradeoff between statistical accuracy and computational stability. Since the convolution smoothed loss function is globally convex and locally strongly convex, we propose an efficient gradient descent algorithm with the Barzilai-Borwein stepsize and a Huber-type initialization.  The proposed algorithm is implemented via \texttt{RcppArmadillo} \citep{ES2014} in the \texttt{R} package \texttt{conquer}. We next focus on large-scale statistical inference (hypothesis testing and confidence estimation) with large $p$ and larger $n$.   We propose a bootstrapped conquer method that has reduced computational complexity when the conquer estimator is used as initialization. Under  appropriate restrictions on dimension, we establish the consistency (or concentration), Bahadur representation, asymptotic normality of the conquer estimator as well as the validity of the bootstrap approximation. In the following, we provide more details on the computational and statistical contributions of this paper.

Theoretically, by allowing $p$ to grow with $n$, the `complexity' of the function classes that we come across  in the analysis  also increases with $n$. 
Conventional asymptotic tools for proving the bootstrap validity are based on weak convergence arguments \citep{VW1996}, which are not directly applicable in the increasing dimension setting, especially with a non-differentiable loss.  In this paper we turn to a more refined and self-contained analysis, and prove a new local restricted strong convexity (RSC) property for the empirical smoothed quantile loss. This validates the key merit of convolution-type smoothing, i.e., local strong convexity.
The smoothing method involves a bandwidth, denoted by $h$.
Theoretically, we show that with sub-exponential random covariates (relaxing the bounded covariates assumption in \cite{FGH2019}), conquer exhibits an $\ell_2$-error in the order of $\sqrt{(p+t)/n}+h^2$ with probability at least $1-2e^{-t}$. 
When $h$ is of order $\{ (p+\log n)/n \}^{\gamma}$ for any $\gamma \in [1/4, 1/2]$, the conquer estimation is first-order equivalent to QR. Under slightly more stringent sub-Gaussian condition on the covariates, we show that the Bahadur-Kiefer linearization error  of conquer is of order $(p+t)/(nh^{1/2})+ h^{3/2} \sqrt{(p+t)/n}+h^4$ with probability at least $1-3e^{-t}$.
Based on such a representation, we establish a Berry-Esseen bound for linear functionals of conquer, which lays the theoretical foundation for testing general linear hypotheses, encompassing  covariate-effect analysis, analysis of variance, and model comparisons, to name a few. It is worth noting that with a properly chosen $h$, the linear functional of conquer is asymptotically normal as long as $p^{8/3}/n\rightarrow 0$, which improves the best known growth condition on $p$ for standard QR \citep{W1989, HS2000, PZ2019}.  We attribute this gain to the effect of smoothing.
Under similar conditions, we further establish   upper bounds on both estimation and Bahadur-Kiefer linearization errors for the bootstrapped conquer estimator.

\begin{figure}[!htp]
  \centering
  \subfigure[Estimation error]{\includegraphics[scale=0.45]{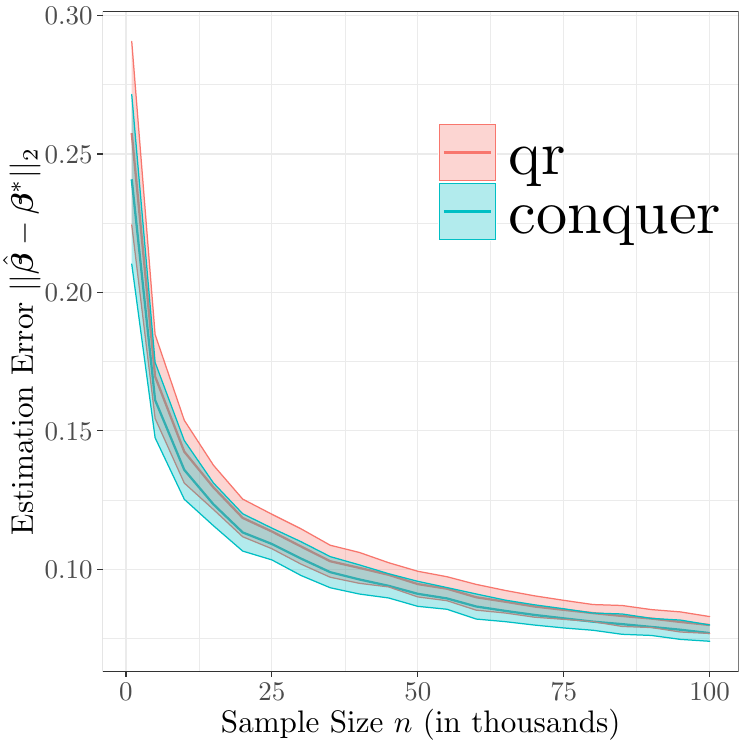}}  \qquad\quad 
  \subfigure[Runtime]{\includegraphics[scale=0.45]{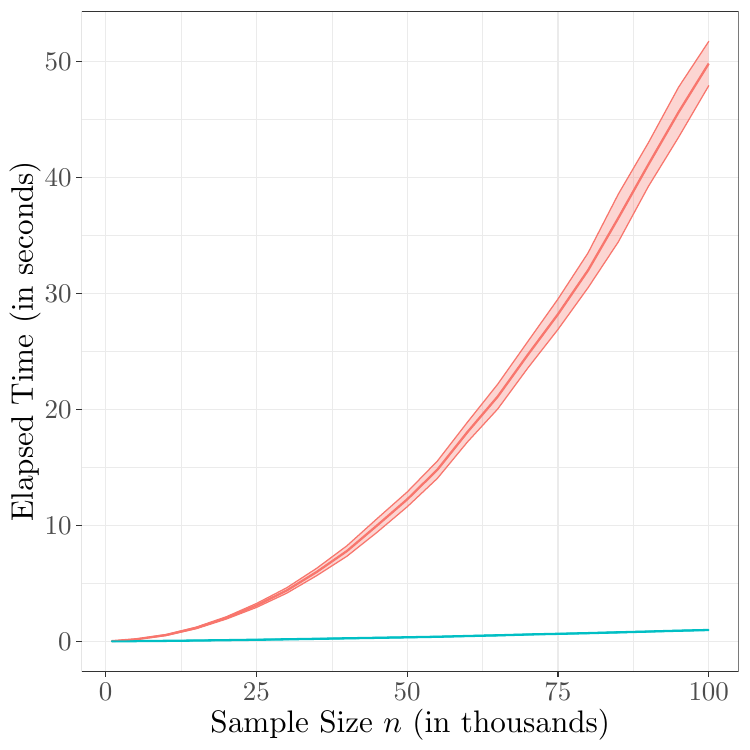}} 
\caption{A numerical comparison between  conquer and QR. The latter is implemented by the \texttt{R} package  \texttt{quantreg} using the ``pfn" method.  Panels (a) and (b) display, respectively, the ``estimation error and its standard deviation versus sample size'' and ``elapsed time and its standard deviation versus sample size" as the size of the problem increases.}
  \label{intro.fig}
\end{figure}

To better appreciate the computational feasibility of conquer for large-scale problems, we compare it with standard QR on large synthetic datasets, where the latter is implemented by the \texttt{R} package \texttt{quantreg} \citep{K2019}  using the Frisch-Newton approach after preprocessing ``pfn''.
We generate independent data vectors $\{y_i, \bx_i\}_{i = 1}^n$ from a linear model $y_i = \beta^*_0 + \langle   \bx_i  , \bbeta^* \rangle + \varepsilon_i$, where $(\beta^*_0, {\bbeta^*}^\T)^\T = (1, \ldots , 1)^\T \in \RR^{p+1}$, $\bx_i \sim \mathcal{N}_p(0, \Ib)$, and the independent errors $\varepsilon_i \sim t_2$, for $i = 1, 2, \dots, n$. 
We report the estimation error and elapsed time for increasing sample sizes $n\in \{1000,5000,10000,\ldots,100000\}$ and  dimension $p = \lfloor n^{1/2} \rfloor$, the largest integer that is less than or equal to $n^{1/2}$. 
Figure~\ref{intro.fig} displays the average estimation error, average elapsed time and their standard deviations based on $100$ Monte Carlo samples.
This experiment shows promise of conquer as a practically useful tool for large-scale quantile regression analysis. More empirical evidence will be given in the latter section.

The rest of the paper is organized as follows.
We start with a brief review of linear quantile regression and  the convolution-type smoothing method in Section~\ref{sec2}.  
Explicit forms of the smoothed check functions are provided for several representative kernel functions in nonparametric statistics.
We introduce the multiplier bootstrap for statistical inference in Section~\ref{subsec:mbi}. 
A Barzilai-Borwein gradient-based algorithm with a Huber-type warm start is detailed in Section~\ref{sec:algorithm}.
In Section~\ref{sec:theory}, we provide a comprehensive theoretical study of conquer from a non-asymptotic viewpoint, which directly leads to array asymptotic results.  Specifically, the bias incurred by smoothing the quantile loss is characterized in Section~\ref{sec:theory:bias}.  
In Section~\ref{sec:estimation error}, we establish the rate of convergence, Bahadur-Kiefer representation, and Berry-Esseen bound for conquer in a large-$p$ and larger-$n$ regime.  Results for its bootstrap counterpart are provided in Section~\ref{sec:boot:theory}.
We conclude the paper with an extensive numerical study in Section~\ref{sec:numerical} to illustrate the finite-sample performance of conquer in large-scale quantile regression analysis. We defer the proofs of all theoretical results as well as the full details of the one-step conquer to online supplementary materials.

\medskip
\noindent
{\sc Notation}: 
For every integer $k\geq 1$, we use $\RR^k$ to denote the the $k$-dimensional Euclidean space. The inner product of any two vectors $\bu=(u_1, \ldots, u_k)^\T, \bv=(v_1, \ldots ,v_k)^\T \in \RR^k$ is defined by $\bu^\T \bv = \langle \bu, \bv \rangle= \sum_{i=1}^k u_i v_i$.
We use $\| \cdot \|_p$ $(1\leq p \leq \infty)$ to denote the $\ell_p$-norm in $\RR^k$: $\| \bu \|_p = ( \sum_{i=1}^k | u_i |^p )^{1/p}$ and $\| \bu \|_\infty = \max_{1\leq i\leq k} |u_i|$. 
Throughout this paper, we use bold capital letters to represent matrices. For $k\geq 2$, $\Ib_k$ represents the identity matrix of size $k$. For any $k\times k$ symmetric matrix $\Ab \in \RR^{k\times k}$, $\| \Ab \|_2$ denotes the operator norm of $\Ab$. If $\Ab$ is positive semidefinite,  we use $\| \cdot \|_{\Ab}$ to denote the vector norm linked to $\Ab$ given by $\| \bu \|_{\Ab} = \| \Ab^{1/2} \bu \|_2$, $\bu \in \RR^k$. 
For $r \geq 0$, define the Euclidean ball and unit sphere in $\RR^k$ as $\BB^k(r) = \{ \bu \in \RR^k : \| \bu \|_2 \leq r\}$ and $\mathbb{S}^{k-1} = \partial \mathbb B^k(1) = \{ \bu \in \RR^k: \| \bu \|_2 =1 \}$, respectively.
For two sequences of non-negative numbers $\{ a_n \}_{n\geq 1}$ and $\{ b_n \}_{n\geq 1}$, $a_n \lesssim b_n$ indicates that there exists a constant $C>0$ independent of $n$ such that $a_n \leq Cb_n$; $a_n \gtrsim b_n$ is equivalent to $b_n \lesssim a_n$; $a_n \asymp b_n$ is equivalent to $a_n \lesssim b_n$ and $b_n \lesssim a_n$.

%%%%%%%%%%%%%%%%%%%%%%%%%%%%%%%%%%%%
%%%%%%%%%%%%%%%%%%%%%%%%%%%%%%%%%%%%
% Smoothing Quantile
%%%%%%%%%%%%%%%%%%%%%%%%%%%%%%%%%%%%
%%%%%%%%%%%%%%%%%%%%%%%%%%%%%%%%%%%%
\section{Smoothed quantile regression}
\label{sec2}
\subsection{The linear quantile regression model}
\label{subsec:problem}
Given a univariate response variable $y\in\RR$ and a $p$-dimensional covariate vector $\bx = (x_1, \ldots, x_p)^\T \in \RR^p$ with $x_1 \equiv 1$, the primary goal here is to learn the effect of $\bx$  on the distribution of $y$.
Let $F_{y|\bx }(\cdot)$ be the conditional distribution function of $y$ given $\bx$. The dependence between $y$ and $\bx$ is then fully characterized by the conditional quantile functions of $y$ given $\bx$, denoted as $F^{-1}_{y|\bx}(\tau)$, for $0<\tau<1$. We consider a linear quantile regression model at a given $\tau \in (0,1)$, that is,
the $\tau$-th conditional quantile function is
\begin{equation}
\label{eq:lqr}
F^{-1}_{y|\bx} (\tau) = \langle \bx, \bbeta^*(\tau) \rangle ,
\end{equation}
where $\bbeta^*(\tau)  = (\beta^*_{1}(\tau)  , \ldots, \beta^*_{ p}(\tau)  )^\T \in \RR^p$ is the true quantile regression coefficient.
%\end{cond}

Let $\{ (y_i, \bx_i) \}_{i=1}^n$ be a random sample from $(y,\bx)$. The standard quantile regression estimator \citep{KB1978} is then given as  
\# \label{qr.loss}
	\hat \bbeta(\tau) \in  \min_{\bbeta \in \RR^p }  \hat Q (\bbeta ) =  \min_{\bbeta \in \RR^p }     \frac{1}{n} \sn \rho_\tau( y_i -  \langle   \bx_i  , \bbeta \rangle )  ,
\#	
where $\rho_\tau(u) = u \{ \tau - \mathbbm{1}(u<0)\}$ is the $\tau$-quantile loss function, also referred to as the check function. Statistical properties of  $\hat \bbeta(\tau)$ have been extensively studied. We refer the reader to  \cite{K2005} and \cite{KCHP2017} for more details.

%%%%%%%%%%%%%%%%%%%%%%%%%%%%%%%%%
%%%%%%%%%%%%%%%%%%%%%%%%%%%%%%%%%
% Smoothing and Approximation Error
%%%%%%%%%%%%%%%%%%%%%%%%%%%%%%%%%
%%%%%%%%%%%%%%%%%%%%%%%%%%%%%%%%%
\subsection{Smoothed estimation equation and convolution-type smoothing}
\label{sec:smoothing}

Let $Q(\bbeta) = \EE \{ \hat Q (\bbeta ) \}$ be the population quantile loss function. Under mild conditions, $Q(\cdot)$  is twice differentiable and strongly convex in a neighborhood of $\bbeta^*$ with Hessian matrix $\Jb := \nabla^2 Q(\bbeta^*) = \EE \{ f_{\varepsilon | \bx} (0) \bx \bx^\T \}$, where $\varepsilon$ is the random noise and $f_{\varepsilon | \bx}(\cdot)$ is the conditional density of $\varepsilon$ given $\bx$. In contrast, the empirical quantile loss $\hat Q(\cdot)$ is not differentiable at $\bbeta^*$, and its ``curvature energy'' is concentrated at a single point. 
This is substantially different from other widely used loss functions that are at least locally strongly convex, such as the squared or logistic loss. 
The non-smoothness property not only brings challenge to theoretical analysis, but more importantly, also prevents gradient-based optimization methods from being efficient.
In his seminal work, \cite{H1998} proposed to directly smooth the check function $\rho_\tau(\cdot)$ to obtain 
\# \label{qr.loss.Horo}
\ell^{\mathrm{Horo}}_{ h}(u) = u \bigl\{ \tau - \cG(-u / h) \bigr\},
\#
where $\cG(\cdot)$ is a smooth function that takes values between 0 and 1, and $h>0$ is a smoothing parameter/bandwidth.
However, Horowitz's smoothing gains smoothness at the cost of convexity, which inevitably raises optimization issues especially when $p$ is large.
On the other hand, by the first-order condition, the population parameter $\bbeta^*$ satisfies the moment condition
\#
	\nabla Q(\bbeta^*) = \EE  \left[  \bigl\{    \mathbbm{1} ( y < \bx^\T \bbeta) - \tau \bigr\} \bx  \right] \Big|_{\bbeta = \bbeta^*}  = \textbf{0}.  \nn
\#
This property motivates a smoothed estimating equation (SEE) estimator \citep{W2006, KS2017}, defined as the solution to the smoothed moment condition
\#
  \frac{1}{n}  \sn \bigl[   \cG\bigl\{  (  \langle \bx_i ,  \bbeta \rangle  - y_i )/h \bigr\}  -\tau \bigr] \bx_i = \textbf{0} . \label{see}
\#

\begin{figure}[H]
  \centering
  \subfigure[Guassian kernel under $\tau = 0.5$.]{\includegraphics[width=0.45\textwidth]{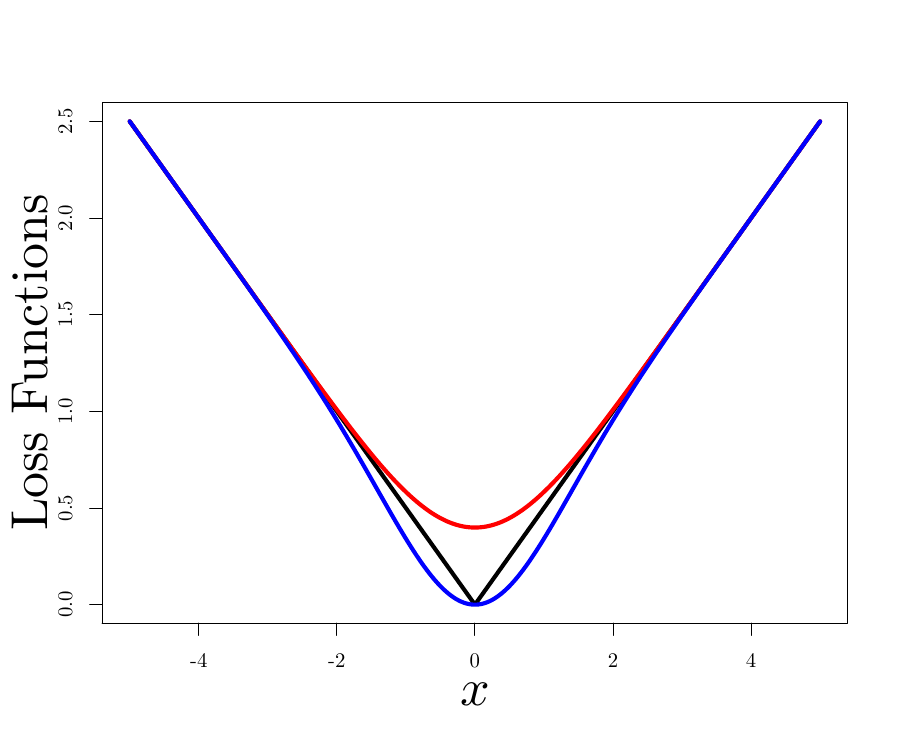}} 
  \subfigure[Uniform kernel under $\tau = 0.7$.]{\includegraphics[width=0.45\textwidth]{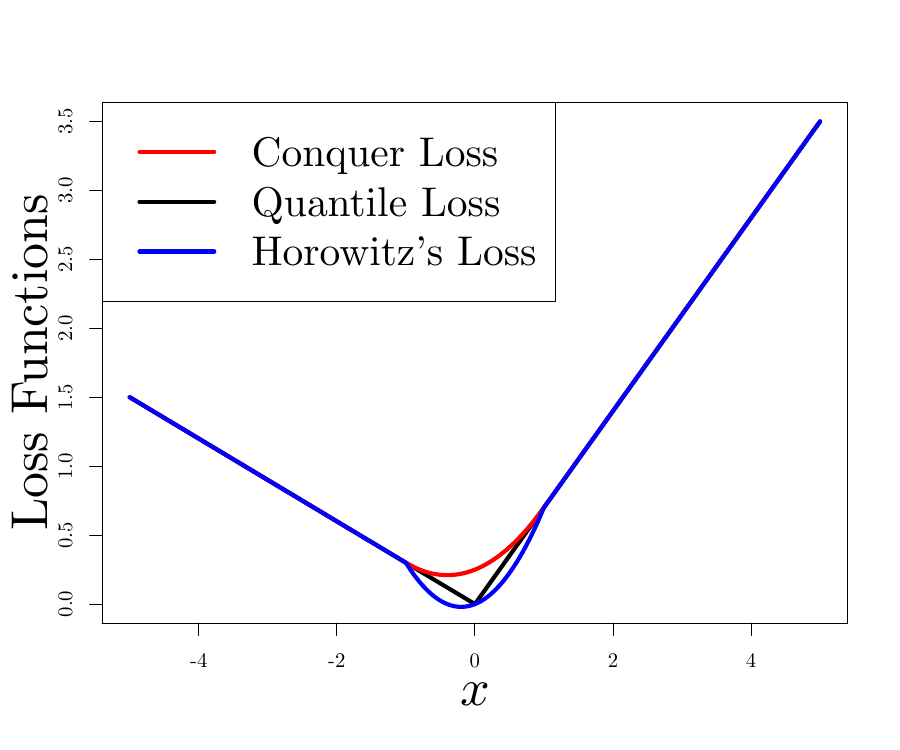}} 
\caption{Visualization of the quantile loss in \eqref{qr.loss}, conquer loss \eqref{convolution.loss}, and Horowitz's smoothed loss \eqref{qr.loss.Horo} with Gaussian and uniform kernels, respectively. 
}
  \label{loss.fig}
\end{figure}

Let $K(\cdot)$ be a kernel function that integrates to one,  and $h>0$ be a bandwidth.  Throughout the paper, we write 
\#
 K_h(u ) = h^{-1} K(u/h), \quad  \cK_h(u) = \cK(u/h) ~\mbox{ and }~ \cK(u) = \int_{-\infty}^u K(v) \, {\rm d} v , \ \ u\in \RR. \label{kernel.def}
\#
From an $M$-estimation viewpoint, the aforementioned SEE estimator can be equivalently defined as a minimizer of the empirical smoothed loss function
\#
\hat Q_{h}(\bbeta) = \frac{1}{n} \sn \ell_{h}(y_i -  \langle   \bx_i ,   \bbeta \rangle  ) ~~\mbox{ with }~~\ell_{h}(u) =  (\rho_\tau * K_h )(u) = \int_{-\infty}^{\infty} \rho_\tau(v) K_h(v- u ) \, {\rm d} v ,  \label{convolution.loss}
\#
where  $*$ denotes the convolution operator.   Therefore, as stated in the Introduction, we refer to the aforementioned smoothing method as \emph{conquer}.
The ensuing conquer estimator is given by
\# \label{sqr}
		\hat \bbeta_h = \hat \bbeta_h(\tau) \in \argmin_{\bbeta \in \RR^p} \hat Q_h(\bbeta) . 
\#
The key difference between the conquer loss \eqref{convolution.loss} and Horowitz's loss \eqref{qr.loss.Horo} is that the former is globally convex, while Horowitz's loss is not. This is illustrated in Figure~\ref{loss.fig}. 

As we shall see later, the ideal choice of bandwidth should adapt to the sample size $n$ and dimension $p$, since the quantile level $\tau$ is prespecified and fixed. Thus, the dependence of $\hat \bbeta_h$ and $\hat Q_h(\cdot)$ on $\tau$ will be assumed without display. Commonly used kernel functions include: (a) uniform kernel $K(u) = (1/2) \mathbbm{1} (|u|\leq 1)$, (b) Gaussian kernel  $K(u) = \phi(u):= (2\pi)^{-1/2} e^{-u^2/2}$, (c)  logistic kernel $K(u) = e^{-u}/(1+ e^{-u})^2 $, (d) Epanechnikov kernel $K(u) = (3 / 4) (1 - u^2) \mathbbm{1}(|u|\leq 1)$, and (e) triangular kernel $K(u) = (1 - |u|) \mathbbm{1}(|u|\leq 1)$. Explicit expressions of the corresponding smoothed loss function $\rho_\tau * K_h$ will be given in Section~\ref{sec:algorithm}.

The convolution-type kernel smoothing yields an objective function $\bbeta \mapsto \hat Q_{h}(\bbeta)$ that is twice continuously differentiable with gradient and hessian matrix
\# 
	\nabla  \hat Q_{h}(\bbeta)  =  \frac{1}{n} \sn  \bigl\{ \cK_h\bigl(   \langle \bx_i, \bbeta\rangle - y_i  \bigr) - \tau \bigr\} \bx_i  ~~\mbox{ and }~~
\nabla^2  \hat Q_{h}(\bbeta)  =  \frac{1}{n} \sn K_h( y_i - \langle \bx_i, \bbeta  \rangle ) \bx_i \bx_i^\T  ,  \label{grad.hess}
\#
respectively, where $\cK_h(\cdot) = \cK(\cdot/h)$ is defined in \eqref{kernel.def}.  Provided that $K$ is non-negative, $\hat Q_{h}(\cdot)$ is a convex function for any $h>0$, and $\hat \bbeta_h = \hat \bbeta_h(\tau) $ satisfies the first-order condition $\nabla  \hat Q_{ h}(\hat \bbeta_h )   = \textbf{0}$. This reveals the connection between the SEE and the conquer methods.
Together, the smoothness and convexity of $\hat Q_{h}(\cdot)$ warrant the superior computation efficiency of first-order gradient based algorithms for solving large-scale smoothed quantile regressions. The computational aspect of conquer will be discussed in Section~\ref{sec:algorithm}.

When the dimension $p$ is fixed, asymptotic properties of the SEE or conquer estimator have been studied by \cite{KS2017} and \cite{FGH2019}. The former used a higher-order kernel to deal with the instrumental variables QR problem (see Section~\ref{sec:related} for further discussions), and the latter showed that the conquer estimator has a lower asymptotic mean squared error than Horowitz's smoothed estimator, and also has a smaller Bahadur linearization error than the standard QR in the almost sure sense. The optimal order of the bandwidth based on the asymptotic mean squared error is unveiled as a function of $n$.
In Section~\ref{sec:theory}, we will establish exponential concentration inequalities and non-asymptotic Bahadur representation for the conquer estimator, while allowing the dimension $p$ to grow with the sample size $n$.  Our results reveal a key feature of the smoothing parameter: the bandwidth should adapt to both the sample size $n$  and dimensionality $p$, so as to achieve a tradeoff between statistical accuracy and computational stability.  
 
\begin{remark}
As discussed in \cite{FGH2019}, another advantage of convolution smoothing is that it facilitates conditional density estimation for quantile regression process. Assume $Q_y(\tau|\bx) = F^{-1}_{y|\bx}(\tau) = \langle \bx, \bbeta^*(\tau) \rangle$ for all $\tau\in [\tau_L, \tau_U] \subseteq (0,1)$. Under mild regularity conditions, $q_y(\tau|\bx):= \partial Q_y(\tau|\bx) /\partial \tau = 1/ f_{y |\bx}(\langle \bx, \bbeta^*(\tau) \rangle)$ exists. The inverse conditional density function plays an important role in, for example, the study of quantile treatment effects through modeling inverse propensity scores \citep{C2010}. By the linear conditional quantile model assumption, $ \partial Q_y(\tau|\bx) /\partial \tau = \langle \bx, \partial
 \bbeta^*(\tau)/\partial \tau \rangle$ for $\tau \in (\tau_L,\tau_U)$. Recall that the conquer estimator $\hat \bbeta_h=\hat \bbeta_h(\tau)$ satisfies the first-order condition $\nabla \hat Q_h(\hat \bbeta_h(\tau) ) =0$. Taking the partial derivative with respect to $\tau$ on both sides, it follows from \eqref{grad.hess} and the chain rule that
 \$
 \frac{\partial \hat \bbeta_h(\tau) }{\partial \tau} =    \big\{ \nabla^2  \hat Q_h\big(\hat \bbeta_h(\tau) \big)    \big\}^{-1}  \frac{1}{n}\sn \bx_i  = \Bigg\{  \frac{1}{n} \sn K_h\big( y_i - \bx_i^\T \hat \bbeta_h(\tau)  \big) \bx_i \bx_i^\T \Bigg\}^{-1}\frac{1}{n}\sn \bx_i  .
 \$
 Consequently, the inverse densities $1/ f_{y_i |\bx_i}(  \bx_i^\T \bbeta^*(\tau) )$ can be directly estimated by $\bx_i^\T  \frac{\partial \hat \bbeta_h(\tau) }{\partial \tau} $. This bypasses the use of any nonparametric method for density estimation with fitted residuals.
\end{remark}

\subsection{Multiplier bootstrap inference}
\label{subsec:mbi}
In this section, we consider a multiplier bootstrap procedure to construct confidence intervals  for conquer.
Independent of the observed sample $\mathcal X_n = \{ (y_i,\bx_i) \}_{i=1}^n$, let $\{w_i\}_{i = 1}^n$ be independent and identically distributed  random variables with $\EE (w_i)=1$ and $\var(w_i)=1$.  Recall that $\hat \bbeta_h = \hat \bbeta_h(\tau) = \min_{\bbeta \in \RR^p} \hat Q_{h}(\bbeta)$ is the conquer estimator. If the minimizer is not unique, we take any of the minima as  $\hat \bbeta_h = (\hat \beta_{h,1},\ldots, \hat \beta_{h,p})^\T$.

The proposed bootstrap method, which dates back to \cite{D1992} and \cite{BB1995},  is based on reweighting the summands of $\hat Q_{h}(\cdot )$ with random weights $w_i$. More specifically, define the weighted quantile loss $\hat Q^{\flat}_h :\RR^p \to \RR $ as  
\#
		\hat Q_{ h}^\flat (\bbeta ) =  \frac{1}{n} \sn  w_i  \ell_{h}(y_i -  \langle \bx_i, \bbeta \rangle ) ,    \label{wsq.loss}
\#
where $\ell_{ h}(u)  = (\rho_\tau * K_h)(u)$ is as in \eqref{convolution.loss}. The ensuing multiplier bootstrap statistic is then given by
\#
 \hat \bbeta_h^\flat =  \hat \bbeta_h^\flat(\tau)  \in  \argmin_{\bbeta \in \Theta  } \hat Q_{h}^{\flat}(\bbeta)  ,  \label{bsqr.1}
\#
where $ \Theta$ is a predetermined subset of $\RR^p$. If the random weights are allowed to take negative values, the weighted quantile loss may be non-convex. We therefore take $\Theta$ as a compact subset, such as $\Theta = \{ \bbeta \in \RR^p: \| \bbeta - \hat \bbeta_h \|_2 \leq R\}$ for some $R\geq 1$, in order to guarantee the existence of local/global optima.

Let $\mathbb E^*$ and $\PP^*$ be the conditional expectation and probability given the observed data $\cX_n$, respectively. 
Observe that $\mathbb E^*\{ \hat  Q_{h}^{\flat}(\bbeta ) \}=  \hat Q_{h}(\bbeta)$ for any $\bbeta \in \RR^{p}$. Consequently, we have
\$
 \argmin_{\bbeta \in \RR^p}  \EE^* \{ \hat  Q_{h}^{\flat}(\bbeta)  \} = \argmin_{\bbeta \in \RR^p}  \hat  Q_{h}(\bbeta)  = \hat \bbeta_h  .
\$
Intuitively, this means that $\hat \bbeta_h^\flat$ is an $M$-estimator of $\hat \bbeta_h$ in the bootstrap world. Since $\hat \bbeta_h$ is also an (approximate) $M$-estimator of $\bbeta^*$,  we expect that the distribution of $\hat \bbeta_h - \bbeta^*$ can be well  approximated with high probability by the conditional distribution of $\hat \bbeta_h^\flat - \hat \bbeta_h$. We will establish the validity of this approximation with explicit rates in Section~\ref{sec:boot:theory}. We refer to \cite{CB2005}  for a general asymptotic theory for weighted bootstrap for estimating equations, where a class of bootstrap weights is considered. Extensions to semiparametric  $M$-estimation can be found in \cite{MK2005} and \cite{CH2010}.

To retain convexity of the loss function, non-negative random weights are preferred, such as $w_i \sim {\rm Exp}(1)$, i.e., exponential distribution with rate 1, and $w_i =1 + e_i$, where $e_i$ are independent Rademacher random variables. {To compute the bootstrap estimate $\hat \bbeta_h^\flat$, we use the same gradient-based algorithm described in Section~\ref{sec:algorithm}, which will have a faster convergence rate thanks to the (provably) good initialization $\hat \bbeta_h$ and the strong convexity of the smoothed loss in a local basin.}

We can construct confidence intervals based on the bootstrap estimates using one of the three classical methods, the percentile method, the pivotal method, and the normal-based method.  
To be specific, for each $q\in (0,1)$  and $1\leq j\leq p$, define the (conditional) $q$-quantile of $\hat \beta^\flat_{h,j}$---the $j^{{\rm th}}$ coordinate of $\hat \bbeta^\flat_h \in \RR^p$---given the observed data as ${\rm c}^\flat_j(q) = \inf  \{ t \in \RR : \PP^* ( \hat \beta^\flat_{h,j}   \leq t ) \geq q   \}$. Then, for a prespecified nominal level $\alpha \in (0,1)$, the corresponding $1-\alpha$ bootstrap percentile and pivotal confidence intervals (CIs) for $\beta^*_j$ ($j=1,\ldots,p$) are, respectively,
\$
   \left[   {\rm c}^\flat_j(\alpha/2 ) , \, {\rm c}^\flat_j( 1 -\alpha/2) \right] ~\mbox{ and }~ \left[  2\hat \beta_{h,j} - {\rm c}^\flat_j( 1 -\alpha/2)  ,  \, 2 \hat \beta_{h,j} - {\rm c}^\flat_j(\alpha/2 ) \right]  .
\$
Numerically,  ${\rm c}^\flat_j(q)$ ($q\in \{ \alpha, 1-\alpha/2\}$) can be calculated with any specified precision by the simulation. In the \texttt{R} package \texttt{conquer}, the default number of bootstrap replications is set as $B=1000$.

In  the next section, we will present a finite-sample theoretical framework for convolution-type smoothed quantile regression, including the concentration inequality and non-asymptotic Bahadur representation for both the conquer estimator \eqref{sqr} and its bootstrap counterpart \eqref{bsqr.1} using Rademacher multipliers. As a by-product, a Berry-Esseen-type inequality  (see Theorem~\ref{thm:clt}) states that, under certain constraints on the (growing) dimensionality and bandwidth, the distribution of any linear projection of $\hat \bbeta_h$ converges to a normal distribution as the sample size increases to infinity. Informally, for any given deterministic vector $\ba \in \RR^p$, the scaled statistic $n^{1/2} \langle \ba, \hat \bbeta^*_h - \bbeta^* \rangle$ is asymptotically normally distributed with asymptotic variance $\sigma_0^2(\ba) := \tau(1-\tau) \,\ba^\T \Jb^{-1} \bSigma \Jb^{-1} \ba$, where $\bSigma$ is the population covariance matrix of the covariates $\bx$.  Another interesting implication from our theoretical analysis is that the unit variance requirement $\var(w_i)=1$ for the random weight is not necessary to ensure (asymptotically) valid bootstrap inference after a proper variance adjustment.  See Remark~\ref{rmk:weights} for more details.

To make inference based on such asymptotic results, we need to consistently estimate the asymptotic variance. \cite{FGH2019} suggested the following estimators 
\#
 \hat \Jb_h := \nabla^2 \hat \cQ_h(\hat \bbeta_h) = \frac{1}{n h}  \sn K(\hat\varepsilon_i /h) \cdot \bx\bx_i^\T   ~~\mbox{ and }~~ \hat \Vb_h := \frac{1}{n } \sn   \bigl\{  \cK_h(-\hat\varepsilon_i ) - \tau \bigr\}^2 \bx_i \bx_i^\T  \label{def:var-cov}
\#
of $\Jb$ and $\tau(1-\tau) \, \bSigma$, respectively, where $\hat \varepsilon_i = y_i - \langle\bx_i, \hat \bbeta_h\rangle$ are fitted residuals. The ensuing $1-\alpha$ normal-based CIs are given by $\hat \beta_{h,j} \pm   \Phi^{-1}(1-\alpha/2) \cdot n^{-1/2}  (\hat \Jb_h^{-1} \hat \Vb_h \hat \Jb_h^{-1} )_{jj}$, $j=1,\ldots,p$.
The normal approximations to the CI may suffer from the sensitivity to the smoothing needed to estimate the conditional densities, namely, the matrix $\Jb = \EE \{ f_{\varepsilon | \bx}(0) \bx \bx^\T\}$. When $p$ is large, inverting the estimated density matrix $\hat \Jb_h$ may be numerically unstable. This is typically true when $\tau$ is in the upper or lower tail. See Section~\ref{sec:normal-boot} for a numerical comparison between normal approximation and bootstrap calibration for confidence construction at high and low quantile levels. As we shall see, the normal-based CIs can be exceedingly wide and thus inaccurate under these situations.
 
\subsection{Connections to instrument variable and nonparametric quantile regression}
\label{sec:related}
 
\subsubsection{Instrument variable quantile regression}
This work focuses on large-scale estimation and inference for  linear quantile regression with many exogenous covariates. However, in many economic applications,  some regressors of interest (e.g., education, prices) are endogenous,  making conventional quantile regression inconsistent for estimating causal quantile effects. To address this problem,  \cite{CH2005} proposed an instrumental variable quantile regression (IVQR) model, which has become a popular tool for estimating quantile effects with endogenous covariates.   Due to the non-convex and non-smooth nature of the IVQR estimation problem,  there is a  burgeoning literature on the computational issues, dating back to \cite{CH2006}.  We refer to \cite{CHW2020}---Chapter 9 of \cite{KCHP2017}---for an overview of IVQR modeling, from identification conditions to estimation and inference.  For more recent progress,  see \cite{KS2017} and \cite{CGKL2019} for smoothed methods,  \cite{CL2018} and \cite{Z2018} for methods based on reformulation as mixed integer optimization (MIO),  \cite{MS2018} for moment-based estimators,  and \cite{KW2021} for a decentralization approach which decomposes the IVQR estimation problem into a set of conventional QR sub-problems.
 
The convolution smoothing method studied in this paper can be directly linked to the SEE approach in \cite{KS2017}.  The latter addressed the more challenging IVQR problem,  and derived both asymptotic mean squared error and normality for the SEE estimator when the dimension is fixed. Our study complements that of \cite{KS2017} in two ways. First, we provide a systematic analysis for smoothed (conventional) QR from an $M$-estimation viewpoint under the growing dimension setting.  Our results provide explicit finite-sample bounds for the estimation error,  Bahadur linearization error as long as their (multiplier) bootstrap counterparts. Asymptotic validity of the multiplier bootstrap is also rigorously established.
Secondly, we propose tailored computational methods for smoothed QR computation,  benefited from the use of non-negative kernels and the resulting local strong convexity.  Compared with generic optimization toolboxes for solving linear programs, the computational efficiency of the gradient-based algorithm for conquer  is considerably improved especially for large-scale problems with many (exogenous) regressors and massive sample size.  A potential application is empirical asset pricing via quantile regression,  extending the existing machine learning tools for average return forecast \citep{GKX2020}.  

In the presence of both exogenous and endogenous covariates, the advantage of smoothing is diluted because the non-convexity issue prevails. The MIO-based IVQR estimation procedure can be implemented by the Gurobi commercial MIO solver, which is free for academic use.  The MIO solver converges fast when the number of endogeous covariates is small and varies in the range of 5 and 20 \citep{Z2018}.  The MIO solver in moderate dimensions typically takes much longer to complete the optimization: optimal solutions may be found in a few seconds, but it can  take much longer to certify optimality via the lower bounds \citep{BKM2016}. \footnote{MIO solvers provide both feasible solutions and lower bounds to the optimal value. As the MIO solver progresses toward the optimal solution, the lower bounds improve and provide an increasingly better guarantee of suboptimality. It is the lower bounds that take so long to converge.} 

Recently, \cite{KW2021} proposed a ``decentralized" approach for IVQR estimation. The idea is to decompose the non-convex program into a set of  convex sub-problems, each being a weighted QR problem. Following \cite{CH2006}, consider the linear IVQR model
\#
	F^{-1}_{y| ( \bx , \bz ) } (\tau) = \bx^\T \bbeta_x^* + \bd^\T \bbeta_d^*  = \sum_{k=1}^{p_x}  x_k \beta^*_{x,k}  + \sum_{j=1}^{p_d} d_j  \beta^*_{d,j} ~~\mbox{ for some }~ \tau \in (0,1),  \nn
\#
where $\bx \in \RR^{p_x}$ is a vector of exogenous covariates, $\bd\in \RR^{p_d}$ is a vector of endogeous treatment variables, and $\bz \in \RR^{p_z}$ is a vector of instruments.  For simplicity, we restrict attention to the ``just-identified" setting where $p_z= p_d$. The corresponding population moment conditions are
\#
	 \EE  \left[      \mathbbm{1}\big( y \leq \bx^\T \bbeta_x^* + \bd^\T \bbeta_d^* \big)  \begin{pmatrix} \bx \\ \bz \end{pmatrix}  \right] = \textbf{0}. \nn
\#
Let $\{ (y_i, \bd_i, \bx_i, \bz_i) \}_{i=1}^n$ be independent observations drawn from the above IVQR model. Assume further that the model is parameterized such that $z_j/d_j>0$ for $j=1,\ldots, p_d$.
Under the parametrization $\btheta = (\btheta_1^\T , \theta_2, \ldots, \theta_J)^\T \in \RR^{p_x+p_d}$, where $\btheta_1 \in \RR^{p_x}$, $\theta_j\in \RR$ for $j=2,\ldots, J := p_d+1$, consider the following weighted empirical QR objection functions:
\#
 \hat Q_{1}(\btheta) & = \frac{1}{n} \sn \rho_\tau ( y_i - \bx_i^\T \btheta_1 - d_{i,1} \theta_2 - \cdots - d_{i, p_d} \theta_J ) ,   \nn \\
\hat Q_{j}(\btheta) & = \frac{1}{n} \sn \rho_\tau ( y_i - \bx_i^\T \btheta_1 - d_{i,1} \theta_2 - \cdots - d_{i, p_d} \theta_J ) \frac{z_{i,j-1}}{d_{i,j-1}}  , \ \ j = 2, \ldots, J.   \nn %\label{IVQR.loss2}
\#
For $j=1,\ldots, J$, let $\btheta_{-j}$ be the subvector of $\btheta$ with its $j$-th component removed. The first step of \cite{KW2021}'s decentralization method is to compute the sample best response (BR)  maps $\hat L_j(\btheta_{-j})$ ($j=1,\ldots, J$) by solving
\#
\hat L_1(\btheta_{-1} )  \in \argmin_{\wt \btheta_1 \in \RR^{p_x}} \hat Q_{1}(\wt \btheta_1 , \btheta_{-1})  ,  \quad   \hat L_j(\btheta_{-j} )   \in \argmin_{\wt \theta_j  \in \RR } \hat Q_{j}(\wt \theta_j , \btheta_{-j})  , \ \ j =2,\ldots, J. \label{IVQR.step1}
\#
In the second step, these estimated BR maps are used to form a sequential dynamical system, to which contraction-based algorithms or root-finding algorithms paired with nesting can be applied. The numerical studies in \cite{KW2021} suggest that estimation procedures based on nesting perform well and are computationally feasible  when $p_d$, the number of endogenous variables, is moderate.
The overall computational complexity of the procedure also depends on the choice of the algorithm for QR fitting in \eqref{IVQR.step1}. When the number of exogenous variables, $p_x$, is large in the range of hundreds to thousands, the proposed framework in this paper, along with the accompanying software \texttt{conquer}, provides a viable option to further reduce the computational cost of the above IVQR method in both estimation and inference. We leave a rigorous theoretical investigation (when $p_d$ is fixed, $p_x= p_x(n) \to \infty$ and $p_x/n \to 0$ as $n\to \infty$) as well as empirical applications with many (exogenous) regressors to future work.

\subsubsection{Nonparametric quantile regression}

Another setting in which one may have a growing number of covariates is nonparametric quantile regression. Suppose that the conditional $\tau$-quantile of $y$ given $\bx \in \cX \subseteq  \RR^p$ is a nonlinear function $f(\bx)= f(\bx; \tau)$, assumed to be sufficiently smooth but does not fall into any parametric family. Two well-documented approaches for fitting nonlinear QR models are based on (i) series approximation \citep{C2006}, and (ii) the use of kernels within a reproducing kernel Hilbert space (RKHS) framework \citep{TLSS2006}.

The idea for QR-series approximation is to first approximate $f$ by its ``projection" on the linear span of $m$ series/basis functions, and then fit the coefficients using the observed data. More specifically, let $ \bz(\bx) = (z_1(\bx), \ldots, z_{m}(\bx) )^\T$ be a vector of series approximating functions of dimension $m$, and define the coefficient vector $\btheta^* = \btheta^*(\tau)=(\theta^*_1, \ldots, \theta^*_m)^\T$ as a solution to 
\$
	 \min_{\btheta \in \RR^m}   \EE\big\{ \rho_\tau ( y - \langle \bz(\bx) , \btheta \rangle ) - \rho_\tau( y - f(\bx ; \tau) ) \big\} .
\$
Commonly used series functions with good approximation properties include B-splines, polynomials, Fourier series and compactly supported wavelets. \cite{C2006} established the consistency and rate of convergence for QR-series estimators at a single quantile index. More recently, \cite{BCCF2019} developed  large sample theory for QR-series coefficient process, including convergence rate and uniform strong approximations.
The choice of the parameter $m$, also known as the order of the series estimator, is crucial for establishing the balance between bias and variance. The former depends on the QR-series approximation error $R(\bx;\tau) :=  f(\bx;\tau) - \langle \bz(\bx), \btheta^* \rangle$. By allowing the dimension  $m= m_n \to \infty$ as $n \to \infty$ at a proper rate, the approximation error vanishes asymptotically, that is, $\sup_{\bx \in \cX} |R(\bx;\tau)| \to 0$. See Lemma~1 of \cite{BCCF2019} for the order of approximation error in the cases of polynomials and B-splines.

Given independent observations $\{ (y_i, \bx_i) \}_{i=1}^n$ from $(y, \bx)$, the QR-series estimator is defined as  $\hat f_{m}(\bx) = \langle \bz(\bx), \hat \btheta_m  \rangle$ with $\hat \btheta_m \in \argmin_{\btheta \in \RR^m} (1/n)\sn \rho_\tau ( y_i - \langle \bz(\bx_i) , \btheta \rangle )$.
With a properly chosen bandwidth, say $h\asymp (m /n)^{2/5}$, a convolution smoothed QR-series estimator is given by $\hat f_{m ,h }(\bx) =  \langle \bz(\bx), \hat \btheta_{m,h} \rangle$ with $\hat \btheta_{m,h} \in \argmin_{\btheta \in \RR^m} (1/n)\sn \ell_h( y_i - \langle \bz(\bx_i) , \btheta \rangle )$, where $\ell_h(\cdot)$ is the convolution smoothed check function  in \eqref{sqr}. Detailed discussions of bandwidth selection will be provided in Sections~\ref{sec:theory} and \ref{sec:numerical}.
When $m$ is large, convolution smoothing admits fast gradient-based algorithms (see Section~\ref{sec:algorithm} for details)  for computing  QR-series estimators over a fine grid of quantiles.
The computational and theoretical aspects of convolution smoothed QR will be investigated in Section~\ref{sec:algorithm} and Section~\ref{sec:theory}, respectively.

For quantile regression in RKHS,  with the smoothing splines estimation as a special case, by the representer theorem \citep{KW1971}, the estimator is defined as a solution to the following optimization problem
\#
	\min_{\bw=(w_1,\ldots, w_n)^\T  } \frac{1}{n} \sn \rho_\tau\Bigg( y_i - \sum_{j=1}^n w_j K(\bx_i, \bx_j ) \Bigg) + \lambda \cdot \bw^\T  \Kb  \bw ,  \label{RKHS-QR}
\#
where $\lambda >0$ is a penality level, $K: \RR^p \times \RR^p \to \RR$ is a positive definite kernel, and $\Kb = ( K(\bx_i , \bx_j ))_{1\leq i, j \leq n}$ is the corresponding $n\times n$ kernel matrix. In general, such an estimator can be computed by quadratic programming methods \citep{TLSS2006}. 
\cite{LLZ2007} derived the rate of convergence for RKHS-QR estimators, and designed a path-following algorithm with a complexity of $\cO(n^3)$. Computationally, we can apply gradient descent methods to solve a convolution smoothed version of the loss function in \eqref{RKHS-QR}. However, the statistically ``optimal" choice of the bandwidth remains unknown for smoothed QR in RKHS. We leave this as future research.

%%%%%%%%%%%%%%%%%%%%%%%%%%%%%%%%
%%%%%%%%%%%%%%%%%%%%%%%%%%%%%%%%
% Algorithm
%%%%%%%%%%%%%%%%%%%%%%%%%%%%%%%%
%%%%%%%%%%%%%%%%%%%%%%%%%%%%%%%%
\section{Computational methods for conquer} 
\label{sec:algorithm}
To solve optimization problems \eqref{sqr} and \eqref{bsqr.1} with non-negative weights, arguably the simplest algorithm is a vanilla gradient descent algorithm (GD). 
For a prespecified $\tau \in (0,1)$ and bandwidth $h>0$, recall that $ \hat Q_h (\bbeta) =  (1/n) \sn \ell_h(y_i - \langle \bx_i, \bbeta \rangle)$.
 Starting with an initial value $\bbeta^0 \in \RR^p$, at iteration $t=0, 1,2,\ldots$, GD computes
\#
  \bbeta^{t+1} =  \bbeta^t - \eta_t \cdot  \nabla \hat Q_h  (  \bbeta^t) =   \bbeta^t -   \frac{\eta_t}{n}  \sn   \big\{   \cK_h (   \langle \bx_i, \bbeta^t \rangle - y_i    ) - \tau  \big\} \bx_i ,
\# 
where $\eta_t>0$ is the stepsize.  In the classical GD method, the stepsize is usually obtained by employing line search techniques. However, line search is computationally intensive  for large-scale settings.
One of the most important issues in GD is to determine a proper update step $\eta_t$ decay schedule.
A common practice in the literature is to use a diminishing stepsize or a best-tuned fixed stepsize. Neither of these two approaches can be efficient, at least compared to the Newton-Frisch algorithm with preprocessing \citep{PK1997}.
Recall that the smoothed loss $\hat Q_h(\cdot)$ is twice differentiable with Hessian $\nabla^2  \hat Q_h (\bbeta)  = (1/n) \sn K_h ( y_i - \langle \bx_i, \bbeta \rangle   ) \bx_i \bx_i^\T$. It is therefore natural to employ the Newton-Raphson method, which at iteration $t$ would read
\#
  \bbeta^{t+1} =  \bbeta^t + \bd^t ~~\mbox{ with }~~ \bd^t :=   -  \bigl\{ \nabla^2  \hat Q_h (  \bbeta^t)  \bigr\}^{-1}  \nabla \hat Q_h  (  \bbeta^t)   .
 \label{NR}
\#
In practice,  the Newton method is often paired with Armoji stepsize: choose a stepsize $\lambda^t  =\max\{ 1, 1/2, 1/4, \ldots \}$ such that $\hat Q_h(\bbeta^t) - \hat Q_h(\bbeta^t + \lambda^t \bd^t) \geq  -c \lambda^t \nabla \hat Q_h(\bbeta^t) \bd^t$, where $c \in (0,1/2)$. Then redefine the current iterate as $\bbeta^{t+1} =  \bbeta^t + \lambda^t \bd^t$.
Since such a backtracking line search requires evaluations of the loss function itself, in the following remark we present the explicit expressions of the convolution smoothed check function for several commonly used kernels.

\begin{remark} \label{rmk1}
  Recall that the check function can be written as $\rho_\tau(u) =|u|/2+   (\tau-1/2)u$, which, after convolution smoothing, becomes $\ell_h (u) =(1/2)\int_{-\infty}^\infty |u + hv | K (v)  \, {\rm d} v +    (\tau-1/2)u$.

\begin{itemize}

\item (Gaussian kernel $K(u)=(2\pi)^{-1/2} e^{-u^2/2}$):  $\ell_{h}(u) = (h/2)\ell^{{\rm G}}(u/h) + (\tau - 1/2 ) u$, where $ \ell^{{\rm G}}(u)  :=   (2/\pi)^{1/2} e^{-u^2/2 } + u \{ 1- 2\Phi(-u) \}$.

\item (Logistic kernel $K(u) = e^{-u}/(1+ e^{-u})^2 $): $\ell_h (u) = (h/2) \ell^{{\rm L}}(u/h) + (\tau-1/2) u$, where $ \ell^{{\rm L}}(u) :=  u + 2\log(1+e^{-u})$.  Logistic kernel smoothed approximation of the check function dates back to \cite{A1982}, which is used as a technical device to simplify the analysis of the asymptotic behavior of a two-stage median regression estimator.

\item  (Uniform kernel  $K(u) = (1/2) \mathbbm{1} (|u|\leq 1)$):  $\ell_h(u) = (h/2)  \ell^{{\rm U}}(u/h) + (\tau - 1/2) u$, where $ \ell^{{\rm U}}(u)  :=  (u^2/2 + 1/2)\mathbbm{1}( |u|\leq 1) +   |u|  \mathbbm{1}(|u|>1)$ is a shifted Huber  loss \citep{H1973}.

\item  (Epanechnikov kernel $K(u) = (3 / 4) (1 - u^2) \mathbbm{1}(|u|\leq 1)$): $\ell_{h}(u) =  (h/2) \ell^{{\rm E}}(u/h) +(\tau - 1/2)u$, where $ \ell^{{\rm E}}(u) :=  (3u^2/4 - u^4/8+ 3/8) \mathbbm{1}(|u|\leq 1) + |u|\mathbbm{1}(|u|>1)$.

\item (Triangular kernel $K(u) = (1 - |u|) \mathbbm{1}(|u|\leq 1)$): $ \ell_{h}(u) = (h/2) \ell^{{\rm T}}(u/h) + ( \tau - 1/2) u$, where $ \ell^{{\rm T}}(u )  :=  ( u^2 - |u|^3/3 + 1/3 ) \mathbbm{1}(|u|\leq 1) + |u|\mathbbm{1}(|u| >1) $.
\end{itemize}

 \end{remark}

%%%%%%%%%%%%%%%%%%%%%%%
%%%%%%%%%%%%%%%%%%%%%%%
% Uniform Kernel --  Coordinate Descent
%%%%%%%%%%%%%%%%%%%%%%%
%%%%%%%%%%%%%%%%%%%%%%%
\subsection{The Barzilai-Borwein stepsize}
\label{subsec:uniform kernel}

In this section,  we propose to solve conquer by means of the gradient descent with a Barzilai-Borwein
update step \citep{BB1988}, which we refer to as the GD-BB algorithm. Motivated by quasi-Newton methods, the BB method has been proven to be very successful in solving nonlinear optimization problems.

 Computing the inverse of the Hessian when $p$ is large is an expensive operation at each Newton step  \eqref{NR}.  Moreover,  in circumstances where $h$ is small or $\tau$ is very close to 0 or 1,  $\nabla^2  \hat Q_h (\cdot)$ may have a large condition number,  thus leading to slow convergence.
For this reason,  many quasi-Newton methods seek a simple approximation of the inverse Hessian matrix, say $(\Jb^t)^{-1}$,  satisfying the secant equation $\Jb^t \bdelta^t = \bg^t$, where 
\#
 \bdelta^t  =  \bbeta^t -  \bbeta^{t-1} ~~\mbox{ and }~~ \bg^t =  \nabla \hat Q_h  ( \bbeta^t) - \nabla \hat Q_h (  \bbeta^{t-1}), \ \ t= 1, 2, \ldots . \label{BB1}
\#
To mitigate the computational cost of inverting a large matrix, the BB method chooses $\eta$ so that $\eta  \nabla \hat Q_h ( \bbeta^t)   =  (\eta^{-1} \Ib_p)^{-1}\nabla \hat Q_h  ( \bbeta^t)$ ``approximates" $(\Jb^t)^{-1} \nabla \hat Q_h  ( \bbeta^t)$. Since $\Jb^t$ satisfies $\Jb^t \bdelta^t = \bg^t$, it is more practical to choose $\eta$ such that $ (1/\eta ) \bdelta^t \approx \bg^t$ or $  \bdelta^t \approx  \eta \bg^t$. Via least squares approximations, one may use  $\eta_{1,t}^{-1}  = \argmin_{ \alpha  }  \|   \alpha  \bdelta^t - \bg^t \|_2^2$ or $\eta_{2,t}  = \argmin_{ \eta }  \|     \bdelta^t - \eta \bg^t \|_2^2$.
The BB stepsizes are then defined as 
\#
	\eta_{1,t} = \frac{\langle \bdelta^t , \bdelta^t\rangle }{\langle \bdelta^t, \bg^t \rangle } ~~\mbox{ and }~~ \eta_{2,t} = \frac{\langle \bdelta^t , \bg^t \rangle }{ \langle \bg^t, \bg^t \rangle } . \label{BB2}
\#
Consequently, the BB iteration takes the form
\#
	 \bbeta^{t+1} =  \bbeta^t - \eta_{\ell ,t}   \nabla \hat Q_h  ( \bbeta^t) , \quad     \ell = 1 {\rm ~or~} 2 . \label{BB3}
\#
Note that the BB step starts at iteration 1, while at iteration 0, we compute $\bbeta^1$ using   standard gradient descent with an initial estimate $\bbeta^0$. The procedure is summarized in Algorithm~\ref{algo:sqr}.
Based on extensive numerical studies, we find that at a fixed $\tau$, the number of iterations is insensitive to varying $(n,p)$ combinations. Moreover, as $h$ increases, the number of iterations declines because the loss function is ``more convex" for  larger $h$.  In Algorithm~\ref{algo:sqr}, the quantity $\delta >0$ is a prespecified tolerance level,  ensuring that the final iterate $ \bbeta^T$ satisfies $\| \nabla \hat Q_h( \bbeta^T) \|_2 \leq \delta$. Provided that $\delta \lesssim \sqrt{p/n}$, the statistical theory developed in Section~\ref{sec:theory} prevails.  In our \texttt{R} package \texttt{conquer}, we set $\delta = 10^{-4}$ as the default value; this value can also be specified by the user.

As $\tau$ approaches 0 or 1, the Hessian matrix becomes more ill-conditioned. As a result, the stepsizes computed in GD-BB may sometimes vibrate drastically, causing instability of the algorithm. Therefore, in practice, we set an upper bound for the stepsizes by taking $\eta_t = \min\{ \eta_{1,t}, \eta_{2,t} , 100\}$, for $t=1, 2,\ldots$.
Another case of an ill-conditioned Hessian arises when we have covariates with very different scales. 
In this case, the stepsize should be different for each covariate, and a constant stepsize will be either too small or too large for one or more covariates, which leads to slow convergence. 
To address this issue, we scale the covariate inputs to have zero mean and unit variance before applying gradient descent.

\begin{algorithm}[!t]
    \caption{ {\small  Gradient descent with Barzilai-Borwein stepsize (GD-BB) for solving conquer.}}
    \label{algo:sqr}
    \textbf{Input:} data vectors $\{(y_i, \bx_i)\}_{i=1}^n$, $\tau\in (0,1)$, bandwidth $h\in (0,1)$, initialization ${\bbeta}^{0}$, and gradient tolerance $\delta$.
    \begin{algorithmic}[1]
    \STATE  Compute $\bbeta^1 \leftarrow  \bbeta^0 -  \nabla \hat Q_h(  \bbeta^0)$
      \FOR{$t=1,2 \ldots $}
          \STATE $\bdelta^t \leftarrow    \bbeta^t -  \bbeta^{t-1}$, $\bg^t  \leftarrow \nabla \hat Q_h ( \bbeta^{t})-\nabla \hat Q_h ( \bbeta^{t-1})$
          \STATE $\eta_{1,t} \leftarrow \langle \bdelta^t, \bdelta^t \rangle / \langle\bdelta^t, \bg^t \rangle$, $\eta_{2,t} \leftarrow \langle \bdelta^t, \bg^t\rangle / \langle \bg^t, \bg^t \rangle$
          \STATE $\eta_t \leftarrow  \min\{\eta_{1,t},\eta_{2,t}, 100\}$ if $\eta_{1,t} >0$ and $\eta_t \leftarrow 1$ otherwise
          \STATE $  \bbeta^{t+1} \leftarrow  \bbeta^t -  \eta_t \nabla \hat Q_h ( \bbeta^{t})$
      \ENDFOR~when $\| \nabla \hat Q_h (   \bbeta^{t}) \|_2\le \delta$
    \end{algorithmic}
\end{algorithm}

\subsection{Warm start via asymmetric Huber regression}
A good initialization helps reduce the number of iterations for GD, and hence facilitates fast convergence.
Recall from Remark~\ref{rmk1} that with a uniform kernel, the smoothed check function is proximal to a Huber loss \citep{H1973}.  Motivated by this subtle proximity, we propose using the asymmetric Huber $M$-estimator as an initial estimate, and then proceed by iteratively applying gradient descent with BB update step.

Let $H_{\tau, \gamma}(u) = |\tau - \mathbbm{1}(u<0)|  \cdot  \{ (u^2/2) \mathbbm{1}(|u|\leq \gamma) + \gamma (|u| - \gamma/2) \mathbbm{1}(|u|> \gamma) \}$ be the asymmetric Huber loss parametrized by $ \gamma>0$. 
The asymmetric Huber $M$-estimator is then defined as 
\#
	\wt \bbeta_\gamma \in \argmin_{\bbeta\in \RR^p}  \hat \cL_\gamma(\bbeta), ~~\mbox{where}~~ \hat \cL_\gamma(\bbeta) = \frac{1}{n} \sn   H_{\tau, \gamma}(y_i - \langle \bx_i, \bbeta \rangle ). \label{huber.reg}
\#
The quantity $ \gamma$ is a shape parameter that controls the amount of robustness. 
The main reason for choosing a fixed (neither diminishing nor diverging) tuning parameter $ \gamma$ in \cite{H1981} is to guarantee robustness towards arbitrary contamination in a neighborhood of the model. This is at the core of the robust statistics idiosyncrasy. In particular, \cite{H1981} proposed $ \gamma=1.35 \sigma$ to gain as much robustness as possible while retaining 95\% asymptotic efficiency for normally distributed data, where $\sigma>0$ is the standard deviation of the random noise.  We estimate $\sigma$ using the median absolute deviation of the residuals at each iteration, i.e., ${\rm MAD}(\{ r_i^t \}_{i=1}^n)= \text{median}(|r_i^t -  \text{median}(r_i^t) |)$.

Noting that the asymmetric Huber loss is twice continuously differentiable, convex, and locally strongly convex, we use the GD-BB method described in the previous section to solve the optimization problem \eqref{huber.reg}.
Starting at iteration 0 with $\bbeta^{0,0} = \textbf{0}$, at iteration $t=0, 1,2,\ldots$, we  compute
\#
	  \bbeta^{0,t+1}  =  \bbeta^{0,t}  - \eta_t  \nabla \hat \cL_\gamma (  \bbeta^{0,t} )  =  \bbeta^{0,t}  + \frac{\eta_t}{n} \sn \psi_{\tau, \gamma}(y_i-\langle\bx_i,   \bbeta^{0,t}  \rangle ) \bx_i
\# 
with $\eta_t>0$ automatically obtained by the BB method, where $\psi_{\tau, \gamma}(u) = | \tau - \mathbbm{1} (u<0)| \cdot H_{\tau, \gamma}'(u) =| \tau - \mathbbm{1} (u<0)| \cdot \min\{\max(- \gamma, u),  \gamma \}$. The final iterate $\bbeta^{0,T'}$ for some $T'>1$ will be used as the initial value in Section~\ref{subsec:uniform kernel}.
We summarize the details in Algorithm~\ref{algo:huber}.

\begin{algorithm}[!t]
    \caption{ {\small  GD-BB method for solving \eqref{huber.reg}. }}
    \label{algo:huber}
    \textbf{Input:} $\{(y_i, \bx_i)\}_{i=1}^n$ and convergence criterion $\delta$.
    \begin{algorithmic}[1]
    \STATE Initialize ${\bbeta}^{0,0} = \textbf{0}$% or $\tilde{\bbeta}^{(0)} \sim \mathcal N(\textbf{0},  p^{-1} \Ib_p)$
    \STATE  Compute $ \gamma^0 = 1.35 \cdot {\rm MAD}(\{ r_i^0 \}_{i=1}^n)$, where $r^0_i \leftarrow y_i - \langle \bx_i,  \bbeta^{0,0} \rangle$, $i=1,\ldots, n$, where $\mathrm{MAD}(\cdot)$ is the median absolute deviation
    \STATE $\bbeta^{0,1} \leftarrow   \bbeta^{0,0} -  \nabla \hat \cL_{\gamma^0}( \bbeta^{0,0})$
      \FOR{$t=1,2 \ldots $}
      	\STATE $\gamma^t = 1.35 \cdot {\rm MAD}(\{ r^t_i\}_{i=1}^n)$, where $r^t_i \leftarrow y_i - \langle \bx_i,   \bbeta^{0,t} \rangle$, $i=1,\ldots, n$
          \STATE $\bdelta^t \leftarrow   \bbeta^{0,t} -  \bbeta^{0, t-1}$, $\bg^t  \leftarrow \nabla  \hat \cL_{\gamma^t}(    \bbeta^{0,t})-\nabla  \hat \cL_{\gamma^t}( \bbeta^{0,t-1})$
          \STATE $\eta_{1,t} \leftarrow \langle \bdelta^t, \bdelta^t \rangle / \langle\bdelta^t, \bg^t \rangle$, $\eta_{2,t} \leftarrow \langle \bdelta^t, \bg^t\rangle / \langle \bg^t, \bg^t \rangle$. 
          \STATE $\eta_t \leftarrow  \min\{\eta_{1,t},\eta_{2,t}, 100\}$ if $\eta_{1,t} >0$ and $\eta_t \leftarrow 1$ otherwise
          \STATE $\bbeta^{0, t+1} \leftarrow \bbeta^{0, t}-  \eta_t  \nabla  \hat \cL_{\gamma^t}(  \bbeta^{0,t})$
      \ENDFOR~when $\| \nabla \hat \cL_{\gamma^t}(  \bbeta^{0,t}) \|_2\le \delta$
    \end{algorithmic}
\end{algorithm}

\begin{remark}
The asymmetric Huber loss $H_{\tau, \gamma}(\cdot)$ approximates the check function $\rho_\tau(\cdot)$  as $\gamma \rightarrow 0$. Therefore, an alternative method for QR computing is to solve asymmetric Huber regression via gradient descent with a shrinking gamma. To evaluate its performance, we implement the above idea by setting $\gamma^{t} = c \cdot \gamma^{t-1}$ for some $c\in (0,1)$ at the $t$-th iteration.  We found that the aforementioned idea is not numerically stable  across several simulated data sets, unless one controls the minimal magnitude of $\gamma$ very carefully.  Furthermore,  the  obtained  solution has a higher estimation error than that of conventional QR and conquer.
\end{remark}

 %%%%%%%%%%%%%%%%%%%%%%%%%%%%%%%%%
% Theoretical Properties
%%%%%%%%%%%%%%%%%%%%%%%%%%%%%%%%%
%%%%%%%%%%%%%%%%%%%%%%%%%%%%%%%%%
\section{Statistical analysis}
\label{sec:theory}

Under the linear quantile regression model in \eqref{eq:lqr},  we write, for convenience, the generic data vector $(y,\bx)$ in a linear model form: given a quantile level $\tau\in (0,1)$ of interest,
\#
 y = \langle \bx, \bbeta^*(\tau) \rangle +  \varepsilon(\tau),  \label{linear.qr}
\#
where the random variable $\varepsilon(\tau)$ satisfies $\PP\{\varepsilon(\tau) \leq 0   | \bx \} = \tau$. Let $ f_{\varepsilon |\bx}(\cdot )$ be the conditional density function of the regression error $\varepsilon = \varepsilon(\tau)$ given $\bx=(x_1, \ldots, x_p)^\T$ ($p\geq 2$).
We first derive upper bounds for the smoothing bias under mild regularity conditions on the conditional density $ f_{\varepsilon |\bx}$ and the kernel function.  
{For any vector $\bu \in \RR^p$, we write $\bu_-\in \RR^{p-1}$ as the sub-vector of $\bu$ with its first component removed.
Recall that $x_1 \equiv 1$,  and $ \bx_- = (x_2, \ldots , x_p)^\T \in \RR^{p-1}$ is assumed to be random.
Without loss of generality, we assume $\bmu_- := \EE(\bx_-) = \textbf{0}$ throughout this section; otherwise,  set $\wt \bx= (1,  (\bx_- - \bmu_-)^\T )^\T$, so that model \eqref{linear.qr} can be written as $ y   = \langle \wt \bx, \wt \bbeta^* \rangle + \varepsilon$,
where $\wt \bbeta^* = ( \wt \beta_1^* ,   \beta^*_2 ,\ldots, \beta^*_p )^\T$ with $ \wt \beta_1^*= \beta^*_1 +  \langle \bmu_-  , \bbeta^*_- \rangle$. The analysis then applies to $\{ (y_i, \wt \bx_i) \}_{i=1}^n$,  and the probabilistic bounds for $\wt \bbeta^*$ naturally lead to those for $\bbeta^*$.}

%%%%%%%%%%%%%%%%%%%%%%%%%%%%%%%%%
%%%%%%%%%%%%%%%%%%%%%%%%%%%%%%%%%
% Finite Sample Bias by Smoothing
%%%%%%%%%%%%%%%%%%%%%%%%%%%%%%%%%
%%%%%%%%%%%%%%%%%%%%%%%%%%%%%%%%%
\subsection{Smoothing bias}
\label{sec:theory:bias}

\begin{cond}[Kernel function] \label{cond.kernel} Let $K(\cdot)$ be a symmetric and non-negative function that integrates to one, that is, $K(u)= K(-u)$, $K(u) \geq 0$ for all $u \in \RR$ and $\int_{-\infty}^\infty K(u) \, {\rm d}u  = 1$. Moreover, $K(\cdot)$ is uniformly bounded with $\kappa_u := \sup_{u\in \RR} K(u) < \infty$. 
\end{cond}

We will use the notation $\kappa_k  = \int_{-\infty}^\infty |u |^k K(u) \,{\rm d} u$ for $k \geq 1$. Furthermore, we define the population smoothed loss function $Q_{  h} (\bbeta) = \EE \{\hat Q_{ h}(\bbeta)\}$, $\bbeta \in \RR^p$ and the pseudo parameter 
\#
	 \bbeta^*_h (\tau) \in \argmin_{\bbeta \in \RR^p } Q_{ h}(\bbeta),  \label{pseudo.parameter}
\#
which is the population minimizer under the smoothed quantile loss. 
For simplicity, we write $\bbeta^* = \bbeta^*(\tau)$ and $\bbeta^*_h = \bbeta_h^*(\tau)$ hereinafter.
In general, $\bbeta^*_h$ differs from $\bbeta^*$, and we refer to $\| \bbeta^*_h - \bbeta^* \|_2$ as the approximation error or smoothing bias.

\begin{cond}[Conditional density]  \label{cond.reg}  There exist $\underbar{$f$} >0$ such that $ f_{\varepsilon |\bx}(0 )\geq \underbar{$f$}$ almost surely (for all $\bx$).
Moreover, there exists a constant $l_0 >0$ such that $|f_{\varepsilon |\bx}(u)- f_{\varepsilon |\bx}(v)| \leq l_0  |u -v |$ for all $u , v\in \RR$ almost surely (over $\bx$).
\end{cond}

\begin{cond}[Random design: moments]  \label{cond.x-mom} The (random) vector $\bx \in \mathbb{R}^p$ of covariates satisfies $m_3 := \sup_{\bu \in \mathbb S^{p-1}}  \EE ( | \langle  \bu, \bSigma^{-1/2} \bx \rangle |^3 )  <\infty$, where   $\bSigma = \EE(\bx \bx^\T)$ is positive definite.
\end{cond}

{Condition~\ref{cond.x-mom} requires that all the one-dimensional marginals of $\bSigma^{-1/2}\bx$ have bounded third absolute moments.
When $\bx_-$ follows a multivariate normal distribution,  Condition \ref{cond.x-mom} holds trivially. 
We refer to Remarks~\ref{rmk:subG} and \ref{rmk:elliptical} below for more examples.} The following result characterizes the smoothing bias from a non-asymptotic viewpoint.

\begin{proposition} \label{prop:bias}
Assume Conditions \ref{cond.kernel}--\ref{cond.x-mom} hold, and let the bandwidth  satisfy $0 < h <   \frac{1}{l_0 \{  \kappa_1 + (m_3 \kappa_2)^{1/2} \}}  \underbar{$f$}$.
Then,   $\bbeta^*_h$ is the unique minimizer of $\bbeta \mapsto Q_{ h}(\bbeta)$ and satisfies
\#
	   \delta_h :=   \| \bbeta^*_h - \bbeta^* \|_{\bSigma }  < \frac{l_0 \kappa_2 h^2 }{ \underbar{$f$} - l_0 \kappa_1 h }.  \label{bias.ubd}
\#
In addition, assume $f_{\varepsilon |\bx}(\cdot)$ is continuously differentiable and satisfies almost surely (over $\bx$) that
$ |  f'_{\varepsilon  |\bx } (u) - f'_{ \varepsilon |\bx } (0) | \leq l_1  | u  |$ for some constant $l_1 >0$. Then
\# \label{bias.leading}
	 \biggl\|  \bSigma^{-1/2}\Jb( \bbeta^*_h - \bbeta^* ) +   \frac{1}{2} \kappa_2 h^2  \cdot \bSigma^{-1/2} \EE   \bigl\{  f_{\varepsilon |\bx}'( 0 )\bx \bigr\}    \biggr\| _2 \leq \frac{	1}{6 } l_1 \kappa_3  h^3 +  \frac{1}{2} l_0  m_3 \delta_h^2   + l_0 \kappa_1 h  \delta_h   ,
\#
where $\Jb = \EE  \{  f_{\varepsilon |\bx}( 0 ) \bx \bx^\T   \}$.
\end{proposition}

To better understand the bounds \eqref{bias.ubd} and \eqref{bias.leading}, note that $\| \bbeta^*_h -\bbeta^* \|^2 _{\bSigma}= \EE \langle \bx, \bbeta^*_h -\bbeta^* \rangle^2$ is the average prediction smoothing error. Interestingly, the upper bound on the right-hand side is dimension-free given $h$ as long as the uniform third moment $m_3$ in Condition~3.3 is dimension-free.  See Remarks~\ref{rmk:subG} and \ref{rmk:elliptical} for examples of multivariate distributions on $\RR^p$ that have dimension-free uniform moment parameter. Another interesting implication is that, when both $f_{\varepsilon |\bx}(0)$ and $f'_{\varepsilon |\bx}(0)$ are independent of $\bx$, i.e., $f_{\varepsilon |\bx}(0)=f_{\varepsilon  }(0)$ and $f'_{\varepsilon |\bx }(0)=f'_{\varepsilon   }(0)$, the leading term in the  bias simplifies to 
\#
 \frac{1}{2} \kappa_2 h^2 \cdot \Jb^{-1} \EE   \bigl\{  f_{\varepsilon |\bx}'(0)\bx \bigr\} = \frac{f'_{\varepsilon  }(0)}{2f_{ \varepsilon  }(0)} \kappa_2 h^2 \cdot \bSigma^{-1} \EE (\bx)   =  \frac{f'_{\varepsilon  }(0)}{2f_{\varepsilon  }(0)} \kappa_2 h^2 \cdot \begin{bmatrix} 
1  \\
 \textbf{0}_{p-1}
\end{bmatrix} . \nn
\#
In other words, the smoothing bias is concentrated primarily on the intercept. 
In the asymptotic setting where $p$ is fixed, and $h=o(1)$ as $n\to \infty$, we refer to Theorem~1 in \cite{FGH2019}  for the expression of asymptotic bias.

\subsection{Finite sample theory}
\label{sec:estimation error}

In this section, we provide two non-asymptotic results, the concentration inequality and the Bahadur-Kiefer representation, for the conquer estimator under random design.

\begin{cond}[Random design: sub-exponential case] \label{cond.predictor}
The predictor $\bx = (x_1,\ldots, x_p)^\T \in \RR^p$ is {\it sub-exponential} with $x_1\equiv 1$ and $\EE(x_j) =0$ for $j=2,\ldots, p$.  That is, there exists $ \upsilon_0 >0$ such that $\PP\{ |\langle \bu  ,   \bw  \rangle |   \geq \upsilon_0     t \} \leq   e^{-t }$ for all $\bu \in  \mathbb{S}^{p-1}$ and $t \geq 0$, where $\bw = \bSigma^{-1/2} \bx $ with  $\bSigma =\EE(\bx \bx^\T)$ being positive definite.
\end{cond}

Condition~\ref{cond.predictor} assumes a sub-exponential condition on the random covariates, which encompasses the bounded case considered by \cite{FGH2019}.
For the standardized predictor $\bw = \bSigma^{-1/2}\bx$, we define the uniform moment parameters (including $m_3$ that first occurred in Condition~\eqref{cond.x-mom})
\#
  m_k = \sup_{\bu \in \mathbb{S}^{p-1}  }  \EE | \langle \bu, \bw  \rangle |^k , \ \  k = 1,2 ,\ldots, \label{def.moment}
\#
with $m_2=1$.  In particular, $m_4$ can be viewed as the uniform kurtosis parameter.
Under Condition~\ref{cond.predictor}, a straightforward calculation shows that  $m_k \leq    \upsilon_0^k k!$, valid for all $k\geq 1$.

\begin{remark} \label{rmk:subG}
{The parameter $\upsilon_0$ is often referred to as the sub-exponential parameter, which along with the sub-Gaussian parameter $\upsilon_1$ defined in Condition~\ref{cond.predictor2} below, plays an important role in non-asymptotic analysis of statistical models with growing dimensions \citep{RV2018,W2019}. For many ``nice" distributions on $\RR^p$,  both $\upsilon_0$ and $\upsilon_1$ are absolute constants and thus are dimension-free. In the following, we list several $p$-dimensional multivariate distributions, all of which have a dimension-free sub-exponential parameter $\upsilon_0$.  
\begin{enumerate} 
\item[(i)] (Multivariate normal).  Let $\bx \sim \cN(\textbf{0} , \bSigma)$ for some  positive definite matrix $\bSigma\in \mathbb{R}^{p\times p}$. 
\item[(ii)]  (Multivariate symmetric Bernoulli).  Let $\bx = (x_1, \ldots, x_p)^\T \sim {\rm Unif} ( \{ -1, 1\}^n)$. That is, $x_1,\ldots, x_p$ are independent, and satisfy $\PP(x_j = 1 ) = \PP(x_j = -1) = 1/2$. 
\item[(iii)] (Uniform distribution on the sphere).   Let $\bx = (x_1, \ldots, x_p)^\T$ be a random vector uniformly distributed on the Euclidean sphere in $\RR^p$ with center at the origin and radius $p^{1/2}$.  
\item[(iv)] (Uniform distribution on the Euclidean ball).  Let $\bx = (x_1, \ldots, x_p)^\T$ be a random vector uniformly distributed on the Euclidean ball $\BB^p(p^{1/2})$ in $\RR^p$.   
\item[(v)] (Uniform distribution on the unit cube).   Let $\bx = (x_1, \ldots, x_p)^\T$ be a random vector uniformly distributed on the unit cube $[-1,1]^p$. That is,  $x_1,\ldots,x_p$ are independent from Unif\,$[-1,1]$.
\item[(vi)]   (Uniform distribution on the $\ell_1$-ball).   Let $\bx = (x_1, \ldots, x_p)^\T$ be a random vector uniformly distributed on the $\ell_1$-ball $\{ \bu \in \RR^p: \| \bu \|_1 \leq r\}$ for $r \asymp p$.
\end{enumerate}
The multivariate distributions from examples (i)--(v) are all sub-Gaussian with a constant parameter,  and hence are also sub-exponential. The distribution from the last example is not sub-Gaussian, but is sub-exponential with a constant parameter.  We refer to Section~3.4 in \cite{RV2018} for a detailed introduction of sub-exponential and sub-Gaussian random vectors, including examples for which the sub-Gaussian parameter does depend on and grow with $p$.}
\end{remark}

{
Another important multivariate distribution is the {\it elliptical distribution}.  We say a  random vector $\bx\in \RR^p$ follows an elliptical distribution,  denoted by $\bx \sim {\rm ED} (\bmu, \bSigma, \xi)$,  if it has a stochastic representation $\bx \stackrel{d}{=} \bmu + \xi \Ab \Ub$, where $\xi$ is a random scalar, $\Ub$ is a random vector uniformly distributed on the unit sphere $\mathbb{S}^{p-1}$ and is independent of $\xi$,  and  $\Ab$ is a deterministic matrix satisfying $\bSigma = \Ab \Ab^\T$. }

\begin{remark} \label{rmk:elliptical}
{
Assume $\bx \sim {\rm ED} (\textbf{0}, \bSigma, \xi)$ for some $\bSigma\in \RR^{p\times p}$ and a random variable $\xi\in \RR$.  With slight abuse of notation,  write $\bx =\xi \Ab \Ub$. Then, for any unit vector $\bu \in \mathbb{S}^{p-1}$, 
\$
	| \langle \bu, \bSigma^{-1/2} \bx \rangle | = | \langle  \Ab^\T  \bSigma^{-1/2} \bu , \Ub \rangle | \cdot | \xi |  \leq \|   \Ab^\T  \bSigma^{-1/2} \bu \|_2 \cdot | \xi| \leq  | \xi | .
\$
This implies that (i) if $\EE |\xi|^3 < \infty$,  Condition~\ref{cond.x-mom} holds, (ii) if $\xi$ is sub-exponential,  Condition~\ref{cond.predictor} holds with a dimension-free $\upsilon_0$, and (iii) if $\xi$ is sub-Gaussian,  Condition~\ref{cond.predictor2} is satisfied with a dimension-free $\upsilon_1$.
}
\end{remark}

 \begin{theorem} \label{thm:concentration}
 Assume Conditions~\ref{cond.kernel}, \ref{cond.reg} and \ref{cond.predictor} hold. For any $t > 0$, the smoothed quantile regression estimator $\hat \bbeta_{h}$ with  $  \underbar{$f$}^{-1} m_3^{1/2} \upsilon_0 \sqrt{(p+t)/n} \lesssim h \lesssim \underbar{$f$} m_3^{-1/2}$ satisfies the bound
 \#
 	\| \hat \bbeta_h  - \bbeta^*  \|_{\bSigma} \leq     \frac{C}{\underbar{$f$}} \Biggl\{  \upsilon_0 \sqrt{  \frac{\log_2(1/h) + p+t}{n}} + l_0 \kappa_2 h^2 \Biggr\} , \label{concentration.ineq}
 \#
 with probability at least $1-2e^{-t}$, where $C>0$ is an absolute constant, and $\log_2(x) := \log\log(x\vee 1)$.
 \end{theorem}
 
Under a high probability statement, the estimation error in \eqref{concentration.ineq} is upper bounded by two terms, $\underbar{$f$}^{-1} l_0\kappa_2 h^2$ and $\underbar{$f$}^{-1} \upsilon_0 \sqrt{(p+t)/n}$,  which can be interpreted as  the bias and statistical rate of convergence,  respectively.   The parameter $t\geq 0$ controls the confidence level through $1-2e^{-t}$. 
The additional factor $\log_2(1/h)$ in the upper bound is a consequence of the peeling argument, which can be removed via a more refined analysis yet under slightly stronger technical conditions; see Section~\ref{proof:thm3.1.refined} in the supplement for details.
Adjusting the proof by changing high probability bounds to $\cO_{\PP}$ statements, it can be shown that $\| \hat \bbeta_h - \bbeta^* \|_{\bSigma} = \cO_{\PP}(\sqrt{p/n}  )$ under the condition $h=\cO((p/n)^{1/4})$ and $\sqrt{p/n}=\cO(h)$. {Next we explain the bandwidth constraint $\sqrt{p/n} \lesssim h \lesssim 1$ required in Theorem~\ref{thm:concentration} and all the other results below.  On one side,  the smoothing parameter should be sufficiently small,  typically $h=h_n \to 0$, so that the smoothing bias is negligible and does not change the target parameter to be estimated. On the other side,  the bandwidth cannot be too small in the sense that we need $h\gtrsim \sqrt{p/n}$. Intuitively, this is because the main motivation for smoothed QR is to seek a tradeoff between statistical rate of convergence and computational precision (unless the data is noiseless).  The standard QR estimator $\hat \bbeta=\hat \bbeta(\tau)$ has a convergence rate $\| \hat \bbeta - \bbeta^* \|_2 = \cO_{\PP}(\sqrt{p/n})$ under the growth condition $p \sup_{\bx \in \cX} \| \bx \|_2^2 \cdot (\log n)^2= o(n)$; see Theorem~1 in \cite{BCCF2019}. Here $\cX \subseteq \RR^p$ is the support of the covariate vector $\bx \in \RR^p$.  Therefore,  smoothing will become redundant if the bandwidth is set at a level below the best possible statistical convergence radius.
 }
 
 {
 Our results provide non-asymptotic bounds via high probability statements, which compliment the classical Big OP ($\cO_{\PP}$) and little op ($o_{\PP}$)  statements frequently used in statistics and econometrics. Probabilistic bounds of this kind can also be extended to analyze high-dimensional models \citep{BC2011,WWL2012} or nonparametric methods \citep{BCCF2019}.
 }

Next, we establish a Bahadur representation for the conquer estimator, which lays the theoretical foundation the ensuing statistical inference.
To this end, we impose a slightly more stringent sub-Gaussian condition on the covariates.

 \begin{cond}[Random design: sub-Gaussian case] \label{cond.predictor2}
 The predictor $\bx = (x_1,\ldots, x_p)^\T \in \RR^p$ is {\it sub-Gaussian} with $x_1\equiv 1$ and $\EE(x_j) =0$ for $j=2,\ldots, p$.  That is, there exists $ \upsilon_1 >0$ such that $\PP\{ |\langle \bu  ,   \bw  \rangle |   \geq \upsilon_1    t \} \leq  2 e^{-t^2/2}$ for all $\bu \in \mathbb{S}^{p-1}$ and $t\geq 0$, where $\bw = \bSigma^{-1/2} \bx$.
 \end{cond}
 
 We are primarily concerned with the cases where $\upsilon_1$ is a dimension-free constant; see the Remarks~\ref{rmk:subG} and \ref{rmk:elliptical}. 
 
%%%%%%%%%%%%%%%%%%%%%%%%%%%%%%%%%%%%%%%
%%%%%%%%%%%%%%%%%%%%%%%%%%%%%%%%%%%%%%%
% Bahadur Representation
%%%%%%%%%%%%%%%%%%%%%%%%%%%%%%%%%%%%%%%
%%%%%%%%%%%%%%%%%%%%%%%%%%%%%%%%%%%%%%%
\begin{theorem} \label{thm:bahadur}
In addition to Conditions~\ref{cond.kernel}, \ref{cond.reg} and \ref{cond.predictor2}, assume $\sup_{u  \in \RR} f_{\varepsilon | \bx} (u) \leq \bar f $ almost surely.
Let $t > 0$, and suppose the sample size $n$ and bandwidth $h$ satisfy $\underbar{$f$}^{-1 }  m_3^{1/2}\upsilon_1   \sqrt{(p+t)/n} \lesssim h \lesssim \underbar{$f$} m_3^{-1/2}$. Then, with probability at least $1-3e^{-t}$, 
 \#
 	\Biggl\|  \bSigma^{-1/2}    \Jb_h(     \hat \bbeta_h  - \bbeta^* )  -  \frac{1}{n}     \sn  \bigl\{  \tau -    \cK_h(- \varepsilon_i  )  \bigr\}  \bSigma^{-1/2} \bx_i   \Biggr\|_2 \leq  C  \Biggl(  \frac{p+ t}{ n  h^{1/2} } +   h^{3/2}\sqrt{\frac{p+t}{n}}  + h^4\Biggr) , \label{bh}
 \#
 where  $\Jb_h = \nabla^2 Q_h(\bbeta^*) =  \EE \bigl\{  K_h(\varepsilon) \bx \bx^\T \bigr\}$, $\cK_h(u) = \int_{-\infty}^{u/h} K(v) \, {\rm d}v$, and $C>0$ is a constant depending only on $(\upsilon_1, \kappa_2, \kappa_u, l_0,\bar f, \underbar{$f$})$. When $\Jb_h$ on the left-hand side of \eqref{bh} is replaced by $\Jb = \EE \{ f_{\varepsilon | \bx} (0) \bx \bx^\T \}$, the upper bound is of order $(p+t)/(n h^{1/2}) + h \sqrt{(p+t)/n} +h^3$.
\end{theorem}

{
With growing dimensions (many regressors), Theorem~\ref{thm:bahadur} is directly comparable to and complements Theorem~2 in \cite{BCCF2019}, although the latter concerns the linear approximation of the quantile regression process. To see the connection, we write 
$n^{1/2} (\hat \bbeta_h - \bbeta^*) = \Jb_h^{-1}\bU_h  + \br_h$, where $\bU_h = n^{-1/2}\sn (1 - \EE)    \{  \tau -    \cK_h(- \varepsilon_i  )    \}    \bx_i$ is a zero-mean random vector, and the remainder $ \br_h$ is such that $\|  \br_h \|_2 \lesssim (p+t)/(nh)^{1/2} + n^{1/2} h^2$ with high probability. Minimizing the right-hand side over $h$ in terms of order leads to a convergence rate $p^{4/5}/n^{3/10}$. For standard QR with fixed design, Theorem~2 in \cite{BCCF2019} implies $n^{1/2} (\hat \bbeta  - \bbeta^*) = \Jb^{-1}\bU   + \br$, where $\bU  = n^{-1/2}\sn \{  \tau -  \mathbbm{1}(\varepsilon_i \leq 0 )   \}    \bx_i$, and $\| \br \|_2 = \cO_{\PP}(p^{3/4} \zeta_p (\log n)^{1/2} n^{-1/4} )$, where $\zeta_p = \sup_{\bx \in \cX } \| \bx \|_2 $.
From an asymptotic perspective, QR estimator has the advantage of being (conditionally) pivotal asymptotically. However, possibly due to the non-smoothness of the check function, the linear approximation error has a slower rate of convergence $(p^5/n)^{1/4}$ even for bounded design, i.e., $\zeta_p \leq B p^{1/2}$ for some constant $B>0$. For the conquer estimator, although the linear term $\bU_h$ is not pivotal, we will show that Rademacher multiplier bootstrap provides accurate approximations both theoretically and numerically.
}

The Bahadur representation can be used to establish the limiting distribution of the estimator or its functionals.
Here we consider a fundamental statistical inference problem for testing the linear hypothesis $H_0:  \langle \ba , \bbeta^* \rangle =  0$, where $\ba \in \RR^p$ is a deterministic vector that defines a linear functional of interest.
It is then natural to consider a test statistic that depends on $n^{1/2}  \langle \ba, \hat \bbeta_h  \rangle$.
Based on the  non-asymptotic result in Theorem~\ref{thm:bahadur}, we establish a Berry-Esseen bound for the linear projection of the conquer estimator.

\begin{theorem} \label{thm:clt}
Assume that the conditions in Theorem~\ref{thm:bahadur} hold, and $\sqrt{(p+\log n)/n}\lesssim h \lesssim 1$. Then,
\#
  \Delta_{n,p}( h) :=\sup_{x\in \RR , \, \ba \in \RR^p } \left| \PP\big(  n^{1/2}   \sigma_{ h}^{-1} \langle \ba, \hat \bbeta_h - \bbeta^*   \rangle   \leq x    \big) - \Phi(x ) \right| \lesssim \frac{p + \log n}{(n h)^{1/2}} + n^{1/2} h^2  , \label{linear.clt}
\#
where $\sigma_{h}^2 = \sigma_{ h}^2(\ba) = \ba^\T \Jb_h ^{-1} \EE \big[  \{ \cK _h(-\varepsilon) - \tau \}^2  \bx \bx^\T   \big]  \Jb_h ^{-1} \ba$,  where $\Phi(\cdot)$ denotes the standard normal distribution function.  Moreover,  
$$
	\sup_{\ba \in \RR^p} \Bigg|\frac{\sigma_{ h}^2(\ba)}{  \ba^\T \Jb_h ^{-1} \bSigma \Jb_h ^{-1} \ba } - \tau(1-\tau)  \Bigg| = O(h)~\mbox{ as } h \to 0.
$$
If, in addition, that $f_{\varepsilon|\bx}(\cdot)$ is twice continuously differentiable and satisfies $|f''_{\varepsilon| \bx}(u) - f''_{\varepsilon | \bx}(v)| \leq l_2(\bx) |u-v|$ for all $u, v\in \RR$ and $\bx\in \RR^p$, and $l_2:\RR^p \to \RR^+$ is such that $\EE \{ l_2^2(\bx) \}\leq C$ for some $C>0$. Then,  
\#
& \sup_{x\in \RR , \, \ba \in \RR^p } \left| \PP\big(     n^{1/2}   \sigma_{ h}^{-1}  \langle \ba, \hat \bbeta_h - \bbeta^*   +   0.5 \kappa_2 h^2  \Jb_h^{-1}  \EE \{ f'_{\varepsilon | \bx}(0) \bx \} \rangle   \leq x    \big) - \Phi(x ) \right|  \nn \\
 & ~~~~~~~~~~~~~~~~~~~~~~~~~~~~~~~~~~~~~~~~~~~~~~~~~~\lesssim \frac{p + \log n}{(n h)^{1/2}} +  (p+\log n)^{1/2} h^{3/2}+ n^{1/2} h^4. \label{linear.clt2}
\#
\end{theorem}
  
Theorem~\ref{thm:clt} shows that for certain choice of bandwidth $h=h_n \to 0$,  all the linear functionals of $\hat \bbeta_h$, after properly standardized, are asymptotically normal as $n, p \to  \infty$ subject to some conditions.  For example,  if $h$ satisfies $h=o(n^{-1/4})$, then the smoothing bias does not affect the asymptotic distribution.   The Berry-Esseen bound \eqref{linear.clt} immediately yields a large-$p$ asymptotic result. Taking $h=h_n = \{ (p+\log n) /n\}^{2/5}$ therein,  the Gaussian approximation error $\Delta_{n,p}(h)$ is of order $(p+ \log n)^{4/5} n ^{-3/10}$. Consequently, $n^{1/2} \langle \ba , \hat \bbeta_h  - \bbeta^*  \rangle$, for any given (deterministic) vector $\ba \in \RR^p$, is asymptotically normally distributed as long as $p^{8/3} /n \to 0$, which improves the best known growth condition on $p$ for quantile regression \citep{W1989}.

\begin{remark}[Optimal bandwidth under AMSE]  \label{rmk:AMSE}
{
From the refined Berry-Esseen bound \eqref{linear.clt2} under an additional smoothness condition on $f_{\varepsilon | \bx}(\cdot)$, we can characterize more precisely the asymptotic bias and variance-covariance matrix of $n^{1/2} (\hat \bbeta_h - \bbeta^*)$, thus leading to the asymptotic mean squared error (AMSE). In fact, the asymptotic distribution of $n^{1/2} (\hat \bbeta_h - \bbeta^*)$ can be approximated by that of 
$$
	\bg_h  \stackrel{{\rm d}}{=}  \cN\Bigg(  \frac{1}{2} \kappa_2  n^{1/2} h^2  \Jb_h^{-1} \EE \{ f'_{ \varepsilon|\bx } (0)  \bx \} , \,  \Jb_h^{-1}   \EE \big[  \{ \cK _h(-\varepsilon) - \tau \}^2  \bx \bx^\T   \big]    \Jb_h^{-1}   \Bigg) .
$$
After some calculations (see Section~\ref{sec:AMSE} of the Appendix), we show that the AMSE of $n^{1/2} (\hat \bbeta_h - \bbeta^*)$  is determined by
\#
   \Jb_h^{-1} \bSigma^{ 1/2} \Bigg[  \tau(1-\tau ) \Ib_p - 2 \bar{\kappa}_1 h \Hb   +  \frac{\kappa_2^2}{4}  n h^4   \EE (\bb) \EE (\bb^\T)   \Bigg] \bSigma^{ 1/2}  \Jb_h^{-1} , \label{asym.cov}
\#
where $\bb=f'_{ \varepsilon|\bx } (0) \bw = f'_{ \varepsilon|\bx } (0) \bSigma^{-1/2}\bx$, $\Hb = \EE \{ f_{\varepsilon |\bx}(0) \bw \bw^\T \}$, and $\bar \kappa_1=  \int_{-\infty}^\infty  u K(u)\cK(u) {\rm d} u >0$.  \cite{KS2017} proposed to choose the bandwidth $h$ by minimizing the trace of AMSE of the smoothed estimator. In this case, it is
\#
 h^* = \argmin_{h > 0 }  {\rm tr}  \big\{   0.25 \kappa_2^2   n h^4  \EE(\bb) \EE(\bb^\T) -  2 \bar{\kappa}_1 h \Hb   \big\} = \Bigg\{ \frac{2 \bar \kappa_1 \EE \{ f_{ \varepsilon|\bx } (0) \bw^\T \bw  \}  }{ \kappa_2^2 n \cdot  \EE(\bb^\T ) \EE(\bb)  } \Bigg\}^{1/3}.\label{AMSE.bandwidth1}
\#
Since $\bx = (1,\bx_-^\T)^\T$ has 1 as its first component, we have $\bSigma^{-1}\EE(\bx) = (1, \textbf{0}_{p-1}^\T)^\T$. As a result, in the special case where $\varepsilon_i$ and $\bx_i$ are independent, the above AMSE-optimal $h^*$ can be reduced to
\#
	h^\star = \Bigg[ \frac{2 \bar \kappa_1  f_{ \varepsilon|\bx } (0)   }{ \{\kappa_2  f'_{ \varepsilon|\bx } (0) \}^2  }  \frac{p}{n}\Bigg]^{1/3}. \label{AMSE.bandwidth2}
\#
When $K(\cdot)$ is the Gaussian kernel,   $\kappa_2 = \EE(g^2) =1 $ and $\bar \kappa_1 = \EE\{ g \Phi(g)\}$, where $g\sim \cN(0,1)$ and $\Phi(\cdot)$ is the standard normal CDF.  The Monte Carlo method computes $\bar \kappa_1 \approx 0.282$.  When $K(\cdot)$ is the logistic kernel (see Remark~\ref{rmk1}),  we have $\kappa_2 = \pi^2/3$ and $\bar \kappa_1 =1/2$.
Although both $h^*$ and $h^\star$ depend on $ f_{ \varepsilon|\bx } (0)$ and $ f'_{ \varepsilon|\bx } (0)$, which can only be estimated by nonparametric estimators with a rule-of-thumb bandwidth, they provide a benchmark bandwidth choice in simulation studies.
}
\end{remark}

\begin{remark}[Large-$p$ asymptotics] \label{rmk2}
A broader view of classical asymptotics recognizes that the parametric dimension of appropriate model sequences may tend to infinity with the sample size; that is $p=p_n \to \infty$ as $n\to \infty$. Results with increasing $p$ are available in \cite{W1989}, \cite{HS2000} and \cite{BCCF2019} when $p= o(n)$, and  in \cite{BC2011},  \cite{WWL2012} and \cite{KCHP2017} for regularized quantile regression when $p\gg n$.
{
In the large-$p$ and larger-$n$ setting---``$p\to \infty$ and $p/n \to 0$", \cite{W1989} shows that $p^3(\log n)^2 / n \to 0$ suffices for a normal approximation.  This growth condition remains the best known one although under weaker assumptions on the (fixed) design \citep{HS2000,BCCF2019}.  To our knowledge, the weakest fixed design assumption is $\max_{1\leq i\leq n} \| \bx_i \|_2^2 = \cO(p)$.
}

For smooth robust regression estimators, asymptotic normality can be proven under less restrictive conditions on $p$.  \cite{H1973} showed that if the loss is twice differentiable, the asymptotic normality for $\langle \ba, \hat \bbeta \rangle$, where $\ba \in \RR^p$, holds if $p^3/n \to 0$ as $n$ increases. \cite{P1985} and \cite{M1989} weakened this condition to $(p \log n)^{3/2} /n \to 0$ and $p^{3/2} \log(n)/ n \to 0$, respectively, when the loss function is four times differentiable.
For Huber loss that has a Lipschitz continuous derivative, \cite{HS2000} obtained the scaling $p^2 \log p = o(n)$ that ensures the asymptotic normality of arbitrary linear combinations of $\hat \bbeta$. Table 1 summarizes our discussion here and shows that the smoothing for conquer helps ensure asymptotic normality of the estimator under weaker conditions on $p$ than what we need for the usual quantile regression estimator.
\end{remark}

\begin{table}[!htp] \label{table1}
\footnotesize
{\scriptsize
\begin{center}
\caption{Summary of scaling conditions required for normal approximation under various loss functions.}
\begin{tabular}{ l | c   | c}
  \hline
Loss function&Design &Scaling condition\\
\hline\hline
Huber loss \citep{H1973} &  Fixed design  &$p^3 =o(n)$  \\ \hline
Four times differentiable loss \citep{P1985}  &  Fixed design (with symmetric error) & $( p \log n )^{3/2}= o(n)$ \\ \hline
Four times differentiable loss \citep{M1989}  &  Fixed design & $p^{3/2}  \log n = o(n)$ \\ \hline
Huber loss \citep{HS2000}   &  Fixed design  &$p^2 \log p =o(n)$  \\ \hline
Huber loss \citep{CZ2019}   &  Sub-Gaussian   &$p^2 =o(n)$  \\ \hline
Quantile loss \citep{W1989, HS2000}   &  Fixed design  & $p^3 (\log n)^2 = o(n)$ \\ \hline
Quantile loss \citep{PZ2019}   &  Sub-Gaussian  & $p^3 (\log n)^2 = o(n)$ \\ \hline
Convolution smoothed quantile loss  &  Sub-Gaussian & $p^{8/3} = o(n)$  \\  \hline
%Three Times Continuously Differentiable Loss \citep{P1986}   &  General Random Design &$p\{\log (p)\}^2 =o(n)$  \\ \hline
\end{tabular}
\label{tablegaussian}
\end{center}
}
\end{table}

\begin{remark}
In this paper, we show that the accuracy of conquer-based inference via the Bahadur representation (and normal approximations) has an error of rate faster than $n^{-1/4}$ yet slower than $n^{-1/2}$; see Theorems~\ref{thm:bahadur} and \ref{thm:clt}. For standard regression quantiles, \cite{P2012} proposed an alternative expansion for the quantile process using the ``Hungarian" construction of Koml\'os, Major and Tusn\'ady. This stochastic approximation yields an error of order $n^{-1/2}$ (up to a factor of $\log n$), and hence provides a theoretical justification for accurate approximations for inference in regression quantile models. 
\end{remark}

%%%%%%%%%%%%%%%%%%%%%%%%%%%%%%%%
%%%%%%%%%%%%%%%%%%%%%%%%%%%%%%%%
% Theoretical Guarantees for Bootstrap Inference
%%%%%%%%%%%%%%%%%%%%%%%%%%%%%%%%
%%%%%%%%%%%%%%%%%%%%%%%%%%%%%%%%
 \subsection{Theoretical guarantees for inference}
 \label{sec:boot:theory}
 
We next investigate the statistical properties of the bootstrapped estimator defined in \eqref{bsqr.1}, with a particular focus on the Rademacher multiplier bootstrap (RMB). To be specific, we use, in this section and the rest of the paper, the random weights $w_i=1+e_i$ for $i=1,\ldots, n$, where  $e_1,\ldots, e_n$ are independent Rademacher random variables, that is, symmetric sign variables with $\PP(e_i=1) = \PP(e_i=-1)=1/2$.
As before, we consider array (non)asymptotics, and the obtained bootstrap approximation errors depend explicitly on $(n,p)$ and $h$.

\begin{theorem} \label{thm:boot.concentration}
Assume Conditions~\ref{cond.kernel}, \ref{cond.reg} and \ref{cond.predictor2} hold. For any given $t\geq 0$,  let the sample size and bandwidth satisfy $\underbar{$f$} ^{-1} m_3^{1/2}\upsilon_1 \sqrt{(p+t)/n }   \lesssim h \lesssim \underbar{$f$} m_3^{-1/2}$.
Then, there exists some ``good" event $\cE(t)$ with $\PP\{ \cE(t) \} \geq 1-3e^{-t}$ such that, 
with $\PP^*$-probability at least $1-2e^{-t}$ conditioned on $\cE(t)$,
\#
	 \| \hat \bbeta^\flat_h  - \bbeta^* \|_{\bSigma} \le   \frac{C^\flat}{\underbar{$f$}} \Biggl\{  \upsilon_1 \sqrt{  \frac{\log_2(1/h) + p+t}{n}} + l_0 \kappa_2 h^2 \Biggr\} ,
\#
where $C^\flat >0$ is a absolute constant.
\end{theorem}

 Analogously to Theorem~\ref{thm:bahadur}, we further provide a Bahadur representation result for the bootstrap estimator $\hat \bbeta^\flat_h $, which paves the way for validating the conquer-RMB method.

\begin{theorem} \label{thm:boot.bahadur}
In addition to Conditions~\ref{cond.kernel}, \ref{cond.reg} and \ref{cond.predictor2}, assume $\sup_{u\in \RR} f_{\varepsilon  | \bx} (u) \leq \bar f$ almost surely (in $\bx$) and $K(\cdot)$ is $l_K$-Lipschitz continuous.  Suppose the sample size satisfies $n\gtrsim q:= p + \log n$, and set the bandwidth as $h \asymp  (q/n)^{2/5}$. Then, there exists a sequence of events $\{ \cF_n\}$ with $\PP(\cF_n) \geq 1- 6 n^{-1}$ such that, with $\PP^*$-probability at least $1-3n^{-1}$ conditioned on $\cF_n$,
\#
     \Bigg\|  \bSigma^{-1/2} \Jb_h  ( \hat \bbeta^\flat_h - \hat \bbeta_h )  - &  \frac{1}{n}   \sn e_i    \bigl\{  \tau -    \cK_h(- \varepsilon_i  )  \bigr\}   \bSigma^{-1/2} \bx_i \Bigg\|_2 \nn \\
   &  \lesssim  \bigg( \frac{q}{n} \bigg)^{4/5}  \bigvee \bigg( \frac{q}{n} \bigg)^{3/5}  \biggl( \frac{p \log n}{n} \bigg)^{1/4} \bigvee \bigg( \frac{q}{n} \bigg)^{3/5} \frac{ p\log n }{n^{1/2}} . \label{boot.br}
\#
\end{theorem}

As suggested by Theorem~\ref{thm:clt} and the discussion below, if we set the order of the bandwidth $h$ as $\{ (p+\log n)/n \}^{2/5}$, the normal approximation to the conquer estimator is asymptotically accurate provided that $p^{8/3} = o(n)$ as $n\to \infty$. For the same $h$, the right-hand side of \eqref{boot.br} is of order $o(n^{-1/2})$ provided that $p^{8/3}(\log n)^{5/3}= o(n)$. Putting these two parts together, we have the following asymptotic bootstrap approximation result. 

\begin{corollary}
Assume the same conditions of Theorem~\ref{thm:boot.bahadur}, and let the bandwidth be of order $h\asymp \{ (p+\log n)/n\}^{2/5}$. If the dimension $p=p_n$ is subject to  $p(\log n)^{5/8} = o(n^{3/8})$, then for any deterministic vector $\ba \in \RR^p$, 
\#
	\sup_{x\in \RR } \, \bigl| \PP\bigl(  n^{1/2} \langle \ba,    \hat \bbeta_h  - \bbeta^*   \rangle   \leq x \bigr)  - \PP^*\bigl( n^{1/2} \langle \ba,       \hat \bbeta^\flat_h - \hat \bbeta_h  \rangle  \leq x \bigr)  \bigr| \stackrel{\PP }{\rightarrow } 0 \mbox{~ as } n \to \infty.  \label{boot.approxi}
\#
\end{corollary}

The proof of \eqref{boot.approxi} follows the same argument as that in the proof of Theorem~\ref{thm:clt}, and therefore is omitted. The additional logarithmic factor in the scaling may be an artifact of the proof technique. For standard quantile regression, \cite{FHH2011} established a fixed-$p$ asymptotic bootstrap approximation result for wild bootstrap under fixed design.

\begin{remark}(Multiplier bootstrap with more general weighting schemes)  \label{rmk:weights}
{
By examining the proof of Theorem~\ref{thm:boot.bahadur}, we see that the assumption  $\EE(e_i^2)=1$ is not necessarily required for the bound \eqref{boot.br} on bootstrap Bahadur linearization error.  To retain the convexity of the bootstrap loss $\hat Q_h^\flat(\bbeta)= (1/n)\sn (1+e_i) \ell_h(y_i-\bx_i^\T \bbeta)$, we restrict our attention to non-negative multipliers $1+e_i\geq 0$. More generally,  assume that $e_1,\ldots, e_n$ are i.i.d. satisfying 
\#
	\EE(e_i)=0, \quad e_i \geq -1 ~\mbox{ and }~ \log \EE e^{\lambda e_i} \leq \lambda^2 \nu/2 ~\mbox{ for all } \lambda\geq 0 \mbox{ and some } \nu>0.  \label{general.weights}
\# 
This means that $e_i$ has sub-Gaussian right tails.  Typical examples satisfying \eqref{general.weights} include: (i) uniform distribution on $[-1,1]$, (ii) symmetric triangular distribution on $[-1,1]$, (iii) shifted folded normal distribution $(\pi/2)^{1/2}|g| -1$ where $g\sim \cN(0,1)$. The proof of the bound \eqref{boot.br} under such a general scheme requires more involved argument; see, for example, the proof of Theorem~2.3 in  \cite{CZ2019} (the unit variance assumption therein can also be relaxed). When $\kappa^2 := \EE(e_i^2) \neq 1$,   although the bootstrap approximation  result \eqref{boot.approxi} will no longer hold,   by a simple variance adjustment it can be shown that 
 $$ 
 \sup_{x\in \RR } \, \bigl| \PP\bigl(  n^{1/2} \langle \ba,    \hat \bbeta_h  - \bbeta^*   \rangle   \leq x \bigr)  - \PP^*\bigl\{  (n/\kappa)^{1/2} \langle \ba,       \hat \bbeta^\flat_h - \hat \bbeta_h  \rangle  \leq x \bigr\}  \bigr| \stackrel{\PP }{\rightarrow } 0 \mbox{~ as } n \to \infty. 
 $$
The pivotal bootstrap confidence intervals can thus be constructed by slightly adapting the method described in Section~\ref{subsec:mbi}.
}

{
Regarding computational complexity, for each bootstrap sample, nonparametric bootstrap and multiplier bootstrap with general  weights refit a quantile regression on the whole sample,  while RMB only solves a conquer optimization on a subsample of size $n/2$ approximately.  This feature makes RMB preferable particularly for large-scale datasets.
}
\end{remark}

{
The numerical performance of the Rademacher multiplier bootstrap inference for conquer will be examined in Section~\ref{subsec:bootstrap inference}. The main advantage of the multiplier bootstrap method is that it does not require estimating the variance-covariance matrices in \eqref{def:var-cov}, which can be quite unstable and thus causes outliers when $\tau$ is close to 0 or 1; see Figure~\ref{norm-boot} for a numerical comparison of the (multiplier) bootstrap percentile method and the normal-based method. When $\tau$ is reasonably bounded away from 0 and 1, say between 0.25 and 0.75, normal calibration also performs well empirically, and is computationally attractive because the estimator is only computed once. 
}

{
The construction of normal-based confidence intervals is based on the estimated variances 
$\hat \sigma^2_{h} (\ba) = \ba^\T \hat \Jb_h^{-1} \hat \Vb_h \hat \Jb_h^{-1} \ba$ for $\ba \in \RR^p$, where $\hat \Jb_h$ and $\hat \Vb_h$ are given in \eqref{def:var-cov}. In view of Theorem~\ref{thm:clt}, the validity of normal calibration relies on the consistency of $\hat \Jb_h$ and $\hat \Vb_h$. In the following, we provide the consistency of  $\hat \Jb_h$ and $\hat \Vb_h$ under the operator norm, again in the regime ``$p/n\to 0$ as $p, n \to \infty$".
}
{
Note that both $\hat \Jb_h$ and $\hat \Vb_h$ depend on the conquer estimator, whose rate of convergence is already established in Theorem~\ref{thm:concentration}. For $\bdelta \in \RR^p$, define matrix-valued functions
\#
 \hat \Jb_h(\bdelta) = \frac{1}{n} \sn K_h( \varepsilon_i -\langle \bx_i, \bdelta \rangle ) \bx_i \bx_i^\T ~~\mbox{ and }~~  \hat \Vb_h(\bdelta) = \frac{1}{n} \sn  \{ \cK_h( \langle \bx_i, \bdelta \rangle  - \varepsilon) - \tau \}^2 \bx_i \bx_i^\T ,
\#
so that $\hat \Jb_h = \hat \Jb_h(\hat \bdelta)$ and $\hat \Vb_h = \hat  \Vb_h(\hat \bdelta)$ with $\hat \bdelta = \hat \bbeta_h - \bbeta^*$. Conditioned on the event $\{ \|\hat \bdelta \|_{\bSigma}\leq r\}$ for some prespecified $r>0$ which determines the convergence rate of $\hat \bbeta_h$, we have
\$
 \| \hat \Jb_h -\Jb_h  \|_2 \leq \sup_{\| \bdelta \|_{\bSigma} \leq r } \|  \hat \Jb_h(\bdelta)  - \Jb_h   \|_2 ~~ \mbox{ and }~~ \| \hat \Vb_h -  \Vb_h \|_2 \leq \sup_{\| \bdelta \|_{\bSigma} \leq r } \|  \hat \Vb_h(\bdelta)  - \Vb_h  \|_2 ,
\$
where $\Vb_h := \EE\big[ \{ \cK_h(-\varepsilon) - \tau\}^2 \bx \bx^\T \big]$.
The problem is thus reduced to controlling the above suprema over a local neighborhood.
\begin{proposition}  \label{thm:var-cov}
In addition to the conditions in Theorem~\ref{thm:bahadur}, assume that the kernel $K(\cdot)$ is $l_K$-Lipschitz continuous. For any given $r\geq 0$, 
\#
\sup_{\| \bdelta \|_{\bSigma} \leq r } \|  \bSigma^{-1/2} \{ \hat \Jb_h(\bdelta)  -  \Jb_h \} \bSigma^{-1/2}  \|_2 \lesssim \sqrt{\frac{p \log n+ t}{n h }}  +    r   \label{J.consistency}
\#
holds with probability at least $1-e^{-t}$, provided that $\max\{ \sqrt{(p+t)/n} , p\log(n)/n \} \lesssim h\lesssim 1$.
The same probabilistic bound also applies to $\sup_{\| \bdelta \|_{\bSigma} \leq r } \|  \bSigma^{-1/2}   \{  \hat \Vb_h(\bdelta)  - \Vb_h \} \bSigma^{-1/2} \|_2$.
\end{proposition}
}

{
Following the discussions below Theorem~\ref{thm:clt}, if we set the bandwidth as $h \asymp \{ (p+\log n)/n \}^{2/5}$, $\| \hat \bbeta_h - \bbeta^* \|_{\bSigma} = \cO_{\PP}( \sqrt{(p+\log n)/n})$ and 
$n^{1/2} \ba^\T (\hat \bbeta_h - \bbeta^*) / \sigma_h(\ba)  \rightarrow \cN(0,1)$ in distribution uniformly over $\ba \in \RR^p$ as $n\to \infty$ under the constraint $p^{8/3} = o(n)$. With the same bandwidth, it follows from Proposition~\ref{thm:var-cov} that
$$
 \max\big( \| \hat \Jb_h -\Jb_h  \|_2  ,  \| \hat \Vb_h   - \Vb_h \|_2 \big)  = \cO_{\PP}  \Big[    \big\{(\log n)^{1/2} p^{3/10} + (\log n)^{3/10} p^{1/2} \big\} n^{-3/10} \Big] = o_{\PP}(1) .
$$ 
This ensures the consistency of variance estimators, that is, $|\hat \sigma_h^2(\ba) /\sigma_h^2(\ba) - 1 | \xrightarrow {\PP} 0$.
}

%%%%%%%%%%%%% %%%%%%%%%%%%% %%%%%%%%%%%%%
% Numerical Studies
%%%%%%%%%%%%% %%%%%%%%%%%%% %%%%%%%%%%%%% 
\section{Numerical studies}
\label{sec:numerical}
In this section, we assess the finite-sample performance of  conquer via extensive numerical studies.  We compare conquer to standard QR \citep{KB1978} and Horowitz's smoothed QR \citep{H1998}.
Both the convolution-type and Horowitz's smoothed methods involve a smoothing parameter $h$.  
In view of Theorem~\ref{thm:clt}, we take $h = \{ (p+\log n) /n\}^{2/5}$ in all of the numerical experiments. 
As we will see from Figure~\ref{sensitivity}, the proposed method is insensitive to the choice of $h$. Therefore, we leave the fine tuning of $h$ as an optional rather than imperative choice. In all the numerical experiments, the convergence criterion in Algorithms~\ref{algo:sqr} and \ref{algo:huber} is taken as $\delta=10^{-4}$.

We first generate the covariates $\bx_i = (x_{i , 1}, \dots, x_{i , p})^\T$ from a multivariate uniform distribution on the  cube $3^{1/2} \cdot [-1, 1]^p$ with covariance matrix $\bSigma=(0.7^{|j - k|})_{1\leq j,k\leq p}$ using the \texttt{R} package \texttt{MultiRNG} \citep{F1999}.   
The random noise $\varepsilon_i$ is generated from two different distributions: (i) Gaussian distribution, $\mathcal{N}(0, 4)$; and (ii) $t$ distribution with two degrees of freedom, $ t_2$.  Let $\bbeta^* = (1, \dots, 1)^\T_{p}$, and $\beta^*_0 = 1$. Given $\tau \in (0, 1)$, we then generate the response $y_i$ from the following homogeneous and heterogeneous models, all of which satisfy the Assumption~\eqref{eq:lqr}:
\begin{enumerate}
\item Homogeneous model: 
\begin{align}
y_i = \beta^*_0 + \langle   \bx_i , \bbeta^* \rangle + \{ \varepsilon_i - F^{-1}_{\varepsilon_i}(\tau) \}, ~~~~i = 1, \dots, n;  \label{model.homo}
\end{align}
\item Linear heterogeneous model: 
\begin{align}
y_i = \beta^*_0 + \langle   \bx_i , \bbeta^* \rangle + (0.5 x_{i, p} + 1) \{ \varepsilon_i - F^{-1}_{\varepsilon_i}(\tau) \} , ~~~~i = 1, \dots, n; \label{model.linear}
\end{align}
\item Quadratic heterogeneous model: 
\begin{align}
y_i = \beta^*_0 + \langle   \bx_i , \bbeta^* \rangle + 0.5\{1+ ( x_{i, p}-1)^2 \} \{ \varepsilon_i - F^{-1}_{\varepsilon_i}(\tau)\}, ~~~~i = 1, \dots, n.   \label{model.quad}
\end{align}
\end{enumerate}

To evaluate the performance of different methods, we calculate the estimation error under the $\ell_2$-norm, i.e., $\|\hat{\bbeta}-\bbeta^*\|_2$, and record the elapsed time.  The details are in Section~\ref{subsec:estimation error}.
In Section~\ref{subsec:bootstrap inference}, we examine the finite-sample performance of the multiplier bootstrap method for constructing confidence intervals  in terms of coverage probability, width of the interval, and computing time.

%%%%%%%%%%%%%%%%%%%%%%%%%%%%%%%%%%%
%%%%%%%%%%%%%%%%%%%%%%%%%%%%%%%%%%%
% Simulation for evaluating estimation error
%%%%%%%%%%%%%%%%%%%%%%%%%%%%%%%%%%%
%%%%%%%%%%%%%%%%%%%%%%%%%%%%%%%%%%%
\subsection{Estimation}
\label{subsec:estimation error}
%First, we evaluate the estimation accuracy of methods, i.e., $\|\hat{\bbeta}-\bbeta^*\|_2$. 
For all the numerical studies in this section, we consider a wide range of the sample size $n$, with the size-dimension ratio fixed at $n / p = 20$. That is, we allow the dimension $p$ to increase as a function of $n$. 
We implement conquer with four different kernel functions as described in Remark~\ref{rmk1}: (i) Gaussian; (ii) uniform; (iii) Epanechnikov; and (iv) triangular. 
The classical quantile regression is implemented via a modified version of the Barrodale and Roberts algorithm \citep{KO1987, KO1994} by setting \texttt{method=} ``\texttt{br}" in the \texttt{R} package \texttt{quantreg}, which is recommended for problems with up to several thousands of observations in \citet{K2019}. For very large problems, the Frisch-Newton approach after preprocessing ``pfn" is preferred.
Since the same size taken to be at most 5000 throughout this section, the two methods, ``br'' and ``pfn", have nearly identical runtime behaviors. 
In some applications where there are a lot of discrete covariates,  it is advantageous to use method ``sfn", a sparse version of Frisch-Newton algorithm that exploits sparse algebra to compute iterates \citep{KN2003}.
Moreover, we implement Horowitz's smoothed quantile regression using the Gaussian kernel, and solve the resulting non-convex optimization via gradient descent with random initialization and stepsize calibrated by backtracking line search (Section 9.3 of \citealp{BV2004}). 
The results, averaged over $500$ replications, are reported.

\begin{figure}[!ht]
  \centering
  \subfigure[Model \eqref{model.homo} with $\mathcal{N}(0, 4)$ error.]{\includegraphics[width=0.32\textwidth]{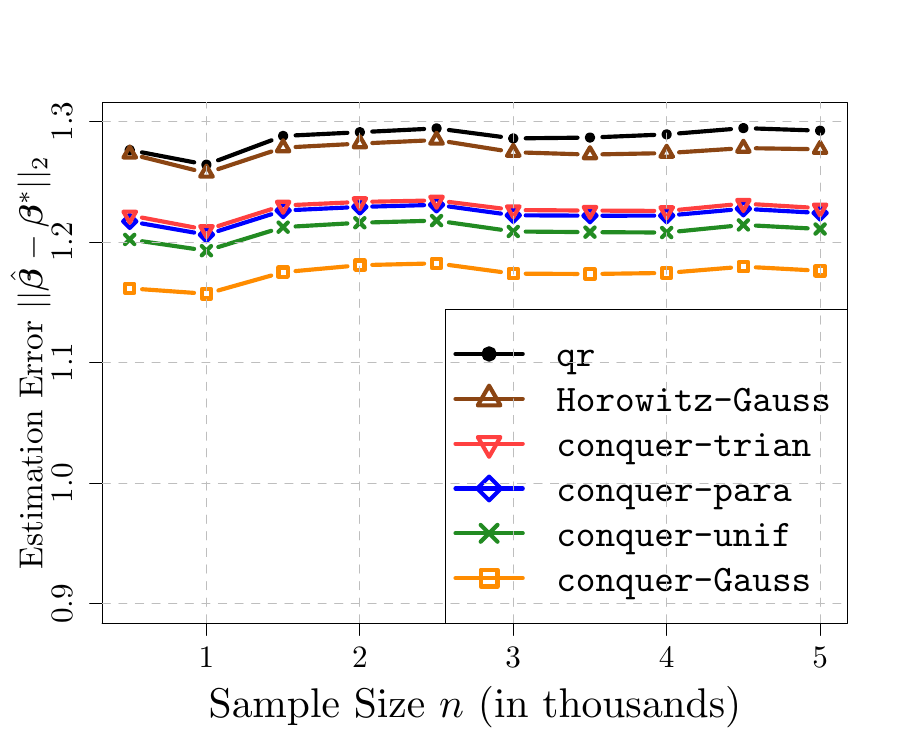}} 
  \subfigure[Model \eqref{model.linear} with $\mathcal{N}(0, 4)$ error.]{\includegraphics[width=0.32\textwidth]{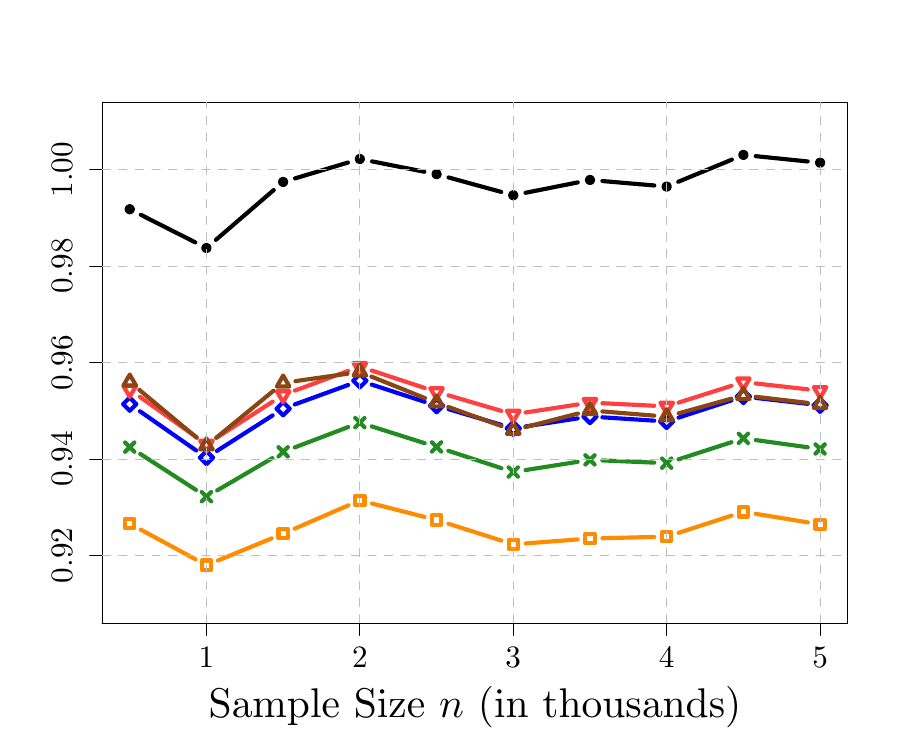}} 
  \subfigure[Model \eqref{model.quad} with $\mathcal{N}(0, 4)$ error.]{\includegraphics[width=0.32\textwidth]{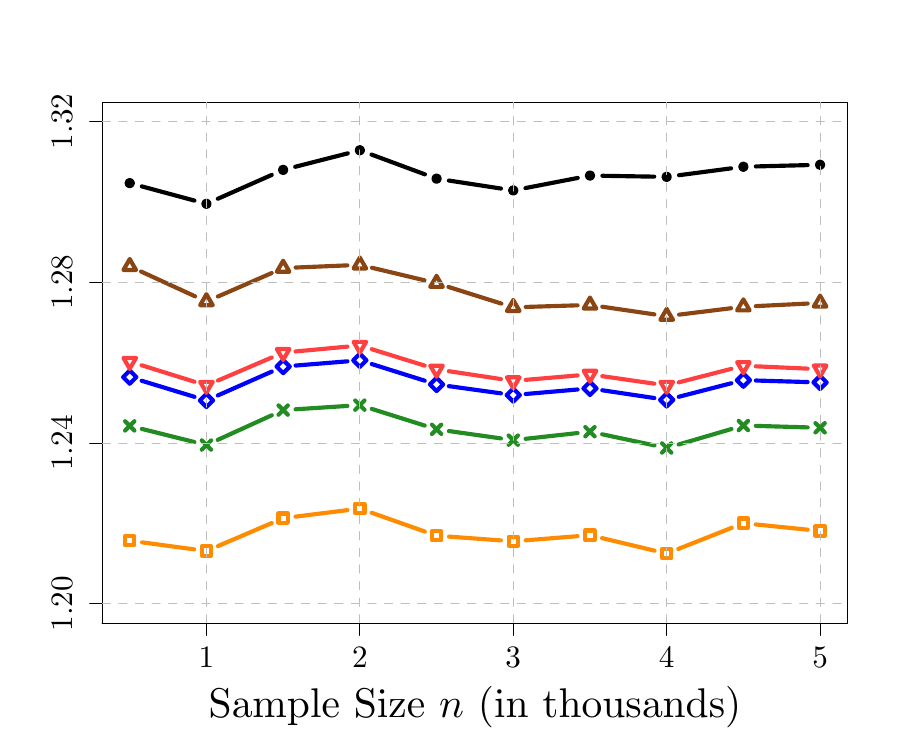}}
  \subfigure[Model \eqref{model.homo} with $t_2$ error.]{\includegraphics[width=0.32\textwidth]{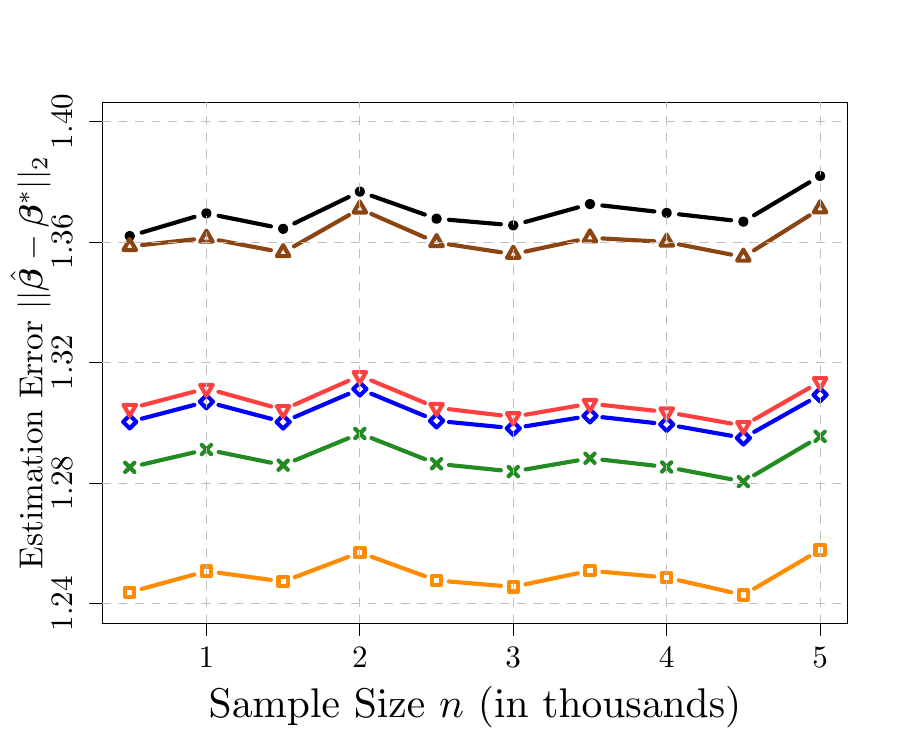}}
  \subfigure[Model \eqref{model.linear} with $t_2$ error.]{\includegraphics[width=0.32\textwidth]{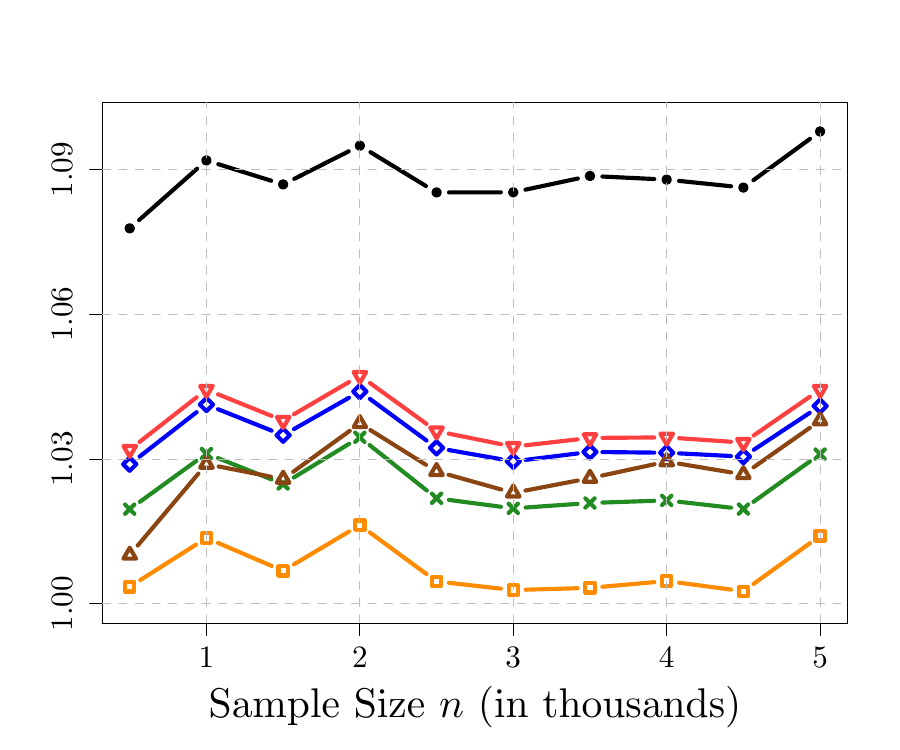}}
  \subfigure[Model \eqref{model.quad} with $t_2$ error.]{\includegraphics[width=0.32\textwidth]{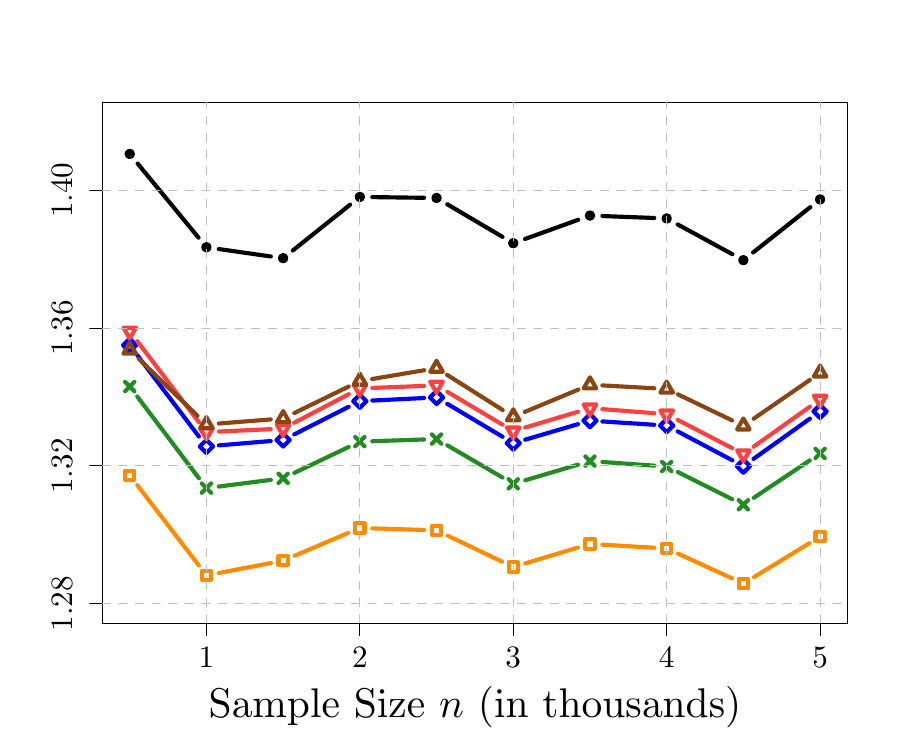}}
\caption{Estimation error under models \eqref{model.homo}--\eqref{model.quad} in Section~\ref{sec:numerical} with $\mathcal{N}(0, 4)$ and $t_2$ errors, $\tau = 0.9$, averaged over 500 data sets for three different methods: (i) quantile regression \texttt{qr}, (ii) Horowitz's method with Gaussian kernel \texttt{Horowitz-Gauss}, and (iii) the conquer method with four different kernel functions \texttt{conquer-trian}, \texttt{conquer-para}, \texttt{conquer-unif}, and \texttt{conquer-Gauss}.}
  \label{est.9}
\end{figure}

Figure~\ref{est.9} depicts estimation error of the different methods under the simulation settings described in Section~\ref{sec:numerical} with $\tau = 0.9$.  
We see  that conquer has a lower estimation error than the classical QR across all scenarios, indicating that smoothing can improve estimation accuracy under the finite-sample setting.  
Moreover, compared to Horowitz's smoothing, conquer has a lower estimation error in most settings.     
Estimation error under various quantile levels $\tau \in \{0.1, 0.3, 0.5, 0.7\}$ with the $\mathcal{N}(0, 4)$ and $t_2$ random noise are also examined.  The results are reported in Figure~\ref{est.1357.normal} and \ref{est.1357} in Appendix~\ref{sec:app.simu}, from which we observe evident advantages of conquer, especially at low and high quantile levels.

%%%%%%%%%%%%%%%%%%%%%%
%%%%%%%%%%%%%%%%%%%%%%
% Computational Aspect
%%%%%%%%%%%%%%%%%%%%%%
%%%%%%%%%%%%%%%%%%%%%%
To assess the computational efficiency, we compute the elapsed time for fitting the different methods.
Figures~\ref{time.9}, \ref{time.1357.normal} and \ref{time.1357} in Appendix~\ref{sec:app.simu} report the runtime for the different methods with growing sample size and dimension under the same settings as in Figures~\ref{est.9}, \ref{est.1357.normal} and \ref{est.1357}, respectively.
We observe that conquer is computationally efficient and stable across all scenarios, and the runtime is insensitive to the choice of kernel functions.  
In contrast, the runtime for classical quantile regression grows rapidly as the sample size and dimension increase.
 Figures~\ref{time.9}, \ref{time.1357.normal} and \ref{time.1357} in Appendix~\ref{sec:app.simu} show that the runtime of Horowitz's smoothing method increases significantly at extreme quantile levels $\tau \in \{0.1,0.9\}$, possibly due to the combination of its non-convex nature and flatter gradient.
{As suggested by a referee, another set of simulations with covariates $\bx_i$'s following multivariate normal disttribution $\mathcal{N}_p(\bm{0}, \bSigma)$ are also conducted, where $\bSigma=(0.7^{|j - k|})_{1\leq j,k\leq p}$. 
When the covariates are unbounded, the assumption \eqref{eq:lqr} is violated under the model \eqref{model.linear}, thus we exclude this setting from our experiments.
The corresponding results of estimation and running time under $t_2$ noise are presented in Figures~\ref{est.unbd}--\ref{time.unbd}, averaged over 500 datasets, and similar performance can be observed.}
In summary, we conclude that conquer significantly improves computational efficiency while retaining high statistical accuracy for fitting large-scale linear quantile regression models.

\begin{figure}[!ht]
  \centering
  \subfigure[Model \eqref{model.homo} with $\mathcal{N}(0, 4)$ error.]{\includegraphics[width=0.32\textwidth]{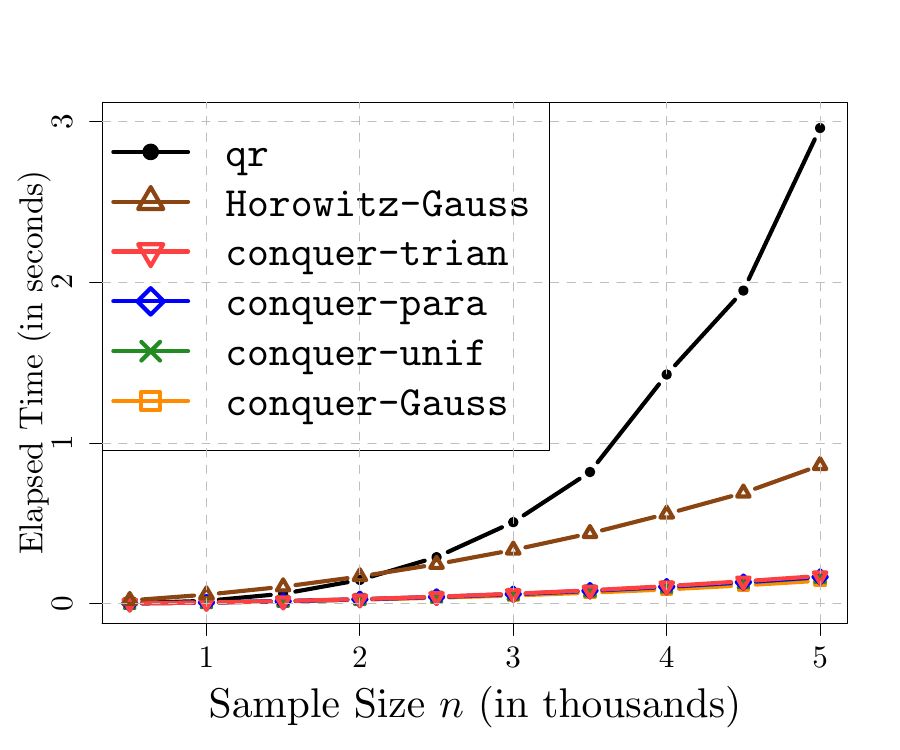}} 
  \subfigure[Model \eqref{model.linear} with $\mathcal{N}(0, 4)$ error.]{\includegraphics[width=0.32\textwidth]{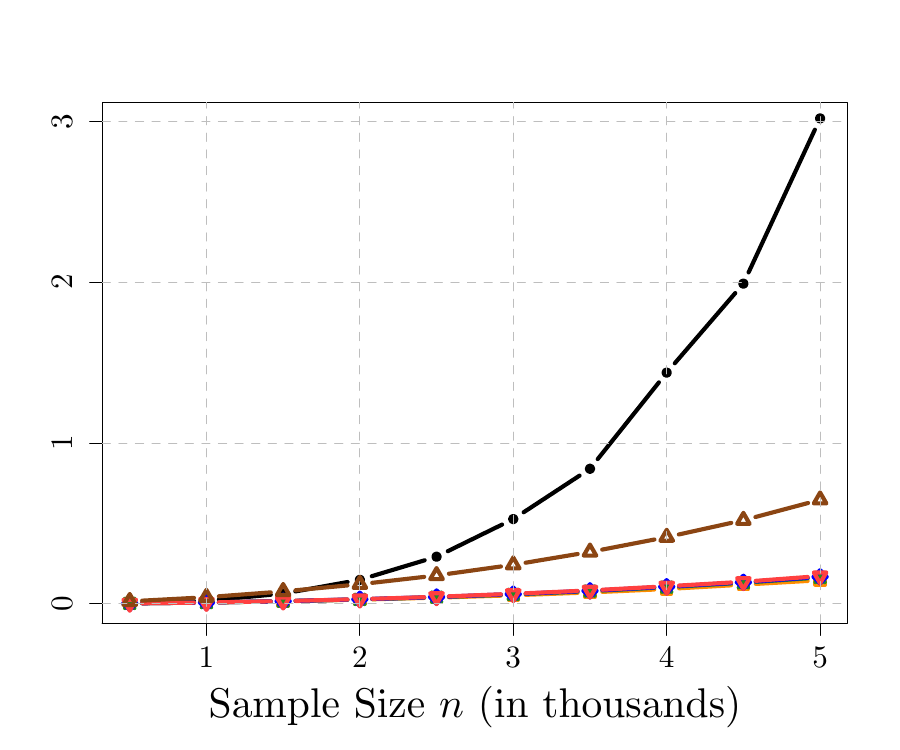}} 
  \subfigure[Model \eqref{model.quad} with $\mathcal{N}(0, 4)$ error.]{\includegraphics[width=0.32\textwidth]{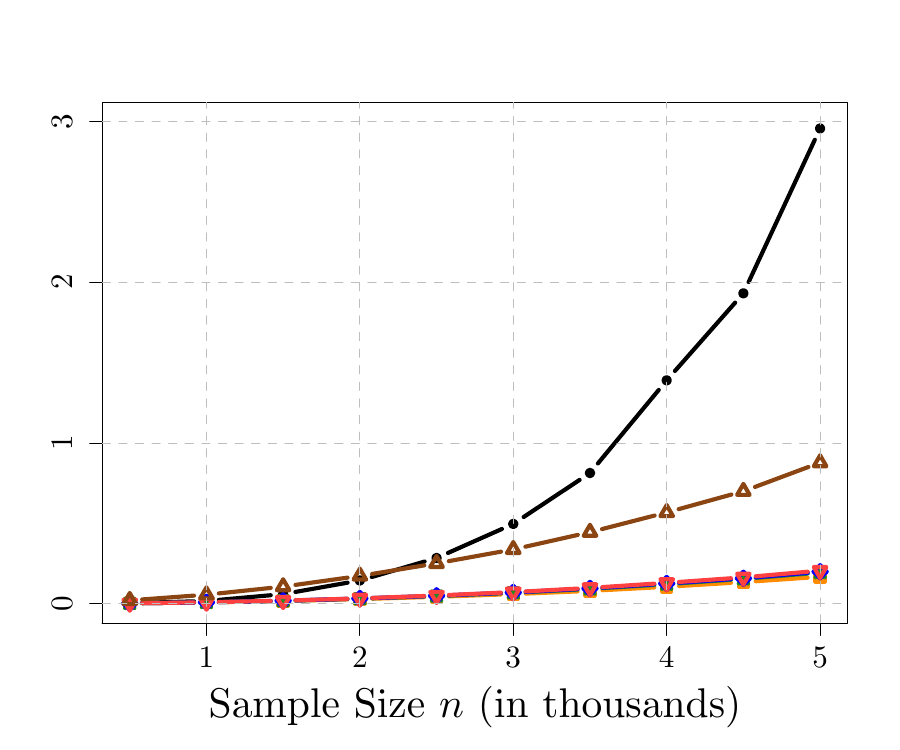}}
  \subfigure[Model \eqref{model.homo} with $t_2$ error.]{\includegraphics[width=0.32\textwidth]{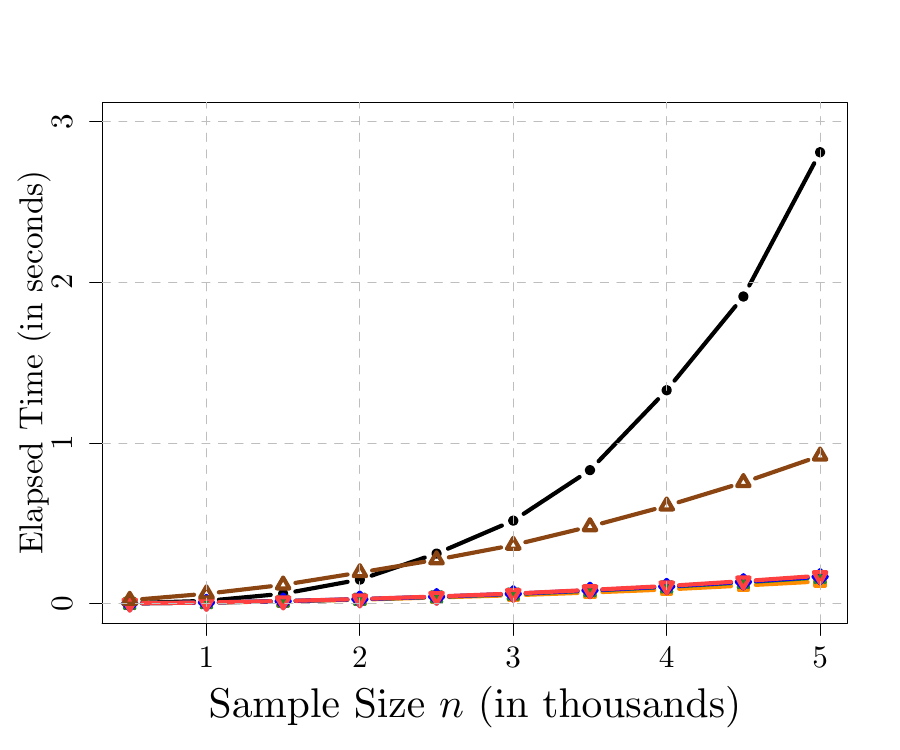}}
  \subfigure[Model \eqref{model.linear} with $t_2$ error.]{\includegraphics[width=0.32\textwidth]{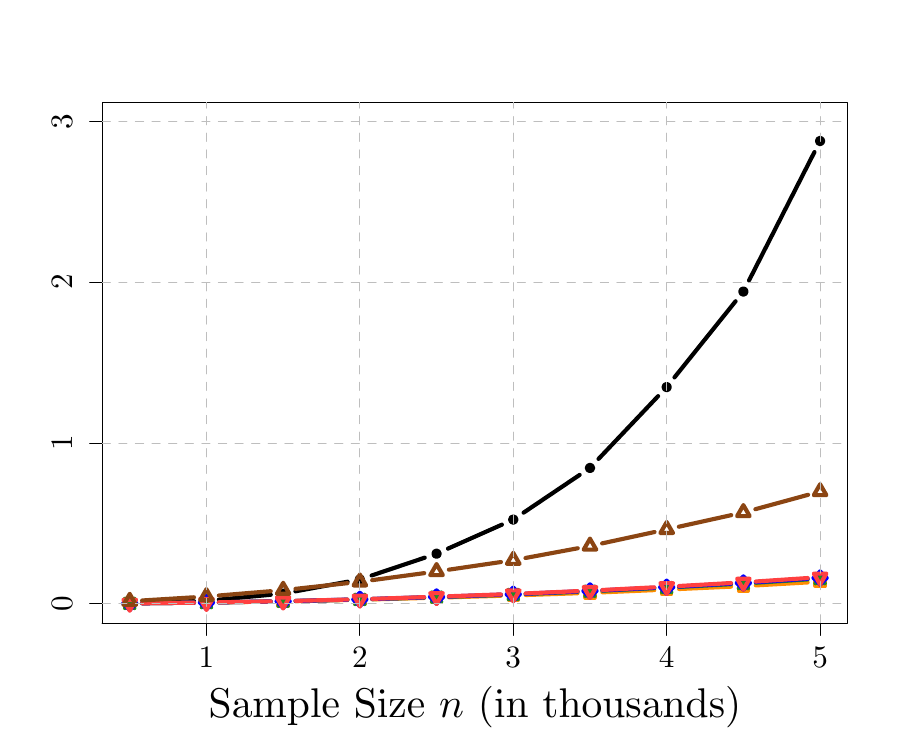}}
  \subfigure[Model \eqref{model.quad} with $t_2$ error.]{\includegraphics[width=0.32\textwidth]{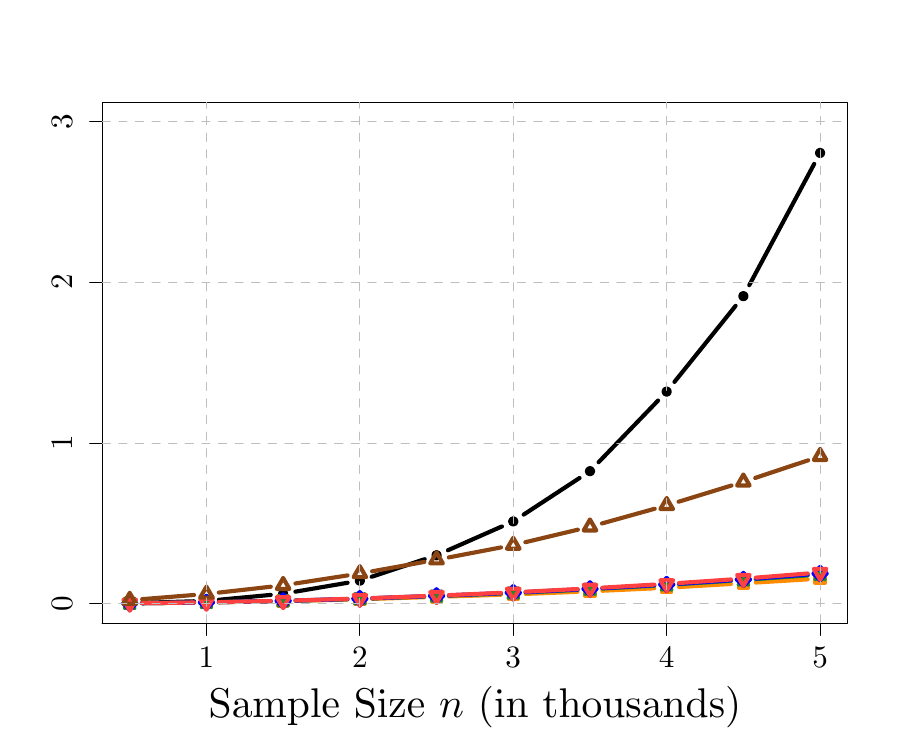}}
\caption{Elapsed time of standard QR, Horowitz's smoothing, and conquer when $\tau = 0.9$. The model settings are the same as those in Figure~\ref{est.9}.}
  \label{time.9}
\end{figure}

\begin{figure}[!ht]
  \centering
  \subfigure[Model \eqref{model.homo} with $\cN(0, 4)$ error.]{\includegraphics[width=0.32\textwidth]{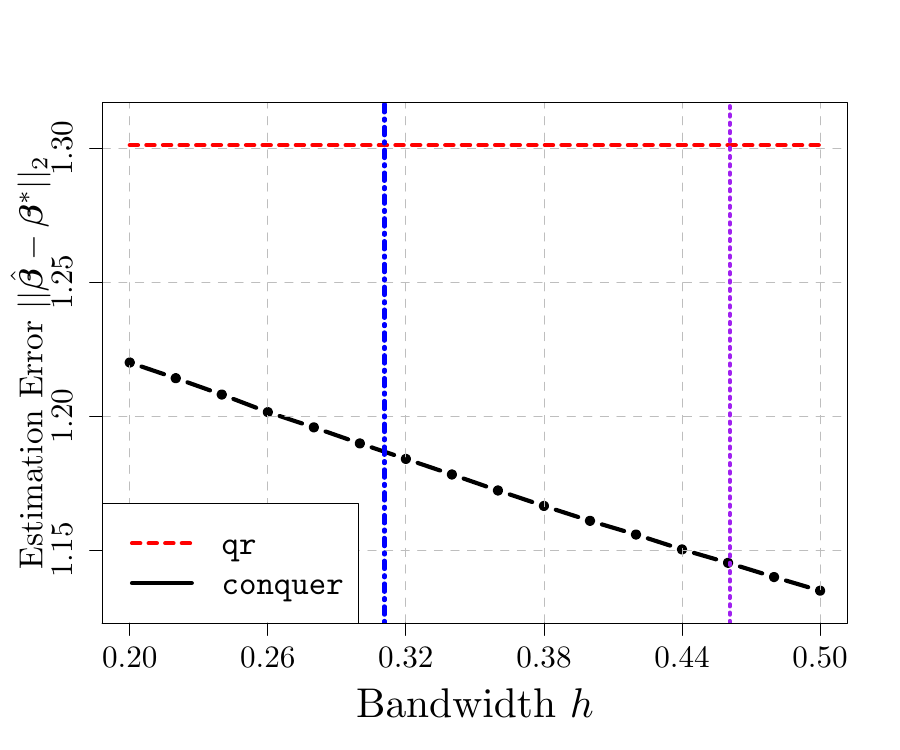}} 
  \subfigure[Model \eqref{model.linear} with $\cN(0, 4)$ error.]{\includegraphics[width=0.32\textwidth]{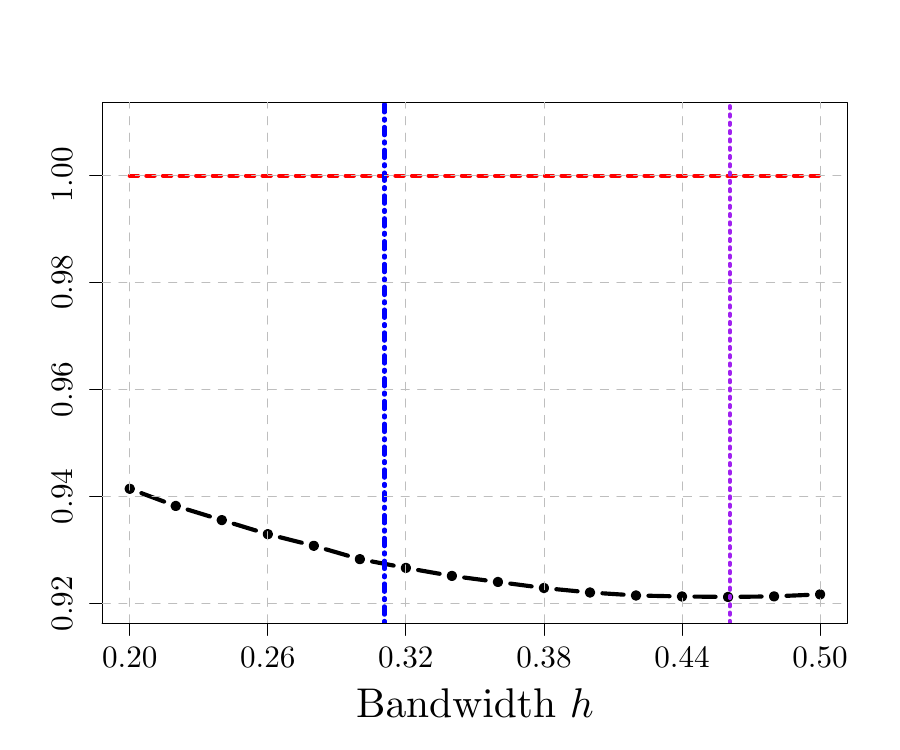}} 
  \subfigure[Model \eqref{model.quad} with $\cN(0, 4)$ error.]{\includegraphics[width=0.32\textwidth]{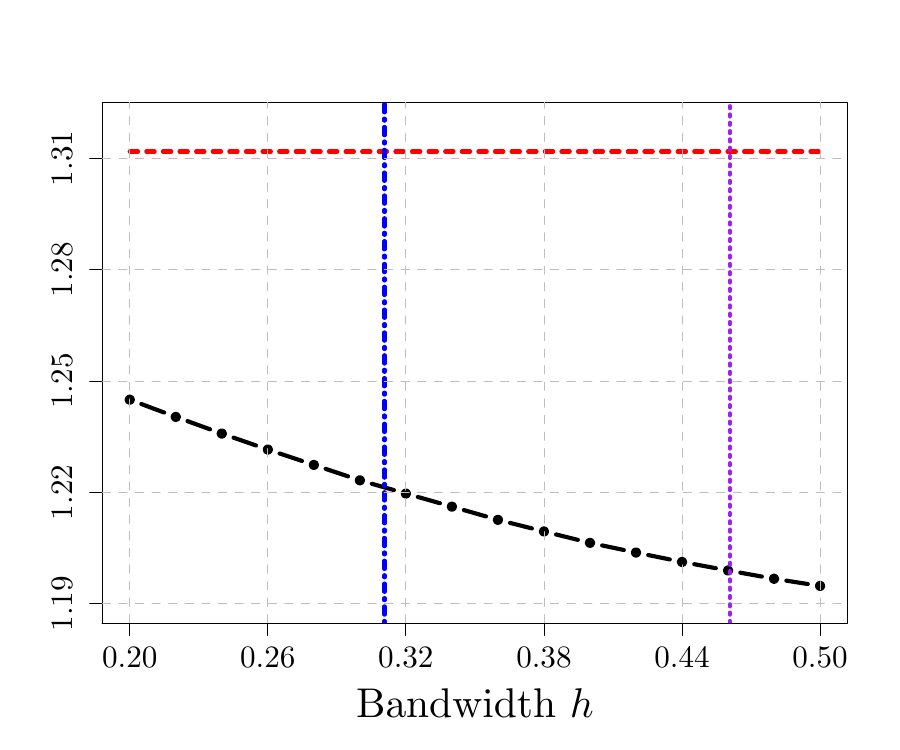}}
  \subfigure[Model \eqref{model.homo} with $t_2$ error.]{\includegraphics[width=0.32\textwidth]{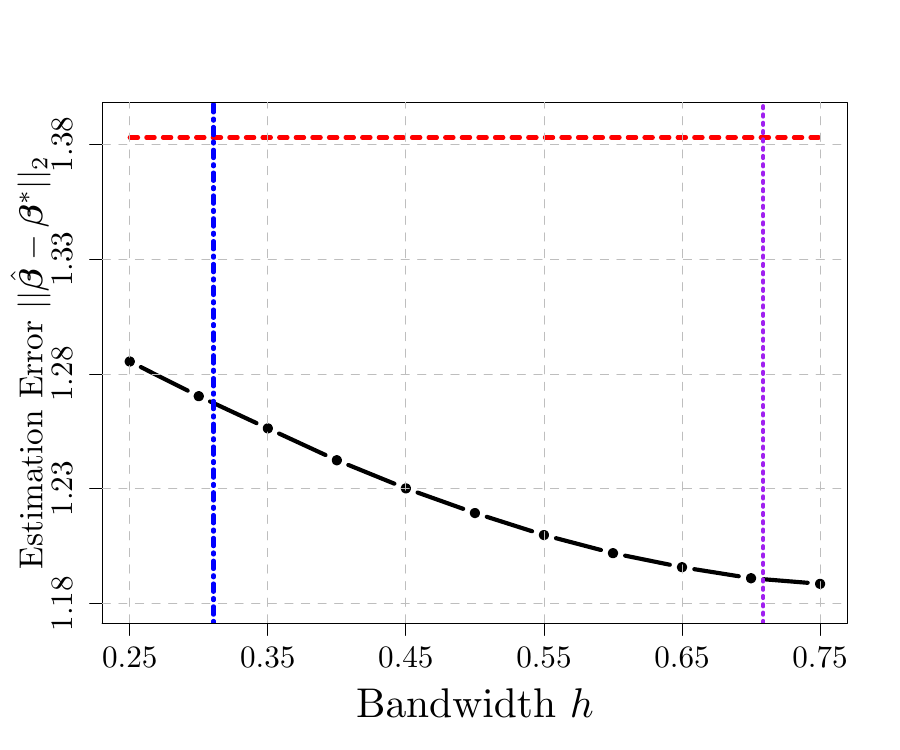}}
  \subfigure[Model \eqref{model.linear} with $t_2$ error.]{\includegraphics[width=0.32\textwidth]{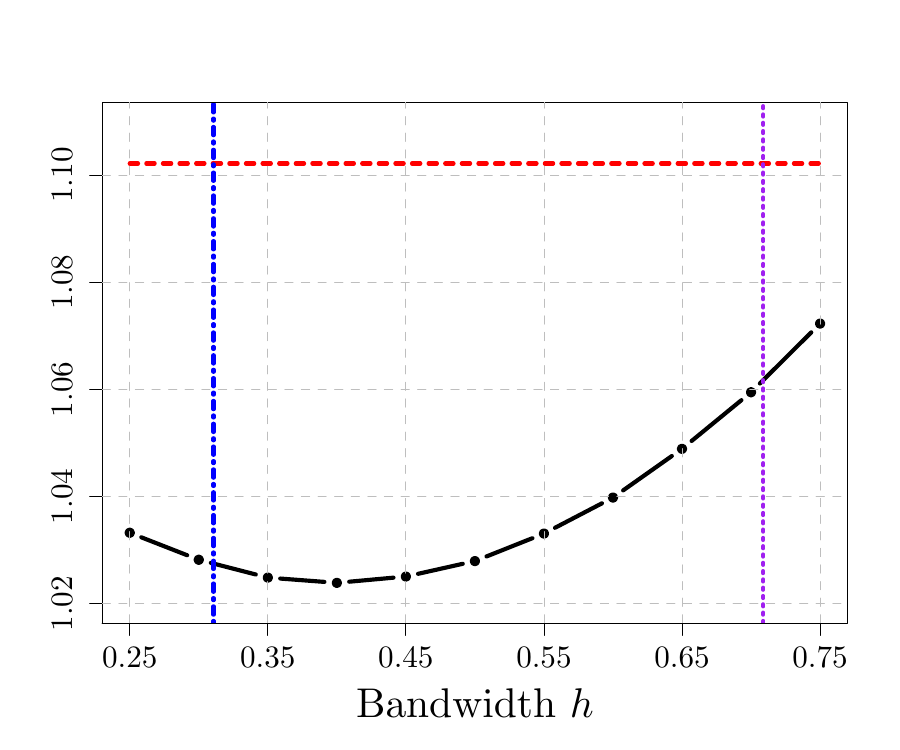}}
  \subfigure[Model \eqref{model.quad} with $t_2$ error.]{\includegraphics[width=0.32\textwidth]{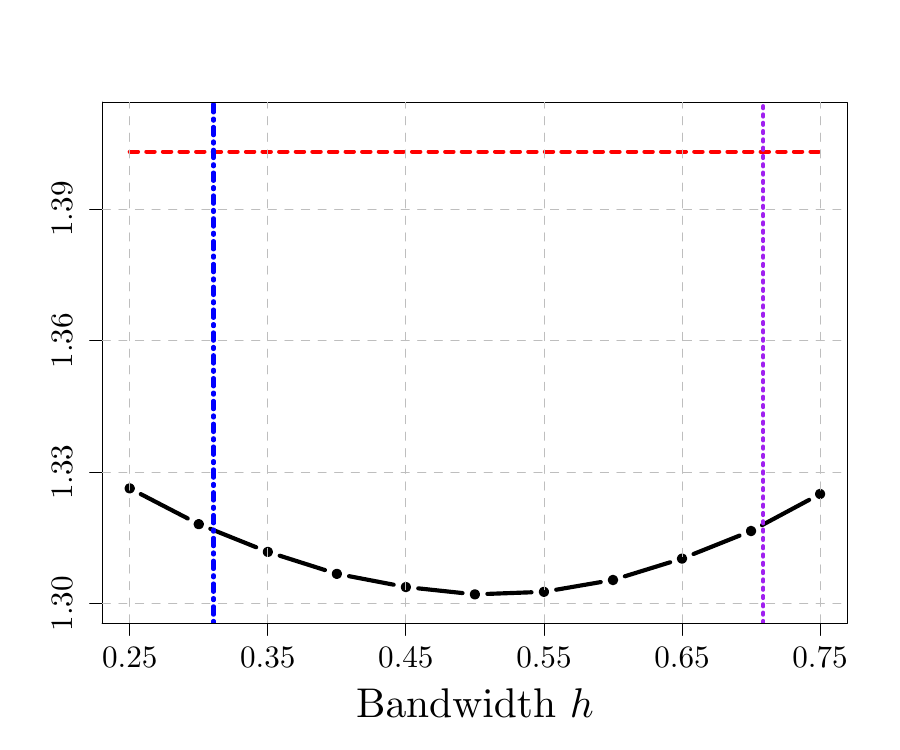}}
\caption{Sensitivity analysis of conquer with a range of bandwidth parameter $h$. Results for $(n, p) = (2000, 100)$, averaged over 500 datasets, with conquer implemented using a Gaussian kernel. The blue vertical dash line represents our default choice $h_{{\rm de}} = \{ (p+\log n) /n\}^{2/5}$, the purple vertical dash line refers to the choice $h_{\mathrm{AMSE}}$ from \eqref{AMSE.bandwidth2} with some oracle knowledge, and the red horizontal dash line represents the estimation error of standard QR.}
  \label{sensitivity}
\end{figure}

Next, we conduct a sensitivity analysis for conquer regarding the smoothing bandwidth $h$. We first set $(n, p) = (2000, 100)$ and consider the simulation settings \eqref{model.homo}--\eqref{model.quad} with $\mathcal{N}(0, 4)$ and $t_2$ noise.
{We perform conquer with $h$ taken from a wide range, including the default value  $h_{{\rm de}} = \{ (p+\log n) /n\}^{2/5} = 0.3107$ guided by Theorem~\ref{thm:clt}, and the bandwidth $h_{\mathrm{AMSE}}$ determined by AMSE in \eqref{AMSE.bandwidth2}, and compare the estimation error with that of QR in Figure~\ref{sensitivity}. For $h_{\mathrm{AMSE}}$, we directly substitute in the oracle values of $f_{ \varepsilon|\bx }(0)$ and $f'_{ \varepsilon|\bx }(0)$ for illustration purpose, and in practice, estimating these two quantities usually involves non-parametric techniques.
We see that the estimation error of conquer is uniformly lower than that of QR over a range of $h$, including our default choice, suggesting that conquer is insensitive to the choice of bandwidth $h$. 
We then consider a growing-scale scenario under model \eqref{model.linear} with $t_2$ error, and compare conquer with QR in Figure~\ref{sensitivity.grow}. 
This can be regarded as an extension of panel (e) in Figure~\ref{sensitivity}.
We observe that the estimation error of conquer using our default bandwidth is  uniformly lower that of QR in various settings, and is comparable to that using AMSE-based bandwidth.}

%%%%%%%%%%%%%%%%%%%%%%%%%%%%%%%%%%%%%%%%%%%%
%%%%%%%%%%%%%%%%%%%%%%%%%%%%%%%%%%%%%%%%%%%%
% Inference
%%%%%%%%%%%%%%%%%%%%%%%%%%%%%%%%%%%%%%%%%%%%
%%%%%%%%%%%%%%%%%%%%%%%%%%%%%%%%%%%%%%%%%%%%
 \subsection{Inference}
 \label{subsec:bootstrap inference}
In this section, we assess the performance of the multiplier bootstrap procedure for constructing confidence interval for each of the regression coefficients obtained from conquer. 
We implement conquer using the Gaussian kernel, and construct three types of confidence intervals: (i) the percentile \texttt{mb-per}; (ii) pivotal \texttt{mb-piv}; (iii) and regular \texttt{mb-norm} confidence intervals, as described in Section~\ref{subsec:mbi}.  
We also refer to the proposed multiplier bootstrap procedure as \texttt{mb-conquer} for simplicity.  
We compare the proposed method to several widely used inference methods for QR.  In particular, we consider confidence intervals by inverting a rank score test, \texttt{rank} (\cite{GJ1992}; Section 3.5 of \cite{K2005}); a bootstrap variant based on pivotal estimating functions, \texttt{pwy} \citep{PWY1994}; and wild bootstrap with Rademacher weights, \texttt{wild} \citep{FHH2011}. 
The three methods \texttt{rank}, \texttt{pwy}, and \texttt{wild} are implemented using the \texttt{R} package \texttt{quantreg}. 
Note that \texttt{rank} is a non-resampling based procedure that relies on prior knowledge on the random noise, i.e., a user needs to specify whether the random noise are independent and identically distributed.  In our simulation studies, we provide \texttt{rank} an unfair advantage by specifying the correct random noise structure.

We set $(n, p) = (800, 20)$, $\tau \in \{0.5, 0.9\}$, and significance level $\alpha = 0.05$. 
All of the resampling methods are implemented with $B = 500$ bootstrap samples.  To measure the reliability, accuracy, and computational efficiency of different methods for constructing confidence intervals, we calculate the average empirical coverage probability, average width of confidence interval, and the average runtime.  The average is taken over all regression coefficients without the intercept.  
Results based on 500 replications are reported in Figure~\ref{inf.t.9}, and Figures~\ref{inf.normal.9}--\ref{inf.normal.5} in Appendix~\ref{sec:app.simu}.

%%%%%%%%%%%%%%%%%%%%%%%%%%%%%%%%%%%%
%%%%%%%%%%%%%%%%%%%%%%%%%%%%%%%%%%%%
% Figure - boxplot for the confidence intervals
%%%%%%%%%%%%%%%%%%%%%%%%%%%%%%%%%%%%
%%%%%%%%%%%%%%%%%%%%%%%%%%%%%%%%%%%%
\begin{figure}[!ht]
  \centering
  \subfigure[Coverage under model \eqref{model.homo}.]{\includegraphics[width=0.32\textwidth]{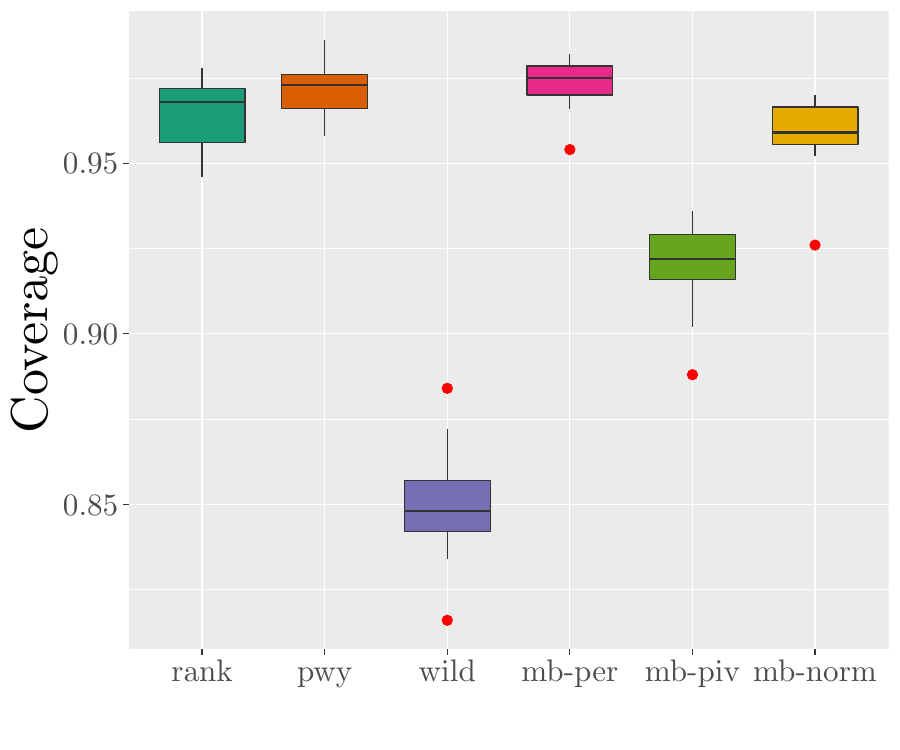}} 
  \subfigure[Coverage under model \eqref{model.linear}.]{\includegraphics[width=0.32\textwidth]{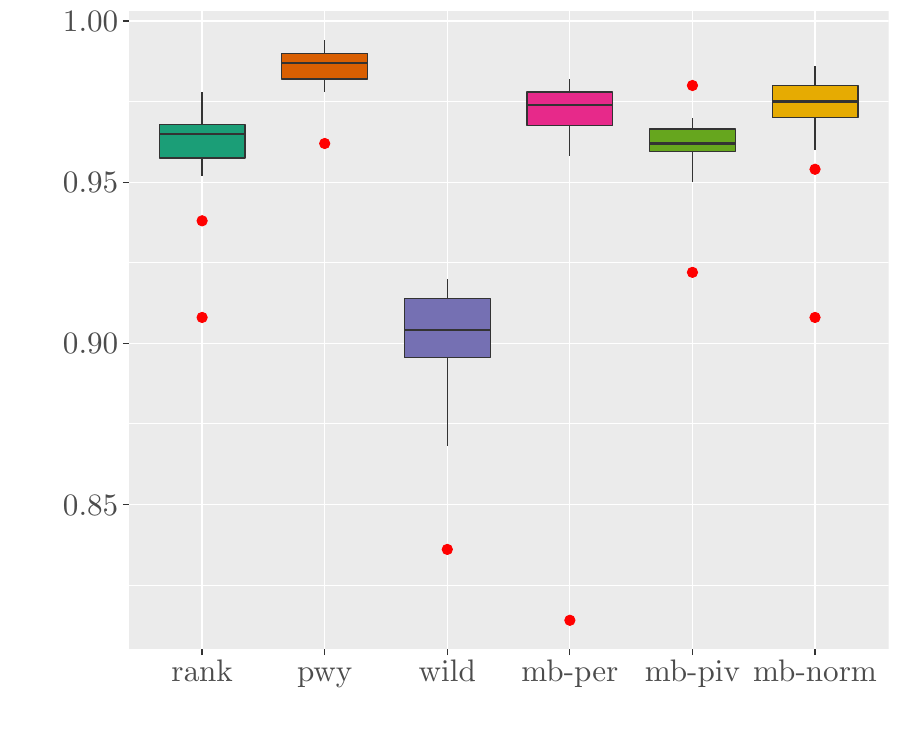}} 
  \subfigure[Coverage under model \eqref{model.quad}.]{\includegraphics[width=0.32\textwidth]{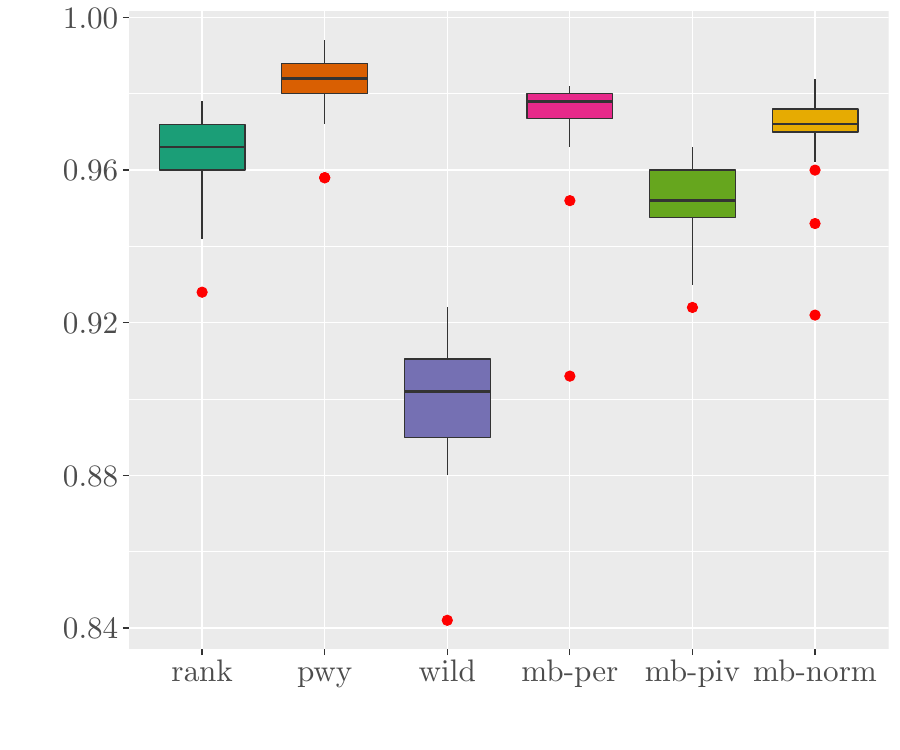}}
  \subfigure[CI width under model \eqref{model.homo}.]{\includegraphics[width=0.32\textwidth]{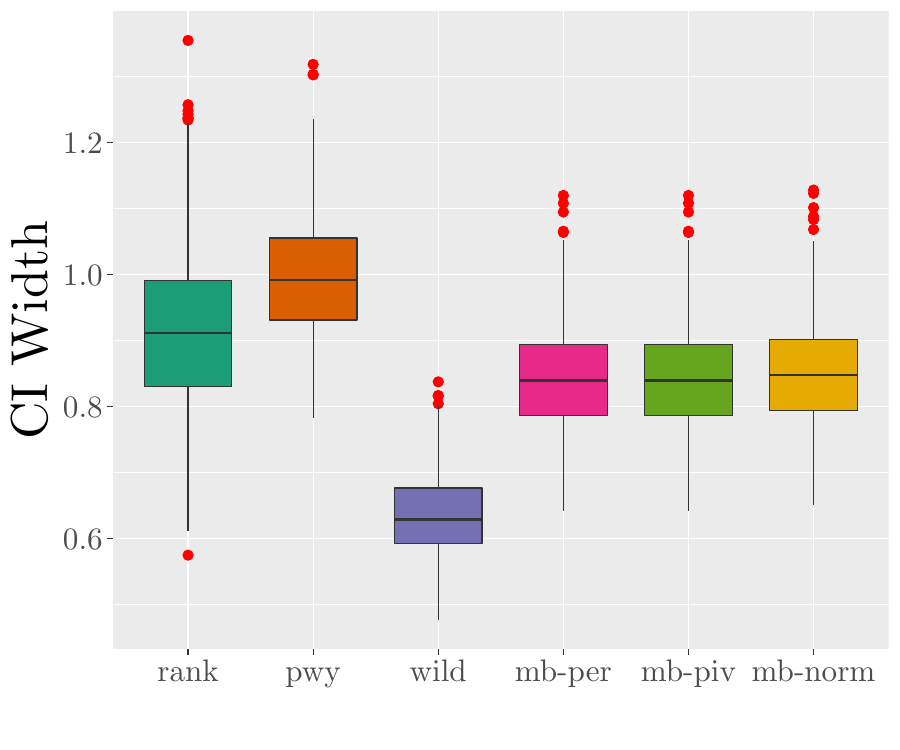}}
  \subfigure[CI width under model \eqref{model.linear}.]{\includegraphics[width=0.32\textwidth]{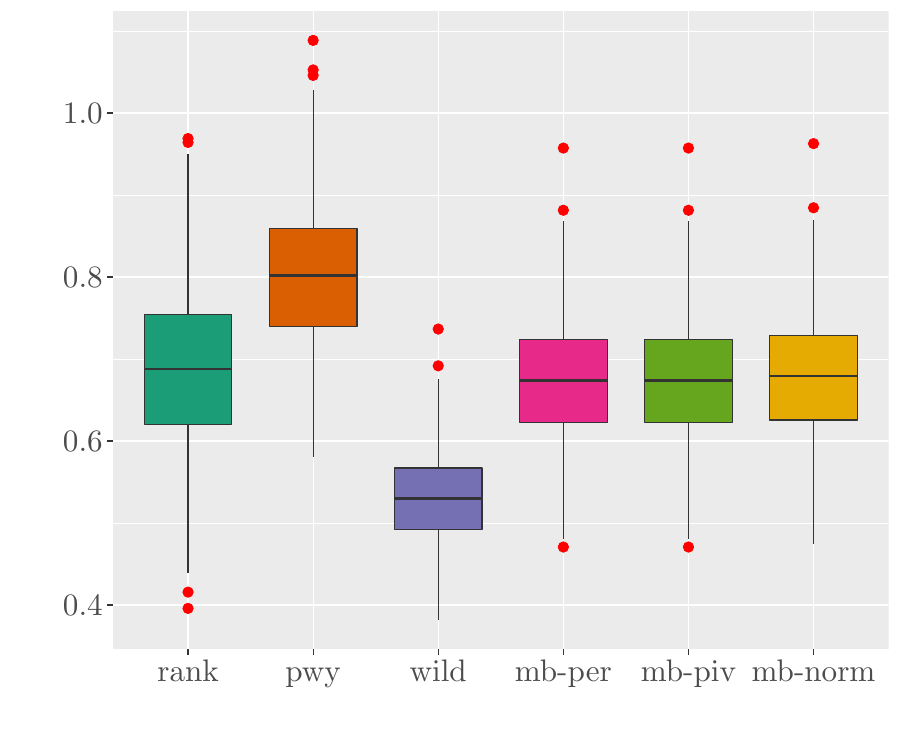}}
  \subfigure[CI width under model \eqref{model.quad}.]{\includegraphics[width=0.32\textwidth]{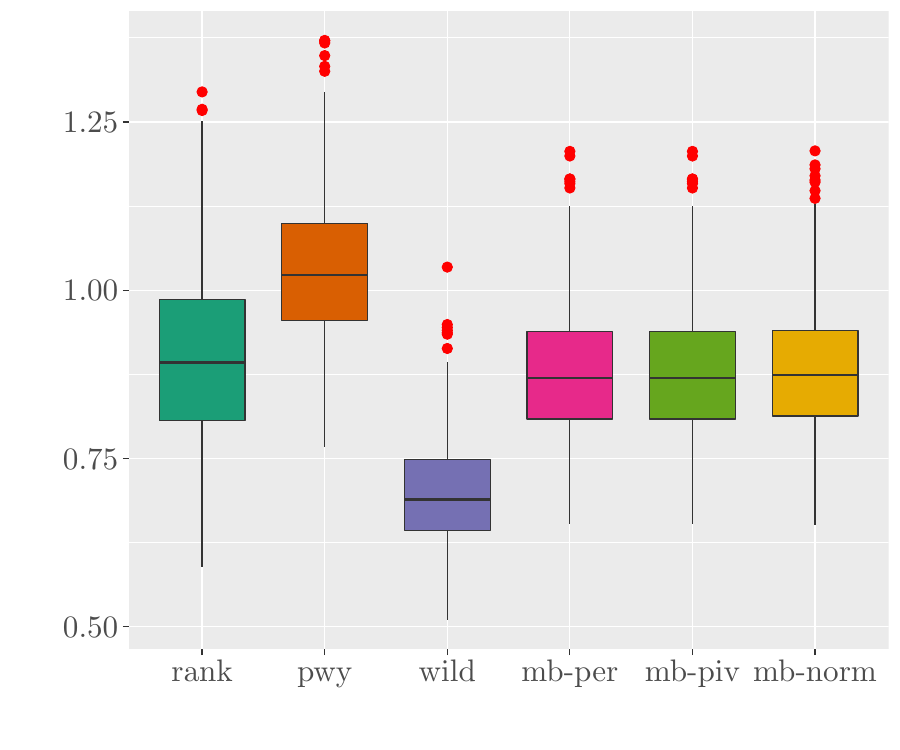}}
  \subfigure[Runtime under model \eqref{model.homo}.]{\includegraphics[width=0.32\textwidth]{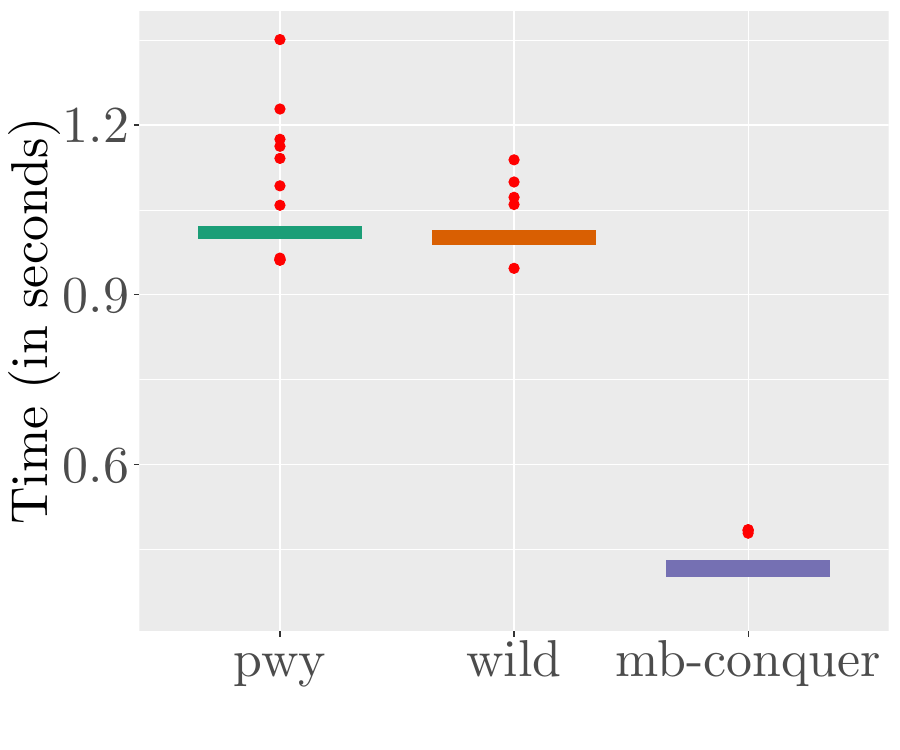}}
  \subfigure[Runtime under model \eqref{model.linear}.]{\includegraphics[width=0.32\textwidth]{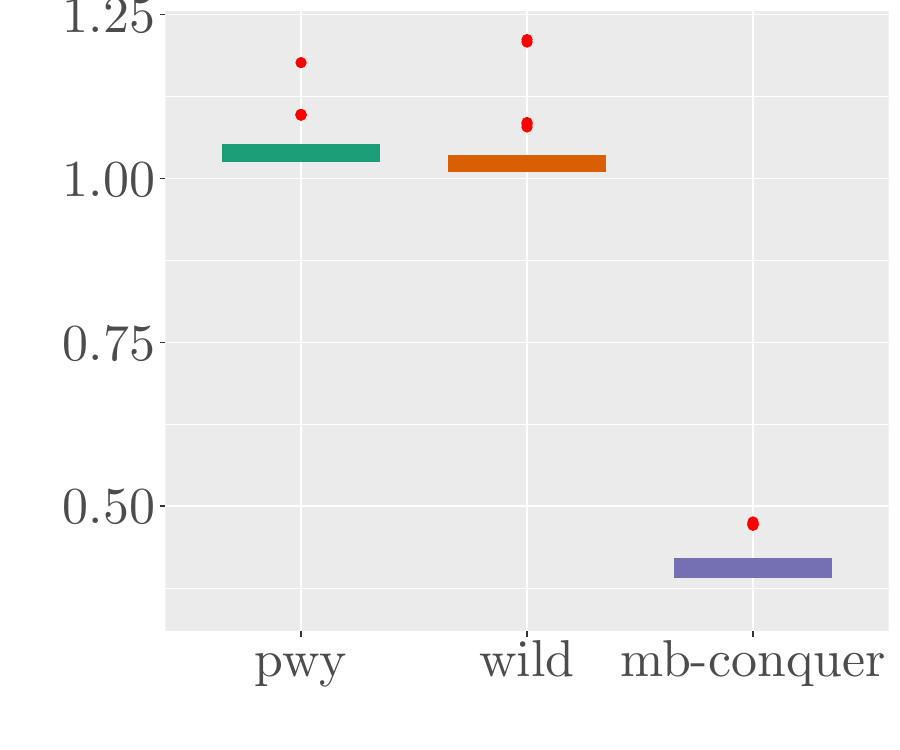}}
  \subfigure[Runtime under model \eqref{model.quad}.]{\includegraphics[width=0.32\textwidth]{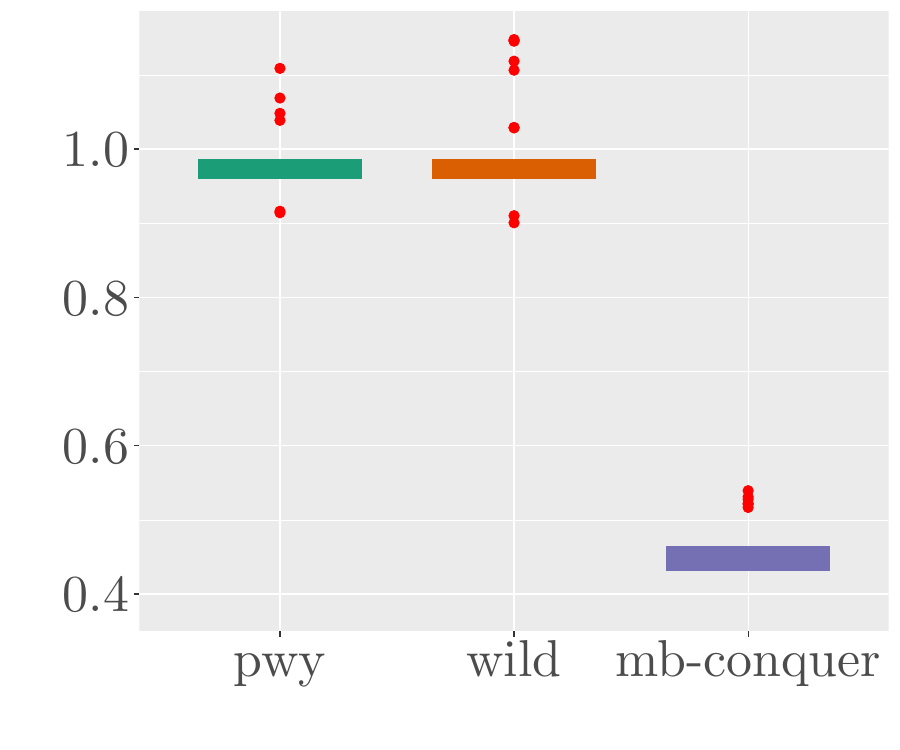}}
\caption{Empirical coverage, confidence interval width, and elapsed time of $six$ methods: \texttt{rank}, \texttt{pwy}, \texttt{wild} and three types of mb-conquer: \texttt{mb-per}, \texttt{mb-piv}, and \texttt{mb-norm} under models \eqref{model.homo}--\eqref{model.quad} with $t_2$ errors. For the running time, \texttt{rank} is not included since it is not a resampling-based method. The quantile level $\tau$ is fixed to be $0.9$, and the results are averaged over 500 data sets.}
  \label{inf.t.9}
\end{figure}

In Figure~\ref{inf.t.9}, and Figures~\ref{inf.normal.9}--\ref{inf.normal.5} in Appendix~\ref{sec:app.simu}, we use the rank-inversion method, \texttt{rank}, as a benchmark since we implement \texttt{rank} using information about the true underlying random noise, which is practically infeasible.
In the case of $\tau=0.9$,  \texttt{pwy} is the most conservative as it produces the widest confidence intervals with slightly inflated coverage probability, and \texttt{wild} gives the narrowest confidence intervals but at the cost of coverage probability.
The proposed methods \texttt{mb-per}, \texttt{mb-piv}, and \texttt{mb-norm} achieve a good balance between reliability (high coverage probability) and accuracy (narrow CI width), and moreover, has the lowest runtime.

\begin{figure}[!ht]
  \centering
  \subfigure[Coverage]{\includegraphics[width=0.32\textwidth]{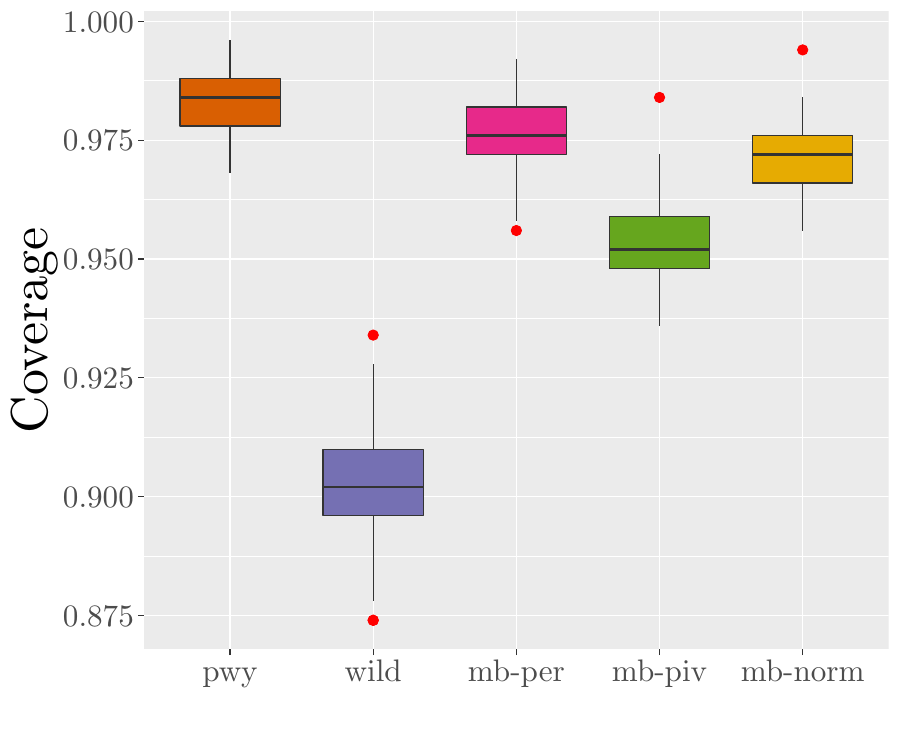}} 
  \subfigure[CI width]{\includegraphics[width=0.32\textwidth]{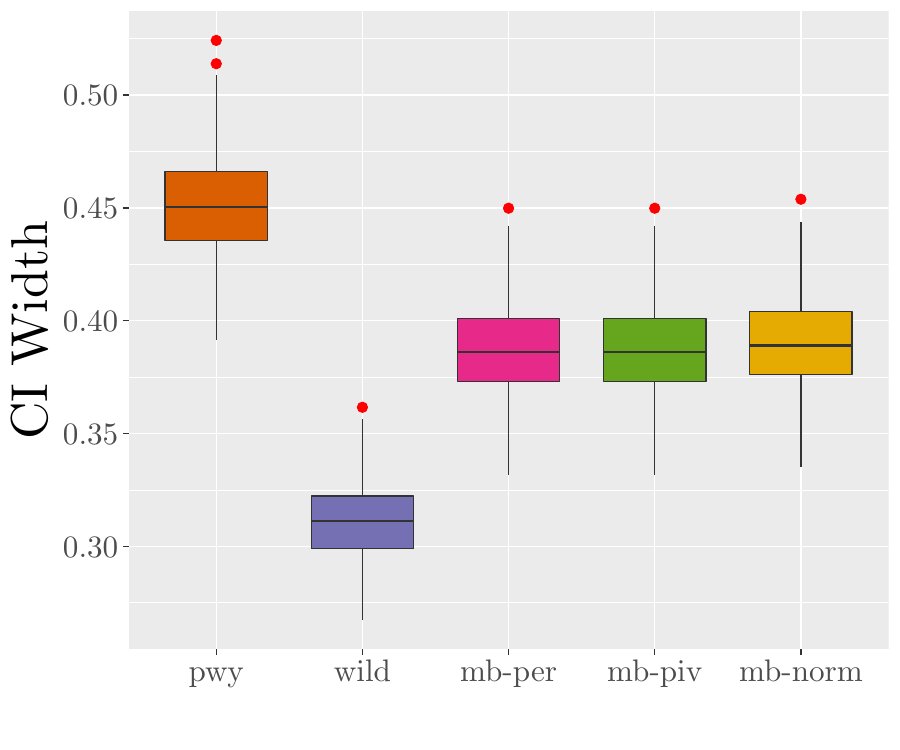}} 
  \subfigure[Runtime]{\includegraphics[width=0.32\textwidth]{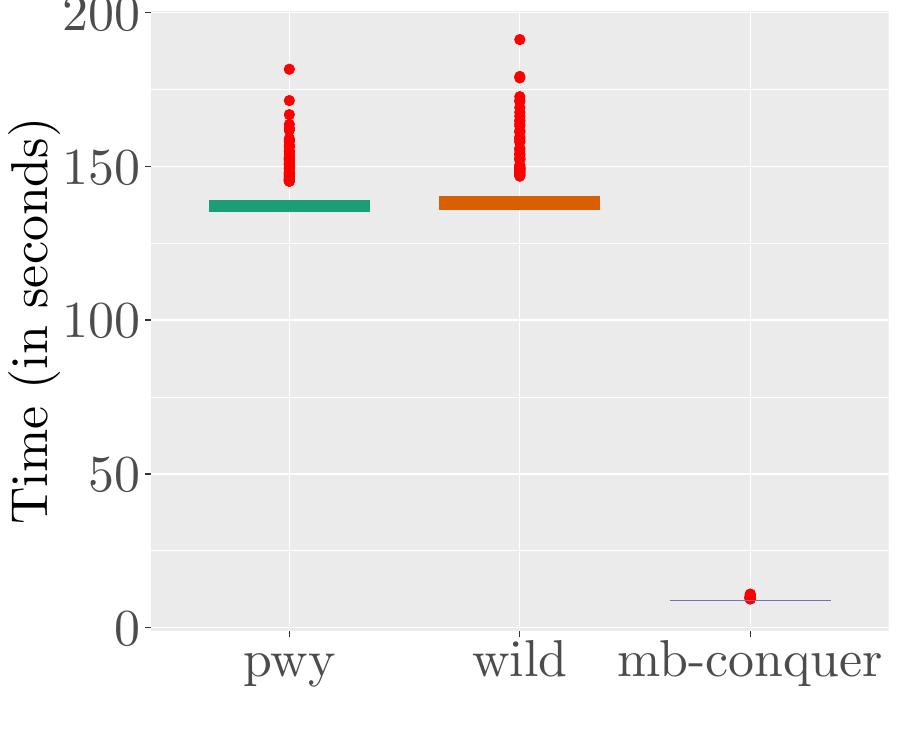}}
\caption{Empirical coverage, confidence interval width and elapsed time of \texttt{pwy}, \texttt{wild} and $3$ types of mb-conquer: \texttt{mb-per}, \texttt{mb-piv}, and \texttt{mb-norm} under quadratic heterogeneous model \eqref{model.quad} with $t_2$ errors. This figure extends the rightmost column of Figure~\ref{inf.t.9} to larger scale: $(n, p) = (4000, 100)$.}
  \label{inf.large.t.9}
\end{figure}

To further highlight the computational gain of the proposed method, we now perform numerical studies with larger $n$ and $p$.  In this case, the rank inversion method \texttt{rank} is computationally infeasible.  For example, when $(n,p)= (5000,250)$, rank inversion takes approximately 80 minutes while conquer with multiplier bootstrap takes 41 seconds for constructing confidence intervals.
We therefore omit \texttt{rank} from the following comparison. We consider the quadratic heterogeneous model \eqref{model.quad} with $(n, p) = (4000, 100)$ and $t_2$ noise.  The results are reported in Figure~\ref{inf.large.t.9}. We see that \texttt{pwy} and \texttt{wild} take up to 200 seconds while \texttt{mb-conquer} takes less than 10 seconds.  In summary, \texttt{mb-conquer} leads to a huge computational gain without sacrificing statistical efficiency.

\subsection{Comparison between normal approximation and bootstrap calibration}
\label{sec:normal-boot}

Finally, we complement the above studies with a comparison between the normal approximation and bootstrap calibration methods for confidence estimation. We consider model \eqref{model.homo} with $(n,p) = (2000,10)$. We use the same $\bbeta^*\in \RR^{p}$ and $\beta^*_0$ as before, and generate random covariates and noise from a multivariate uniform distribution and $t_{1.5}$-distribution, respectively. For each of the $p$ regression coefficients, we apply the proposed bootstrap percentile method and the normal-based method \citep{FGH2019} to construct pointwise confidence intervals at quantile indices close to 0 and 1, that is, $\tau \in \{ 0.05, 0.1, 0.9, 0.95\}$.

\begin{figure}[!ht]
  \centering
  \subfigure[Boxplots of empirical coverage]{\includegraphics[width=0.32\textwidth]{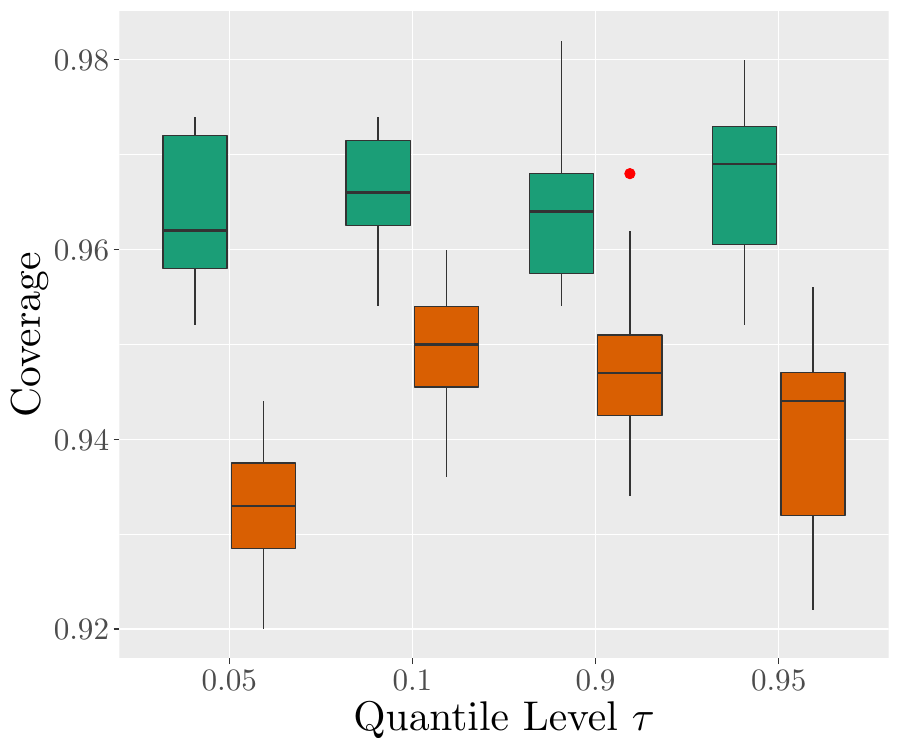}}~~~~~~~~~~~~~~~~~
  \subfigure[Boxplots of CI width]{\includegraphics[width=0.32\textwidth]{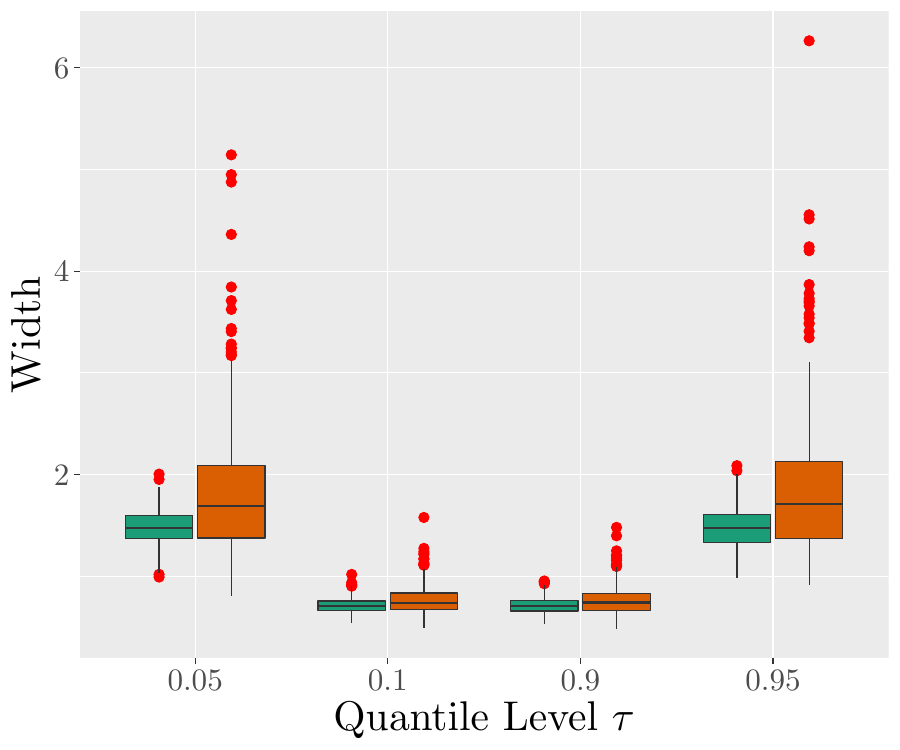}} 
\caption{Boxplots of the empirical coverage and CI width for the bootstrap percentile method \protect\includegraphics[height=0.88em]{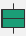} and normal-based method \protect\includegraphics[height=0.88em]{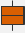} under model \eqref{model.homo} with $t_{1.5}$ error.}
  \label{norm-boot}
\end{figure}

Boxplots of the empirical coverage and CI width for the two methods are reported in Figure~\ref{norm-boot}. Considering that extreme quantile regressions are notoriously hard to estimate, the bootstrap method can produce much more reliable (high coverage) and accurate (narrow width) confidence intervals than the normal-based counterpart. Therefore, for applications in which extreme quantiles are of particular interest, such as  the problem of forecasting the conditional value-at-risk of a financial institution \citep{CU2001}, the bootstrap provides a more reliable approach for quantifying the uncertainty of the estimates.

%%%%%%%%%%%%%%%%%%%%%%%%%%%%%%%%%%%%%%%%%
%%%%%%%%%%%%%%%%%%%%%%%%%%%%%%%%%%%%%%%%%
% Bibliography
%%%%%%%%%%%%%%%%%%%%%%%%%1%%%%%%%%%%%%%%%%
%%%%%%%%%%%%%%%%%%%%%%%%%%%%%%%%%%%%%%%%%
\section{Discussion}
In this paper, we provide a comprehensive study on the statistical properties of {\it conquer}, namely, convolution-type smoothed quantile regression, under the non-asymptotic setting in which $p$ is allowed to increase as a function of $n$ while $p/n$ being small. 
When a non-negative kernel is used,  the smoothed objective function is  convex,  twice continuously differentiable, and locally strongly convex in a neighborhood of $\bbeta^*$ (with high probability). An efficient gradient-based algorithm is proposed to compute the conquer estimator,  which is scalable to very large-scale problems. 
{
For traditional QR computation with linear programming, interior point algorithms are typically used to get solutions with high precision (low duality gap) \citep{PK1997}. When applied to large-scale datasets, this may be inefficient for two reasons: (i) it takes a lot more time to reach a duality gap of the order of machine precision, and (ii) such a generic algorithm, which is less tailored to problem structure, tends to be very slow or even run out of  memory (unless with a high performance computing cluster).  In this regard, convolution smoothing offers a balanced tradeoff between statistical accuracy and computational complexity.   
} 
 
 {
 In the context of nonparametric density or regression estimation,  it is known that when higher-order kernels are used (and if the density or regression function has enough derivatives), the bias is proportional to $h^\nu$ for some $\nu\geq 4$ which is of better order than $h^2$. While a higher-order kernel has negative parts, the resulting smoothed loss is non-convex and thus brings the computational issue once again. 
 Motivated by the two-stage procedure proposed by \cite{B1975} whose original idea is to improve an initial estimator that is already consistent but not efficient, we further propose a one-step conquer estimator using higher-order kernels but without the need for solving a large-scale non-convex optimization. With increasing degrees of smoothness, the one-step conquer is asymptotically normal under a milder dimension constraint of roughly $p^2/n \to 0$.  Due to space limitations,  the details of this method are relegated to Section~\ref{sec:one-step} in the supplementary material. 
 }
 
 {In high-dimensional settings in which $p \gg n$, various authors have studied the regularized quantile regression under the sparsity assumption that most of the regression coefficients are zero \citep{BC2011,WWL2012,ZPH2015}.  Due to the vast literature in regularized quantile regression, we refer the reader to  Chapters 15--16 of \cite{KCHP2017} for an overview.  }
The computation of $\ell_1$-penalized QR is based on either reformulation as linear programs or alternating direction method of multiplier algorithms \citep{Gu2018}. 
Since the conquer loss is convex and twice differentiable, we expect that gradient-based algorithms, such as coordinate gradient descent or proximal gradient descent, will enjoy superior computational efficiency for solving regularized conquer without sacrificing statistical accuracy.

%%%%%%%%%%%%%%%%%%%%%%%%%%%%%%%%%%%%%%%%%
%%%%%%%%%%%%%%%%%%%%%%%%%%%%%%%%%%%%%%%%%
% Bibliography
%%%%%%%%%%%%%%%%%%%%%%%%%1%%%%%%%%%%%%%%%%
%%%%%%%%%%%%%%%%%%%%%%%%%%%%%%%%%%%%%%%%%

\newpage 
\appendix

%%%%%%%%%%%%% %%%%%%%%%%%%% %%%%%%%%%%%%%
% Proofs
%%%%%%%%%%%%% %%%%%%%%%%%%% %%%%%%%%%%%%%

\section{Asymptotic mean squared error}
\label{sec:AMSE}

In this section, we provide details for the statements in Remark~\ref{rmk:AMSE}. 
Throughout, assume the conditional density $f_{\varepsilon|\bx}(\cdot)$ satisfies (i)
$|f_{\varepsilon| \bx}(u)-f_{\varepsilon| \bx}(v)| \leq l_0(\bx) |u-v|$, and (ii) $|f''_{\varepsilon| \bx}(u) - f''_{\varepsilon | \bx}(v)| \leq l_2(\bx) |u-v|$ for all $u, v\in \RR$ and $\bx\in \RR^p$, where $\ell_0(\cdot)$ and $\ell_2(\cdot)$ are such that $\max\{\EE \ell_0^2(\bx)  ,  \EE   \ell_2^2(\bx) \} \leq C$ for some $C>0$. The refined Berry-Esseen bound \eqref{linear.clt2} in Theorem~\ref{thm:clt} indicates that
\#
	n^{1/2}\frac{\big\langle \ba ,  \hat \bbeta_h - \bbeta^* + 0.5 \kappa_2 h^2  \Jb_h^{-1} \EE \{f'_{\varepsilon | \bx}(0) \bx \}   \big\rangle }{\sqrt{\ba^\T   \Jb_h^{-1}  \EE \big[ \{ \cK_h(-\varepsilon) - \tau  \}^2 \bx \bx^\T \big] \Jb_h^{-1}   \ba }}  \xrightarrow {{\rm d}} \cN(0, 1) \nn
\#
uniformly over $\ba \in \RR^p$ as $n\to \infty$, provided $\sqrt{(p+\log n)/n} \lesssim h\lesssim \{ (p+\log n)/n\}^{1/4}$ and $p+\log n = o((nh)^{1/2})$. With dimension $p$ fixed, the asymptotic mean squared error (AMSE) of $n^{1/2}(\hat \bbeta_h - \bbeta^*)$ is thus determined by that of
$$
	\bg_h  \stackrel{{\rm d}}{=}  \cN\Bigg(  \frac{1}{2} \kappa_2  n^{1/2} h^2  \Jb_h^{-1} \EE \{ f'_{ \varepsilon|\bx } (0)  \bx \} , \,  \Jb_h^{-1}   \EE \big[  \{ \cK _h(-\varepsilon) - \tau \}^2  \bx \bx^\T   \big]    \Jb_h^{-1}   \Bigg) .
$$
To prove \eqref{asym.cov}, it suffices to derive the asymptotic expression of $\EE(\xi^2 \bx \bx^\T)$, where $\xi  = \tau - \cK_h(-\varepsilon )$.

Using integration by parts and a change of variable, we have
\$
& \EE (  \xi^2 |\bx ) = \int_{-\infty}^\infty \{ \tau - \cK(-t/h) \}^2 {\rm d} F_{ \varepsilon |\bx }(t)   \\
& = \tau^2 - \frac{2}{h}\int_{-\infty}^\infty   \{ \tau - \cK(-t/h)\} K(t/h) F_{ \varepsilon|\bx }(  t)  \, {\rm d}t   = \tau^2 - 2\int_{-\infty}^{\infty} \{ \cK(u) - 1 + \tau \} K(u) F_{\varepsilon|\bx}(hu)\mathrm{d}u .
\$
Note that $\int_{-\infty}^\infty K(u) \cK(u) {\rm d}u = 1/2$ and $F_{\varepsilon|\bx}(0) = \tau$. Then, applying a Taylor series expansion of $F_{ \varepsilon|\bx }(\cdot)$ around $0$ yields
\#
& \EE (  \xi^2 |\bx ) = \tau^2 - 2\int_{-\infty}^{\infty} \{ \cK(u) - 1 + \tau \} K(u) \Bigg[  \tau + hu f_{\varepsilon|\bx}(0) +  h u \int_0^1 \{ f_{\varepsilon|\bx}(h u w) - f_{\varepsilon|\bx}(0) \} {\rm d} w    \Bigg] \mathrm{d}u,  \nn \\
&= \tau^2 - 2\tau \int_{-\infty}^{\infty} \{ \cK(u) - 1 + \tau \} K(u)\mathrm{d}u - 2h f_{\varepsilon|\bx}(0) \int_{-\infty}^{\infty} u \{ \cK(u) - 1 + \tau \} K(u)\mathrm{d}u \nn \\
&~~~~~~~ - 2  h   \int_{-\infty}^{\infty} \int_{0}^1 u \{ \cK(u) - 1 + \tau \} K(u) \{ f_{\varepsilon|\bx}(h u w) - f_{\varepsilon|\bx}(0) \} {\rm d} w    \mathrm{d}u \nn \\
&= \tau(1 - \tau) - 2  \bar \kappa_1  f_{\varepsilon|\bx}(0) h -  2  h   \int_{-\infty}^{\infty} \int_{0}^1 u \{ \cK(u) - 1 + \tau \} K(u) \{ f_{\varepsilon|\bx}(h u w) - f_{\varepsilon|\bx}(0) \} {\rm d} w    \mathrm{d}u ,   \nn
\#
where $\bar \kappa_1  := \int_{-\infty}^\infty  u K(u)\cK(u) {\rm d} u   >0$. Consequently, 
\$
	 \big\| \EE  (\xi^2 \bx \bx^\T )   - \tau(1-\tau)  \bSigma + 2 \bar \kappa_1 h   \Jb   \big\|_2  \lesssim  \sup_{\bu \in \mathbb{S}^{p-1}} \EE \{  l_0(\bx) \langle \bx, \bu \rangle^2 \} \cdot  h^2
	\lesssim   \big\{ m_4 \cdot  \EE l_0^2(\bx) \big\}^{1/2} h^2 .
\$
Putting together the pieces we conclude that
\$
\EE \big( \bg_h \bg_h^\T \big)  & = \Jb_h^{-1} \EE (\xi^2 \bx \bx^\T) \Jb_h^{-1}  + \frac{1}{4} \kappa_2^2  n h^4 \Jb_h^{-1} \EE \{  f'_{\varepsilon |\bx} (0) \bx \} \EE \{ f'_{\varepsilon |\bx} (0) \bx^\T \} \Jb_h^{-1} \\
& =  \Jb_h^{-1} \Bigg[   \tau(1-\tau )\bSigma - 2 \bar{\kappa}_1 h \Jb   +  \frac{1}{4} \kappa_2^2 n  h^4   \EE \big\{ f'_{ \varepsilon|\bx } (0) \bx  \big\}  \EE \big\{ f'_{ \varepsilon|\bx } (0) \bx^\T   \big\}+  \cO (h^2)  \Bigg] \Jb_h^{-1} . 
%& =  \Jb^{-1} \Bigg\{ \frac{1}{n} \cov(\xi_1 \bx_1) + \EE(\xi_1 \bx_1) \EE(\xi_1 \bx_1)^\T \Bigg\} \Jb^{-1}  \\
%& =  \Jb^{-1} \Bigg\{ \frac{1}{n} \EE(\xi_1^2 \bx_1 \bx_1^\T) + \bigg(1+\frac{1}{n}\bigg)\EE(\xi_1 \bx_1) \EE(\xi_1 \bx_1)^\T \Bigg\} \Jb^{-1} \\
%& =  \Jb^{-1} \Bigg[ \frac{1}{n} \big\{  \tau(1-\tau )\bSigma - 2 \bar{\kappa}_1 h \Jb \big\}  +  \frac{\kappa_2^2}{4} h^4   \EE \big\{ f'_{ \varepsilon|\bx } (0) \bx_i \big\}  \EE \big\{ f'_{ \varepsilon|\bx } (0) \bx_i \big\}^\T  +  \cO\big(h^2 n^{-1} + h^6\big)  \Bigg] \Jb^{-1}
\$
Under the prescribed bandwidth constraint, the leading term of the MSE of $\bg_h$, and thus of the AMSE of $n^{1/2}(\hat \bbeta_h - \bbeta^*)$, is  
\#
  \Jb_h^{-1} \bSigma^{ 1/2} \Bigg[  \tau(1-\tau ) \Ib_p - 2 \bar{\kappa}_1 h \Hb   +  \frac{1}{4}  \kappa_2^2 n h^4   \EE (\bb) \EE (\bb^\T)   \Bigg] \bSigma^{ 1/2} \Jb_h^{-1} , \nn
\#
where $\bb=f'_{ \varepsilon|\bx } (0) \bw = f'_{ \varepsilon|\bx } (0) \bSigma^{-1/2}\bx$ and $\Hb  = \EE \{ f_{\varepsilon |\bx}(0) \bw \bw^\T \}$. This proves \eqref{asym.cov}.

\section{One-step conquer with higher-order kernels}
\label{sec:one-step}

As noted in Section~\ref{sec:theory:bias}, the smoothing bias is of order $h^2$ when a non-negative kernel is used. The ensuing empirical loss $\bbeta \mapsto    (1/n) \sn (\rho_\tau * K_h)(y_i -  \langle \bx_i , \bbeta \rangle )$ is not only twice-differentiable and convex, but also (provably) strongly convex in a local vicinity of $\bbeta^*$ with high probability. 
Kernel smoothing is ubiquitous in nonparametric statistics. The order of a kernel, $\nu$, is defined as the order of the first non-zero moment. The order of a symmetric kernel is always even. 
A kernel is called {\it high-order} if $\nu>2$, which inevitably has negative parts and thus is no longer a probability density.
Thus far we have focused on conquer with second-order kernels, and the resulting estimator achieves an $\ell_2$-error of the order $\sqrt{p/n} + h^2$.

Let $G(\cdot)$ be a higher-order symmetric kernel with order $\nu\geq 4$, and $b >0$ be a bandwidth. Again, via convolution smoothing, we may consider a bias-reduced estimator that minimizes the empirical loss $\bbeta \mapsto \hat Q_b^G(\bbeta):= (1/n) \sn (\rho_\tau * G_b )(y_i -  \langle \bx_i , \bbeta \rangle )$. This, however, leads to a non-convex optimization. Without further assumptions, finding a global minimum is computationally intractable: finding an $\epsilon$-suboptimal  point for a $k$-times continuously differentiable loss function requires at least $\Omega\{ (1/\epsilon)^{p/k}\}$ evaluations of the function and its first $k$ derivatives, ignoring problem-dependent constants; see Section~1.6 in \cite{NY1983}. Instead, various gradient-based methods have been developed for computing {\it stationary points}, which are points $\bbeta$ with sufficiently small gradient  $\| \nabla \hat Q_b^G(\bbeta) \|_2 \leq \epsilon$, where $\epsilon \geq 0$ is optimization error. However, the equation $\nabla \hat Q_b^G(\bbeta)  = \textbf{0}$ does not necessarily have a unique solution, whose statistical guarantees remain unknown.

Motivated by the classical one-step estimator \citep{B1975}, we further propose a one-step conquer estimator using high-order kernels, which bypasses solving a large-scale non-convex optimization. To begin with, we choose two symmetric kernel functions, $K: \RR \mapsto [0,\infty)$ with order two and $G(\cdot)$ with order $\nu\geq 4$, and let $h, b >0$ be two bandwidths.
First, compute an initial conquer estimator $\overbar \bbeta \in \argmin_{\bbeta \in \RR^p} \hat Q_{h}^K(\bbeta)$, where $\hat Q_{h}^K(\bbeta)  = (1/n) \sn (\rho_\tau * K_h) (y_i -  \langle \bx_i ,  \bbeta \rangle )$. Denote by $\bar \varepsilon_i = y_i - \langle \bx_i, \overbar \bbeta \rangle$ for $i=1,\ldots, n$ the fitted residuals. Next, with slight abuse of notation, we define the one-step conquer estimator $\hat \bbeta$ as a solution to the equation $  \nabla^2  \hat Q_b^G (\overbar \bbeta)  ( \hat \bbeta  - \overbar \bbeta) = - \nabla \hat Q_b^G(\overbar \bbeta)$, or equivalently,
 \#
 \Biggl\{ \frac{1}{n } \sn  G_b  ( \bar  \varepsilon_i  )   \bx_i \bx_i^\T \Biggr\} ( \hat \bbeta  - \overbar \bbeta) = \frac{1}{n} \sn   \bigl\{ \cG_b (  \bar \varepsilon_i   ) + \tau - 1 \bigr\} \bx_i . \label{os.conquer1}
 \#
where  $\hat Q_b^G(\bbeta)  = (1/n) \sn (\rho_\tau * G_b ) (y_i - \langle \bx_i , \bbeta \rangle)$ and $\cG_b(u) = \int_{-\infty}^{u/b} G(v) \, {\rm d}v$.
Provided that $\nabla^2  \hat Q_b^G (\overbar \bbeta)$ is positive definite, the one-step conquer estimate $\hat \bbeta$ essentially performs a Newton-type step based on $\overbar \bbeta$:
\#
\hat \bbeta = \overbar  \bbeta -  \bigl\{ \nabla^2  \hat Q_b^G (\overbar \bbeta) \bigr\}^{-1}  \nabla \hat Q_b^G(\overbar \bbeta). \label{os.conquer2}
\#
In this case, $\hat \bbeta$ can be computed by the conjugate gradient method \citep{HS1952}.

Theoretical properties of the one-step estimator $\hat \bbeta$ defined in \eqref{os.conquer1}, including the Bahadur representation and asymptotic normality with explicit Berry-Esseen bound, will be provided in on-line supplementary materials. For practical implementation, we consider higher-order Gaussian-based kernels. For $r=1,2,\ldots,$ the $(2r)$-th order Gaussian kernels are 
$$
 G ( u;2r) = \frac{(-1)^r \phi^{(2r-1)}(u)}{2^{r-1} (r-1)! u} = \sum_{\ell=0}^{r-1} \frac{(-1)^\ell}{2^\ell \ell!} \phi^{(2\ell)} (u) ;
$$
see Section~2 of \cite{WS1990}. Integrating $G(\cdot ; 2r)$ yields
$$
    \cG (v;2r) = \int_{-\infty}^v  G( u; 2r) \, {\rm d}u = \sum_{\ell=0}^{r-1} \frac{(-1)^\ell}{2^\ell \ell!} \phi^{(2\ell -1)}(v).
$$
In fact, both $G(\cdot; 2r)$ and $\cG(\cdot; 2r)$ have simpler forms $G(u;2r) = p_r(u) \phi(u)$ and $\cG(u; 2r) = \Phi(u)+ P_r(u) \phi(u)$, where $p_r(\cdot)$ and $P_r(\cdot)$ are polynomials in $u$. For example, $p_1(u) = 1$, $P_1(u)=0$, $p_2(u) =(-u^2+3)/2$, $P_2(u)=u/2$, $p_3(u)= (u^4 - 10 u^2 + 15)/8$, and $P_3(u)=(-u^3 +7u)/8$.
We refer to \cite{O2020} for more details when $r$ is large.

\section{Proofs}

Recall that $\bx=(x_1,\ldots , x_p)^\T$ is such that $x_1 \equiv 1$, $\EE(x_j) =0$ for $j=2,\ldots, p$, and $\bSigma = \EE(\bx \bx^\T)$ is positive definite. In this case, 
$$
\bSigma = \begin{bmatrix}
   1      &  \textbf{0}_{p-1}^\T \\
   \textbf{0}_{p-1}      & \Sb
\end{bmatrix} ~~\mbox{ with }~~\Sb = \EE(\bx_- \bx_-^\T ) ~~\mbox{ and }~~ \bw = \begin{bmatrix}
   1       \\
   \Sb^{-1/2} \bx_-
\end{bmatrix},
$$
where $\textbf{0}_{k}$ is the zero vector in $\RR^k$ ($k\geq 2$).
For every $r\geq 0$, define the ellipse $\Theta(r) = \{ \bu \in \RR^p : \| \bu \|_{\bSigma} \leq r \}$ and its boundary $\partial \Theta(r) = \{ \bu \in \RR^p : \| \bu \|_{\bSigma} = r \}$. 

%Throughout the proof, we write $\bOmega = \bSigma^{-1}$ and $\| \bu \|_{\bOmega} =  \| \bSigma^{-1/2} \bu \|_2 $ for $\bu \in \RR^p$. By H\"older's inequality, $|\langle \bu, \bv \rangle| \leq \| \bu \|_{\bSigma} \cdot \| \bv \|_{\bOmega}$ for any $\bu, \bv \in \RR^p$.

\subsection{Proof of Proposition~\ref{prop:bias}}

To begin with, define $\bdelta_h = \bbeta^*_h - \bbeta^* \in \RR^p$ and $\delta_h = \| \bdelta_h \|_{\bSigma}$.  By the convexity of $\bbeta \mapsto Q_h(\bbeta)$ and the first-order optimality condition $\nabla Q_h(\bbeta^*_h) = \textbf{0}$, we have
\#
 0\leq   \langle \nabla Q_h(\bbeta^*_h )  - \nabla Q_h(\bbeta^*) , \bbeta^*_h - \bbeta^* \rangle  =  \langle    - \nabla Q_h(\bbeta^*) ,  \bdelta_h  \rangle \leq \|  \bSigma^{-1/2} \nabla Q_h(\bbeta^*) \|_2 \cdot   \| \bdelta_h \|_{\bSigma}  ,  \label{pop.foc}
\#
where the last step follows from H\"older's inequality. Note that $\nabla Q_h(\bbeta^*) = \EE \{  \cK(-\varepsilon /h) - \tau \} \bx $. By integration by parts and a Taylor series expansion,
\#
   \EE \bigl\{   \cK(-\varepsilon /h) | \bx \bigr\} &  = \int_{-\infty}^\infty  \cK(-t/h) \, {\rm d} F_{ \varepsilon|\bx }(   t) \nn \\
& = -\frac{1}{h}\int_{-\infty}^\infty  K(-t/h) F_{ \varepsilon|\bx }(  t)  \, {\rm d}t = \int_{-\infty}^\infty  K(u) F_{\varepsilon|\bx}(  -h u)  \, {\rm d}u  \nn \\
& = \tau    + \int_{-\infty}^\infty  K(u) \int_{0}^{-hu} \bigl\{  f_{ \varepsilon|\bx } (  t) - f_{\varepsilon|\bx} (0) \bigr\} \, {\rm d} t  \, {\rm d}u ,\nn
\#
from which it follows that $| \EE \{   \cK(-\varepsilon /h) | \bx \} - \tau | \leq  0.5 l_0 \kappa_2    h^2$. 
Consequently,
\#
 &  \|  \bSigma^{-1/2}  \nabla Q_h(\bbeta^*) \|_2 = \sup_{\bu \in \mathbb S^{p-1}}  \EE  \bigl\{ \cK(-\varepsilon/h) -\tau \bigr\} \langle  \bu  ,  \bSigma^{-1/2} \bx \rangle \leq 0.5 l_0 \kappa_2 h^2 .  \label{score.mean.ubd}
\#

Turning to the left-hand side of \eqref{pop.foc}, applying the mean value theorem for vector-valued functions implies
\#
 \nabla Q_h( \bbeta^*_h )  - \nabla Q_h(\bbeta^*)  =   \int_0^1   \nabla^2 Q_h(  \bbeta^* + t \bdelta_h )  \, {\rm d} t  \,  \bdelta_h  , \label{pop.grad.diff}
\# 
where  $\nabla^2 Q_h( \bbeta  )  =  \EE \bigl\{ K_h(  y - \langle \bx, \bbeta \rangle ) \bx \bx^\T \bigr\}$ for $\bbeta \in \RR^p$. With $\bdelta = \bbeta - \bbeta^*$, note that
\#
 \EE \bigl\{ K_h(  y - \langle \bx, \bbeta \rangle ) | \bx \bigr\}  = \frac{1}{h} \int_{-\infty}^\infty K\Bigg(  \frac{u- \langle \bx, \bdelta  \rangle  }{h} \Bigg) f_{\varepsilon|\bx} (u) \, {\rm d} u =  \int_{-\infty}^\infty K (v )  f_{\varepsilon|\bx} ( \langle \bx ,  \bdelta \rangle + h v ) \, {\rm d} v . \nn
\#
By the Lipschitz continuity of $f_{ \varepsilon|\bx }(\cdot)$,
\#
	\EE \bigl\{ K_h(  y - \langle \bx, \bbeta \rangle ) | \bx \bigr\}   =  f_{ \varepsilon|\bx }( 0 ) + R_h(\bdelta)  \label{Kh.mean.ubd1}
\#
with $R_h(\bdelta)$ satisfying $|R_h(\bdelta)|  \leq   l_0  \bigl( |\langle \bx, \bdelta \rangle | + \kappa_1 h \bigr)$.
Together,   \eqref{pop.grad.diff}, \eqref{Kh.mean.ubd1} and the assumption $f_{\varepsilon | \bx }(0) \geq \underbar{$f$}  >0$ (almost surely) yield 
\#
& 	\langle \nabla Q_h(  \bbeta_h^* )  - \nabla Q_h(\bbeta^*) , \ \bbeta_h^*  - \bbeta^* \rangle \nn \\
& \geq  \underbar{$f$}  \cdot \|  \bdelta_h  \|_{\bSigma}^2    - 0.5 l_0  \EE |  \langle \bx, \bdelta_h  \rangle  |^3 - l_0 \kappa_1 h \cdot \| \bdelta_h  \|_{\bSigma}^2  \geq   \underbar{$f$} \cdot  \delta_h^2     -   0.5 l_0 m_3 \cdot  \delta_h^3  - l_0 \kappa_1 h \cdot \delta_h^2  .  \label{pop.foc.lbd}
\#

Combining \eqref{pop.foc} with the upper and lower bounds \eqref{score.mean.ubd} and \eqref{pop.foc.lbd},  we find that $\delta_h \geq 0$ satisfies $0.5 l_0 m_3  \cdot \delta_h^2 - (\underbar{$f$} - l_0 \kappa_1 h ) \delta_h + 0.5 l_0 \kappa_2 h^2 \geq 0$. 
Provided that $l_0 \{ \kappa_1 +  ( m_3  \kappa_2)^{1/2} \} h >   \underbar{$f$}$,  solving this inequality  yields
\#
	\delta_h %\leq \frac{\underbar{$f$} - l_0 \kappa_1 h - \sqrt{(\underbar{$f$} - l_0 \kappa_1 h)^2 - l_0^2 m_3 \kappa_2 h^2}}{l_0 m_3 } 
\leq   \frac{l_0   \kappa_2 h^2 }{  \underbar{$f$} - l_0 \kappa_1 h +  \Delta_h^{1/2} } ~~\mbox{ or }~~\delta_h %\leq \frac{\underbar{$f$} - l_0 \kappa_1 h - \sqrt{(\underbar{$f$} - l_0 \kappa_1 h)^2 - l_0^2 m_3 \kappa_2 h^2}}{l_0 m_3 } 
\geq   \frac{ \underbar{$f$}  -  l_0 \kappa_1 h +  \Delta_h^{1/2} }{   l_0 m_3  } .  \label{two-sided.sol}
\#
where $\Delta_h := (\underbar{$f$} - l_0 \kappa_1 h)^2 - l_0^2 m_3 \kappa_2 h^2 > 0$.
 It remains to rule out the second bound in \eqref{two-sided.sol}.

Assume $\delta_h$ satisfies the second bound  in \eqref{two-sided.sol},  so that $\delta_h > l_0 (m_3 \kappa_2)^{1/2} h /(l_0 m_3) = (\kappa_2/m_3)^{1/2} h =: r_0$.  Then, there exists some $\eta \in (0,1)$ such that $\wt \bbeta := (1-\eta ) \bbeta^* +\eta  \bbeta^*_h$ satisfies $\| \wt \bbeta - \bbeta^* \|_{\bSigma} = \eta   \delta_h = r_0$.  By the convexity of $\bbeta\mapsto Q_h(\bbeta)$ and Lemma~C.1 in the supplementary material of \cite{SZF2019},
\#
  \langle \nabla Q_h( \wt \bbeta )  - \nabla Q_h(\bbeta^*) , \wt \bbeta  - \bbeta^* \rangle  \leq \eta  \cdot  \langle \nabla Q_h(\bbeta^*_h )  - \nabla Q_h(\bbeta^*) , \bbeta^*_h - \bbeta^* \rangle  =  \langle    - \nabla Q_h(\bbeta^*) , \wt \bbeta  - \bbeta^* \rangle . \nn
\#
 Repeating the above analysis, we find that the right-hand side of the above inequality $\leq 0.5 l_0 \kappa_2 h^2 \cdot r_0$, and the left-hand side $\geq \underbar{$f$} \cdot r_0^2     -   0.5 l_0 m_3 \cdot r_0^3  - l_0 \kappa_1 h \cdot r_0^2= \{ \underbar{$f$} - l_0 \kappa_1 h - 0.5 l_0 (m_3 \kappa_2)^{1/2}  h\} r_0^2$. Canceling  out the common factor $r_0 $ from both sides, we obtain
 \$
 r_0 \leq \frac{0.5 l_0 \kappa_2 h^2 }{ \underbar{$f$} - l_0 \kappa_1 h - 0.5 l_0 (m_3 \kappa_2)^{1/2}} <  \frac{0.5 l_0 \kappa_2 h^2 }{  0.5 l_0 (m_3 \kappa_2)^{1/2} h} = (\kappa_2/m_3)^{1/2}  h = r_0 ,
 \$
 which leads to a contradiction.  Consequently,   $\delta_h$ must satisfy the first bound in \eqref{two-sided.sol},  which in turn implies the claimed result.

Next, to investigate the leading term in the bias, define 
\#
	\bDelta  =  \bSigma^{-1/2}  \big\{   \nabla Q_h(\bbeta^*_h) - \nabla Q_h(\bbeta^*)  - \Jb  ( \bbeta^*_h -\bbeta^* ) \big\} ~~\mbox{ and }~\Hb = \bSigma^{-1/2} \Jb \bSigma^{-1/2} = \EE \{ f_{\varepsilon |\bx}(0) \bw \bw^\T\} , \nn
\#
where $\bw = \bSigma^{-1/2} \bx$.
Again, by the mean value theorem for vector-valued functions,
\#
\bDelta &  =  \Bigg\{    \bSigma^{-1/2}  \int_0^1  \nabla^2 Q_h(  \bbeta^* + t \bdelta_h ) \, {\rm d}t  \, \bSigma^{-1/2}  - \Hb \Bigg\}  \bSigma^{1/2}  \bdelta_h  ~\mbox{ with }~ \bdelta_h = \bbeta^*_h - \bbeta^*  .  \label{residual} 
\#
The Lipschitz continuity of $f_{\varepsilon |\bx}(\cdot)$ ensures that
\#
	&  \Bigg\|  \bSigma^{-1/2}  \int_0^1  \nabla^2 Q_h( \bbeta^* + t \bdelta_h ) \, {\rm d}t  \, \bSigma^{-1/2}  - \Hb \Bigg\|_2  \nn \\
&=\Bigg\|   \EE   \int_0^1  \int_{-\infty}^\infty K(u)  \bigl\{ f_{\varepsilon|\bx }(       t \langle  \bx,  \bdelta_h \rangle - hu ) - f_{ \varepsilon|\bx }( 0) \bigr\} \, {\rm d} u \, {\rm d} t \,   \bw  \bw^\T \Bigg\|_2 \nn \\
& \leq  l_0 \sup_{\bu \in \mathbb{S}^{p-1}} \EE    \int_0^1  \int_{-\infty}^\infty K(u) \bigl(  t |   \langle \bx,  \bdelta_h \rangle | + h |u |   \bigr)  \, {\rm d}u \, {\rm d} t  \, \langle \bw , \bu \rangle^2 \nn \\
& \leq 0.5 l_0  \sup_{\bu \in \mathbb{S}^{p-1}}  \EE \bigl( |\langle \bx,  \bdelta_h \rangle | \langle \bw , \bu \rangle^2 \bigr) +    l_0  \kappa_1 h \nn \\
& \leq  0.5 l_0  m_3 \| \bdelta_h \|_{\bSigma} +  l_0 \kappa_1 h. \nn
\#
This bound, together with \eqref{residual}, implies
\#
	\| \bDelta \|_2 \leq  l_0   \bigl( 0.5 m_3 \| \bdelta_h \|_{\bSigma}   + \kappa_1 h \bigr) \|  \bdelta_h \|_{\bSigma}  . \label{residual.bound}
\#
Moreover, applying a second-order Taylor series expansion to $f_{\varepsilon|\bx}(\cdot)$ yields
\#
 & \EE \bigl\{  \cK(-\varepsilon /h) | \bx \bigr\}  - \tau \nn \\
  & =  \int_{-\infty}^\infty  K(u) \int_{0}^{-h u} \bigl\{  f_{ \varepsilon|\bx } (   t) - f_{\varepsilon|\bx} ( 0) \bigr\} \, {\rm d} t  \, {\rm d}u \nn \\
& =  0.5  \kappa_2  h^2 \cdot  f_{ \varepsilon|\bx}'(0 ) +   \int_{-\infty}^\infty  \int_{0}^{-hu}   \int_0^t   K(u) \bigl\{  f_{ \varepsilon|\bx}'(v) - f_{ \varepsilon|\bx }'( 0 ) \bigr\}   \, {\rm d} v \, {\rm d} t  \, {\rm d}u  . \nn
\#
For $\nabla Q_h(\bbeta^*) = \EE \{ \cK(-\varepsilon /h) - \tau \} \bx $, it follows that
\#
	 \bigg\|  \bSigma^{-1/2}  \nabla Q_h(\bbeta^*)     - \frac{1}{2} \kappa_2 h^2 \cdot \bSigma^{-1/2}  \EE   \bigl\{  f_{\varepsilon|\bx }'( 0)\bx \bigr\}  \bigg\|_2 \leq \frac{	1}{6 } l_1 \kappa_3 h^3 . \label{gradient.approxi}
\#
Combining \eqref{residual.bound} and \eqref{gradient.approxi} completes the proof of \eqref{bias.leading}.  \qed

%%%%%%%%%%%%%%%%%%%%%%%%%%%%%%%%%%%
%%%%%%%%%%%%%%%%%%%%%%%%%%%%%%%%%%%
% Proof of Theorem 3.1 (New Version)
%%%%%%%%%%%%%%%%%%%%%%%%%%%%%%%%%%%
%%%%%%%%%%%%%%%%%%%%%%%%%%%%%%%%%%%
\subsection{Proof of Theorem~\ref{thm:concentration}}
\label{proof:thm3.1}
For every $\bdelta \in \RR^p$, define $\hat D_h(\bdelta) = \hat Q_h(\bbeta^* + \bdelta) - \hat Q_h(\bbeta^*)$, $D_h(\bdelta) = Q_h(\bbeta^* + \bdelta) -  Q_h(\bbeta^*)$,  as well as first-order Taylor series remainder terms $\hat R_h(\bdelta )=\hat D_h(\bdelta) - \langle \nabla  \hat Q_h(\bbeta^*) , \bdelta \rangle$ and $R_h(\bdelta ) = D_h(\bdelta)  -  \langle \nabla  Q_h(\bbeta^*), \bdelta \rangle$.
With these notations, we have
\#
\hat D_h(\bdelta) &= \langle  \nabla   Q_h(\bbeta^*), \bdelta \rangle +  R_h(\bdelta )  + \{ \hat D_h(\bdelta )  - D_h(\bdelta )  \} \nn\\
& \geq R_h(\bdelta ) -  \| \bSigma^{-1/2}  \nabla Q_h(\bbeta^*) \|_2 \cdot \| \bdelta \|_{\bSigma} -\{   D_h(\bdelta )  - \hat  D_h(\bdelta )  \} \nn \\  
& \geq  R_h(\bdelta ) - 0.5 l_0 \kappa_2 h^2 \cdot \| \bdelta \|_{\bSigma} - \underbrace{ \{   D_h(\bdelta )  -  \hat  D_h(\bdelta )  \}   }_{{\rm sampling~error}} , \label{Dh.lbd1}
\#
where we used \eqref{score.mean.ubd} in the last step.
Following the same argument that leads to \eqref{pop.foc.lbd},   it can be  shown that
\#
R_h(\bdelta ) \geq \frac{1}{2} \big(  \underbar{$f$}    - l_0 \kappa_1 h  - 0.5 l_0 m_3 \cdot    \| \bdelta \|_{\bSigma}   \big) \cdot   \| \bdelta \|_{\bSigma}^2 ~\mbox{ for all } \bdelta \in \RR^p. \label{Dh.lbd2}
\#
Take $r_0 = (2\kappa_2 / m_3)^{1/2}  h$ as an intermediate convergence radius. For any $\bdelta \in \partial  \Theta(r_0)$, i.e., $\| \bdelta \|_{\bSigma} = r_0$, the last two displays imply
\#
 R_h(\bdelta ) - 0.5 l_0 \kappa_2 h^2 \cdot \| \bdelta \|_{\bSigma}  \geq   \frac{1}{2}
\big\{ \underbar{$f$}    - l_0 \kappa_1 h - l_0 (2\kappa_2  m_3)^{1/2}    h  \big\}   r_0^2   ~\mbox{ for all } \bdelta \in \partial  \Theta(r_0). \label{Dh.lbd3}
\#

To control the last term on the right-hand side of \eqref{Dh.lbd1},  the following lemma provides some type of uniform law of large numbers for the zero-mean stochastic process  $\{ \hat D_h(\bdelta) - D_h(\bdelta), \bdelta \in \Theta(r) \}$, $r>0$.  This implies a form of restricted
strong convexity (RSC) of the empirical loss $\hat Q_h(\cdot)$.

\begin{lemma} \label{lem:loss.difference}
Given any $r \geq 0$, the bound
\#
   \sup_{\bdelta \in \Theta(r)} \{  D_h(\bdelta )  - \hat  D_h(\bdelta)  \} \leq  3  \bar \tau \upsilon_0  r\cdot  \Bigg( \sqrt{\frac{u}{n}} + \frac{u}{n} \Bigg)    \label{difference.loss.bd1}
\#
holds with probability at least $1-e^{4 p - u}$ for any $u\geq 0$, where $\bar \tau = \max(\tau, 1-\tau)$. In addition, given any $r_u > r_l >0$, with probability at least $1- \lceil e \log (\frac{r_u}{r_l} ) \rceil e^{4 p - u}$ for any $u\geq 0$,
\#
  D_h(\bdelta) - \hat D_h(\bdelta)   \leq  4.25  \bar \tau \upsilon_0    \| \bdelta \|_{\bSigma}      \Bigg( \sqrt{\frac{u}{n}} + \frac{u}{n} \Bigg)  ~\mbox{ holds for all $\bdelta$ satisfying $r_l  \leq \|\bdelta \|_{\bSigma} \leq r_u$} . \label{difference.loss.bd2}
\#
\end{lemma}
 
First, applying \eqref{difference.loss.bd1} with $r=r_0$ and $u= 4p + t$ yields that, with probability at least $1-e^{-t}$,
\#
 \sup_{\bdelta \in \Theta(r_0) } \{   D_h(\bdelta) - \hat D_h(\bdelta) \}  \leq 3  \bar \tau  \upsilon_0  r_0  \Bigg( \sqrt{\frac{4 p + t}{n}} + \frac{4 p + t}{n} \Bigg) .  \label{Dh.unif.bd1}
\#
Combining this bound with \eqref{Dh.lbd1} and \eqref{Dh.lbd3} implies that with probability at least $1-e^{-t}$,  $\hat D_h(\bdelta) >0$ for all $\bdelta \in \partial \Theta(r_0)$ as long as the bandwidth is subject to $\underbar{$f$}^{-1} m_3^{1/2} \upsilon_0 \sqrt{(p+t)/n} \lesssim  h \lesssim \underbar{$f$}m_3^{-1/2}$.  On the other hand, by the optimality of $\hat \bbeta_h$,  $\hat \bdelta := \hat \bbeta_h - \bbeta^*$ satisfies $\hat D_h(\hat \bdelta )\leq 0$. Consequently, the convexity of $\hat Q_h(\cdot)$ ensures that $\| \hat \bdelta \|_{\bSigma} \leq r_0$. See, e.g., Lemma~9.21 in \cite{W2019}.

Next, we refine the convergence rate of $\hat \bbeta_h$ from $r_0$ to the claimed one under the above event.  Consider the ring-shaped set $\Theta(r_l, r_0) = \{ \bdelta \in \RR^p: r_l \leq \| \bdelta \|_{\bSigma} \leq r_0\}$ with $r_l = r_0 h$.  If $\hat \bdelta \notin \Theta(r_l, r_0)$,  we must have $\hat \bdelta \in \Theta(r_0 h)$, and thus the claimed bound follows immediately.  Hereinafter, we assume $\hat \bdelta \in \Theta(r_l, r_0)$. Using the second inequality in Lemma~\ref{lem:loss.difference} with $(r_l,r_u)=(r_0 h, r_0)$ and $u=\sqrt{ \log( e \log h^{-1} ) +4 p + t  }$, we find that with probability at least $1-e^{-t}$,
\#
    D_h(\bdelta) - \hat D_h(\bdelta)    \leq    \| \bdelta \|_{\bSigma} \cdot \underbrace{ 4.25  \bar \tau \upsilon_0 \Bigg\{ \sqrt{\frac{  \log (e\log h^{-1}) + 4 p +  t}{n}}   + \frac{  \log (e\log h^{-1}) + 4 p +  t}{n} \Bigg\} }_{=: r_1 } \label{Dh.unif.bd2}
\#
holds for all $\bdelta\in \Theta(r_l, r_0)$, hence including $\hat \bdelta$.  Applying this along with the earlier bounds \eqref{Dh.lbd1}, \eqref{Dh.lbd2} and the fact $\hat D_h(\hat \bdelta)\leq 0$,   we obtain that 
\$
 (\underbar{$f$} - l_0 \kappa_1 h ) \|  \hat \bdelta \|_{\bSigma}^2 \leq  ( 2 r_1 + l_0 \kappa_2 h^2 )  \|  \hat\bdelta \|_{\bSigma}  + 0.5 l_0 m_3  \| \hat \bdelta \|_{\bSigma}^3 \leq 2  (  r_1 + l_0 \kappa_2 h^2    )  \|  \hat\bdelta \|_{\bSigma}. 
\$
Canceling out $\|  \hat\bdelta \|_{\bSigma}$ proves the claimed bound. \qed

\subsubsection{Proof of Lemma~\ref{lem:loss.difference}}

For each sample $\bz_i= (\bx_i, \varepsilon_i)$,  define the loss difference
$d_h(\bdelta ; \bz_i) = \ell_h(  \varepsilon_i - \langle \bx_i,  \bdelta \rangle ) -  \ell_h( \varepsilon_i   )$, so that $\hat D_h(\bdelta) =  (1/n) \sn d_h(\bdelta ; \bz_i)$.  By the Lipschitz continuity of $u\mapsto \ell_h(u)$,  $d_h(\bdelta ; \bz_i)$ is $\bar \tau$-Lipschitz continuous in $\langle\bx_i, \bdelta\rangle$. That is, for any $\bz_i$ and $\bdelta , \bdelta' \in \RR^p$,  $|d_h(\bdelta ; \bz_i)- d_h(\bdelta' ; \bz_i)| \leq \bar \tau | \langle\bx_i, \bdelta\rangle-\langle\bx_i, \bdelta'\rangle|$.

For any given $r>0$ and some $\epsilon \in (0,1)$ to be determined, define the random variable $\Delta_\epsilon(r)  = n (1-\epsilon) \sup_{\bdelta \in \Theta(r)} \{  D_h(\bdelta) - \hat D_h(\bdelta) \}/(2\bar \tau r)$, where  $D_h(\bdelta) = \EE  \hat D_h(\bdelta)$.  By Chernoff's inequality, for any $u\geq 0$,
\#
	\PP\{ \Delta_\epsilon(r) \geq u \} \leq  \exp \Bigg[  - \sup_{\lambda \geq 0 }   \big\{  \lambda u - \log \EE e^{\lambda  \Delta_\epsilon(r) }  \big\}  \Bigg] . \label{chernoff.bd}
\#
To control the moment generating function $ \EE e^{\lambda  \Delta_\epsilon(r) }$,  by Rademacher symmetrization we have 
\$
\EE e^{\lambda  \Delta_\epsilon(r) } \leq  \EE  \exp\Bigg\{   2 \lambda (1-\epsilon) \sup_{\bdelta \in \Theta(r)}   \frac{1}{2\bar \tau r}   \sn e_i  d_h(\bdelta ; \bz_i)   \Bigg\} , 
\$
where $e_1, \ldots, e_n$ are independent Rademacher random variables.  Recall that $d_h(\bdelta ; \bz_i)$ is $\bar \tau$-Lipschitz continuous in $\langle \bx_i, \bdelta\rangle$,  and $d_h(\bdelta ; \bz_i)=0$ if $\langle \bx_i, \bdelta \rangle = 0$. Applying the Ledoux-Talagrand contraction inequality (see Theorem~4.12 and inequality (4.10) in \cite{LT1991}) yields
\$
& \EE  \exp\Bigg\{   2 \lambda (1-\epsilon) \sup_{\bdelta \in \Theta(r)}  \frac{1}{2 \bar \tau r}   \sn e_i  d_h(\bdelta ; \bz_i)  \Bigg\}  \\
&  \leq  \EE  \exp\Bigg\{  \frac{\lambda }{r} (1-\epsilon) \sup_{\bdelta \in \Theta(r)}     \sn e_i  \langle \bx_i, \bdelta \rangle  \Bigg\}  \leq \EE \exp \Bigg\{  \lambda  (1-\epsilon)\bigg\|  \sn e_i \bw_i \bigg\|_2 \Bigg\} ,
\$
where $\bw_i = \bSigma^{-1/2} \bx_i$.
For this $\epsilon \in (0,1)$, there exists an $\epsilon$-net $\{ \bu_1,  \ldots , \bu_{N_\epsilon} \}$ of $\mathbb{S}^{p-1}$ with cardinality $N_\epsilon \leq (1+2/\epsilon)^p$ such that $\| \sn e_i \bw_i  \|_2 \leq (1-\epsilon)^{-1} \max_{1\leq j\leq N_\epsilon} \sn e_i \bu_j^\T\bw_i $. This implies
\$
 \EE \exp \Bigg\{  \lambda  (1-\epsilon)  \bigg\|  \sn e_i \bw_i \bigg\|_2 \Bigg\} 
  \leq \sum_{j=1}^{N_\epsilon}   \EE \exp \Bigg\{ \lambda  \sn e_i \bu_j^\T\bw_i  \Bigg\} .
\$
Write $S_j = \sn e_i \bu_j^\T\bw_i$, which is a sum of zero-mean random variables.
Note that $e_i \in \{-1, 1\}$ is symmetric, and  Condition~\ref{cond.predictor} ensures that for any $k\geq 3$,  $\EE |\bu_j^\T\bw_i |^k \leq   \upsilon_0^k k \int_0^\infty t^{k-1} e^{-t } {\rm d}t =  \upsilon_0^k k!$. Hence, for every $0\leq c < 1/\upsilon_0$,
\$
& \EE  e^{c e_i \bu_j^\T\bw_i}   = 1 + \frac{c^2}{2} \EE (e_i \bu_j^\T\bw_i )^2 + \sum_{\ell=3}^\infty  \frac{c^\ell}{\ell !} \EE(e_i \bu_j^\T\bw_i)^\ell  \\
&   \leq  1 + \frac{c^2}{2} + \sum_{\ell=2}^\infty  \frac{c^{2\ell}}{(2\ell)!}   \upsilon_0^{2\ell } (2\ell)!  \leq   1 + \frac{c^2}{2} +  \sum_{\ell=2}^\infty  (c^2\upsilon_0^2 )^\ell     \leq 1 + \frac{\upsilon_0^2}{2} \sum_{\ell\geq 2} c^\ell (\sqrt{2}\upsilon_0)^{\ell-2}.
\$
 It follows that for every $0<\lambda  < 1/(\sqrt{2 }\upsilon_0)$ and $j=1,\ldots, N_\epsilon$,
 \$
  \log \EE e^{\lambda  S_j } \leq  \frac{n\upsilon_0^2 \lambda^2 }{2(1- \sqrt{2} \upsilon_0 \lambda)} ~\mbox{ and thus }~ \log  \EE e^{\lambda  \Delta_\epsilon(r) }  \leq  \log N_\epsilon +  \frac{n \upsilon_0^2 \lambda^2 }{2(1- \sqrt{2} \upsilon_0 \lambda)}.
 \$
For any $u\geq 0$,  note that
\$
\sup_{\lambda \geq 0}  \big\{  \lambda u - \log \EE e^{\lambda  \Delta_\epsilon(r) }  \big\} \geq -\log N_\epsilon  + \sup_{\lambda \in (0,  (\sqrt{2}\upsilon_0)^{-1})}  \Bigg\{  \lambda  u   - \frac{n \upsilon_0^2 \lambda^2}{2(1-\sqrt{2}\upsilon_0 \lambda )}\Bigg\}
\$
Substituting this into \eqref{chernoff.bd}, and  following the proof of Bernstein's inequality (see, e.g., Theorem~2.10 in \cite{BLM2013}),  it can be shown that  with probability at least $1-\exp \{ p\log(1+2/\epsilon) - u \}$,
\#
	 \sup_{\bdelta \in \Theta(r)}\{   D_h(\bdelta) - \hat D_h(\bdelta) \}    \leq \frac{2\sqrt{2} }{1-\epsilon} \bar \tau \upsilon_0 r  \cdot \Bigg(      \sqrt{\frac{u}{n}} +  \frac{u}{n}   \Bigg) . \label{Deltar.tail.bound}
\#
This proves \eqref{difference.loss.bd1} by taking $\epsilon=2/(e^4 -1)$.

%Putting together the pieces, we conclude that
%\$
% \EE e^{\lambda  \Delta(r) }  \leq \inf_{ \epsilon \in (0,1)}  \sum_{j=1}^{N_\epsilon}  e^{1.1 (1-\epsilon)^{-2} \upsilon_1^2 \lambda^2  } \leq \inf_{\epsilon\in (0,1)}  \exp\big\{ p\log(1+2/\epsilon) + 1.1 (1-\epsilon)^{-2} \upsilon_0^2 \lambda^2 \big\}.
%\$
%Substituting this bound into \eqref{chernoff.bd}, and with a properly chosen $\epsilon$,  we obtain that
%\#
%	\PP \{ \Delta(r) \geq u \} \leq \inf_{\lambda \geq 0 }   \exp (  2\upsilon_1^2 \lambda^2 + 2.2 p - \lambda t  )   = \exp\big\{ 2.2 p - u^2/ (8\upsilon_0^2) \big\}.  \label{Deltar.tail.bound}
%\#

Next, we prove the uniform bound \eqref{difference.loss.bd2}, which holds for all $\bdelta\in  \Theta(r_l, r_u):=\{ \bv \in \RR^p: r_l \leq \|\bv \|_{\bSigma} \leq r_u \}$, via a peeling argument. For some $\gamma>1$ (to be specified) and positive integers $k =1 ,\ldots, N := \lceil \log (\frac{r_u}{r_l} )/\log(\gamma) \rceil$, define the sets $\Theta_k=\{ \bv \in \RR^p:  \gamma^{k-1} r_l \leq \|\bv \|_{\bSigma} \leq \gamma^{k } r_l \}$, so that $\Theta(r_l, r_u) \subseteq \cup_{k=1}^N \Theta_k$.
Then,
\$
 & 	\PP\Bigg\{   \exists \bdelta \in \Theta(r_l, r_u) \mbox{ s.t. }    D_h(\bdelta) - \hat D_h(\bdelta)   >  \frac{2\sqrt{2}\gamma}{1-\epsilon} \bar \tau \upsilon_0   \| \bdelta \|_{\bSigma} \cdot   \Bigg(      \sqrt{\frac{u}{n}} +  \frac{u}{n}   \Bigg)   \Bigg\} \\
& \leq \sum_{k=1}^N  	\PP\Bigg\{   \exists \bdelta \in \Theta_k  \mbox{ s.t. }    D_h(\bdelta) - \hat D_h(\bdelta)   >\frac{2\sqrt{2}\gamma }{1-\epsilon}  \bar \tau  \upsilon_0  \gamma^{k-1} r_l   \cdot \Bigg(      \sqrt{\frac{u}{n}} +  \frac{u}{n}   \Bigg)   \Bigg\} \\
& \leq   \sum_{k=1}^N  	\PP\Bigg\{  \sup_{\bdelta  \in \Theta(\gamma^k r_l) }    D_h(\bdelta) - \hat D_h(\bdelta)    >  \frac{2\sqrt{2}}{1-\epsilon} \bar \tau  \upsilon_0 \gamma^{k } r_l \cdot  \Bigg(      \sqrt{\frac{u}{n}} +  \frac{u}{n}   \Bigg)    \Bigg\} \\
& \stackrel{({\rm i})}{\leq } \sum_{k=1}^N   \exp\big\{ p \log(1+2/\epsilon) - u \big\} \leq  \lceil \log (\tfrac{r_u}{r_l})/\log(\gamma) \rceil  \exp\big\{ p \log(1+2/\epsilon) - u  \big\} ,
\$
where inequality (i) is obtained by repeatedly using \eqref{Deltar.tail.bound} with $r= \gamma^{k } r_l$ for $k=1,\ldots, N$.  Taking $\epsilon=2/(e^4-1)$ and $\gamma = e^{1/e}$ yields the claimed bound \eqref{difference.loss.bd2}.  \qed

\subsection{An alternative proof to Theorem~\ref{thm:concentration}}
\label{proof:thm3.1.refined}

From the proof of Theorem~\ref{thm:concentration} in Section~\ref{proof:thm3.1},  we see that the use of peeling argument will create an additional term $\log_2(1/h)$ in the upper bound, although it is further bounded by $\log(\log n)$ (under the prescribed constraint on $h$), a very slowly growing function of $n$. In this section, we argue that this extra term is an artifact of the proof technique, and can be avoided through a more careful analysis regarding the (local) restricted strong convexity of the empirical loss $\hat Q_h(\cdot)$.

The key is to refine the convergence rate of $\hat \bbeta_h$ from $r_0 \asymp  h$ to the claimed one, conditioned on $\hat \bdelta = \hat \bbeta_h - \bbeta^* \in \Theta(r_0)$.
The first-order optimality condition implies $\nabla \hat Q_h(\hat \bbeta_h ) = 0$, and hence 
\#
 & \langle \nabla \hat Q_h(\hat \bbeta_h )  - \nabla \hat Q_h(\bbeta^* ) , \hat \bdelta \rangle \nn \\
&= \langle - \nabla \hat Q_h(\bbeta^* ) , \hat \bdelta \rangle \leq   \Big(  \| \bSigma^{-1/2} \{ \nabla \hat Q_h(\bbeta^* )  -  \nabla  Q_h(\bbeta^* )  \} \|_2 + \underbrace{ \|  \bSigma^{-1/2}\nabla  Q_h(\bbeta^* )  \|_2  }_{\leq 0.5 l_0 \kappa_2 \cdot h^2 } \Big) \| \hat \bdelta \|_{\bSigma}  . \label{new.foc}
\# 
Define the symmetrized Bregman divergence associated with the convex function $\hat Q_h(\cdot)$:
\#  
	D  (\bbeta_1, \bbeta_2) = \langle \nabla \hat Q_h  (\bbeta_1)-  \nabla \hat Q _h (\bbeta_2), \bbeta_1 - \bbeta_2 \rangle \geq 0  ,~~\bbeta_1, \bbeta_2  \in \RR^p. \label{Bregman.div}
\#
Then the left-hand side of \eqref{new.foc} reads $D  (\bbeta^* + \hat \bdelta , \bbeta^*)$. Starting from \eqref{new.foc}, we need to bound $\| \bSigma^{-1/2} \{ \nabla \hat Q_h(\bbeta^* )  -  \nabla  Q_h(\bbeta^* )  \} \|_2$ from above, and derive a lower bound for $D  (\bbeta^* +   \bdelta , \bbeta^*)$ uniformly over $\bdelta \in \RR^p$ in a local neighborhood of the origin. The following two lemmas serve for this purpose.

%%%%%%%%%%%%%%%%%%%%%%%%%%%%%%%%%%%%%%%%%%%%
%%%%%%%%%%%%%%%%%%%%%%%%%%%%%%%%%%%%%%%%%%%%
% Lemma Score
%%%%%%%%%%%%%%%%%%%%%%%%%%%%%%%%%%%%%%%%%%%%
%%%%%%%%%%%%%%%%%%%%%%%%%%%%%%%%%%%%%%%%%%%%
\begin{lemma} \label{lem:score}
Assume Conditions~\ref{cond.kernel}--\ref{cond.predictor} hold. For any $t\geq 0$,
\#
 \|   \bSigma^{-1/2}  \{ \nabla  \hat Q_h  (\bbeta^* )  - \nabla   Q_h  (\bbeta^* )\}   \|_2  \leq  1.46 \upsilon_0  \Biggl(     C_\tau \sqrt{   \frac{  4p +2 t }{n}} +  \bar \tau \frac{ 4p+2t }{n} \Biggr)       \label{score.ubd}
\#
with probability at least $1-e^{-t}$, where $C_\tau^2 = \tau(1-\tau) + (1+\tau) l_0 \kappa_2 h^2$.
\end{lemma}

In addition to Condition~\ref{cond.reg}, assume $f_{\varepsilon |\bx}(0) \leq \bar f $ for some $\bar f\geq \underbar{$f$} >0$. Then, for all $0< h \leq  \underbar{$f$}/(4l_0)$,
\#
	  \frac{7}{8} \underbar{$f$} \leq \inf_{|u|\leq h/2} f_{\varepsilon | \bx} (u)  \leq \sup_{|u|\leq h/2} f_{\varepsilon | \bx} (u) \leq \frac{9}{8} \bar{f} , \label{density.local.bound}
\#
almost surely (over $\bx$). Moreover, for $\delta \in (0,1]$, define $\iota_\delta \geq 0$ as
\#
	\iota_\delta = \inf\bigl\{  \iota >0 :  \EE \bigl\{     \langle \bu , \bw \rangle^2 \mathbbm{1}\big( |\langle \bu , \bw \rangle| > \iota  \big)  \bigr\} \leq \delta ~\mbox{ for all } \bu \in \mathbb{S}^{p-1}  \bigr \} , \label{def.eta}
\#
where $\bw = \bSigma^{-1/2} \bx$ is the standardized covariate vector satisfying $\EE(\bw \bw^\T) = \Ib_p$, and hence $\EE \langle \bu , \bw \rangle^2 = 1$ for any $\bu \in \mathbb{S}^{p-1} $.
It can be shown that $\iota_\delta$ depends only on $\delta$ and $\upsilon_0$ in Condition~\ref{cond.predictor}, and the map $\delta \mapsto \iota_\delta$ is non-increasing with $\iota_\delta \downarrow 0$ as $\delta \uparrow 1$ and $\iota_1=0$.  By Markov's inequality,  for any $\iota>0$ it holds $\sup_{\bu \in \mathbb{S}^{p-1}}\EE  \{     \langle \bu , \bw \rangle^2 \mathbbm{1}  ( |\langle \bu , \bw \rangle| > \iota   )  \} \leq \iota^{-2} m_4$. 
Hence, a rather crude upper bound for $\iota_\delta$ is $\iota_\delta \leq (m_4/\delta)^{1/2}$.

%%%%%%%%%%%%%%%%%%%%%%%%%%%%%%%%%%%%%%%%%%%%
%%%%%%%%%%%%%%%%%%%%%%%%%%%%%%%%%%%%%%%%%%%%
% Lemma RSC
%%%%%%%%%%%%%%%%%%%%%%%%%%%%%%%%%%%%%%%%%%%%
%%%%%%%%%%%%%%%%%%%%%%%%%%%%%%%%%%%%%%%%%%%%
\begin{lemma} 
\label{lem:rsc}
Assume the kernel $K(\cdot)$ is such that $\kappa_l = \min_{|u|\leq 1} K(u)>0$, and let the bandwidth satisfy $0 < h \leq   \underbar{$f$} / (4l_0 )$. Given any $0<r \leq  h/(4 \iota_{1/4})$ with $\iota_{1/4}$ defined in \eqref{def.eta}, 
\#
	\inf_{\bdelta  \in  \Theta(r) }   \frac{  D(  \bbeta^*+ \bdelta , \bbeta^* )  }{ \kappa_l  \| \bdelta \|_{\bSigma}^2 }  \geq 
	 c_0\underbar{$f$}   -  c_1 \Bigg\{ \bar f^{1/2} \frac{ 1 }{r} \sqrt{\frac{ph}{n}} +  (\bar f m_4)^{1/2}\sqrt{\frac{t}{nh}}  +   \frac{  h t}{r^2 n} \Bigg\}   \label{rsc.ineq}
\#
with probability at least $1-e^{-t}$ for any $t\geq 0$, where $c_1 >0$ is an absolute constant, and $c_0=21/32$.
\end{lemma}

Let $\cG(t)$ be the ``good" event that the bounds \eqref{score.ubd} and \eqref{rsc.ineq} are satisfied.
Together, Lemma~\ref{lem:score}, Lemma~\ref{lem:rsc} and \eqref{new.foc} imply that, conditioned on $\{ \hat \bdelta \in \Theta(r_0)\} \cap \cG(t)$ with $r_0 \leq   h/(4\iota_{1/4})$,
\$
c_0 \kappa_l  \underbar{$f$} \cdot \| \hat \bdelta \|_{\bSigma}^2 &  \leq   \left\{ 1.46 \upsilon_0  \Biggl(     C_\tau \sqrt{   \frac{  4p +2 t }{n}} +  \bar \tau \frac{ 4p+2t }{n} \Biggr)      + 0.5 l_0\kappa_2 h^2  \right\}\| \hat \bdelta \|_{\bSigma} \nn \\
& ~~~~~~+c_1 \kappa_l  \left\{ \bar f^{1/2} \frac{ 1 }{ r_0 } \sqrt{\frac{ph}{n}} +  (\bar f m_4)^{1/2}\sqrt{\frac{t}{nh}}  +   \frac{  h t}{ r_0^2  n} \right\} r_0 \cdot \| \hat \bdelta \|_{\bSigma}  .
\$
Consequently, we obtain the bound
\$
c_0\underbar{$f$}   \cdot \| \hat \bdelta \|_{\bSigma} &  \leq   1.46  \kappa_l^{-1} \upsilon_0  \Biggl(     C_\tau \sqrt{   \frac{  4p +2 t }{n}} +  \bar \tau \frac{ 4p+2t }{n} \Biggr)      + 0.5 \kappa_l^{-1}  l_0\kappa_2 h^2 \nn \\
& ~~~~~~+  c_1   \Bigg\{ \bar f^{1/2}   \sqrt{\frac{ph}{n}} +   (\bar f m_4)^{1/2}  r_0 \sqrt{\frac{   t}{n h  }}  +     \frac{  h  t}{ r_0  n} \Bigg\}  
\$
without having the additional $\log_2(1/h)$ term. For example, we may take $r_0 = h/(8 m_4^{1/2})$.\qed

%%%%%%%%%%%%%%%%%%%%%%%%%%%%%%%%%%
%%%%%%%%%%%%%%%%%%%%%%%%%%%%%%%%%%
% Proof of Lemma 5.1 
%%%%%%%%%%%%%%%%%%%%%%%%%%%%%%%%%%
%%%%%%%%%%%%%%%%%%%%%%%%%%%%%%%%%%

\subsubsection{Proof of Lemma~\ref{lem:score}}

Define $\xi_i = \cK(-\varepsilon_i /h) - \tau$ for $i=1,\ldots, n$, so that $\bSigma^{-1/2}\{  \nabla \hat Q_h (\bbeta^*) - \nabla Q_h (\bbeta^*)  \}= (1/n) \sn \{\xi_i \bw_i -  \EE(\xi_i \bw_i)\} \in \RR^p$, where $\bw_i = \bSigma^{-1/2} \bx_i$. Using a covering argument, for any $\epsilon\in (0,1)$, there exists an $\epsilon$-net $\cN_\epsilon$ of the unit sphere with cardinality $|\cN_{\epsilon}|\le (1+2/\epsilon)^{p}$ such that 
$$
 \| \bSigma^{-1/2} \{ \nabla \hat Q_h (\bbeta^*) - \nabla Q_h (\bbeta^*)  \}   \|_2 \leq (1-\epsilon)^{-1} \max_{\bu\in \cN_{\epsilon}} \bigl \langle \bu,   \bSigma^{-1/2} \{ \nabla \hat Q_h (\bbeta^*) - \nabla Q_h (\bbeta^*)  \}  \bigr\rangle.
$$
For each  unit vector $\bu\in\cN_{\epsilon}$, define centered random variables $\gamma_{\bu,i} = \langle \bu,   \xi_i \bw_i -\EE(\xi_i \bw_i ) \rangle$. It can be shown that $|\xi_i| \leq \bar \tau := \max(1-\tau,\tau)$ and $\EE(\xi_i^2 | \bx_i ) \leq C_\tau  = \tau(1-\tau) + (1+\tau) l_0 \kappa_2 h^2$ (see the proof of Theorem~\ref{thm:clt} below). 
Hence, for $k=2,3,\ldots$,
\#
\EE \bigl( | \langle \bu,   \xi_i \bw_i \rangle|^k\bigr) &\le  \bar \tau^{k-2}\, \EE  \bigl\{  | \langle \bu,  \bw_i \rangle|^k \cdot\EE(\xi_i^2 | \bx_i ) \bigr\}    \nn \\
&\le  C_\tau^2   \bar \tau^{k-2}    \upsilon^k_0 \int_0^{\infty}\PP (| \langle \bu, \bw_i\rangle|\ge \upsilon_0 t) kt^{k-1}\,\mathrm{d}t \nn\\
&\le C_\tau^2    \bar \tau^{k-2}       \upsilon^k_0  k \int_0^{\infty} t^{k-1}e^{-t}\, \mathrm{d}t\nn\\
&=    k! \cdot    C_\tau^2  \bar \tau^{k-2}  \cdot   \upsilon^k_0 \nn \\
&\le \frac{k!}{2} \cdot    ( C_\tau \upsilon_0)^2     \cdot  \left( 2 \bar \tau  \upsilon_0\right)^{k-2}. \nn
\#
Consequently, it follows from Bernstein's inequality that for every $u\geq 0$,
\#
 \frac{1}{n} \sn \gamma_{\bu,i}  \leq   \upsilon_0 \Biggl(  C_\tau  \sqrt{  \frac{2 u}{n}} +  \frac{ 2 \bar \tau  u}{n} \Biggr) \nn
\#
with probability at least $1-e^{-u}$.  

Finally, applying a union bound over  $\bu \in \cN_{\epsilon}$ yields
\#
 \| \bSigma^{-1/2} \{  \nabla \hat Q_h (\bbeta^*) - \nabla Q_h (\bbeta^*)   \} \|_2  \leq \frac{\upsilon_0}{1-\epsilon}\Biggl(   C_\tau \sqrt{  \frac{2 u}{n}} +  \frac{ 2 \bar \tau  u}{n}\Biggr)  \nn
\#
with probability at least $1- e^{\log(1+2/\epsilon) p -u}$. 
Taking $\epsilon = 2/(e^2-1)$  and $u=  2p+ t$ ($t\geq 0$) proves the claimed result. \qed

\subsubsection{Proof of Lemma~\ref{lem:rsc}}

Recall that the empirical loss $ \hat Q_h(\cdot)$ in \eqref{convolution.loss} is convex and twice continuously differentiable with $\nabla  \hat Q_h(\bbeta) =  (1/n) \sn \{ \cK_h   (  \langle \bx_i, \bbeta \rangle - y_i )  - \tau \} \bx_i$ and $\nabla^2 \hat Q_h(\bbeta) = (1/n) \sn K_h(\langle \bx_i, \bbeta \rangle - y_i) \bx_i \bx_i^\T$. For the symmetrized Bregman divergence $D: \RR^p \times \RR^p \to [0,\infty)$ defined in \eqref{Bregman.div}, we have
\#
   D(\bbeta^*+\bdelta, \bbeta^* )   =  \frac{1}{n} \sn\left\{ \cK \left(  \frac{\langle \bw_i, \bv \rangle - \varepsilon_i}{h}\right) - \cK\left( \frac{    -\varepsilon_i}{h} \right) \right\} \langle \bw_i, \bv \rangle ,\label{def.D}
\#
where $\bw_i = \bSigma^{-1/2} \bx_i$ and $\bv = \bSigma^{ 1/2}\bdelta$.  
Define the events  $\cE_i =  \{ |   \varepsilon_i | \leq h/2  \} \cap \{  |\langle \bw_i,  \bv \rangle |    \leq     \| \bv\|_2 \cdot h/(2r )   \}$ for $i=1,\ldots, n$. For any $\bv \in \BB^p(r)$, note that $|\varepsilon_i - \langle \bw_i, \bv  \rangle| \leq h$ on $\cE_i$, implying
\#
 D(\bbeta^*+ \bdelta , \bbeta^*)  \geq  \frac{\kappa_l}{n h} \sn \langle \bw_i ,  \bv \rangle^2 \mathbbm{1}_{ \cE_i } , \label{eq:dblower1}
\#
where $ \mathbbm{1}_{\cE_i }$ is the indicator function of $\cE_i $ and $\kappa_l = \min_{|u|\leq 1} K(u)$. It then suffices to bound the right-hand side of the above inequality from below uniformly over $\bv \in \BB^p(r)$. 

For $R>0$, define the function  $\varphi_R(u)=u^2 \mathbbm{1}(|u|\leq R/2) + \{ u\sign(u) -R \}^2 \mathbbm{1}( R/2 < |u| \leq R)$, which is $R$-Lipschitz continuous and satisfies
\#
  u^2 \mathbbm{1}(|u| \leq R / 2)	\leq \varphi_R(u) \leq  u^2  \mathbbm{1}(|u| \leq R ) . \label{phi.bound} 
\#
Moreover,  note that $\varphi_{cR}(cu) = c^2 \varphi_{R}(u)$ for any $c > 0$ and $\varphi_0(u) = 0$.
Hence,
\$
\langle \bw_i ,  \bv \rangle^2 \mathbbm{1}_{ \cE_i } \geq \varphi_{\| \bv \|_2 h/(2r)} (\| \bv \|_2 \langle  \bw_i ,  \bv /\| \bv \|_2\rangle )   \cdot \omega_i  = \| \bv \|_2^2 \cdot   \varphi_{  h/(2r)} (   \langle  \bw_i ,  \bv /\| \bv \|_2\rangle )   \cdot \omega_i,
\$
where $\omega_i :=  \mathbbm{1}( | \varepsilon_i | \leq h/2 )$. By a change of variable, the problem is reduced to bounding
\#
	D_0( \bv) :=   \frac{1}{n h  } \sn \omega_i  \cdot  \varphi_{  h/(2r)  } (\langle \bw_i, \bv \rangle )   \label{def.D0}
\# 
from below uniformly over $\bv \in \mathbb{S}^{p-1}$.

In the following, we bound the expectation $\EE \{   D_0(\bv ) \}$ and the random fluctuation $  D_0(\bv  ) - \EE   \{ D_0(\bv )\} $, separately, starting with the former. By \eqref{density.local.bound},
\#
	\frac{7}{8}\underbar{$f$}  h \leq  \EE (   \omega_i | \bx_i ) =  \int_{-h/2}^{h/2} f_{\varepsilon_i | \bx_i} (u) \, {\rm d}u   \leq \frac{9}{8}\bar f  h.   \label{indicator.bound}
\#
Moreover, define $\xi_{\bv} = \langle \bw, \bv \rangle  $ such that $\EE (\xi_{\bv}^2) =1$. 
By~\eqref{phi.bound} and \eqref{indicator.bound}, 
\#
 &	\EE \bigl\{ \omega_i  \cdot \varphi_{    h/(2r)} (\langle \bw_i , \bv \rangle )  \bigr\} 
  \geq   \frac{7}{8}\underbar{$f$}   h \cdot \EE \varphi_{  h/(2r)   } (\langle \bw_i , \bv \rangle )    
  \geq  \frac{7}{8}\underbar{$f$}   h    \cdot   \left[ 1 -  \EE  \xi_{\bv}^2  \mathbbm{1}  \{ | \xi_{\bv} | >  h/(4r) \}  \right]  . \nn
\#
With $r \leq   h/(4 \iota_{1/4})$ and $\iota_{1/4}$ defined in \eqref{def.eta}, it follows that
\#
	\inf_{ \bv  \in \mathbb{S}^{p-1} }  \EE  \{D_0(\bv) \}  \geq   \frac{7}{8}\underbar{$f$}    \cdot    \left\{ 1 - \sup_{\bu \in \mathbb{S}^{p-1}} \EE  \langle \bw , \bu\rangle^2  \mathbbm{1} ( |\langle \bw , \bu\rangle | \geq \iota_{1/4} )   \right\}   \geq   \frac{21}{32}  \underbar{$f$}  . \label{D0.mean.lbd1}
\#

Turning to the random fluctuation, we will use Theorem~7.3 in \cite{B2003}  (a refined Talagrand's inequality) to bound
\#
	\Delta  = \sup_{  \bv \in\mathbb{S}^{p-1}  }   \bigl\{  D^-_0(\bv) - \EE D^-_0(\bv)  \bigr\} , \label{Deltar.def}
\#
where $D^-_0(\bv) := -D_0(\bv) $.
Note that $0\leq \varphi_R(u) \leq (R/2)^2$ for all $u\in \RR$ and $\omega_i \in \{0, 1\}$. Hence,  $ \chi_i :=  (\omega_i/h) \cdot  \varphi_{  h/(2r) } (\langle \bw_i, \bv  \rangle )    \geq 0$ is bounded by $h/(4r)^2 $.  Moreover, it follows from  \eqref{indicator.bound} that $\EE  (\chi_i^2) \leq   9\bar f  m_4 / (8h) $.
We then apply Theorem~7.3 in \cite{B2003} and obtain  that, for any $t>0$,
\#
	\Delta  & \leq \EE (\Delta)  +   \{\EE( \Delta )\}^{1/2}   \frac{1}{2 r}\sqrt{\frac{ h t}{   n }} +  \frac{3}{2}  (\bar f   m_4 )^{1/2}   \sqrt{    \frac{  t  }{   n  h  }} +  \frac{h}{(4 r)^2} \frac{  t}{3  n  }      \label{Deltar.concentration}
\#
with probability at least $1-e^{-t}$.

It remains to bound the expected value $\EE (\Delta) $. Recall that $\omega_i   = \mathbbm{1} ( | \varepsilon_i | \leq h/2 ) \in \{0, 1\}$, and hence
$ \omega_i  \varphi_{ h/(2r) } (\langle \bw_i, \bv  \rangle ) = \omega_i^2  \varphi_{ h/(2r) } (\langle \bw_i, \bv  \rangle ) =  \varphi_{   \omega_i h/(2r)  } (\langle \omega_i\bw_i  , \bv    \rangle ) = \varphi_{   h/(2r)  } (\langle \omega_i\bw_i  , \bv    \rangle ) $.
Define 
\#
	 \bar \bw_i =\omega_i \bw_i ~~\mbox{ and }~~ \cE(\bv;  \bar \bw_i  )  =   \varphi_{   h/(2r)  } (\langle  \bar \bw_i  , \bv    \rangle )  , \ \ \bv \in \mathbb{S}^{p-1} . \nn
\#
By Rademacher symmetrization,
\#
 \EE (\Delta)  \leq 2 \EE \left\{  \sup_{ \bv  \in  \mathbb{S}^{p-1}  }    \frac{1}{n h} \sn  e_i  \cdot  \cE(\bv;  \bar \bw_i  )    \right\}  , \nn
\# 
where $e_1,\ldots, e_n$ are independent Rademacher random variables.
Since $\varphi_R(\cdot)$ is $R$-Lipschitz,  $\cE(\bv; \bz_i)$ is an $(h/2r)$-Lipschitz function in $\langle  \bar \bw_i  , \bv  \rangle$, i.e.,  for any $ \bar \bw_i $ and parameters $\bv , \bv'\in \mathbb{S}^{p-1}$, 
\#
	\bigl| \cE(\bv;  \bar \bw_i ) - \cE(\bv';  \bar \bw_i ) \bigr|  \leq  \frac{h}{2r}  \bigl| \langle   \bar \bw_i   , \bv  \rangle - \langle  \bar \bw_i  , \bv' \rangle \bigr| \label{eq:lips}.
\#
Moreover, observe that $\cE(\bv;   \bar \bw_i  ) = 0$ for any $\bv$ such that $ \langle  \bar \bw_i   , \bv  \rangle = 0$. With the above preparations, we are   ready to use Talagrand's contraction principle to bound $\EE( \Delta)$.
 Define the subset $T\subseteq \RR^n$ as
$$
	T = \bigl\{ \bt = (t_1, \ldots , t_n)^\T : t_i = \langle  \bar \bw_i  , \bv    \rangle,  i=1,\ldots, n, \, \bv \in \mathbb{S}^{p-1} \bigr\},
$$
and contractions $\phi_i:\RR \to \RR$ as $\phi_i(t) = (2r/h) \cdot   \varphi_{ h /(2r)} (t)$. By \eqref{eq:lips}, $|\phi(t)-\phi(s)|\leq |t-s|$ for all $t, s \in \RR$. 
Applying Talagrand's contraction principle (see, e.g., Theorem~4.12 and (4.20) in \cite{LT1991}), we have
\#
 \EE (\Delta)   & \leq 2 \EE \left\{  \sup_{ \bv  \in \mathbb{S}^{p-1} }    \frac{1}{n h } \sn  e_i   \cdot   \cE(\bv;   \bar \bw_i  )    \right\}  =   \EE \left\{ \sup_{\bt \in T} \frac{1}{nr } \sn e_i \phi_i(t_i) \right\}   \leq   \EE \left( \sup_{\bt \in T} \frac{1}{n r} \sn e_i  t_i  \right)  \nn \\
&= \EE \left(  \sup_{ \bv  \in \mathbb{S}^{p-1}  }    \frac{1}{n r} \sn  e_i      \langle   \bar \bw_i  , \bv  \rangle  \right)   \leq     \EE  \Bigg\|  \frac{1}{n r } \sn e_i  \bar \bw_i   \Bigg\|_2 \leq      \bar f^{1/2} \frac{3}{2 r } \sqrt{   \frac{  p h }{2n } }. \nn
\#
This, combined with \eqref{Deltar.def} and \eqref{Deltar.concentration}, yields that
\#
	\Delta  \leq    \frac{3}{2} \bar f^{1/2}      \Bigg(  1.25  \frac{1}{r}\sqrt{\frac{  p h }{ 2n }} +   m_4^{1/2} \sqrt{\frac{   t}{  n h }} \, \Bigg)      + (1+1/48)  \frac{  h  t}{ r^2  n }  \nn
\#
with probability at least $1-e^{-t}$. Combining this with \eqref{def.D}, \eqref{def.D0} and \eqref{D0.mean.lbd1} proves \eqref{rsc.ineq}.  \qed

%%%%%%%%%%%%%%%%%%%%%%%%%%%%%%%%%
%%%%%%%%%%%%%%%%%%%%%%%%%%%%%%%%%
% Proof of Theorem 3.2 
%%%%%%%%%%%%%%%%%%%%%%%%%%%%%%%%%
%%%%%%%%%%%%%%%%%%%%%%%%%%%%%%%%%
\subsection{Proof of Theorem~\ref{thm:bahadur}}
We keep the notation used in the proof of Theorem~\ref{thm:concentration}, and for any $t\geq 0$, let $r=r(n,p,t) \asymp \sqrt{(p+t)/n} + h^2 >0$ be such that $\PP\{\hat \bbeta_h \in \bbeta^* + \Theta(r) \} \geq 1-2e^{-t}$, provided $\sqrt{(p+t)/n} \lesssim h \lesssim 1$.
Define the vector-valued random process
\#
	\Delta(\bdelta) =   \bSigma^{-1/2}  \bigl\{   \nabla \hat Q_h(\bbeta^*+ \bdelta) - \nabla \hat Q_h(\bbeta^*)  - \Jb_h \bdelta \bigr\}   , \ \ \bdelta \in \RR^p, \label{br.remainder.def}
\#
where $\Jb_h =  \nabla^2 Q_h(\bbeta^*)$ is the population Hessian at $\bbeta^*$.  Since $\hat \bbeta_h$ falls in a local neighborhood of $\bbeta^*$ with high probability, it suffices to bound the local fluctuation $\sup_{\bdelta \in   \Theta(r) } \| \Delta(\bdelta) \|_2$.   By the triangle inequality, 
\#
\sup_{\bdelta \in  \Theta(r) }  \| \Delta(\bdelta)   \|_2  \le  \sup_{\bdelta \in  \Theta(r) }   \|  \EE \Delta(\bdelta)  \|_2+ \sup_{\bdelta \in \Theta(r) }  \| \Delta(\bdelta) - \EE \Delta(\bdelta )    \|_2 :=I_1+I_2. 
\label{bahaduri1i2}
\#
We now provide upper bounds for $I_1$ and $I_2$, respectively.

%%%%%%%%%%%%%%%%%%%%%%%%%%%%%%%%%%%%%%%%%%
%%%%%%%%%%%%%%%%%%%%%%%%%%%%%%%%%%%%%%%%%%
% Upper bound for I1
%%%%%%%%%%%%%%%%%%%%%%%%%%%%%%%%%%%%%%%%%%
%%%%%%%%%%%%%%%%%%%%%%%%%%%%%%%%%%%%%%%%%%
\noindent
\textsc{Upper bound for $I_1$:}
By the mean value theorem for vector-valued functions,
\#
	\EE  \Delta(\bdelta)  & =   \bSigma^{-1/2} \biggl\langle   \int_0^1 \nabla^2  Q_h\bigl(   \bbeta^* + t \bdelta  \bigr) {\rm d}t ,  \bdelta  \biggr\rangle - \bSigma^{-1/2} \Jb_h \bdelta \nn \\
	& =\biggl\langle  \bSigma^{-1/2} \int_0^1 \nabla^2  Q_h\bigl(   \bbeta^* + t  \bdelta  \bigr) {\rm d}t \, \bSigma^{-1/2} - \Hb_h   ,   \bSigma^{1/2}  \bdelta \biggr\rangle , \nn
\#
where $ \Hb_h := \bSigma^{-1/2} \Jb_h \bSigma^{-1/2} = \EE \{ K_h(\varepsilon) \bw \bw^\T\}$.
By law of iterative expectation and a change of variable, 
\#
   \bSigma^{-1/2}  \nabla^2   Q_h(\bbeta^* + t  \bdelta  )  \bSigma^{-1/2}  & = \EE \bigl\{ K_h(  t \langle \bx, \bdelta \rangle - \varepsilon ) \bw  \bw^\T \bigr\}   =   \EE   \Biggl\{  \int_{-\infty}^\infty K(u)  f_{\varepsilon | \bx}(t \langle \bx  , \bdelta \rangle  - h u ) \, {\rm d} u \cdot   \bw \bw^\T \Biggr\}. \nn
\#
Write $\bv = \bSigma^{ 1/2} \bdelta$ for $\bdelta \in  \Theta(r)$, so that $\| \bv \|_2 \leq r$ and
$$
	  \bSigma^{-1/2} \nabla^2  Q_h \bigl(    \bbeta^* + t\bdelta \bigr)  \bSigma^{-1/2}  =   \EE \Biggl\{ \int_{-\infty}^\infty K(u)  f_{\varepsilon | \bx}( t \langle \bw  ,\bv  \rangle - h u ) \, {\rm d} u \cdot   \bw\bw^\T  \Biggr\} . 
$$	
By the Lipschitz continuity of $f_{\varepsilon |\bx}(\cdot)$,  
\$
	& \bigl\|  \bSigma^{-1/2}  \nabla^2  Q_h\bigl(    \bbeta^* + t \bdelta \bigr)  \bSigma^{-1/2}    -  \Hb_h   \bigr\|_2 \nn \\
	&  = \Biggl\|     \EE  \int K(u) \bigl\{ f_{\varepsilon|\bx } ( t \langle \bw , \bv \rangle   - h u)    - f_{\varepsilon|\bx } ( -hu ) \bigr\} \, {\rm d} u  \, \cdot \bw  \bw^\T   \Biggr\|_2  \nn \\
%	&  \leq \Biggl\|   \Jb^{-1/2}   \EE  \biggl[ \bigl( l_0 t \langle \bx , \bdelta \rangle   - l_0 \kappa_1 h     \bigr)  \, \bx  \bx^\T  \biggr]  \Jb^{-1/2}   \Biggr\|_2  \nn \\
	   & \leq l_0 t \sup_{\| \bu \|_2 =1} \EE \bigl( \langle \bw, \bu \rangle^2 |\langle \bw, \bv \rangle | \bigr)   \leq    l_0 m_3 r t  ,
\$
where the third inequality holds by the Cauchy-Schwarz inequality.
Consequently, 
\#
	\sup_{\bdelta \in  \Theta(r) }   \| \EE  \Delta ( \bdelta )   \|_2 \leq  
    0.5  l_0 m_3  r^2 .  \label{mean.Delta.ubd}
\#

%%%%%%%%%%%%%%%%%%%%%%%%%%%%%%%%%%%%%%%%%%
%%%%%%%%%%%%%%%%%%%%%%%%%%%%%%%%%%%%%%%%%%
% Upper bound for I2
%%%%%%%%%%%%%%%%%%%%%%%%%%%%%%%%%%%%%%%%%%
%%%%%%%%%%%%%%%%%%%%%%%%%%%%%%%%%%%%%%%%%%
\noindent
\textsc{Upper bound for $I_2$:}
Next, we provide an upper bound for the supremum of the zero-mean stochastic process $\Delta(\bdelta) - \EE  \Delta(\bdelta)$ under $\ell_2$-norm. Define the centered gradient process $G(\bbeta) = \bSigma^{-1/2} \{ \nabla \hat Q_h(\bbeta) - \nabla Q_h(\bbeta) \}$, so that $\Delta(\bdelta) - \EE  \Delta(\bdelta)  =   G(\bbeta^*+\bdelta) - G(\bbeta^*)$.
Again, by a change of variable $\bv =  \bSigma^{ 1/2} \bdelta$,  
\#
 	\sup_{\bdelta \in  \Theta(r)}   \| \Delta(\bbeta) - \EE \Delta(\bbeta) \|_2 
& \leq \sup_{ \bdelta \in   \Theta(r) }  \|    G(\bbeta^*+ \bdelta) - G(\bbeta^*)     \|_2 \nn \\
& =   \sup_{ \| \bv \|_2 \leq r }   \| \underbrace{      G(\bbeta^* + \bSigma^{-1/2} \bv ) - G(\bbeta^*)  }_{ =:  \Delta_0(\bv) }   \|_2  .\nn 
\#
We will employ Theorem~A.3 in \cite{S2013} to bound the supremum $\sup_{ \| \bv \|_2 \leq r } \| \Delta_0(\bv) \|_2$, where $\Delta_0(\cdot)$ defined above satisfies $\Delta_0(\textbf{0}) = \textbf{0}$ and $\EE \{\Delta_0(\bv)\} =\textbf{0}$.
%Recall that $\bw_i = \bSigma^{-1/2} \bx_i \in \RR^{p}$ are isotropic random vectors, i.e., $\EE(\bw_i \bw_i^\T) = \Ib_p$. 
Taking the gradient with respect to $\bv$ yields
\#
	\nabla \Delta_0 (\bv) = \frac{1}{n} \sn \bigl\{ K_{i,\bv} \bw_i \bw_i^\T - \EE\bigl( K_{i,\bv} \bw_i \bw_i^\T \bigr) \bigr\} , \nn
\#
where $K_{i,\bv} : =  K_h(\langle \bw_i, \bv \rangle - \varepsilon_i )$ satisfies $0\leq K_{i,\bv} \leq  \kappa_u h^{-1}$.
For any $\bu, \bu' \in \mathbb{S}^{p-1}$ and $\lambda \in \RR$, using the elementary inequality $|e^u- 1 - u| \leq u^2 e^{|u|}/2$, we obtain 
\#
	 & \EE \exp\bigl\{   \lambda n^{1/2}  \langle \bu,   \nabla  \Delta_0 (\bv)  \bu' \rangle \bigr/\upsilon_1^2 \big\}  \nn \\
	& \leq   \prod_{i=1}^n  \Biggl\{  1 +   \frac{\lambda^2 }{2\upsilon_1^4 n} e^{ \frac{  \bar f |\lambda | }{ \upsilon_1^2 \sqrt{n}} \EE |  \langle \bw_i, \bu \rangle \langle \bw_i,  \bu' \rangle |}\EE    \bigl\{   K_{i,\bv} \langle \bw_i, \bu \rangle \langle \bw_i,  \bu' \rangle  - \EE (  K_{i,\bv}   \langle \bw_i, \bu \rangle \langle \bw_i,  \bu' \rangle)  \bigr\}^2    e^{\frac{ \kappa_u|\lambda | }{ h \sqrt{n}} | \langle \bw_i, \bu \rangle \langle \bw_i,  \bu' \rangle |/\upsilon_1^2 }   \Biggr\} \nn \\
	& \leq   \prod_{i=1}^n  \Biggl\{  1 +   \frac{\lambda^2 }{2 \upsilon_1^4 n} e^{ \frac{  \bar f |\lambda | }{ \sqrt{n}} }\EE    \bigl\{   K_{i,\bv} \langle \bw_i, \bu \rangle \langle \bw_i,  \bu' \rangle  - \EE (  K_{i,\bv}   \langle \bw_i, \bu \rangle \langle \bw_i,  \bu' \rangle)  \bigr\}^2    e^{\frac{ \kappa_u|\lambda | }{ h \sqrt{n}} | \langle \bw_i, \bu \rangle \langle \bw_i,  \bu' \rangle |/\upsilon_1^2}   \Biggr\} , \label{mgf.ubd1}
\#
where we used the bound $\EE |  \langle \bw_i, \bu \rangle \langle \bw_i, \bu' \rangle |  \leq (\EE \langle \bw_i, \bu \rangle^2 )^{1/2} (\EE \langle \bw_i,  \bu' \rangle^2 )^{1/2} = 1$ in the second inequality.
Moreover, the first and second conditional moments of $K_{i,\bv}$ can be rewritten as follows:
\# \,
\EE \bigl(   K_{i,\bv} | \bx_i \bigr)     = \frac{1}{h} \int_{-\infty}^\infty K \Biggl(  \frac{\langle \bw_i,  \bv  \rangle - t}{h} \Biggr) f_{\varepsilon_i  | \bx_i  } (t) \, {\rm d} t   = \int_{-\infty}^\infty K(u) f_{\varepsilon_i   | \bx_i  } (\langle \bw_i,  \bv \rangle - hu  ) \, {\rm d} u; \nn \\
\EE \bigl(   K^2_{i,\bv} | \bx_i \bigr)     = \frac{1}{h^2} \int_{-\infty}^\infty K^2 \Biggl(  \frac{\langle \bw_i,  \bv\rangle - t}{h} \Biggr) f_{\varepsilon_i  | \bx_i  } (t) \, {\rm d} t   = \frac{1}{h} \int_{-\infty}^\infty K^2(u) f_{\varepsilon_i  | \bx_i  } (\langle \bw_i,  \bv \rangle - hu  ) \, {\rm d} u, \nn
\#
from which it follows that $|\EE (   K_{i,\bv} | \bx_i ) | \leq \bar f $ and $\EE (   K_{i,\bv}^2 | \bx_i )\leq   \kappa_u \bar f   h^{-1}$ almost surely.

By the Cauchy-Schwarz inequality and the inequality $ab\le  a^2/2+b^2/2$, $a,b\in\RR$, we have 
\#
 &  \EE   ( \langle \bw_i, \bu \rangle \langle \bw_i,  \bu' \rangle )^2   e^{t  | \langle \bw_i, \bu \rangle \langle \bw_i,  \bu' \rangle | }   \nn \\
& \leq \EE   ( \langle \bw_i, \bu \rangle \langle \bw_i,  \bu' \rangle )^2   e^{ \frac{ t }{ 2 }   \langle \bw_i, \bu \rangle^2 +\frac{  t }{ 2 }  \langle \bw_i,  \bu' \rangle^2 }  \nn \\
& \leq    \Bigl(  \EE \langle \bw_i, \bu \rangle^4   e^{  t  \langle \bw_i, \bu \rangle^2   } \Bigr)^{1/2} \Bigl(  \EE \langle \bw_i,  \bu' \rangle^4   e^{ t \langle \bw_i,  \bu' \rangle^2   } \Bigr)^{1/2} , ~\mbox{ valid for any } t>0 .\nn
\#
Given a unit vector $\bu$, let $\chi = \langle \bw , \bu \rangle^2/(2\upsilon_1)^2$ so that under Condition~\ref{cond.predictor2}, 
$\PP( \chi \geq u)  \leq 2 e^{-2u}$ for any $u\geq 0$. It follows that
$\EE (e^\chi) = 1 + \int_0^\infty e^u \PP(\chi \geq u) {\rm d}u \leq 1 + 2 \int_0^\infty e^{-u} {\rm d}u = 3$, and 
\#
\EE (\chi^2 e^{\chi} ) = \int_0^\infty (u^2 + 2u) e^{u} \PP(\chi \geq u) {\rm d} u \leq 2 \int_0^\infty (u^2 + 2u) e^{-u} {\rm d} u = 8. \nn
\#
Taking the supremum over $\bu   \in \mathbb{S}^{p-1}$, we have 
$$
	\sup_{\bu \in \mathbb{S}^{p-1} } \EE e^{\langle \bw, \bu \rangle^2 / (2\upsilon_1)^2  } \leq 3 
	~\mbox{ and } ~  \sup_{\bu \in \mathbb{S}^{p-1} } \EE \langle \bw, \bu \rangle^4 e^{\langle \bw, \bu \rangle^2 / (2\upsilon_1)^2  } \leq 8(2\upsilon_1)^4 .
$$
Substituting the above bounds into \eqref{mgf.ubd1} yields that, for any $|\lambda | \leq \min\{   h/( 4  \kappa_u  ),  1/ \bar f  \}n^{1/2} $,
\#
 & \EE \exp\bigl\{ \lambda n^{1/2} \langle \bu,   \nabla \Delta_0 (\bv)   \bu'  \rangle/\upsilon_1^2 \bigr\}   \nn \\
& \leq   \prod_{i=1}^n  \Biggl[  1 +   \frac{ e\lambda^2    }{2\upsilon_1^4 n}\EE    \bigl\{   K_{i,\bv} \langle \bw_i, \bu \rangle \langle \bw_i,  \bu' \rangle  - \EE  (K_{i,\bv}   \langle \bw_i, \bu \rangle \langle \bw_i,  \bu' \rangle)  \bigr\}^2    e^{ | \langle \bw_i, \bu \rangle \langle \bw_i,  \bu' \rangle | /(4\upsilon_1^2)}   \Biggr]  \nn \\
& \leq  \prod_{i=1}^n  \Biggl[  1 +   \frac{ e\lambda^2    }{ \upsilon_1^4 n}\EE    \bigl(   K_{i,\bv} \langle \bw_i, \bu \rangle \langle \bw_i,  \bu' \rangle  \bigr)^2      e^{ | \langle \bw_i, \bu \rangle \langle \bw_i,  \bu' \rangle | / (4\upsilon_1^2) }   \nn \\
 &~~~~~~~~~~~~~~~~~~~~~~~~~~~~~~~ + \frac{e \lambda^2}{\upsilon_1^4 n} \bigl\{ \EE \bigl(  K_{i,\bv}   \langle \bw_i, \bu \rangle \langle \bw_i,  \bu' \rangle \bigr)\bigr\}^2 \EE  e^{ | \langle \bw_i, \bu \rangle \langle \bw_i,  \bu' \rangle |/ (4\upsilon_1^2)  }   \Biggr]  \nn \\
 & \leq  \prod_{i=1}^n  \Biggl(    1 + C_0^2  \frac{\lambda^2}{2 n h}  \Biggr) \leq \exp\bigl\{ C_0^2  \lambda^2 /(2h) \bigr\} ,\nn
\#
where $C_0>0$ depends only on $(\kappa_u, \bar f )$.
We have thus verified condition (A.4) in \cite{S2013} with ${\rm g}= \min\{  h /( 4\kappa_u  ), 1 / \bar f   \} (n/2)^{1/2}$ and $\nu_0= C_0 h^{-1/2}$.  Applying Theorem~A.3 therein to the process $\{ \Delta_0(\br)/\upsilon_1^2, \bv \in \BB^p(r) \}$, we obtain that with probability at least $1-e^{-t}$,
\#
  \sup_{ \| \bv \|_2 \leq r }   \|  \Delta_0(\bv)   \|_2 \leq 6 C_0  \upsilon_1^2 r  \sqrt{\frac{4 p + 2  t}{n h }}  \nn
\#
as long as $h \geq 8 \kappa_u  \sqrt{(2p+t)/n}$ and $n\geq 4 \bar f^2 (2p+t)$.

Joint with \eqref{bahaduri1i2} and \eqref{mean.Delta.ubd}, this implies that with probability at least $1-e^{-t}$,
\#\label{br.remainder.ubd}
\sup_{\bdelta \in  \Theta(r) }   \|  \Delta ( \bdelta )   \|_2 \leq  6 C_0 \upsilon_1^2 r  \sqrt{\frac{4 p + 2  t}{n h }}   + 0.5  l_0     m_3 r^2  . 
\#
Recall from the beginning of the proof that $\hat \bdelta = \hat \bbeta_h -  \bbeta ^* \in \Theta(r)$ with probability at least $1-2e^{-t}$ with  $r=r(n,p,t) \asymp \sqrt{(p+t)/n} +h^2$. Combined with \eqref{br.remainder.ubd}, we  conclude that with probability at least $1-3e^{-t}$,
$\| \Delta(\hat \bdelta ) \|_2 \lesssim (p+t)/(h^{1/2}n) +  h^{3/2} \sqrt{(p+t)/n} + h^4$, as claimed.
 \qed

%%%%%%%%%%%%%%%%%%%%%%%%%%%%%%%%%%
%%%%%%%%%%%%%%%%%%%%%%%%%%%%%%%%%%
% Proof of CLT
%%%%%%%%%%%%%%%%%%%%%%%%%%%%%%%%%%
%%%%%%%%%%%%%%%%%%%%%%%%%%%%%%%%%%
\subsection{Proof of Theorem~\ref{thm:clt}}
Let $\ba \in \RR^p$ be an arbitrary vector defining a linear functional of interest.  Given $h=h_n>0$, define $S_n =  n^{-1/2} \sn \gamma_i \xi_i$ and its centered version $S_n^0= S_n - \EE (S_n)$, where $\xi_i  =   \tau - \cK(-\varepsilon_i /h)$ and $\gamma_i = \langle \Jb_h^{-1} \ba, \bx_i \rangle$.  By the Lipschitz continuity of $f_{\varepsilon | \bx}(\cdot)$ and the fundamental theorem of calculus,  it can be shown that $| \EE(\xi_i | \bx_i)  |\leq  0.5 l_0  \kappa_2 h^2$, from which it follows by the law of iterated expectation that $|\EE (\gamma_i \xi_i) | \leq 0.5 l_0 \kappa_2 \| \Jb_h^{-1} \ba \|_{\bSigma}  \cdot h^2$.

Let $\eta_n = (p+\log n)/n$.  Then, applying  \eqref{br.remainder.def} and \eqref{br.remainder.ubd}  with $t=\log n$ and the triangle inequality, we obtain that under the constraint $\eta_n ^{1/2} \lesssim h \lesssim 1$,
\#
&	\bigl| n^{1/2} \langle \ba, \hat \bbeta_h - \bbeta^* \rangle   - S_n^0 \bigr| \nn \\ 
	&= n^{1/2}  \biggl|   \bigg\langle  \bSigma^{ 1/2}  \Jb_h^{-1} \ba,  \bSigma^{-1/2}   \Jb_h(\hat \bbeta_h - \bbeta^* ) -   \bSigma^{-1/2} 	 \frac{1}{n} \sn \bigl\{ \tau - \cK(-\varepsilon_i/h) \bigr\} \bx_i  \bigg\rangle \biggr|   + | \EE (S_n) |  \nn \\
	& \leq  c_1 \|  \Jb_h^{-1} \ba \|_{\bSigma}  \cdot   n^{1/2}  \bigl(   h^{-1/2} \delta_n  + h^2  \bigr) \label{bh.error}
\#
with probability at least $1 - 3 n^{-1}$ for some constant $c_1>0$.

For the centered partial sum $S_n^0 = S_n - \EE (S_n) = n^{-1/2} \sn (1-\EE) \gamma_i \xi_i$, we have $\var(S_n^0) = \var(S_n) = \EE (\gamma \xi)^2 - \{\EE (\gamma \xi)\}^2$, where $\gamma = \langle\Jb_h^{-1} \ba, \bx\rangle$ and $\xi = \tau - \cK(-\varepsilon/h)$.  By the Berry-Esseen inequality (see, e.g., \cite{T2011}),
\#
	\sup_{x\in \RR} \bigl| \PP\bigl\{ S_n^0 \leq \var(S_n)^{1/2} x \bigr\} - \Phi(x) \bigr| \leq  \frac{\EE \{ |\gamma \xi - \EE (\gamma \xi) |^3\} }{2 [\EE ( \gamma \xi )^2 -  \{\EE( \gamma \xi)\}^2 ]^{3/2} \sqrt{n}  } . \label{be.bound1}
\#
We have shown that $|\EE(\gamma \xi) |  \lesssim  \| \Jb_h^{-1} \ba \|_{\bSigma} \cdot h^2$.  
To bound $\EE(\xi^2|\bx)$,  by a change of variable and integration by parts, 
\#
\EE\{\cK^2(-\varepsilon/h) |\bx\} &= 2\int_{-\infty}^{\infty} K(v) \cK (v) F_{\varepsilon|\bx}(-vh)\mathrm{d}v\nn \\
&=2 \tau  \underbrace{  \int_{-\infty}^{\infty} K(v) \cK(v) \mathrm{d}v }_{=1/2}	-2h f_{\varepsilon|\bx}(0)  \underbrace{  \int_{-\infty}^{\infty}v K(v) \cK(v) \mathrm{d}v }_{= \int_{0}^\infty \cK(v)\{ 1-\cK(v) \} {\rm d}v   >0 } \nn \\
&\quad+ 2\int_{-\infty}^{\infty}\int_{0}^{-vh} \{f_{\varepsilon| \bx}(t) -f_{\varepsilon| \bx}(0) \} K(v) \cK(v) \mathrm{d}t\mathrm{d}v\nn \\ 
&\le \tau + l_0 \kappa_2 h^2\nn ,
\#
where $\kappa_2$ and $l_0$ are the constants from Conditions \ref{cond.kernel} and \ref{cond.reg}.
It then follows that $$
    \tau(1-\tau) - C h - (1+\tau) l_0 \kappa_2 h^2  \leq      \EE(\xi^2 | \bx_i ) \leq   \tau(1-\tau) + (1+\tau) l_0 \kappa_2 h^2 ,
$$
where $C = 2 \bar f \cdot \int_{0}^\infty \cK(v)\{ 1-\cK(v) \} {\rm d}v$.
For all sufficiently small $h$, $\var(S_n)= \{ \tau (1-\tau) + O(h)\} \| \Jb_h^{-1} \ba \|_{\bSigma}^2$. On the other hand, $\EE (|\gamma \xi |^3) \leq \max(\tau, 1-\tau) \EE (\xi^2 |\gamma  |^3) \leq m_3 \{ \tau(1-\tau) + O(h^2) \}   \| \Jb_h^{-1} \ba \|_{\bSigma}^3$. Substituting these bounds into \eqref{be.bound1} yields
\#
	\sup_{x\in \RR} \bigl| \PP\bigl\{ S_n^0 \leq \var(S_n)^{1/2} x \bigr\} - \Phi(x) \bigr|  \leq c_2 n^{-1/2} \label{be.bound2}
\#
for some constant $c_2>0$. 

Recall that $\sigma_{h}^2= \EE \{ \cK_h(-\varepsilon) - \tau \}^2 \langle \Jb_h^{-1}\ba , \bx \rangle^2 = \EE(\gamma \xi)^2$, and thus $|\var(S_n) - \sigma_{h}^2 | = (\EE \gamma \xi)^2 \lesssim \| \Jb_h^{-1} \ba \|_{\bSigma}^2\cdot  h^4$.
By an application of Lemma A.7 in the supplement of \citet{SZ2015}, for sufficiently small $h$, we have
\#
	 \sup_{x\in \RR } \bigl| \Phi( x/\var(S_n)^{1/2}  )   -\Phi(x/\sigma_h ) \bigr|   \leq c_3 h^4. \label{gaussian.comparison}
\#

Before proceeding, we note that the constants $c_1$--$c_3$ appeared above are all independent of $\ba$ and $(n,p)$.
Let $G\sim N(0,1)$.  Putting together the above derivations, for any $x\in \RR$ and $\ba \in \RR^p$,  we obtain 
\#
	& \PP\bigl(  n^{1/2}  \langle \ba, \hat \bbeta_h -\bbeta^* \rangle  \leq x \bigr) \nn \\
	& \leq  \PP\bigl\{    S_n^0    \leq x  +  c_1 \|  \Jb_h^{-1} \ba \|_{\bSigma}  \cdot  n^{1/2}  \bigl(   h^{-1/2} \eta_n  + h^2  \bigr) \bigr\} + 3 n^{-1} \nn \\
	& \leq  \PP\bigl\{   \var(S_n)^{1/2} G \leq x+ c_1 \|  \Jb_h^{-1} \ba \|_{\bSigma}  \cdot n^{1/2}  \bigl(   h^{-1/2} \eta_n  + h^2  \bigr) \bigr\} +c_2n^{-1/2} + 3 n^{-1} \nn \\
	&\leq  \PP\bigl\{   \sigma_h  G \leq x+ c_1 \|  \Jb_h^{-1} \ba \|_{\bSigma} \cdot n^{1/2}  \bigl(   h^{-1/2} \eta_n  + h^2  \bigr) \bigr\} + c_2 n^{-1/2} + c_3 h^4 +  3 n^{-1} \nn \\
	& \leq  \PP\bigl(   \sigma_h  G \leq x  \bigr) + c_1(2\pi)^{-1/2}  \|  \Jb_h^{-1} \ba \|_{\bSigma} \cdot  n^{1/2}  \bigl(   h^{-1/2} \eta_n  + h^2  \bigr)/\sigma_{ h} + c_2 n^{-1/2} + c_3 h^4 +  3 n^{-1} , \nn
\#
where the first, second, and third inequalities holds by \eqref{bh.error}, \eqref{be.bound2}, and \eqref{gaussian.comparison}, respectively, and the last inequality follows from the fact that for any $a\leq b$ and $\sigma>0$, $\Phi(b/\sigma ) - \Phi(a/\sigma) \leq  (2\pi)^{-1/2}(b-a)/\sigma$.
Recall that $\sigma_{ h}^2 = \{ \tau(1-\tau) +O(h) \} \|\Jb_h^{-1}\ba\|_{\bSigma}^2$,   $\|  \Jb_h^{-1} \ba \|_{\bSigma} / \sigma_{ h} \lesssim \{ \tau(1-\tau)\}^{-1/2}$ for all sufficiently small $h$.
A similar argument leads to a series of reverse inequalities. 
The above bounds are independent of $x$ and $\ba$, and therefore hold uniformly over all  $x$ and $\ba$.

Putting together the pieces, we conclude that under the bandwidth requirement $\eta_n^{1/2} \lesssim h \lesssim 1$,
\#
	\sup_{x\in \RR , \, \ba \in \RR^p } \bigl| \PP\bigl(  n^{1/2}  \langle \ba, \hat \bbeta_h -\bbeta^* \rangle  \leq \sigma_{h} x \bigr)  - \Phi(x ) \bigr|  \lesssim  \frac{p+ \log n}{ (n h)^{1/2}} + n^{1/2} h^2 , \nn
\#
as claimed.

Under the additional smoothness condition on $f_{\varepsilon | \bx}(\cdot)$, using a higher-order Taylor series expansion of $F_{\varepsilon|\bx}(\cdot)$ gives
\#
&~~~~ \EE \bigl\{  \cK(-\varepsilon_i /h) | \bx_i \bigr\} =   \int_{-\infty}^\infty  K(u) F_{\varepsilon|\bx}( h u)  \, {\rm d}u  \nn \\
&= \int_{-\infty}^\infty  K(u) \Bigg[  F_{\varepsilon|\bx}(0) + hu f_{ \varepsilon|\bx }(0) + \frac{h^2 }{2} u^2 f'_{ \varepsilon|\bx }(0) + \frac{h^3 }{3!} u^3 f''_{ \varepsilon|\bx } (0)   \nn \\
&~~~~~~~~~~~~~~~~~~~~~~~~~~~~~~~~~~~~~~~~~~ + \frac{ h^3}{2!}  u^3\int_0^1 (1-w)^2 \big\{ f''_{ \varepsilon|\bx } (h u w) - f''_{ \varepsilon|\bx }(0) \big\} {\rm d} w \Bigg]  {\rm d}u  \nn \\
&= \tau + \frac{\kappa_2  }{2} f'_{ \varepsilon|\bx } (0) h^2 +   O(h^4) \kappa_4 l_2(\bx) , \nn
\#
from which it follows that 
\$
	 \Bigg|  \EE(\gamma_i \xi_i)  + \frac{\kappa_2}{2}   \langle \Jb_h^{-1} \ba,   \EE  \{ f'_{\varepsilon | \bx}(0) \bx_i  \} \rangle     \cdot   h^2 \Bigg| \lesssim  \kappa_4   \EE | \langle \Jb_h^{-1} \ba, \bx_i \rangle l_2(\bx) |  \cdot  h^4 \lesssim  \kappa_4 \| \Jb_n^{-1} \ba \|_{\bSigma} \big\{ \EE l_2^2(\bx) \big\}^{1/2} h^4 .
\$
Together, the above bound and \eqref{br.remainder.ubd} with $t=\log n$ imply
\#
\bigg| n^{1/2} \bigg\langle \ba, \hat \bbeta_h - \bbeta^* &	 + \frac{\kappa_2}{2} \Jb_h^{-1}  \EE  \{ f'_{\varepsilon | \bx}(0) \bx_i  \} \bigg\rangle   - S_n^0 \bigg|  \nn \\
& ~~~~~~~~~~~~~~ \lesssim   \|  \Jb_h^{-1} \ba \|_{\bSigma}  \cdot   n^{1/2}  \Biggl(  \frac{p+\log n}{n h^{1/2}} + h^{3/2} \sqrt{\frac{p+ \log n}{n}} +  h^4  \Biggr)  \label{bh.error2}
\#
with probability at least $1-3n^{-1}$. Repeating the above analysis, with \eqref{bh.error} replaced by \eqref{bh.error2}, proves the refined Berry-Esseen bound \eqref{linear.clt2}.  \qed

\subsection{Proof of Theorem~\ref{thm:boot.concentration}}

Keep the notation used in the proof of Theorem~\ref{thm:concentration}. 
With non-negative multipliers $w_i$'s, the weighted loss $\hat Q_h^\flat(\cdot)$ given in \eqref{wsq.loss} is also convex.  For $\bdelta \in \RR^p$, define the bootstrap counterpart of $\hat D_h(\cdot)$ as $\hat D^\flat_h(\bdelta) = \hat Q_h^\flat(\bbeta^* + \bdelta) - \hat Q_h^\flat(\bbeta^*)$. It is easy to see that $\EE^*\{ \hat D^\flat_h(\bdelta) \} = \hat D_h(\bdelta) $.  Similarly to \eqref{Dh.lbd1}, we have
\#
\hat D^\flat_h(\bdelta) & = \hat D_h(\bdelta) +           \{ \hat D^\flat_h(\bdelta)  - \hat D_h(\bdelta) \} \nn \\
& \geq  R_h(\bdelta) - 0.5 l_0 \kappa_2 h^2 \cdot \| \bdelta \|_{\bSigma} -  \{  D_h(\bdelta)  -  \hat  D_h(\bdelta)\}   -     \{ \hat D_h(\bdelta) - \hat D^\flat_h(\bdelta) \}  \nn \\
& \geq   \frac{1}{2} \big( \underbar{$f$} - l_0 \kappa_1 h - 0.5 l_0 m_3  \| \bdelta \|_{\bSigma}  - l_0 \kappa_2 h^2 /  \| \bdelta \|_{\bSigma} \big)  \| \bdelta \|_{\bSigma}^2  \nn \\
&~~~~~~~~~~~~~~~~~~~~~~~~~~~~~- \underbrace{ \{  D_h(\bdelta)  -  \hat  D_h(\bdelta)\}      }_{{{\rm sampling~error}}} -  \underbrace{   \{   \hat D_h(\bdelta) - \hat D^\flat_h(\bdelta)  \}  }_{{\rm bootstrap~error}} ,  \ \ \bdelta \in \RR^p. \label{boot.Dh.lbd1}
\#
As before,  let $r_0 = (2\kappa_2/m_3)^{1/2}   h$ be an intermediate convergence radius, and set $r_l = r_0 h$.  For any $\bdelta \in \partial \Theta(r_0)$, it follows that
\#
\hat D^\flat_h(\bdelta) \geq  \frac{1}{2} \big\{  \underbar{$f$} - l_0 \kappa_1 h -  l_0 (2\kappa_2/m_3)^{1/2}   h \big\} r_0^2  -  \{  D_h(\bdelta)  -  \hat  D_h(\bdelta)\}       -  \{     \hat   D_h(\bdelta) - \hat D^\flat_h(\bdelta)\}  . \label{boot.Dh.lbd2}
\#

As a bootstrap counterpart of Lemma~\ref{lem:boot.loss.difference}, the following lemma provides high probability bounds for the bootstrap error $\hat   D_h(\bdelta) - \hat D^\flat_h(\bdelta)$ uniformly over $\bdelta \in \RR^p$ in some compact subset.

\begin{lemma} \label{lem:boot.loss.difference}
For each $t \geq 0$,  there exists an event $\cE_1(t)$ with  $\PP\{ \cE_1(t) \} \geq 1-e^{-t}$ such that conditioned on $\cE_1(t)$,
\#
	\PP^*  \Bigg\{  \sup_{\bdelta \in \Theta(r)} \{\hat   D_h(\bdelta) - \hat D^\flat_h(\bdelta)  \}\geq   C \bar \tau  \upsilon_1  r  \sqrt{\frac{ u}{n}} \, \Bigg\} \leq e^{2p-u}    \nn
\#
for any $u\geq 0$ as long as $n\gtrsim p+t$.
Moreover,  for any $r_u>r_l>0$, with $\PP^*$-probability (over $\{e_i\}_{i=1}^n$) at least $1-\lceil e \log(\frac{r_u}{r_l} ) \rceil e^{2p-u}$ conditioned on $\cE_1(t)$,
\#
\hat   D_h(\bdelta) - \hat D^\flat_h(\bdelta) \leq  C'  \bar \tau  \upsilon_1  \| \bdelta \|_{\bSigma}  \sqrt{\frac{ u}{n}} \nn
\#
holds for all $\bdelta \in \RR^p$ satisfying $r_l \leq \| \bdelta \|_{\bSigma} \leq r_u$ as long as $n\gtrsim p+t$.  Here both $C, C'>0$ are absolute constants.
\end{lemma}

Let $\cE_1(t)$ be the event in Lemma~\ref{lem:boot.loss.difference} that occurs with probability at least $1-e^{-t}$.
Applying the first inequality with $r=r_0$ and $u=2p+t$,  we obtain that conditioned on  $\cE_1(t)$,
\#
\hat   D_h(\bdelta) - \hat D^\flat_h(\bdelta) \leq C \bar \tau \upsilon_1  r_0 \sqrt{\frac{2p+t}{n}} .  \nn
\#
Let $\cE_2(t)$ be the event that the bounds \eqref{Dh.unif.bd1} and \eqref{Dh.unif.bd2} hold. Then, the ``good event" $\cE(t):= \cE_1(t) \cap \cE_2(t)$ satisfies $\PP\{ \cE(t) \} \geq 1-3e^{-t}$.  By \eqref{boot.Dh.lbd2} and the above bound,  we find that with $\PP^*$-probability at least $1-e^{-t}$ conditioned on $\cE(t)$,  $\hat D_h^\flat(\bdelta) >0$ for all $\bdelta \in \partial \Theta(r_0)$ as long as the bandwidth satisfies $\underbar{$f$}^{-1} m_3^{1/2} \upsilon_1 \sqrt{(p+t)/n} \lesssim  h \lesssim \underbar{$f$}m_3^{-1/2}$.  
Let $\hat \bdelta^\flat = \hat \bbeta_h^\flat - \bbeta^*$. Then,  $\hat D^\flat_h(\hat \bdelta^\flat )\leq 0$ by the optimality of $\hat \bbeta^\flat_h$.  This, combined with the convexity of $\hat Q^\flat_h(\cdot)$, ensures that $ \hat \bdelta^\flat \in \Theta(r_0)$.  

Next, consider the decomposition $\Theta(r_0) = \Theta(r_l) \cap \Theta(r_l, r_0)$, where $\Theta(r_l, r_0)=\{ \bdelta \in \RR^d: r_l \leq \| \bdelta \|_{\bSigma} \leq r_0 \}$ and $r_l= r_0 h = (2\kappa_2/m_3)^{1/2}   h^2$.  If $ \hat  \bdelta^\flat \in \Theta(r_l)$, the claimed bound holds trivially. Throughout the rest, we assume $ \hat  \bdelta^\flat \in \Theta(r_l, r_0)$. Taking $u=\log(e\log h^{-1}) +2p + t$ in the second inequality in Lemma~\ref{lem:boot.loss.difference} yields that, with $\PP^*$-probability at least $1-e^{-t}$ conditioned on $\cE(t)$, 
\$
  \hat D_h( \hat  \bdelta^\flat )  - \hat D_h^\flat( \hat \bdelta^\flat)  \leq  \| \hat  \bdelta^\flat \|_{\bSigma}  \cdot \underbrace{  C'  \bar \tau  \upsilon_1    \sqrt{\frac{ \log(e\log h^{-1}) +2p + t}{n}} }_{=: r_1^\flat} .
\$
Substituting the above bound and \eqref{Dh.unif.bd2} into \eqref{boot.Dh.lbd1}, we conclude that
\$
 	( \underbar{$f$} - l_0 \kappa_1 h ) \| \hat \bdelta^\flat \|_{\bSigma}^2 &  \leq   ( 2 r_1 + 2 r_1^\flat +  l_0 \kappa_2 h^2 )    \| \hat \bdelta^\flat \|_{\bSigma}   + 0.5 l_0 m_3 \| \hat \bdelta^\flat \|_{\bSigma}^3  \leq 2 (   r_1 +   r_1^\flat +  l_0 \kappa_2 h^2 )    \| \hat \bdelta^\flat \|_{\bSigma} .\nn
\$
Canceling out a factor of $\| \hat \bdelta^\flat \|_{\bSigma}$ from both sides yields the claimed bound. \qed

\subsubsection{Proof of Lemma~\ref{lem:boot.loss.difference}}

We will use a similar argument as in the proof of Lemma~\ref{lem:loss.difference}.  Consider the bootstrap loss difference
$\hat D_h^\flat(\bdelta) - \hat D_h(\bdelta) = (1/n) \sn e_i d_h(\bdelta ; \bz_i)$, where $d_h(\bdelta ; \bz_i)= \ell_h(\varepsilon_i - \langle \bx_i, \bdelta \rangle ) - \ell_h(\varepsilon_i)$ with $\bz_i=(\bx_i, \varepsilon_i)$, and $e_i$'s are independent Rademacher random variables that are independent of $\{ \bz_i \}_{i=1}^n$.

For any $r>0$ and any $\epsilon \in (0,1)$,  applying a conditional version of Chernoff's inequality to $\Delta^\flat_\epsilon(r) :=n (1-\epsilon)\sup_{\bdelta \in \Theta(r)}  \{  \hat D_h(\bdelta) - \hat D_h^\flat(\bdelta) \} / (\bar \tau r)$ yields
\#
\PP^* \big\{ \Delta^\flat_\epsilon(r)  \geq u \big\} \leq \exp\Bigg[ - \sup_{\lambda \geq 0}  \big\{ \lambda u - \log \EE^* e^{\lambda \Delta^\flat_\epsilon (r)} \big\} \Bigg].  \nn
\#
Again, by the Ledoux-Talagrand contraction principle and discretization via an $\epsilon$-net,
\#
& \EE^* e^{\lambda \Delta^\flat_\epsilon (r)}   \leq \EE^* \exp\Bigg\{ \frac{\lambda}{r} (1-\epsilon) \sup_{\bdelta \in \Theta(r)  } \sn e_i \langle \bx_i, \bdelta \rangle    \Bigg\}   \nn \\
& \leq   \EE^*  \exp \Bigg\{ \lambda(1-\epsilon) \bigg\| \sn e_i \bw_i \bigg\|_2 \Bigg\} \leq \EE^*  \exp \Bigg\{ \lambda \max_{1\leq j\leq N_\epsilon }  \sn e_i  \bu_j^\T \bw_i  \Bigg\} \nn \\
& \leq \sum_{j=1}^{N_\epsilon}  \prod_{i=1}^n \EE^* e^{\lambda e_i  \bu_j^\T \bw_i   } \leq  \sum_{j=1}^{N_\epsilon}  e^{ ( \lambda^2/2) \sn  \langle \bu_j, \bw_i\rangle^2 } , \nn
\#
where we used the Rademacher moment bound $\EE e^{\lambda e_i} \leq e^{\lambda^2/2}$ for any $\lambda \in \RR$. Here $\{ \bu_j\}_{j=1}^{N_\epsilon}$ are unit vectors, and $N_\epsilon \leq (1+2/\epsilon)^p$.  Consequently,
\#
\log \EE^* e^{\lambda \Delta^\flat_\epsilon (r)} \leq \log N_\epsilon + \frac{\lambda^2}{2} \max_{1\leq j\leq N_\epsilon } \sn \langle \bu_j , \bw_i\rangle^2  ,  \nn
\#
so that for any $u\geq 0$,
\#
 \sup_{\lambda \geq 0}  \big\{ \lambda u - \log \EE^* e^{\lambda \Delta^\flat_\epsilon (r)} \big\}  & \geq -\log N_\epsilon +  \sup_{\lambda \geq 0}  \Bigg(  \lambda u -  \frac{1}{2} \lambda^2 \max_{1\leq j\leq N_\epsilon } \sn \langle \bu_j , \bw_i\rangle^2 \Bigg) \nn \\
& =   -\log N_\epsilon +  \frac{u^2 }{2\max_{1\leq j\leq N_\epsilon } \sn \langle \bu_j , \bw_i\rangle^2}.  \nn
\#
Substituting this into the earlier Chernoff's inequality, and by a change of variable, we obtain that for any $u \geq 0$, 
\$
  \sup_{\bdelta \in \Theta(r)}  \hat D_h(\bdelta) - \hat D_h^\flat(\bdelta) \leq  \frac{  \bar \tau}{1-\epsilon}  \sqrt{ \max_{1\leq j\leq N_\epsilon } \frac{1}{n} \sn \langle \bu_j , \bw_i\rangle^2}  \cdot r \sqrt{\frac{ 2 u}{n}}  
\$
holds with $\PP^*$-probability (over $\{e_i\}_{i=1}^n$) at least $1- e^{p\log(1+2/\epsilon) - u}$.

The above bound holds for any given $r>0$.   Proceed via a peeling argument, we obtain that for any prespecified $\gamma>1$ and radii $r_u>r_l>0$, 
\#
\hat D_h(\bdelta) - \hat D_h^\flat(\bdelta)  \leq  \frac{ \gamma \bar \tau }{1-\epsilon}  \sqrt{ \max_{1\leq j\leq N_\epsilon } \frac{1}{n} \sn \langle \bu_j , \bw_i\rangle^2}  \cdot \| \bdelta \|_{\bSigma}  \sqrt{\frac{ 2 u}{n}} ~\mbox{ holds for all }  r_l \leq \| \bdelta \|_{\bSigma} \leq r_u   \label{boot.different.unif.ubd1}
\#
with $\PP^*$-probability (over $\{e_i\}_{i=1}^n$) at least $1-   \lceil \log(\frac{r_u}{r_l}) / \log(\gamma) \rceil e^{p\log(1+2/\epsilon) - u}$.  

Next, we bound the data-dependent quantity $\max_{1\leq j\leq N_\epsilon} (1/n) \sn \langle \bu_j , \bw_i\rangle^2$ with $N_\epsilon \leq (1+2/\epsilon)^p$ and $\bu_j \in \mathbb{S}^{p-1}$.  Note that $\EE \langle \bu_j , \bw_i\rangle^2=1$.  By the sub-Gaussianity of $\bw_i \in \RR^p$ ensured by Condition~\ref{cond.predictor2},  for integers $k=2,3,\ldots$,
\#
	\EE \bigl(   \langle \bu_j , \bw_i \rangle^2 \bigr)^k & \leq  \upsilon_1^{2k}  \cdot 2k \int_0^\infty \PP\bigl(  | \langle \bu_j , \bw_i \rangle | \geq \upsilon_1 u \bigr) u^{2k-1} \, {\rm d} u \nn \\
	& \leq    \upsilon_1^{2k} \cdot  4k \int_0^\infty u^{2k-1} e^{-u^2/2 } \, {\rm d} u \nn \\
	& =   2^k \upsilon_1^{2k} \cdot  2k \int_0^\infty v^{ k-1} e^{-v } \, {\rm d}v = 2^{k+1} k!  \cdot      \upsilon_1^{2k}   . \nn
\#
In particular, $\EE(  \langle \bu, \bw_i \rangle^4)\leq     16 \upsilon_1^4$ and $\EE (    \langle \bu, \bw_i \rangle^2  )^k \leq \frac{k!}{2} \cdot    (4  \upsilon_1^2)^2 \cdot   (2   \upsilon_1^2 )^{k-2}$ for $k\geq 3$. With the above calculations,  applying Bernstein's inequality (see, e.g., Theorem~2.10 in \cite{BLM2013}),  and taking the union bound over $j=1,\ldots, N_\epsilon$, we obtain that for any $v\geq 0$,
\#
	\PP \Bigg(   \max_{1\leq j\leq N_\epsilon }\frac{1}{n} \sn \langle \bu_j, \bw_i \rangle^2 \geq  1  + 4 \upsilon_1^2   \sqrt{\frac{2    v}{n}} +  2\upsilon_1^2    \frac{ v}{n} \Bigg) \leq  \exp\big\{ p\log(1+2/\epsilon) -v \big\}  .\label{boot.different.unif.ubd2}
\#

Finally, we take $\epsilon=2/(e^2-1)$, $\gamma=e^{1/e }$ and  $v= 2p+t$ ($t\geq 0)$ in \eqref{boot.different.unif.ubd1} and \eqref{boot.different.unif.ubd2}.  Let $\cE_1(t)$ be the event that \eqref{boot.different.unif.ubd2} holds. Then, $\PP\{ \cE_1(t) \} \geq 1-e^{-t}$, and with $\PP^*$-probability (over $\{e_i\}_{i=1}^n$) at least $1- \lceil e \log(\frac{r_u}{r_l})  \rceil e^{2p - u}$ conditioned on $\cE_1(t)$,
\$
\hat D_h(\bdelta) - \hat D_h^\flat(\bdelta)  \leq 3 \bar \tau    \sqrt{1 + 4\upsilon_1^2 \sqrt{(4p+2t)/n} +  \upsilon_1^2 (4 p+2 t)/n }  \cdot \| \bdelta \|_{\bSigma}  \sqrt{\frac{ u}{n}}
\$
holds for all $\bdelta \in \RR^p$ satisfying $r_l \leq \| \bdelta \|_{\bSigma} \leq r_u$.  This proves the claimed bound under the scaling $n\gtrsim p+t$.   \qed

\subsection{Proof of Theorem~\ref{thm:boot.bahadur}}

The proof is based on an argument similar to that used in the proof of Theorem~\ref{thm:bahadur}. To begin with, define the random process
\#
	\Delta^\flat (\bdelta) =  \bSigma^{-1/2}  	  \bigl\{  \nabla \hat Q_h^\flat (\bbeta^* + \bdelta) -  \nabla \hat Q_h^\flat (\bbeta^*)  - \Jb_h \bdelta \bigr\}   , ~~ \bdelta \in \RR^p.
\#
For a prespecified $r>0$, a key step is to bound the local fluctuation $\sup_{\bdelta \in \Theta(r)} \| \Delta^\flat (\bdelta ) \|_2$. Since $\EE(w_i) = \EE(1+e_i)=1$, we have $\EE^* \{  \nabla \hat Q_h^\flat (\bbeta )  \} =  \nabla \hat Q_h  (\bbeta )$.
Define the (conditionally) centered process 
$$ 
	G^\flat (\bbeta ) = \bSigma^{-1/2} \big\{ \nabla \hat Q_h^\flat (\bbeta ) - \nabla \hat Q_h  (\bbeta ) \big\} = \frac{1}{n} \sn   \big\{ \cK_h (   \bx_i^\T \bbeta  -y_i  ) - \tau  \big\}  e_i \bw_i, 
$$
so that $\Delta^\flat(\bdelta)$ be be written as
\#
	\Delta^\flat(\bdelta)      =  \underbrace{   \bigl\{    G^\flat (\bbeta^*+\bdelta ) -  G^\flat (\bbeta^*)  \bigr\} }_{{\rm bootstrap~error}} +  \underbrace{  \Delta(\bdelta)  }_{{\rm sampling~error}}  , \nn
\#
where $\Delta(\bdelta)$ is given in \eqref{br.remainder.def}. By the triangle inequality,
\#
\sup_{\bdelta \in   \Theta(r)}  \| \Delta^\flat (\bdelta )   \|_2 \leq \sup_{\bdelta \in  \Theta(r)}  \| 	   G^\flat (\bbeta^* + \bdelta ) -  G^\flat (\bbeta^*)    \|_2 + \sup_{\bdelta \in   \Theta(r)}   \| \Delta  (\bdelta )   \|_2 . \label{Deltab.decomposition}
\#

Let $\cE_3(t)$ denote the event that the bound \eqref{br.remainder.ubd} holds.
It suffices to deal with the first term on the right-hand side of \eqref{Deltab.decomposition}.  Using a change of variable $\bv = \bSigma^{1/2} \bdelta \in \mathbb{B}^p(r)$ for $\bdelta \in   \Theta(r)$, we have $y_i - \bx_i^\T \bbeta  = \varepsilon_i -   \bw_i^\T  \bv$ and
\#
 \sup_{\bdelta \in   \Theta(r)}   \| 	    G^\flat (\bbeta^*+ \bdelta ) -  G^\flat (\bbeta^*)   \|_2  %& \leq \underbar{$f$}^{-1/2} \sup_{ \| \bv \|_2 \leq r } \bigl\|  \underbrace{  	\bSigma^{-1/2} \bigl\{    G^\flat (\bbeta^* + \bSigma^{-1/2} \bv ) -  G^\flat (\bbeta^*)  \bigr\} }_{:= \Delta_0^\flat(\bv) } \bigr\|_2  \nn \\
 %& \leq   \sup_{   \bu  , \bv \in \BB^p(1) }  \langle G^\flat ( \bbeta^* +  \bSigma^{-1/2} r \bv  ) -  G^\flat( \bbeta^* )   ,     \bu   \rangle \nn \\
 \leq    \sup_{ \bv \in \mathbb{B}^p(r) }    \| 	 \underbrace{   G^\flat (\bbeta^* + \bSigma^{-1/2} \bv ) -  G^\flat (\bbeta^*)    }_{  =:   \, \Delta^\flat_0(   \bv)  }  \|_2  , \label{Deltab.ub1}
\#
where $\Delta^\flat_0( \bv) =(1/n)\sn  e_i   \bw_i  \{  \cK_h(    \bw_i^\T  \bv - \varepsilon_i   ) - \cK_h(-\varepsilon_i ) \}$ satisfies $\Delta^\flat_0(\textbf{0})= \textbf{0}$ and $\EE^*\{ \Delta^\flat_0( \bv) \}=\textbf{0}$. Note that 
$\nabla \Delta^\flat_0( \bv) = (1/n) \sn e_i K_{i,\bv} \bw_i \bw_i^\T$, where $K_{i,\bv}  = K_h (\bw_i^\T  \bv - \varepsilon_i)$. For any $\lambda \in \RR$ and $\bu, \bu' \in \mathbb{S}^{p-1}$, we have
\$
 & \EE^* \exp \bigl\{  \lambda n^{1/2} \bu^\T \Delta^\flat_0( \bv) \bu'  \bigr\}    = \prod_{i=1}^n \EE^* \exp\big\{ \lambda n^{-1/2} e_i K_{i,\bv}  \bw_i^\T \bu \cdot \bw_i^\T   \bu'  \big\} \\
& \leq  \prod_{i=1}^n    \exp\Bigg\{ \frac{\lambda^2}{2 n} K^2_{i,\bv} ( \bw_i^\T \bu \cdot  \bw_i^\T \bu' )^2  \Bigg\} =  \exp\Bigg\{  \frac{\lambda^2}{2 n}  \sn K^2_{i,\bv} ( \bw_i^\T \bu \cdot  \bw_i^\T \bu' )^2  \Bigg\} . 
\$
Note that $K_{i,\bv} \leq \kappa_u h^{-1}$, and by H\"older's inequality,
\$
\frac{1}{n }  \sn K^2_{i,\bv} ( \bw_i^\T \bu \cdot  \bw_i^\T \bu' )^2 \leq  \frac{\kappa_u}{h} \cdot \Biggl\{ \frac{1}{n }  \sn K_{i,\bv} ( \bw_i^\T \bu )^4\Biggr\}^{1/2} \Bigg\{  \frac{1}{n }  \sn K_{i,\bv} ( \bw_i^\T \bu' )^4 \Bigg\}^{1/2}.
\$
Define the function $\Lambda_{r,h}(\cdot, \cdot) : \RR^p \times \RR^p \mapsto [0,\infty)$ 
\#
	\Lambda_{r,h}(\bu, \bv) = \frac{1}{n} \sn  K_h(r \bw_i^\T \bv - \varepsilon_i) ( \bw_i^\T \bu)^4 ~\mbox{ for } \bu, \bv \in \mathbb{S}^{p-1},  \label{def:Lambda}
\#
and write $\| \Lambda_{r,h} \|_\infty = \sup_{\bu, \bv \in \mathbb{S}^{p-1}} \Lambda_{r,h}(\bu, \bv ) $. 
With this notation, it follows that for any $\lambda \in \RR$ and $\bu, \bv \in \mathbb{S}^{p-1}$,
\$
  \sup_{ \bdelta  \in \mathbb{B}^p(r)}  \EE^* \exp \bigl\{  \lambda n^{1/2} \bu^\T \Delta^\flat_0( \bdelta) \bv  \bigr\}  \leq \exp\Bigg\{  \frac{ \lambda^2}{ 2 } \kappa_u h^{-1} \| \Lambda_{r,h} \|_\infty \Bigg\} .
\$
Thus, applying a conditional version of Theorem~A.3 in \cite{S2013} yields
\#
 \sup_{\bv \in \mathbb{B}^p(r) }  \| \Delta^\flat_0(\bv) \|_2 \leq  6 \kappa_u^{1/2} \| \Lambda_{r,h} \|_\infty^{1/2}  \cdot  r \sqrt{\frac{4p+ 2t}{n h }}
\#
with $\PP^*$-probability at least $1-e^{-t}$.

Next, we bound the data-dependent quantity $\| \Lambda_{r,h} \|_\infty $.
For any $\epsilon_1, \epsilon_2 \in (0,1)$, there exist $\epsilon_1$- and $\epsilon_2$-nets $\{ \bu_1, \ldots,  \bu_{d_1}\}$ and $\{ \bv_1,\ldots, \bv_{d_2}\}$ of $\mathbb{S}^{p-1}$ with $d_1 \leq (1+ 2/\epsilon_1)^p$ and $d_2 \leq (1+2/\epsilon_2)^p$. Given $\bu, \bv \in \mathbb{S}^{p-1}$, there exist  some $1\leq j\leq d_1$ and $1\leq k\leq d_2$ such that $\| \bu- \bu_j \|_2 \leq \epsilon_1$ and $\| \bv - \bv_k \|_2 \leq \epsilon_2$. At $(\bu,\bv)$, consider the decomposition 
\#
 \Lambda_{r,h}(\bu , \bv ) = \Lambda_{r,h}(\bu , \bv) -\Lambda_{r,h}(\bu , \bv_k )  +   \Lambda_{r,h}(\bu , \bv_k ) .\nn
\#
%\$
%	&  | \Lambda_{r,h}(\bu , \bv) -\Lambda_{r,h}(\bu_j, \bv_k )  | =   | \Lambda_{r,h}(\bu , \bv) - \Lambda_{r,h}(\bu , \bv_k)- \Lambda_{r,h}(\bu , \bv_k) - \Lambda_{r,h}(\bu_j, \bv_k ) |  \\
%& \leq  \Bigg|  \frac{1}{n} \sn   \big\{  K_h(r \bw_i^\T \bv - \varepsilon_i) - K_h(r \bw_i^\T \bv_k - \varepsilon_i) \big\} ( \bw_i^\T \bu)^4 \Bigg| + 
%\Bigg|  \frac{1}{n} \sn  K_h(r \bw_i^\T \bv_k - \varepsilon_i)   \big\{ ( \bw_i^\T \bu)^4 -  ( \bw_i^\T \bu_j)^4 \big\} \Bigg|  \\ 
%& ~~~~~~+  |  \Lambda_{r,h} (\bu_j, \bv_k) -  \EE \Lambda_{r,h} (u_j , v_k)  | + \Bigg|  \frac{1}{n} \sn  (1-\EE)  K_h(r \bw_i^\T \bv_k - \varepsilon_i)   ( \bw_i^\T \bu_j)^4  \Bigg|  \\
%& = : I_1 + I_2 + I_3 + I_4. 
%\$
For $ \Lambda_{r,h}(\bu , \bv) -\Lambda_{r,h}(\bu , \bv_k )$,  the Lipschitz continuity of $K(\cdot)$ ensures that
\#
	& |  \Lambda_{r,h}(\bu , \bv) -\Lambda_{r,h}(\bu , \bv_k ) |  \nn \\
	& \leq \frac{l_K r}{ n h^2 } \sn |\bw_i^\T (\bv -\bv_k ) |  ( \bw_i^\T \bu)^4 \leq    \frac{ l_K    r \epsilon_2}{h^2}  \cdot \max\nolimits_{1\leq i\leq n} \| \bw_i \|_2 \cdot   \frac{1}{n} \sn  ( \bw_i^\T \bu)^4 .   \label{decom1}
\#
For $\Lambda_{r,h}(\bu , \bv_k ) $, by the triangle inequality for the $\ell_4$-norm we have
\$
 \Lambda_{r,h}^{1/4}(\bu , \bv_k )   &= \Bigg\{   \frac{1}{n} \sn  K_h(r \bw_i^\T \bv_k - \varepsilon_i)   ( \bw_i^\T \bu)^4\Bigg\}^{1/4} \\
 & \leq  \Bigg\{   \frac{1}{n} \sn  K_h(r \bw_i^\T \bv_k - \varepsilon_i)   ( \bw_i^\T \bu_j)^4\Bigg\}^{1/4} + \Bigg\{   \frac{1}{n} \sn  K_h(r \bw_i^\T \bv_k - \varepsilon_i)   \langle \bw_i , \bu - \bu_j \rangle^4\Bigg\}^{1/4} \\
 & \leq \Lambda^{1/4}_{r,h}(\bu_j,\bv_k)  +  \epsilon_1 \cdot \sup\nolimits_{\bu \in \mathbb{S}^{p-1}} \Lambda_{r,h}^{1/4}(\bu , \bv_k ) ,
\$
which in turn implies
\#
\sup\nolimits_{\bu \in \mathbb{S}^{p-1}} \Lambda_{r,h}(\bu , \bv_k )  \leq (1-\epsilon_1)^{-4} \max\nolimits_{1\leq j\leq d_1}  \Lambda_{r,h}(\bu_j,\bv_k) . \label{decom2}
\#
In view of \eqref{decom1} and \eqref{decom2}, it suffices to bound $\max\nolimits_{1\leq i\leq n} \| \bw_i \|_2$, $\sup_{\bu\in \mathbb{S}^{ p-1}} (1/n) \sn (\bw_i^\T \bu)^4$ and $\max_{(j,k) \in[d_1] \times [d_2]} \Lambda_{r,h}(\bu_j,\bv_k)$.

\begin{lemma} \label{lem:max.l2norm}
For any $t\geq 0$, $\max_{1\leq i\leq n} \| \bw_i \|_2^2 \leq C_1 \upsilon_1^2 (p+ \log n + t )$ with probability at least $1-e^{-t}$, where $C_1>0$ is an absolute constant.
\end{lemma}

%By \eqref{max.predictor.bound}, the bound
%\$
%\max\nolimits_{1\leq i\leq n} \| \bw_i \|_2 \lesssim  (p + \log n + t)^{1/2}  ~\mbox{ holds on event $\cE_{\max}(t)$}.
%\$
For the supremum $\sup_{\bu\in \mathbb{S}^{p-1}} (1/n) \sn (\bw_i^\T \bu)^4$,  similarly to \eqref{decom2} it can be shown that
\$
\sup\nolimits_{\bu \in \mathbb{S}^{p-1}}  \frac{1}{n} \sn (\bw_i^\T \bu)^4 \leq (1-\epsilon_1)^{-4} \max\nolimits_{1\leq j\leq d_1 } \frac{1}{n} \sn (\bw_i^\T \bu_j)^4.
\$
Fix $j$ and $k$, Condition~\ref{cond.predictor2} implies
\$
 \EE e^{ \{ (\bw_i^\T \bu_j)^4/( 36\upsilon_1^4 ) \}^{1/2}}  & =  \EE e^{ (\bw_i^\T \bu_j)^2/ (6\upsilon_1^2 ) } = 1 + \int_0^\infty e^u \PP\big\{  |\bw_i^\T \bu_j |  \geq  \upsilon_1 (6u)^{1/2} \big\} {\rm d} u \\
& \leq 1 + 2 \int_0^\infty    e^{u- 3u }{\rm d} u = 1 + 1 = 2 .
\$
Therefore, $\| (\bw_i^\T \bu_j)^4 \|_{\psi_{1/2}} \leq 36 \upsilon_1^4$, where $\| \cdot \|_{\psi_r}$ denotes the $\psi_r$-norm $(r>0)$; see Definition~2.1 in \cite{Adam2011}. Since $0\leq K_h(r \bw_i^\T \bv_k - \varepsilon_i) \leq \kappa_u h^{-1}$, it is easy to see that $\|K_h(r \bw_i^\T \bv_k - \varepsilon_i)  (\bw_i^\T \bu_j)^4 \|_{\psi_{1/2}} \leq 36 \kappa_u \upsilon_1^4 h^{-1}$. 
Moreover, note that $\EE  (\bw_i^\T \bu_j)^4 \leq m_4$ and 
\$
\EE \big\{  K_h(r \bw_i^\T \bv_k - \varepsilon_i) (\bw_i^\T \bu_j)^4  \big\} = \EE \big[ \EE \big\{ K_h(r \bw_i^\T \bv_k - \varepsilon_i) | \bx_i \big\}  (\bw_i^\T \bu_j)^4 \big] \leq \bar f m_4 .
\$
Hence, for any $u, v \geq 3$, applying inequality (3.6) (and those above it) with $s=1/2$ in \cite{Adam2011} and the union bound, we obtain that  
\$
 \PP \left\{ \max_{1\leq j\leq d_1} \frac{1}{n} \sn (\bw_i^\T \bu_j)^4 \geq  m_4 + C_2 \upsilon_1^4 \, \Bigg(  \sqrt{\frac{u}{n}} +  \frac{u^2}{n} \Bigg) \right\} \leq d_1 e^{-u}
\$
and
\$
 \PP \left\{   \max_{(j,k)\in[d_1]\times [d_2]}  \Lambda(\bu_j,\bv_k) \geq   \bar f m_4 + C_2 \kappa_u \upsilon_1^4   \, \Bigg(  \sqrt{\frac{v}{n h^2 }} +  \frac{v^2}{n h} \Bigg) \right\} \leq  d_1 d_2 e^{-v} .
\$
Taking  $\epsilon_1 = 1-2^{-1/4}$, $\epsilon_2 =n^{-2}$, $u=p\log(1+2/\epsilon_1)+t$ and $v= p \log\{ (1+2/\epsilon_1)   (1+2/\epsilon_2) \}  +t$ in the above bounds, it follows that with probability at least $1-2e^{-t}$,
\$
	\max_{1\leq j\leq d_1} \frac{1}{n} \sn (\bw_i^\T \bu_j)^4 \leq  m_4 + C_3 \upsilon_1^4  \Bigg\{ \sqrt{\frac{p+t}{n}} + \frac{(p+t)^2}{n} \Bigg\}
\$
and 
\#
  \max_{(j,k)\in[d_1]\times [d_2]}  \Lambda(\bu_j,\bv_k) \leq    \bar f m_4 + C_4 \kappa_u \upsilon_1^4   \Bigg\{   \sqrt{\frac{ p \log n +t}{n h^2 }} +  \frac{(p\log n + t)^2}{n h} \Bigg\} .\nn
\#
Denote by $\cE_3(t)$ the event that the above two bounds and the bound in Lemma~\ref{lem:max.l2norm} are satisfied.
Then $\PP\{ \cE_3(t) \} \geq 1- 3e^{-t}$. Conditioned on $\cE_3(t)$ with $t=\log n$, we have
%For any $r\lesssim h$, conditioned on $\cE_{\max}(t) \cap \cE_4(t)$ with $t=\log n$, we have
\#
  \| \Lambda_{r,h}  \|_\infty & \leq 2 \max_{(j,k)\in[d_1]\times [d_2]}  \Lambda(\bu_j,\bv_k) +   \frac{ 2 l_K r   }{   n^2  h^2   } \cdot \max\nolimits_{1\leq i\leq n} \| \bw_i \|_2 \cdot  \sup_{1\leq j\leq d_1} \frac{1}{n} \sn  ( \bw_i^\T \bu_j) ^4    \nn \\
& \leq 2 \bar f m_4 +  C_5 \frac{\upsilon_1^4}{h} \Bigg\{    \sqrt{\frac{ p  \log n}{n   }} +  \frac{(p \log n)^2}{n } \Bigg\}
+ C_6  \frac{ \upsilon_1^5 r}{(n h)^2}   (p+\log n)^{1/2} . \label{Lambda.ubd}
\#

With the above preparations, we are ready to prove the Bahadur representation for the bootstrap estimate. 
Let $\cE(t)=\cE_1(t) \cap \cE_2(t)$ be the event from Theorem~\ref{thm:boot.concentration} and its proof.
%Let $\cE_1(t)$ and $\cE_2(t)=\cE_{\max}(t) \cap \cE_{\smo}(t) \cap \cE_{\lsc}(t)$ be the events defined in the proof of Lemmas~\ref{lem:boot.score} and \ref{lem:boot.rsc}. 
In the rest of the proof, we take $t=\log n$, and set the bandwidth $h\asymp (q/n)^{2/5}$ with $q= p + \log n$. 
Recall that $\hat \bdelta = \hat \bbeta_h - \bbeta^*$ and $\hat \bdelta^\flat = \hat \bbeta^\flat_h - \bbeta^*$.
Conditioned on $\cE_1(t)\cap \cE_2(t)$, 
$
 \| \hat \bdelta  \|_{\bSigma} \leq r_{\est} \asymp  \sqrt{q/n} 
$, and $ \| \hat \bdelta^\flat   \|_{\bSigma} \leq r_{\est}^\flat \asymp \sqrt{q/n}$   with $\PP^*$-probability at least $1-2n^{-1}$. Conditioned further on $\cE_3(t)$,
$$ \| \bSigma^{-1/2}  \Jb_h ( \hat \bbeta_h - \bbeta^*) + \bSigma^{-1/2} \nabla \hat Q_h(\bbeta^* ) \|_{2}  = \| \Delta(\hat \bdelta  ) \|_2 \leq  \sup_{\bdelta \in  \Theta(r_{\est})} \| \Delta(\bdelta) \|_2   \lesssim   \bigg( \frac{q}{n} \bigg)^{4/5},$$
 and with $\PP^*$-probability at least $1-3n^{-1}$,
\$
 &  \| \bSigma^{-1/2} \Jb_h ( \hat \bbeta_h^\flat - \bbeta^*) +  \bSigma^{-1/2}  \nabla \hat Q_h^\flat(\bbeta^* ) \|_2  = \| \Delta^\flat(\hat \bdelta^\flat) \|_2    \\ 
 &\leq  \sup_{\bdelta \in  \Theta(r^\flat_{\est})} \| \Delta^\flat(\bdelta) \|_2  \lesssim \bigg( \frac{q}{n} \bigg)^{4/5}  \bigvee \bigg( \frac{q}{n} \bigg)^{3/5} \bigg( \frac{p \log n}{n} \bigg)^{1/4} \bigvee \bigg( \frac{q}{n} \bigg)^{3/5} \frac{ p  \log n }{n^{1/2}}   .
\$
Together, the above two bounds proves the claimed result. \qed

\subsubsection{Proof of Lemma~\ref{lem:max.l2norm}}
 
 Note that, after centering, $\bw_i = (1, \bw_{i,-}^\T)^\T$, where $\bw_{i,-}\in \RR^{p-1}$ is a zero-mean sub-Gaussian random vector. Under Condition~\ref{cond.predictor2}, there exists some constant $\upsilon_2 \asymp \upsilon_1$ depending only on $\upsilon_1$ such that $\EE \exp\{\balpha^\T (\bw-\be_1 ) \} \leq \exp(  \| \balpha \|_2^2  \upsilon_2^2/2 )$ for all $\balpha \in \RR^{p}$, where $\be_1 = (1,0,\ldots, 0)^\T$. Then, applying Theorem~2.1 in \cite{HKZ2012} with $\Sigma=A= \Ib_p$ yields that, for any $u\geq 0$,
\$
 	 \|  \bw_i \|_2^2 \leq \upsilon_2^2 \big(   p + 2\sqrt{p u } + 2 u \big) + 1 + 2(u/p)^{1/2} \leq \upsilon_2^2 \big(   2 p +  3 u \big) + 1 + 2(u/p)^{1/2}
\$
holds with probability at least $1-e^{-u}$. Taking the union bound over $i=1,\ldots, n$ and setting $u=\log n  + t$ prove the claimed bound. \qed
 
\subsection{Proof of Proposition~\ref{thm:var-cov}}
 
 Consider the change of variable $\bv = \bSigma^{1/2} \bdelta$, so that $\bdelta \in \Theta(r)$ is equivalent to $\bv \in \BB^p(r)$.
Write $\bw_i = \bSigma^{-1/2} \bx_i \in \RR^p$, which are isotropic random vectors, and define
\#
	 \hat \Hb_h(\bv ) = \bSigma^{-1/2}  \hat \Jb_h(\bdelta) \bSigma^{-1/2} = \frac{1}{n} \sn K_h(\varepsilon_i - \bw_i^\T \bv ) \bw_i \bw_i^\T  , \quad \Hb_h(\bv ) = \EE \bigl\{ \hat \Hb_h(\bv) \bigr\}.
\#
For any $\epsilon_1 \in (0,r)$, there exists an $\epsilon$-net $\cN_1 :=\{ \bv_1, \ldots, \bv_{d_1 }\} \subseteq \BB^p(r)$ with $d_1 \leq (1+ 2r/\epsilon_1)^p$ satisfying that, for every $\bv \in \BB^p(r)$, there exists some $1\leq j\leq d_1$ such that $\| \bv - \bv_j \|_2 \leq \epsilon_1$. Hence,
 \#
	&	 \| \hat \Hb_h(\bv ) -    \Hb_h( \textbf{0} )    \|_2  \nn \\
	& \leq   \|\hat \Hb_h(\bv ) - \hat \Hb_h(\bv_j)   \|_2 +    \| \hat \Hb_h(\bv_j)  -  \Hb_h(\bv_j)     \|_2 +    \| \Hb_h(\bv_j) -   \Hb_h( \textbf{0} )   \|_2 \nn \\
	& =: I_1(\bv ) +  I_{2,j} + I_{3,j}. \nn
\# 
For $I_1(\bv ) $, note that $K_h(u) = (1/h)K(u/h)$ is Lipschitz continuous, i.e. $|K_h(u) -K_h(v)| \leq  l_K h^{-2} |u-v|$ for all $u,v \in \RR$.
It follows that
\#
I_1(\bv )  & \leq \sup_{\bu, \bu'  \in \mathbb{S}^{p-1}} \frac{1}{n} \sn   | K_h(\varepsilon_i - \bw_i^\T \bv) - K_h(\varepsilon_i -\bw_i^\T \bv_j )   | \cdot  | \bw_i^\T \bu \cdot \bw_i^\T \bu'    |  \nn \\
  & \leq  l_K h^{-2}\sup_{\bu , \bu' \in \mathbb{S}^{p-1}} \frac{1}{n} \sn   | \bw_i^\T ( \bv -\bv_j ) \cdot   \bw_i^\T \bu  \cdot  \bw_i^\T \bu'   | \nn \\
 & \leq   l_K h^{-2} \epsilon_1  \underbrace{  \sup_{\bu \in  \mathbb{S}^{p-1} }  \frac{1}{n} \sn | \bw_i^\T \bu |^3 }_{=: M_{n,3}}  . \label{I1.ubd1}
\#
Next, we use the standard covering argument to bound $M_{n,3}$. Given $\epsilon_2 \in(0,1)$, let $\cN_2$ be an $\epsilon_2$-net  of the unit sphere $\mathbb{S}^{p-1}$ with $d_2 := | \cN_2 | \leq (1+2/\epsilon_2)^p$ such that for every $\bu \in \mathbb{S}^{p-1}$, there exists some $\bu' \in \cN_2$ satisfying $\| \bu - \bu' \|_2 \leq \epsilon_2$. Define the (standardized) design matrix $\Wb_n =  n^{-1/3}(\bw_1,\ldots, \bw_n)^\T  \in \RR^{n\times p}$, so that $M_{n,3} = \sup_{\bu \in \mathbb{S}^{p-1}} \| \Wb_n \bu \|_3^3$.
 By the triangle inequality,
\begin{align}
	  \|  \Wb_n  \bu  \|_3 &   \leq  \|  \Wb_n \bu'   \|_3 +  \| \Wb_n ( \bu - \bu')  \|_3 \nn \\
	& =    \| \Wb_n \bu'  \|_3 +\bigg(   \frac{1}{n} \sn   | \bw_i^\T (  \bu - \bu' ) |^3 \bigg)^{1/3}  \leq    \| \Wb_n \bu'  \|_3 + \epsilon_2 \cdot  M_{n,3}^{1/3} . \nn
\end{align}
Taking the maximum over $\bu' \in \cN_2$, and then taking the supremum over $\bu \in \mathbb{S}^{p-1}$, we arrive at
\#
	M_{n,3}  \leq   (1-\epsilon_2)^{-3}  \cdot  N_{n, 3}   := (1-\epsilon_2)^{-3}  \cdot \max_{\bu' \in \cN_2 } \frac{1}{n}\sn | \bw_i^\T \bu' |^3 .  \label{discretization}
\#	 

For every $\bu' \in \cN_2$, note that
\$
 \EE e^{ \{ |\bw_i^\T \bu' |^3/( 6^{3/2}\upsilon_1^3) \}^{2/3}  } = 1 + \int_0^\infty e^{u} \PP\big\{  |\bw_i^\T \bu' | \geq \upsilon_1 (6u)^{1/2}  \big\} {\rm d} u  \leq 1 + 2 \int_0^\infty e^{-2u} {\rm d} u = 2 , 
\$
implying $\|   |\bw_i^\T \bu' |^3 \|_{\psi_{2/3}} \leq 6^{3/2} \upsilon_1^3$. Hence, by inequality (3.6) in \cite{Adam2011} with $s=2/3$, we obtain that for any $z\geq 3$,
$$
 \frac{1}{n} \sn |\bw_i^\T \bu' |^3 \leq   \EE  |\bw^\T \bu' |^3  + C_1 \upsilon_1^3 \Biggl(  \sqrt{\frac{z }{n}}  + \frac{z^{3/2}}{n} \Bigg)
$$
with probability at least $1-e^{ -z}$.
Taking the union bound over all vectors $\bu'$ in $\cN_2$ yields that, with probability at least $1-d_2 e^{ -z} \geq 1- e^{ p\log(1+2/\epsilon_2) - z}$,
\$
	N_{n,3}   \leq m_3   + C_1 \upsilon_1^3 \Biggl(  \sqrt{\frac{z }{n}}  + \frac{z^{3/2}}{n} \Bigg)
\$
where $m_3 = \sup_{\bu \in \mathbb{S}^{p-1} } \EE |\bw^\T \bu|^3$. Reorganizing the terms, we get
\begin{align}
	N_{n,3}   \leq  m_3 + C_1  \upsilon_1^3 \Bigg[ \sqrt{\frac{p\log (1+2/\epsilon_2) +\log 2 + t }{n}} + \frac{\{ p\log(1+2/\epsilon_2) + \log 2 + t \}^{3/2}}{n} \Bigg]  \label{discrete.concentration}
\end{align}
with probability at least $1- e^{- t}/2$. Taking $\epsilon_2 =1/8$ in \eqref{discretization} and \eqref{discrete.concentration} implies
\#
 M_{n,3} \leq  1.5 m_3  +  1.5 C_1 \upsilon_1^3   \Bigg\{ \sqrt{\frac{3p+ 1 +t }{n}} + \frac{ ( 3p + 1+ t )^{3/2}}{n} \Bigg\} . \nn
\#
Under the sample size scaling $n\gtrsim p+t$, plugging the above bound into \eqref{I1.ubd1} yields 
\#
  \sup_{\bv \in \BB^p(r)}  I_1(\bv )   \lesssim \upsilon_1^3 (p + t)^{1/2}  h^{-2} \epsilon_1 \label{I1.ubd2}
\#
with probability at least $1- e^{-t}/2$.

Turning to $I_{2,j}$, note that $\hat \Hb_h(\bv_j) - \Hb_h(\bv_j )  = (1/n) \sn (1-\EE) \phi_{ij} \bw_i \bw_i^\T$, where $\phi_{ij} =K_h(\varepsilon_i - \bw_i^\T  \bv_j )$ satisfy $|\phi_{ij}| \leq \kappa_u h^{-1}$ and
\# 
\EE \bigl(  \phi_{ij}^2 | \bx_i \bigr)     = \frac{1}{h^2} \int_{-\infty}^\infty K^2 \Biggl( \frac{  \langle  \bw_i , \bv \rangle - t}{h} \Biggr) f_{\varepsilon_i  | \bx_i  } (t) \, {\rm d} t   = \frac{1}{h} \int_{-\infty}^\infty K^2(u) f_{\varepsilon_i  | \bx_i  } (\bw_i^\T \bv - b u  ) \, {\rm d} u  \leq  \bar f \kappa_u h^{-1}  \nn 
\#
almost surely.
Given $\epsilon_3\in (0,1/2)$, there exits an $\epsilon_3$-net $\cN_3$ of the sphere $\mathbb{S}^{p-1}$ with $|\cN_3 |\leq (1+2/\epsilon_3)^p$ such that $\| \hat  \Hb_h(\bv_j) -\Hb_h(\bv_j )\|_2 \leq (1-2\epsilon_3)^{-1} \max_{\bu \in \cN_3} |  \bu^\T \{ \hat \Hb_h(\bv_j) -\Hb_h(\bv_j ) \}\bu   |$. Given $\bu \in \cN_3$ and $k=2,3,\ldots$, we bound the higher order moments of $\phi_{ij}   (\bw_i^\T  \bu )^2$ by
\#
	\EE |\phi_{ij}   (\bw_i^\T  \bu )^2 |^k &  \leq  \bar f \kappa_u h^{-1} \cdot   (\kappa_u h^{-1})^{k-2}   \upsilon_1^{2k}   \cdot 2k \int_0^\infty \PP\bigl( |\bw_i^\T \bu  | \geq\upsilon_1 u \bigr) u^{2k-1} {\rm d} u  \nn \\
	& \leq  \bar f \kappa_u h^{-1}  \cdot   ( \kappa_u h^{-1})^{k-2}   \upsilon_1^{2k} \cdot 4k  \int_0^\infty  u^{2k-1} e^{-u^2/2} {\rm d} u   \nn \\
	& \leq  \bar f \kappa_u h^{-1}  \cdot   (\kappa_u h^{-1})^{k-2}   \upsilon_1^{2k}  \cdot  2^{k+1} k! .\nn
\#
In particular, $\EE  \phi_{ij}^2 (\bw_i^\T \bu)^4 \leq (4\upsilon_1^2)^2  \bar f  \kappa_u h^{-1} $, and $\EE |\phi_{ij} (\bw_i^\T \bu)^2 |^k\leq \frac{k!}{2} \cdot ( 4 \upsilon_1^2)^2  \bar f  \kappa_u h^{-1} \cdot (2\upsilon_1^2 \kappa_u  h^{-1} )^{k-2}$ for $k\geq 3$. Applying  Bernstein's inequality and the union bound, we find that for any $u\geq 0$,
\#
 &   \| \hat \Hb_h (\bv_j) -  \Hb_h (\bv_j )   \|_2 \nn \\
 & \leq \frac{1}{1-2\epsilon_3 } \max_{\bu \in \cN_3} \Biggl|   \frac{1}{n} \sn (1-\EE) \phi_{ij}  (\bw_i^\T \bu)^2 \Biggr| \leq \frac{ 2\upsilon_1^2}{1- 2\epsilon_3} \Biggl(   2 \sqrt{2   \bar f \kappa_u  \frac{u}{n h }} +\kappa_u  \frac{u}{n h} \Biggr) \nn 
\#
with probability at least $1-2 (1+2/\epsilon_3)^p e^{-u} = 1 - (1/2) e^{\log(4) + p\log(1+2/\epsilon_3) - u}$. Setting $\epsilon_3=2/(e^3-1)$ and $u= \log(4) + 3p + v$, it follows that with probability at least $1 -e^{-v}/2$,
\#
	I_{2,j}\lesssim  \upsilon_1^2 \Bigg(   \sqrt{  \frac{ p + v  }{n h} }  + \frac{ p + v }{n h}   \Bigg) . \nn
\#
Once again,  taking the union bound over $j=1,\ldots, d_1$ and setting $v= p \log(1 + 2r/\epsilon_1)+ t$, we obtain that with probability at least $1-d_1 e^{-v} \geq 1-e^{-t}/2$, 
\#
\max_{1\leq j\leq d_1 } I_{2,j} \lesssim \sqrt{\frac{p\log(3 e r/\epsilon_1)+t}{n h }} +   \frac{p\log(3e r/\epsilon_1)+t}{n h}  . \label{I2.ubd}
\#

For $I_{3,j}$, following the proof of \eqref{residual.bound} it can similarly shown that $I_{3,j} \leq  0.5 l_0  m_3  r$. Combining this with  \eqref{I1.ubd2} and \eqref{I2.ubd}, and taking $\epsilon_1 =     r/n \in (0, r)$ in the beginning of the proof, we conclude that with probability at least $1-  e^{-t}$,
\#
 & \sup_{\bdelta  \in \Theta(r)} \|  \bSigma^{-1/2}  \{  \hat \Jb_h(\bdelta) - \Jb_h \}\bSigma^{-1/2}   \|_2  \lesssim \sqrt{\frac{p \log n+ t}{n h }} +  \frac{p\log n +t}{n h } + \frac{(p + t)^{1/2}  r }{n h^2} +      r   \nn
\#
as long as $n\gtrsim p+t$. This leads to \eqref{J.consistency} under the prescribed bandwidth constraint. 

To derive the same bound for $\hat \Vb_h(\bdelta) - \Vb_h$, notice that $u \mapsto \{\cK_h(u) - \tau \}^2$ is a $(2\bar \tau \kappa_u/h)$-Lipschitz continuous function, where $\bar \tau =\max(\tau, 1-\tau)$. Moreover, for every $\bdelta \in  \RR^p$, the random variable $\{\cK_h(\langle \bx_i, \bdelta \rangle - \varepsilon_i) - \tau \}^2$ takes values in $[0 , \bar \tau^2]$. We can thus apply the same argument to bound $\sup_{\bdelta \in \Theta(r)} \| \bSigma^{-1/2} \{ \hat \Vb_h(\bdelta) - \Vb_h \} \bSigma^{-1/2} \|_2$.
\qed
 
 %%%%%%%%%%%%% %%%%%%%%%%%%% %%%%%%%%%%%%%
% Higher-order
%%%%%%%%%%%%% %%%%%%%%%%%%% %%%%%%%%%%%%%
\section{Theoretical properties of one-step conquer}

In this section, we provide theoretical properties of the one-step conquer estimator $\hat \bbeta$, defined in Section~\ref{os.conquer1}. The key message is that, when higher-order kernels are used (and if the conditional density $f_{\varepsilon|\bx}(\cdot)$ has enough derivatives), the asymptotic normality of the one-step estimator holds under weaker growth conditions on $p$. For example, the scaling condition $p=o(n^{3/8})$ that is required for the conquer estimator can be  reduced to  roughly $p=o(n^{7/16})$ for the one-step conquer estimator using a kernel of order 4.

\begin{condition} \label{cond:higher-order.kernel}
Let $G(\cdot)$ be a symmetric kernel of order $\nu>2$, that is, $\int_{-\infty}^\infty u^k G(u) \, \dd u=0$ for $k=1,\ldots, \nu-1$ and $\int_{-\infty}^\infty u^\nu G(u) \, {\rm d} u \neq 0$. Moreover, $g_k  = \int_{-\infty}^\infty |u^k G(u) | \, \dd u <\infty$ for $1\leq k\leq \nu$, $G$ is uniformly bounded with $g_u= \sup_{u\in \RR} |G(u) | <\infty$ and  is $l_G$-Lipschitz continuous for some $l_G>0$.
\end{condition}
 
 As before, we write $\cG_b(u) = \cG(u/b)$ and $\cG(u) = \int_{-\infty}^u G(v)\,{\rm d v}$ for $u\in \RR$ and $b>0$.
 The use of a higher-order kernel does not necessarily reduce bias unless the  conditional density $ f_{\varepsilon | \bx}(\cdot)$ of $\varepsilon$ given $\bx$ is sufficiently smooth. Therefore, we further impose the following smoothness conditions on $f_{\varepsilon | \bx}(\cdot)$.
 
%\begin{condition} \label{cond:density}
%There exist $\bar f  \geq \underbar{$f$} >0$ such that $ \underbar{$f$}  \leq f_{\varepsilon |\bx}(0 ) \leq \bar f $ almost surely (for all $\bx$).
%Moreover, there exists a constant $l_0 >0$ such that $|f_{\varepsilon |\bx}(u)- f_{\varepsilon |\bx}(v)| \leq l_0  |u - v|$ for all $u , v \in \RR$ almost surely.
%\end{condition} 

 \begin{condition} \label{cond:smoothness}
 Let $\nu\geq 4$ be the integer in Condition~\ref{cond:higher-order.kernel}.
 The conditional density $ f_{\varepsilon | \bx}(\cdot )$ is $(\nu-1)$-times differentiable, and satisfies $| f_{\varepsilon | \bx}^{(\nu-2)} (u) - f_{\varepsilon | \bx}^{(\nu-2)} (0) | \leq l_{\nu-2} |u|$ for all $u\in \RR$ almost surely (over the random vector $\bx$), where $l_{\nu-2}>0$ is a constant. Also,   there exists some constant $C_G>0$ such that
$ \int_{-\infty}^\infty | u^{\nu-1} G(u)| \cdot \sup_{|t| \leq |u|}  |  f_{\varepsilon | \bx}^{(\nu-1)} (t)- f_{\varepsilon | \bx}^{(\nu-1)} (0)   | \, \dd u \leq C_G$ almost surely.
 \end{condition}

% \begin{condition}[Random design] \label{cond.covariate}
%The predictor $\bx \in \RR^p$ is {\it sub-Gaussian}: there exists $ \upsilon_1>0$ such that $\PP( | \bw^\T \bu |    \geq \upsilon_1    t  ) \leq  2e^{-t ^2/2}$ for every unit vector $\bu \in \mathbb{S}^{p-1}$ and $t \geq 0$, where $\bw = \bSigma^{-1/2} \bx$ and $\bSigma = \EE(\bx \bx^\T)$ is positive definite. In addition, write  $ m_k = \sup_{\bu \in \mathbb{S}^{p-1}  }  \EE |  \bw^\T \bu  |^k$ for $k\geq 1$.
%\end{condition}
 
% Define the $p\times p$ matrix $\Jb = \EE \{ \bx \bx^\T f_{\varepsilon | \bx} (0) \}  $, and the empirical and population smoothed loss
% \#
% 	\hat Q^G_b( \bbeta) = \frac{1}{n} \sn (\rho_\tau * G_b) (y_i - \bx_i^\T \bbeta) ~~\mbox{ and }~~ Q^G_b( \bbeta) = \EE \bigl\{ \hat Q^G_b( \bbeta) \bigr\},
% \#
% where $h>0$ is a bandwidth. 

Notably, we have
 \#
 	\nabla Q^G_b( \bbeta) = \EE \bigl\{  \cG_b\big(   \langle \bx , \bbeta \rangle  - y \big)  - \tau \bigr\} \bx  ~~\mbox{ and }~~
 	\nabla^2 Q^G_b (\bbeta ) = \EE \big\{ G_b( y- \langle \bx , \bbeta \rangle ) \bx \bx^\T \big\} ,
 \#
 representing the population score and Hessian of $Q^G_b (\cdot) = \EE \hat Q^G_b(\cdot)$. As $b \to 0$, we expect $\nabla Q^G_b( \bbeta^*)$ and  $\nabla^2 Q^G_b (\bbeta^*)$ to converge to $\textbf{0}$ (zero vector in $\RR^p$) and $\Jb = \EE \{  f_{\varepsilon | \bx} (0) \bx \bx^\T \}$, respectively.
 %For any symmetric matrix $\Ab \in \RR^{p\times p}$, with slight abuse of notation we use $\| \cdot \|_{\bOmega}$ to denote the relative operator norm, defined as $\| \Ab \|_{\bOmega} =\| \bSigma^{-1/2} \Ab \bSigma^{-1/2} \|_2$. 
 The following proposition validates this claim by providing explicit error bounds.

%which extends Lemma~1 in \cite{FGH2019} with explicit error bounds, validates this claim. 

\begin{proposition} \label{prop:higher-order.bias} 
Let $b\in (0, 1)$ be a bandwidth. Under Conditions~\ref{cond:higher-order.kernel} and \ref{cond:smoothness}, we have
 \#
 \|  \bSigma^{-1/2} \nabla Q^G_b(\bbeta^*) \|_2 & \leq  l_{\nu-2}  g_{ \nu} \, b^\nu / \nu!  
~~\mbox{ and }~~
 \|     \bSigma^{-1/2} \nabla^2  Q^G_b(\bbeta^*) \bSigma^{-1/2} - \Hb  \|_2   \leq  C_G  \,    b^{\nu-1} / (\nu-1)! , \nn
 \#
 where $\Hb = \bSigma^{-1/2} \Jb \bSigma^{-1/2} = \EE \{ f_{\varepsilon |\bx} (0)  \bw \bw^\T\}$ with $\bw = \bSigma^{-1/2} \bx$.
\end{proposition}
 
Proposition~\ref{prop:Hessian.uniform} shows that when a higher-order kernel is used, the bias is significantly reduced in the sense that 
$\| \nabla Q^G_b(\bbeta^*)  \|_2 = \cO( b^\nu)$  and $\| \nabla^2  Q^G_b(\bbeta^*) - \Jb \|_2 = \cO(b^{\nu-1})$, where $\nu\geq 4$ is an even integer. Notably, even if the kernel $G$ has negative parts, the population Hessian $\nabla^2  Q^G_b(\bbeta^*) $ preserves the positive definiteness of $\Jb$ as long as the bandwidth $b$ is sufficiently small.  
 
To construct the one-step conquer estimator,  two key quantities are the sample Hessian $\nabla^2 \hat Q^G_b( \cdot ) $ and sample gradient $\nabla \hat Q^G_b(\cdot)$, both evaluated at $\overbar \bbeta$, a consistent initial estimate. In the next two propositions, we establish  uniform convergence results of the Hessian and gradient of the empirical smoothed loss to their population counterparts. As a direct consequence, $\nabla^2 \hat Q^G_b( \overbar \bbeta)$ is positive definite with high probability, provided that $\overbar \bbeta$ is consistent (i.e., in a local vicinity of $\bbeta^*$). 
To be more specific, for $r>0$, we define the local neighborhood 
\#
	\Theta^*(r) = \bigl\{ \bbeta \in \RR^p : \| \bbeta - \bbeta^* \|_{\bSigma} \leq r \bigr\} . \label{local.region}
\#

 \begin{proposition} \label{prop:Hessian.uniform}
 Conditions~\ref{cond:higher-order.kernel}, \ref{cond:smoothness} and \ref{cond.predictor2} ensure that, with probability at least $1-e^{-t}$,
\#
 & \sup_{\bbeta \in \Theta^*(r)} \| \bSigma^{-1/2}  \{ \nabla^2 \hat Q^G_b(\bbeta ) -  \nabla^2  Q^G_b(\bbeta )  \}\bSigma^{-1/2}   \|_2 \lesssim \sqrt{\frac{p \log n + t}{n b }} +  \frac{p\log n+t}{n  b } + \frac{(p + t)^{1/2}   r }{n b^2}  \nn
\#
as long as $n\gtrsim p+t$.
 \end{proposition}
 
\begin{proposition} \label{prop:gradient.uniform}
 Conditions~\ref{cond:higher-order.kernel}, \ref{cond:smoothness} and \ref{cond.predictor2}   ensure that, with probability at least $1- e^{-t}$,
 \#
	& \sup_{\bbeta \in \Theta^*(r)}   \|   \bSigma^{-1/2}  \{   \nabla \hat Q^G_b(\bbeta ) - \nabla \hat Q^G_b(\bbeta^* )  \} -  \bSigma^{-1/2}  \Jb  ( \bbeta - \bbeta^*)    \|_2  \lesssim r \,\Biggl(  \sqrt{\frac{p+ t}{n b }}  +   r +b^{\nu-1}   \Bigg)  \label{gradient.uniform.bound}
 \#
 as long as $\sqrt{(p+t) /n } \lesssim b $.
 \end{proposition}

 With the above preparations, we are ready to present the Bahadur representation for the one-step conquer estimator $\hat \bbeta$.

\begin{theorem} \label{thm:one-step.conquer}
Assume Conditions~\ref{cond.kernel}, \ref{cond.reg} and \ref{cond.predictor2} in the main text and Conditions~\ref{cond:higher-order.kernel} and \ref{cond:smoothness} hold. For any $t>0$, let the sample size $n$, dimension $p$ and the bandwidths $h , b>0$ satisfy $n\gtrsim p(\log n)^2 +t$, $\sqrt{(p+t)/n} \lesssim h \lesssim \{ (p+t)/n \}^{1/4}$ and $\sqrt{(p+t)/n }   \lesssim b \lesssim  \{ (p+t)/n \}^{1/(2\nu)}$. Then, the one-step conquer estimator $\hat \bbeta$ satisfies the bound
\# \label{os.conquer.br}
\Bigg\|   \bSigma^{-1/2} \Jb( \hat \bbeta - \bbeta^*   ) -   \frac{ 1}{n}  \sn  \bigl\{  \tau - \cG_b(-\varepsilon_i  ) \bigr\}   \bSigma^{-1/2}  \bx_i   \Bigg\|_2 \lesssim  \Biggl\{  \underbrace{  \sqrt{ (p\log n + t)/(n b ) } }_{{\rm variance~ term}}+  \underbrace{  b^{\nu-1} }_{{\rm bias~term}}  \Biggr\}   \sqrt{\frac{p+t}{n}} 
\#
with probability at least $1-5e^{-t}$.
\end{theorem} 
 
 Theorem~\ref{thm:one-step.conquer} shows that using a higher-order kernel ($\nu\geq 4$) allows one to choose larger bandwidth, thereby reducing the ``variance" and the total Bahadur linearization error.
Similarly to Theorem~\ref{thm:clt} in the main text, the following asymptotic normal approximation result for linear projections of one-step conquer is a direct consequence of Theorem~\ref{thm:one-step.conquer}.
 
\begin{theorem}  \label{thm:higher-order.clt}
Assume Conditions~\ref{cond.kernel}, \ref{cond.reg} and \ref{cond.predictor2} in the main text and Conditions~\ref{cond:higher-order.kernel} and \ref{cond:smoothness} hold. Let the bandwidths satisfy $( q/n  )^{1/2} \lesssim h \lesssim ( q/n )^{1/4}$ and $(q/n )^{1/2} \lesssim b \lesssim (q/n)^{1/(2\nu)}$, where $q:=p + \log n$.
Then, 
\#
 \sup_{x\in \RR , \, \ba \in \RR^p } \bigl| \PP\bigl(  n^{1/2} \langle \ba, \hat \bbeta  - \bbeta^*   \rangle/\sigma_0 \leq  x \bigr) - \Phi(x ) \bigr| \lesssim   \sqrt{\frac{ (p+\log n)p \log n}{n b}}  + n^{1/2} b^{\nu}  , \label{linear.clt.higher}
\#
where $\sigma_0^2 = \sigma_0^2(\ba) =  \tau (1-\tau ) \|   \Jb^{-1} \ba \|_{ \bSigma}^2$. In particular, with a choice of bandwidth $b \asymp (q/n)^{2/(2\nu +1)}$,
\$
 \sup_{x\in \RR , \, \ba \in \RR^p } \bigl| \PP\bigl(  n^{1/2} \langle \ba, \hat \bbeta  - \bbeta^*   \rangle \leq \sigma_0  x \bigr) - \Phi(x) \bigr|  \to  0   \nn
\$
as $n \to \infty$ under the scaling $p^{4\nu/(2\nu -1)} (\log n)^{(2\nu+1 )/(2\nu-1)} =o(n)$.
\end{theorem}

Let $G(\cdot)$ be a kernel of order $\nu=4$. In view of Theorem~\ref{thm:higher-order.clt}, we take $h \asymp \{(p+\log n)/n\}^{2/5}$ as in the main text and $b = \{(p+\log n)/n\}^{2/9}$, thereby obtaining that $n^{1/2}\langle \ba , \hat \bbeta -\bbeta \rangle$, for an arbitrary $\ba \in \RR^p$, is asymptotically normally distributed as long as $p  (\log n)^{9/16} = o(n^{7/16})$ as $n\to \infty$.

\subsection{Proof of Proposition~\ref{prop:higher-order.bias}}
We start from the gradient $\bSigma^{-1/2}\nabla Q^G_b(\bbeta^*)=\EE  \{ \cG_b ( - \varepsilon ) - \tau \}\bw$ with $\bw=\bSigma^{-1/2} \bx$. Let $\EE_{\bx}$ be the conditional expectation given $\bx$. By integration by parts,
 \#
 \EE_{\bx}\cG_b (  -\varepsilon  ) & = \int_{-\infty}^\infty \cG (-t/b ) \, {\rm d} F_{\varepsilon | \bx} (t) =   \int_{-\infty}^\infty G(u)  F_{\varepsilon | \bx} (- b u )  \, {\rm d} u. \label{cond.mean1}
 \#
Applying a Taylor series expansion with integral remainder on $F_{\varepsilon | \bx} (- b u )$ yields 
\$
F_{\varepsilon | \bx} (- b u ) &= F_{\varepsilon | \bx}(0) + \sum_{\ell =1}^{\nu-1} F_{\varepsilon | \bx}^{(\ell)} (0) \frac{(- b u)^\ell}{\ell!} + \frac{(-b u)^{\nu-1}}{(\nu-2)!} \int_0^1 (1-w)^{\nu-2} \big\{  F_{\varepsilon | \bx}^{(\nu-1)} (- b u w)  - F_{\varepsilon | \bx}^{(\nu-1)} (0) \bigr\} {\rm d} w \\
&=  \tau  + \sum_{\ell =0}^{\nu-2} f_{\varepsilon | \bx}^{(\ell )} (0) \frac{(- b u)^{\ell+1}}{(\ell+1)!} + \frac{(- b u)^{\nu-1}}{(\nu-2)!} \int_0^1 (1-w)^{\nu-2} \big\{ f_{\varepsilon | \bx}^{(\nu-2)} (- b u w)  - f_{\varepsilon | \bx}^{(\nu-2)} (0) \bigr\} {\rm d} w .
\$
Recall that $G$ is a kernel of order $\nu\geq 4$ (an even integer) and $g_\nu= \int_{-\infty}^\infty |u^\nu  G (u) |\, \dd u <\infty$. Substituting the above expansion into \eqref{cond.mean1}, we obtain
 \#
\EE_{\bx}\cG_b (  -\varepsilon   ) = \tau - \frac{b^{\nu-1}}{(\nu-2)!}   \int_{-\infty}^\infty \int_0^1 u^{\nu-1} G(u) (1-w)^{\nu-2} \big\{ f_{\varepsilon | \bx}^{(\nu-2)} (- b u w)  - f_{\varepsilon | \bx}^{(\nu-2)} (0) \bigr\} {\rm d} w {\rm d} u .\nn
 \#
Furthermore, by the Lipschitz continuity of $f_{\varepsilon | \bx}^{(\nu-2)}(\cdot)$ around 0,
\#
     | \EE_{\bx}\cG_b (  -\varepsilon   ) - \tau  |  &  \leq    \frac{l_{\nu-2} b^{\nu}  }{(\nu-2)!}   \int_{-\infty}^\infty \int_0^1 | u^{\nu } G(u)|  (1-w)^{\nu-2} w  \,{\rm d} w {\rm d} u \nn \\
 & =   B(2,\nu-1)   l_{\nu-2}  g_\nu\,  b^{\nu } / (\nu-2)! , \nn
\#
where $B(x,y) := \int_0^1 t^{x-1} (1-t)^{y-1} {\rm d}t$ denotes the beta function. In particular, $B(2,\nu-1)= \Gamma(2)\Gamma(\nu-1)/\Gamma(\nu+1)=  (\nu-2)!/\nu!$. Putting together the pieces yields
\#
   \|  \bSigma^{-1/2}  \nabla Q^G_b(\bbeta^*)   \|_2  & = \sup_{\bu \in \mathbb{S}^{p-1}} \EE\,  \EE_{\bx} \big\{ \cG_b(-\varepsilon  ) - \tau \big\} \bw^\T \bu \leq  l_{\nu-2}  g_\nu\, b^\nu / \nu!   . \nn
\#
 
 Turning to the Hessian, note that
 \#
 & \|  \bSigma^{-1/2} \{   \nabla^2  Q^G_b(\bbeta^*)  -  \Jb \} \bSigma^{-1/2} \|_2  =  \Bigg\|   \EE  \int_{-\infty}^\infty G(u)   \big\{ f_{\varepsilon | \bx} (- b u) - f_{\varepsilon | \bx}(0) \big\} {\rm d}u \, \bw \bw^\T      \Bigg\|_2 .\nn
 \#
 Applying a similar Taylor series expansion as above,  we have
 \#
 f_{\varepsilon | \bx} (t) &= f_{\varepsilon | \bx}(0) + \sum_{\ell =1}^{\nu-1} f_{\varepsilon | \bx}^{(\ell)} (0) \frac{t^\ell}{\ell!} + \frac{t^{\nu-1}}{(\nu-2)!} \int_0^1 (1-w)^{\nu-2} \big\{  f_{\varepsilon | \bx}^{(\nu-1)} (t w )  - f_{\varepsilon | \bx}^{(\nu-1)} (0) \bigr\} {\rm d} w . \label{taylor}
 \#
Under Conditions~\ref{cond:higher-order.kernel} and \ref{cond:smoothness}, it follows that
 \$
 &  \| \bSigma^{-1/2}\{   \nabla^2  Q^G_b(\bbeta^*)     - \Jb \}\bSigma^{-1/2}  \|_{\bOmega}  \\
 & \leq  \frac{b^{\nu-1}}{(\nu-2)!} \sup_{\bu, \bdelta \in \mathbb{S}^{p-1}}    \EE  \int_{-\infty}^\infty \int_0^1   u^{\nu-1}   G(u) (1-w)^{\nu-2}\big\{  f_{\varepsilon | \bx}^{(\nu-1)} (- b  u w)  - f_{\varepsilon | \bx}^{(\nu-1)} (0) \bigr\}    {\rm d} w\,   {\rm d}u \,  \big\langle \bw  , \bu \big\rangle \big\langle \bw , \bdelta  \big\rangle   \\
 & \leq   \frac{b^{\nu-1}}{(\nu-1)!} \sup_{\bu, \bdelta \in \mathbb{S}^{p-1}}    \EE  \int_{-\infty}^\infty  |  u^{\nu-1}   G(u) |  \sup_{|t| \leq b |u|} \big| f_{\varepsilon | \bx}^{(\nu-1)} (t)  - f_{\varepsilon | \bx}^{(\nu-1)} (0) \bigr|    \,   {\rm d}u \cdot    |    \langle \bw , \bu  \rangle \langle \bw , \bdelta \big\rangle | \\
 & \leq \frac{C_G b^{\nu-1 }}{(\nu-1)!}   \sup_{ \bu  \in \mathbb{S}^{p-1}}   \EE \big\langle \bw , \bu \big\rangle^2  = \frac{C_G}{(\nu-1)!}         b^{\nu-1 } .
 \$
 This completes the proof. \qed

\subsection{Proof of Proposition~\ref{prop:Hessian.uniform}}

The proof is almost identical to that of Proposition~\ref{thm:var-cov}, and thus is omitted. \qed

 \subsection{Proof of Proposition~\ref{prop:gradient.uniform}}
 
 Define the stochastic process $\Delta_b(\bbeta) = \bSigma^{-1/2}    \{ \nabla \hat Q^G_b(\bbeta ) - \nabla \hat Q^G_b(\bbeta^* ) -   \Jb ( \bbeta - \bbeta^*) \}$. By the triangle inequality, 
\$
	\sup_{\bbeta \in \Theta^*(r)}   \| \Delta_b(\bbeta)   \|_2 \leq  \sup_{\bbeta \in \Theta^*(r)}   \| \EE  \Delta_b(\bbeta)   \|_2 + \sup_{\bbeta \in \Theta^*(r)}   \| \Delta_b(\bbeta) - \EE  \Delta_b(\bbeta)  \|_2
\$
Recall that $\Hb = \bSigma^{-1/2} \Jb \bSigma^{-1/2} = \EE \{ f_{\varepsilon |\bx}(0) \bw \bw^\T\}$.
 For the first term on the right-hand side, using the mean value theorem for vector-valued functions yields
 \$
 \EE  \Delta_b(\bbeta) = \Biggl\{ \bSigma^{-1/2}  \int_0^1 \nabla^2 Q^G_b\bigl( (1-s) \bbeta^* + s \bbeta  \bigr) \, \dd s \,  \bSigma^{-1/2} - \Jb_0  \Biggr\}  \bSigma^{1/2} ( \bbeta -\bbeta^* ) . 
 \$
By a change of variable $\bdelta = \bSigma^{1/2}(\bbeta - \bbeta^*)$,
 \$
  \nabla^2 Q^G_b\bigl( (1-s) \bbeta^* + s \bbeta  \bigr)  = \EE  \int_{-\infty}^\infty  G(u) f_{\varepsilon | \bx} (s \bw^\T \bdelta - b u ) \, \dd u \cdot \bx \bx^\T .
 \$
For every $s \in [0,1]$ and $u\in \RR$, it ensures from $f_{\varepsilon | \bx}(\cdot)$ being Lipschitz  that 
$$
  | f_{\varepsilon | \bx} (s \bw^\T \bdelta - b u ) - f_{\varepsilon | \bx} ( - b u )   | \leq l_0  s \cdot  | \bw^\T \bdelta | .
$$
Moreover, by the Taylor series expansion \eqref{taylor}, 
\$
  f_{\varepsilon | \bx} (  - b u ) &  =   f_{\varepsilon | \bx}(0) + \sum_{\ell=1}^{\nu-1}   f_{\varepsilon | \bx}^{(\ell)}(0) \frac{(  - b u )^\ell}{\ell!} + \frac{( - b u)^{\nu-1}}{(\nu-2)!} \int_0^1 (1-w)^{\nu-2} \bigl\{   f_{\varepsilon | \bx}^{(\nu-1)}(  -  b u w ) -f_{\varepsilon | \bx}^{(\nu-1)}(0)\bigr\} \, \dd w.
\$
Consequently,
\$
	& \Bigg\|   \bSigma^{-1/2}  \int_0^1  \nabla^2 Q^G_b\bigl( (1-s) \bbeta^* + s\bbeta  \bigr)\, \dd  s \,  \bSigma^{-1/2}   -  \Hb \Bigg\|_2   \\ 
&\leq \Bigg\|  \frac{b^{\nu-1}}{(\nu-2)!}  \EE \int_{-\infty}^\infty \int_0^1  u^{\nu-1} G(u) (1-w)^{\nu-2} \bigl\{   f_{\varepsilon | \bx}^{(\nu-1)}(  - b u w ) -f_{\varepsilon | \bx}^{(\nu-1)}(0)\bigr\} \, \dd w \dd u  \cdot    \bw \bw^\T  \Bigg\|_2 \\
&~~~~ +  0.5 l_0 \cdot \big\|  \EE \,  |\bw^\T  \bdelta |  \cdot    \bw \bw^\T  \big \|_2   \\
& \leq  \frac{C_G b^{\nu-1}}{(\nu-1)!} \sup_{\bu \in \mathbb{S}^{p-1}}  \EE \big\langle \bw , \bu \big\rangle^2 + \frac{l_0}{2} \sup_{\bu \in \mathbb{S}^{p-1}}  \EE \big(  | \bw^\T \bdelta | \langle \bw , \bu  \rangle^2 \big) \\ 
& \leq       \Bigg(   \frac{ C_G }{(\nu-1)! }  b^{\nu-1} + 0.5 l_0 m_3\| \bdelta \|_2 \Bigg) .
\$
Taking the supremum over $\bbeta \in \Theta^*(r)$, or equivalently $\bdelta \in \BB^p(r)$, yields
\#
	\sup_{\bbeta \in \Theta^*(r)}   \| \EE \Delta_b(\bbeta)   \| \leq      \Bigg( \frac{ C_G }{(\nu-1)!  }  b^{\nu-1} + 0.5 l_0 m_3 r \Bigg) r  \lesssim (b^{\nu-1} + r ) r .  \nn
\# 

For the stochastic term $ \sup_{\bbeta \in \Theta^*(r)}   \| \Delta_b(\bbeta) - \EE  \Delta_b(\bbeta)   \|_2$, following the proof of Theorem~\ref{thm:bahadur}, it can be similarly shown that with probability at least $1-e^{-t}$,
\#
\sup_{\bbeta \in \Theta^*(r)}   \| \Delta_b(\bbeta) - \EE  \Delta_b(\bbeta)   \|_2 \lesssim r  \sqrt{\frac{p+t}{n b}} \nn
\#  
as long as $\sqrt{(p+t)/n} \lesssim b$. 

Combining the last two displays completes the proof of \eqref{gradient.uniform.bound}. \qed

\subsection{Proof of Theorem~\ref{thm:one-step.conquer}}

\noindent
{\sc Step 1} (Consistency of the initial estimate). First, note that the consistency of the initial estimator $\overbar  \bbeta$---namely, $\overbar  \bbeta$ lies in a local neighborhood of $\bbeta^*$ with high probability, is a direct consequence of Theorem~\ref{thm:concentration}. Given a non-negative kernel $K(\cdot)$ and for any $t>0$, the initial estimator $\overbar \bbeta$ satisfies
\#
	 \| \overbar \bbeta - \bbeta^* \|_{\bSigma} \leq r_1 \asymp \sqrt{\frac{p+t}{n}}  \label{initial.concentration}
\#
with probability at least $1-2e^{-t}$ as long as $ (\frac{p+t}{n})^{1/2} \lesssim h \lesssim (\frac{p+t}{n})^{1/4}$.  Let $\cE_{{\rm init}}(t)$ be the event that \eqref{initial.concentration} holds.
Provided that the sample Hessian $\nabla^2 \hat Q^G_b(\overbar \bbeta)$ is  invertible, we have
\#
 \Jb (  \hat \bbeta -  \overbar \bbeta  )  & = -  \Jb  \bigl\{  \nabla^2 \hat Q^G_b(\overbar \bbeta)   \bigr\}^{-1}   \nabla \hat Q^G_b(\overbar \bbeta ) \nn \\
 & =   -    \Jb  \bigl\{  \nabla^2 \hat Q^G_b(\overbar \bbeta)   \bigr\}^{-1} \bSigma^{ 1/2}  \cdot  \bSigma^{-1/2}   \bigl\{  \nabla \hat Q^G_b(\overbar \bbeta )  - \nabla \hat Q^G_b(\bbeta^* ) -  \Jb (\overbar \bbeta - \bbeta^* )  \bigr\} \nn \\
 & ~~~~~ -  \Jb  \bigl\{  \nabla^2 \hat Q^G_b(\overbar \bbeta)   \bigr\}^{-1} \bSigma^{ 1/2} \cdot  \bSigma^{-1/2} \bigl\{   \Jb  (\overbar \bbeta - \bbeta^* ) +  \nabla \hat Q^G_b(\bbeta^* )\bigr\}  ,\nn
\#
or equivalently,
\#
\bSigma^{-1/2} \Jb (  \hat \bbeta -   \bbeta^*  )  &=   -   \Hb  \hat \Hb_b^{-1}\cdot  \bSigma^{-1/2}   \bigl\{  \nabla \hat Q^G_b(\overbar \bbeta )  - \nabla \hat Q^G_b(\bbeta^* ) -  \Jb (\overbar \bbeta - \bbeta^* )  \bigr\}  \nn\\
& ~~~~~ + \big( \Ib_p - \Hb \hat \Hb_b^{-1} \big) \Hb \cdot \bSigma^{1/2}  (\overbar \bbeta - \bbeta^* ) - \Hb \hat \Hb_b^{-1} \cdot \bSigma^{-1/2} \nabla \hat Q^G_b(\bbeta^*) , \nn
\#
where $\Hb = \EE \{ f_{\varepsilon |\bx}(0) \bw \bw^\T \} = \bSigma^{-1/2} \Jb \bSigma^{-1/2}$ and  $\hat \Hb_b := \bSigma^{-1/2} \nabla^2 \hat Q^G_b(\overbar \bbeta) \bSigma^{-1/2} $.  It follows that
\#
 &   \| \bSigma^{-1/2}  \Jb (  \hat \bbeta -   \bbeta^*  )   +  \bSigma^{-1/2} \nabla  \hat Q^G_b(\bbeta^*)   \|_2  \nn \\
& \leq  \|  ( \Ib_p - \Hb  \hat \Hb_b^{-1} ) \Hb \|_2  \cdot       \| \overbar \bbeta - \bbeta^*\|_{\bSigma} +  \| \Ib_p - \Hb \Hb_b^{-1} \|_2 \cdot \| \nabla  \hat Q^G_b(\bbeta^*) \|_{\bOmega} \nn \\
&~~~~~ + \|  \Hb \hat \Hb_b^{-1} \|_2 \cdot   \|   \nabla \hat Q^G_b(\overbar \bbeta )  - \nabla \hat Q^G_b(\bbeta^* ) -  \Jb (\overbar \bbeta - \bbeta^* )  \|_{\bOmega} .  \label{os.conquer.ubd1}
\#
In view of \eqref{os.conquer.ubd1}, it remains to bound the following three quantities:
$$
   \|  \Ib_p -   \Hb  \hat \Hb_b^{-1}    \|_2 , ~~\|  \nabla \hat Q^G_b(\bbeta^* )   \|_{\bOmega}
~\mbox{ and }~
 \|    \nabla \hat Q^G_b(\overbar  \bbeta )  - \nabla \hat Q^G_b(\bbeta^* )   -\Jb (\overbar  \bbeta - \bbeta^* )   \|_{\bOmega}.
$$

\medskip
\noindent
{\sc Step 2} (Consistency of the sample Hessian $\nabla^2 \hat Q^G_b(\overbar  \bbeta)$). Recall that $\hat \Hb_b  = \bSigma^{-1/2} \nabla^2 \hat Q^G_b(\overbar  \bbeta) \bSigma^{-1/2} $ and $\Hb   = \bSigma^{-1/2} \Jb \bSigma^{-1/2}$. By the triangle inequality,
\#
	&  \| \hat \Hb_b   - \Hb    \|_2    \leq  \|  \bSigma^{-1/2}  \{ \nabla^2 \hat Q^G_b(\overbar  \bbeta)  - \nabla^2   Q^G_b( \bbeta^*) \} \bSigma^{-1/2}    \|_2+   \|  \bSigma^{-1/2}  \nabla^2  Q^G_b(  \bbeta^*)  \bSigma^{-1/2}  - \Hb \|_2.   \nn
\#
Let the bandwidth $b$ satisfy $ \max\{  \frac{p\log n + t}{n} , n^{-1/2}\} \lesssim b \lesssim 1$. Conditioned on $\cE_{{\rm init}}(t)$, applying Propositions~\ref{prop:higher-order.bias} and \ref{prop:Hessian.uniform} with $r=r_1$ yields that, with probability $1- e^{-t}$,
\#
	&   \|  \hat \Hb_b - \Hb   \|_2    \leq r_2  \asymp \sqrt{\frac{p\log n + t}{ nb }} + b^{\nu-1} .  \nn
\#
Under Condition~\ref{cond.reg},  $0< \underbar{$f$} \leq \lambda_{\min}(\Hb) \leq \lambda_{\max}(\Hb) \leq  \bar f$, so that $\| \Hb^{-1} \|_2 \leq \underbar{$f$}^{-1}$. For sufficiently large $n$ and small $b$, this implies $\|  \hat \Hb_b  \Hb^{-1} - \Ib_p \|_2 \leq \underbar{$f$}^{-1} r_2 < 1$, and hence
\#
	  \|  \Hb \hat \Hb_b^{-1}  - \Ib_p \|_2 \leq \frac{r_2}{ \underbar{$f$} - r_2 }  , ~~~  \| \Hb\hat \Hb_b^{-1}  \|_2 \leq   \frac{\underbar{$f$}}{ \underbar{$f$} - r_2 }.  \label{Hessian.rate}
\#

\medskip
\noindent
{\sc Step 3} (Controlling the scores). For $\| \bSigma^{-1/2} \nabla \hat Q^G_b(\bbeta^* )\|_2$, it follows from Lemma~\ref{lem:score} and Proposition~\ref{prop:higher-order.bias} that with probability at least $1-e^{-t}$,
\#
	   \|   \bSigma^{-1/2} \nabla \hat Q^G_b(\bbeta^* )    \|_2  \lesssim \sqrt{\frac{p+ t}{n}} +  b^\nu . \label{score.concentration1}
\# 
Turning to  $\|  \bSigma^{-1/2} \{   \nabla \hat Q^G_b(\overbar  \bbeta )  - \nabla \hat Q^G_b(\bbeta^* ) \}  -  \bSigma^{-1/2} \Jb (\overbar  \bbeta - \bbeta^* )  \|_2$, applying the concentration bound \eqref{initial.concentration} and Proposition~\ref{prop:gradient.uniform} we obtain that, with probability at least $1- e^{-t}$ conditioned on $\cE_{{\rm init}}(t)$,
\#
 \|  \bSigma^{-1/2}\{    \nabla \hat Q^G_b(\overbar  \bbeta )  - \nabla \hat Q^G_b(\bbeta^* )  \}  -    \bSigma^{-1/2}\Jb (\overbar  \bbeta - \bbeta^* )   \|_2 \lesssim  \frac{p+ t}{n b^{1/2}}     + b^{\nu -1} \sqrt{\frac{p+t}{n}} \label{score.concentration2}
\#
as long as $ (\frac{p+t}{n})^{1/2} \lesssim b \lesssim 1$.

Finally, combining the bounds \eqref{initial.concentration}--\eqref{score.concentration2}, we conclude that with probability at least $1-5e^{-t}$, 
\#
 \|  \bSigma^{-1/2}   \Jb (\hat  \bbeta - \bbeta^* ) +\bSigma^{-1/2} \nabla \hat Q^G_b(\bbeta^* )   \|_2
 \lesssim   \Biggl(  \sqrt{\frac{p\log n+ t}{n b}}+ b^{\nu-1} \Biggr) \Bigg(  \sqrt{\frac{p+t}{n}} + b^\nu \Bigg)  , \nn
\#
provided that $ \max\bigl\{ \frac{p\log n + t}{n} ,  (\frac{p + t}{n})^{1/2} \bigr\} \lesssim b \lesssim 1$. Under the sample size scaling  $n\gtrsim p \,(\log n)^2+t$ and bandwidth constraint $b \lesssim (\frac{p+t}{n})^{1/(2\nu)}$, this leads to the claimed bound \eqref{os.conquer.br}. \qed

 \section{Additional simulation results} \label{sec:app.simu}
\begin{figure}[!htp]
  \centering
  \subfigure[Model \eqref{model.homo} under $\tau = 0.1$.]{\includegraphics[width=0.32\textwidth]{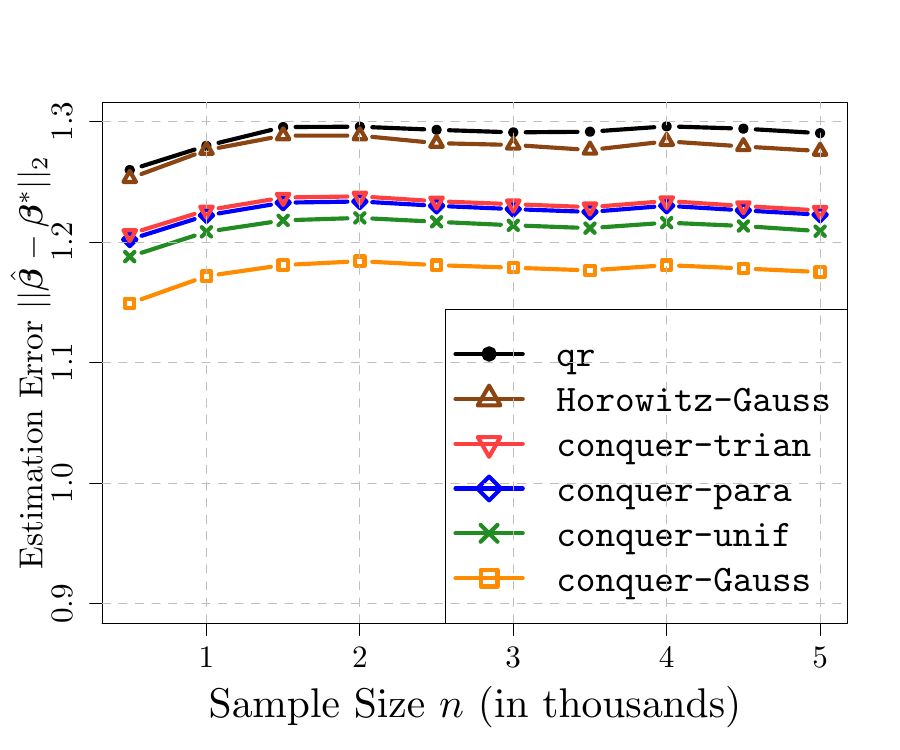}} 
  \subfigure[Model \eqref{model.linear} under $\tau = 0.1$.]{\includegraphics[width=0.32\textwidth]{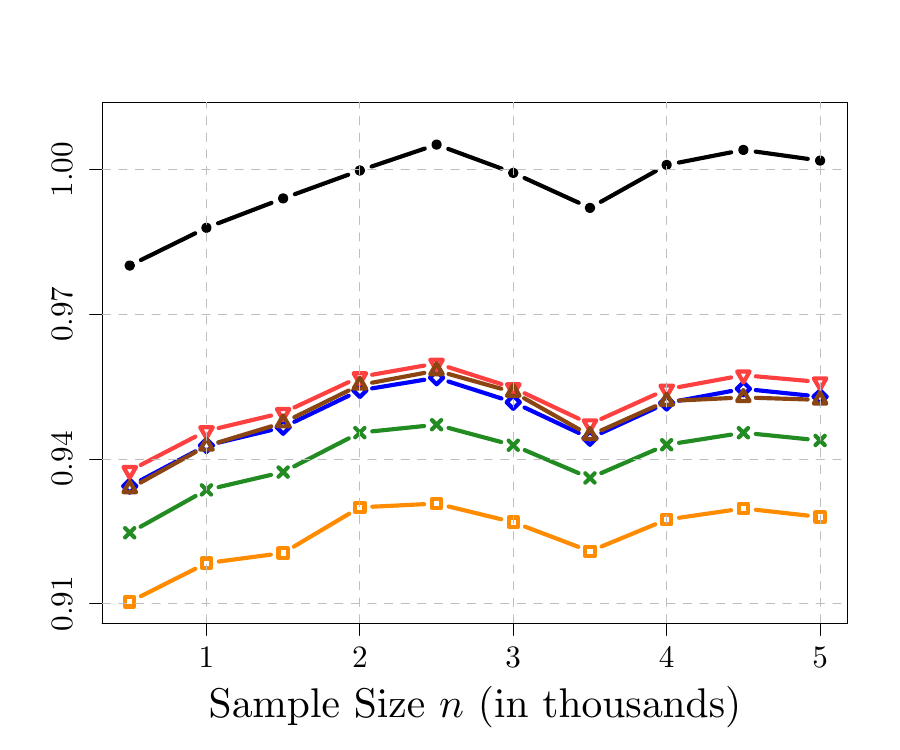}} 
  \subfigure[Model \eqref{model.quad} under $\tau = 0.1$.]{\includegraphics[width=0.32\textwidth]{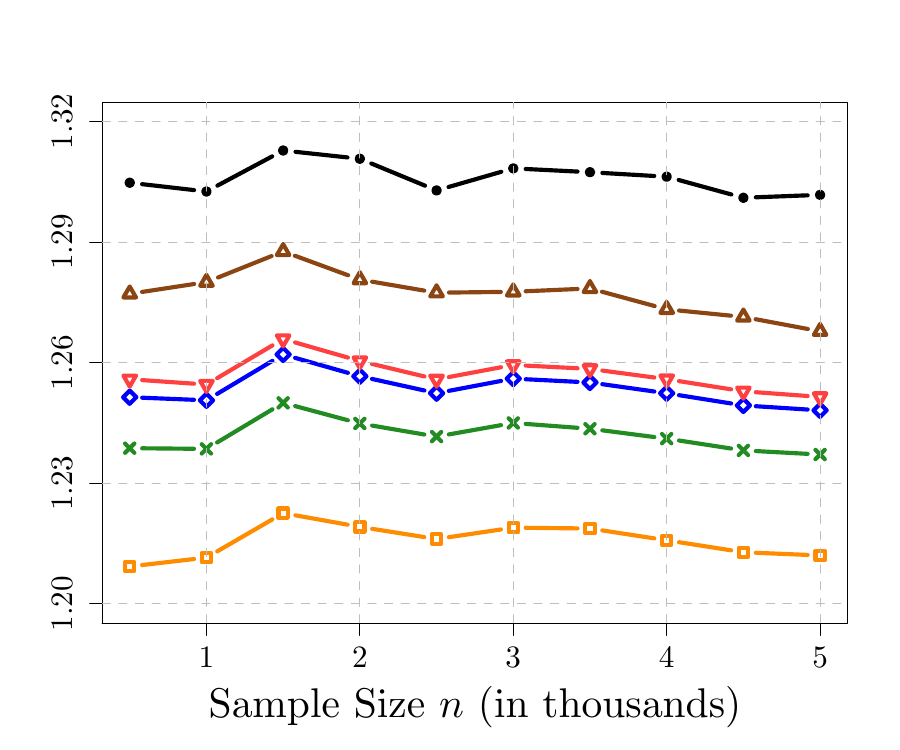}}
  \subfigure[Model \eqref{model.homo} under $\tau = 0.3$.]{\includegraphics[width=0.32\textwidth]{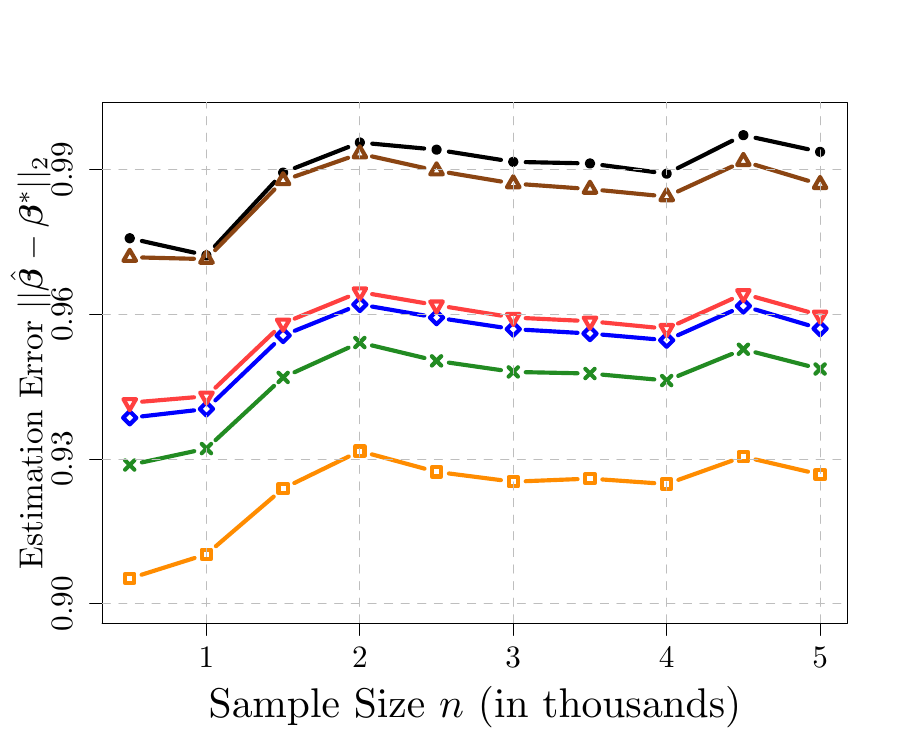}}
  \subfigure[Model \eqref{model.linear} under $\tau = 0.3$.]{\includegraphics[width=0.32\textwidth]{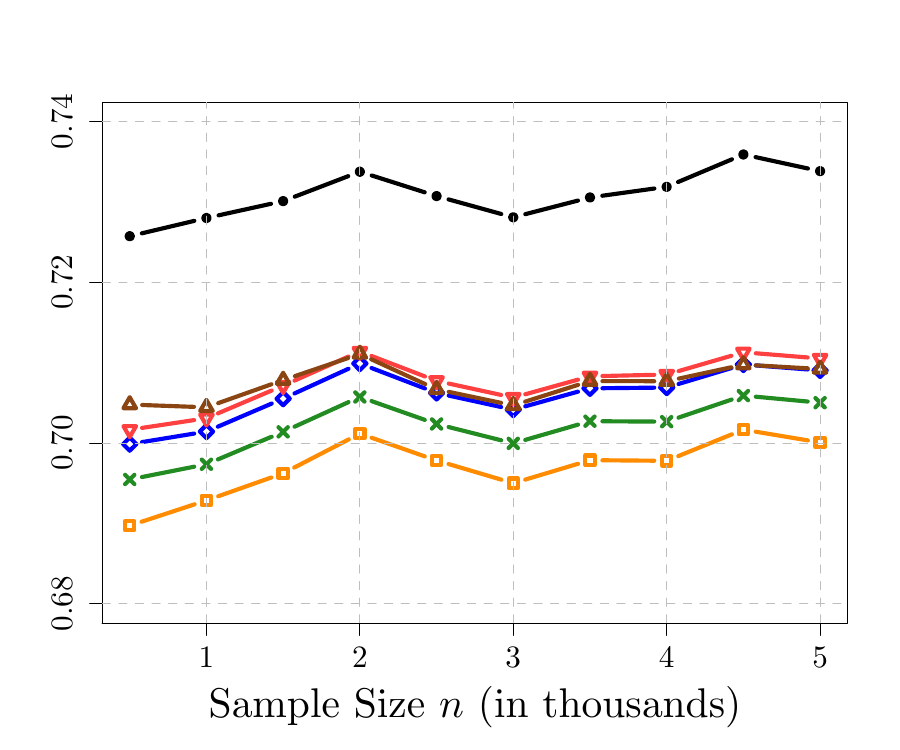}}
  \subfigure[Model \eqref{model.quad} under $\tau = 0.3$.]{\includegraphics[width=0.32\textwidth]{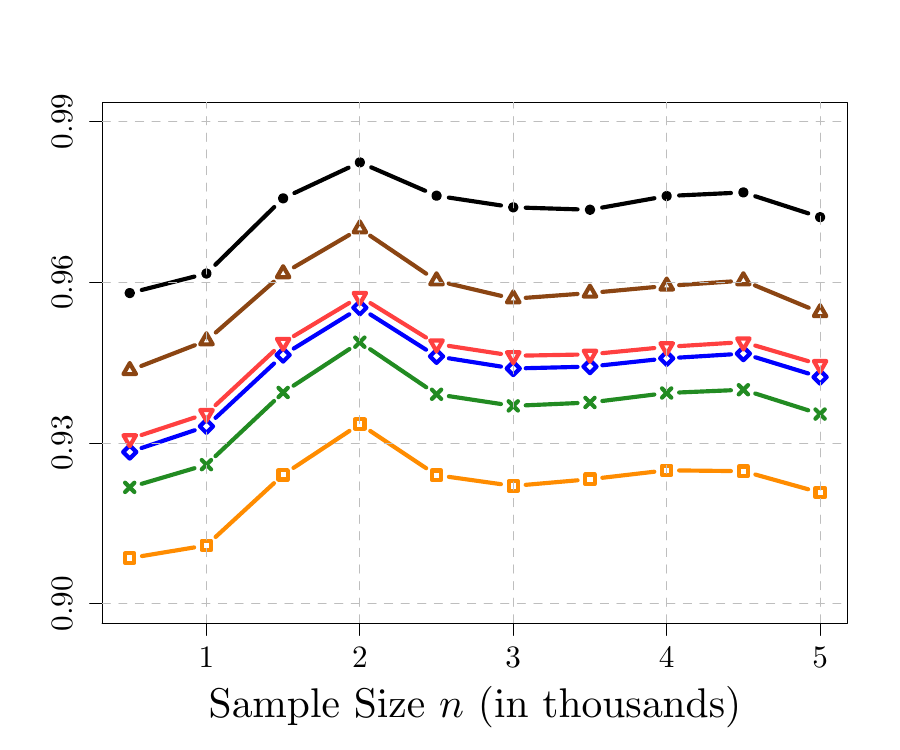}}
  \subfigure[Model \eqref{model.homo} under $\tau = 0.5$.]{\includegraphics[width=0.32\textwidth]{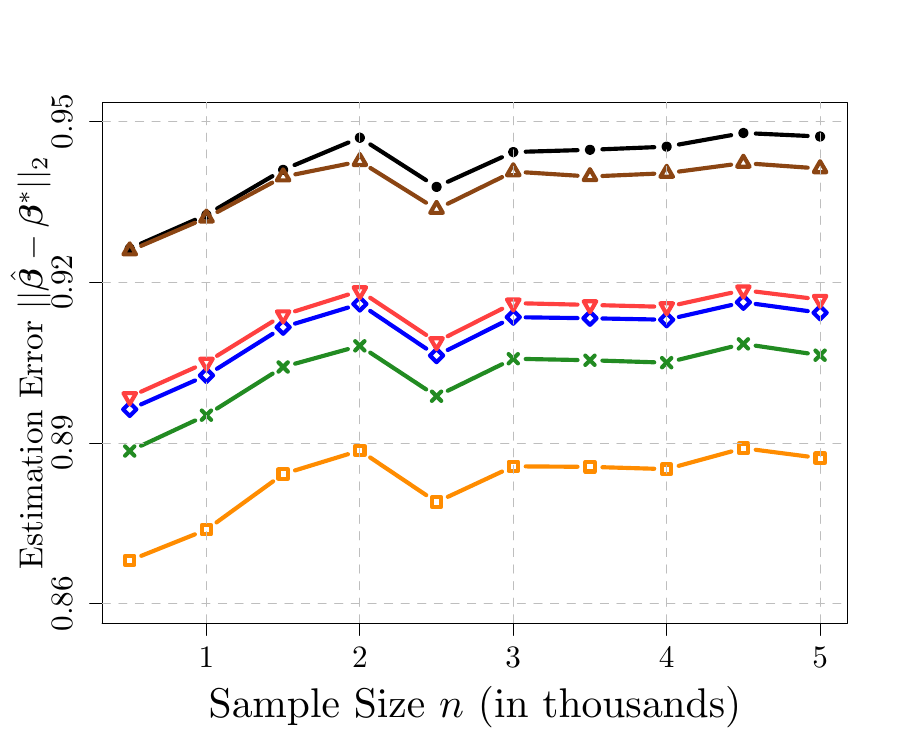}}
  \subfigure[Model \eqref{model.linear} under $\tau = 0.5$.]{\includegraphics[width=0.32\textwidth]{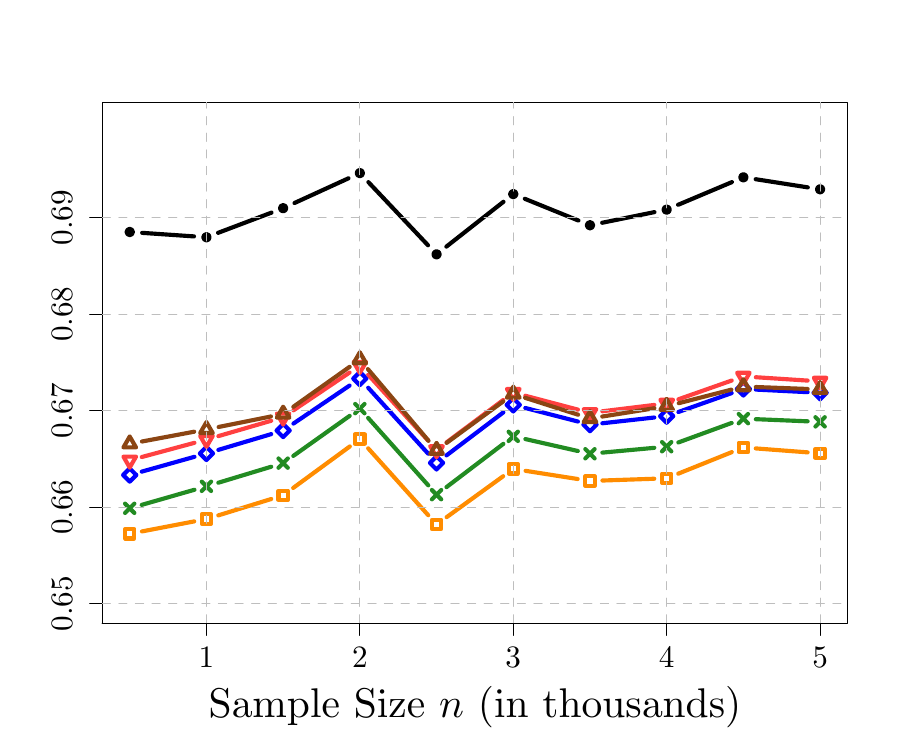}}
  \subfigure[Model \eqref{model.quad} under $\tau = 0.5$.]{\includegraphics[width=0.32\textwidth]{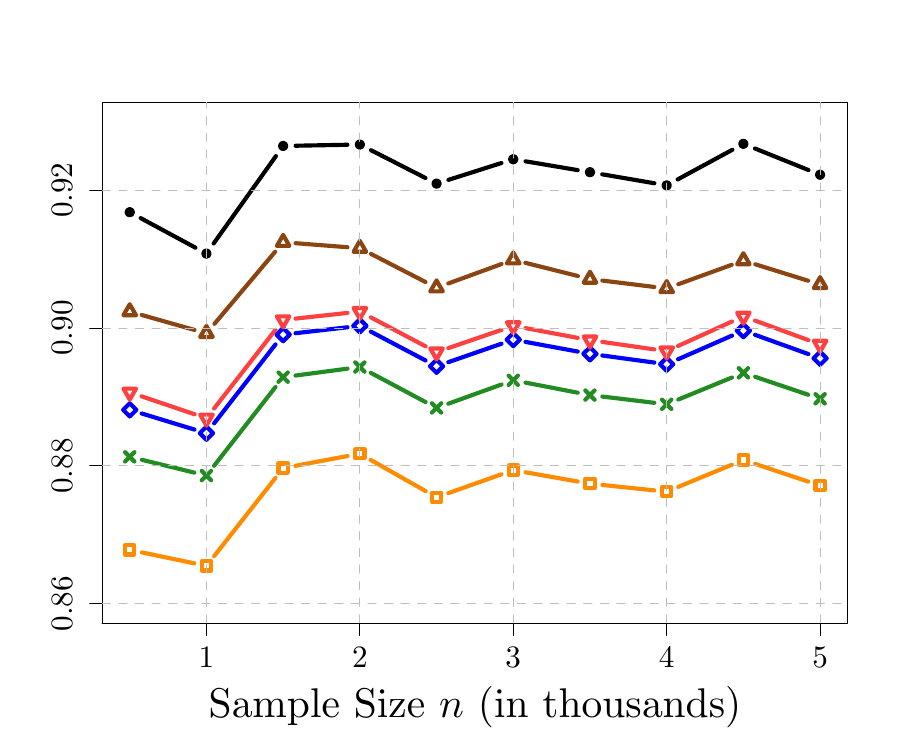}}
  \subfigure[Model \eqref{model.homo} under $\tau = 0.7$.]{\includegraphics[width=0.32\textwidth]{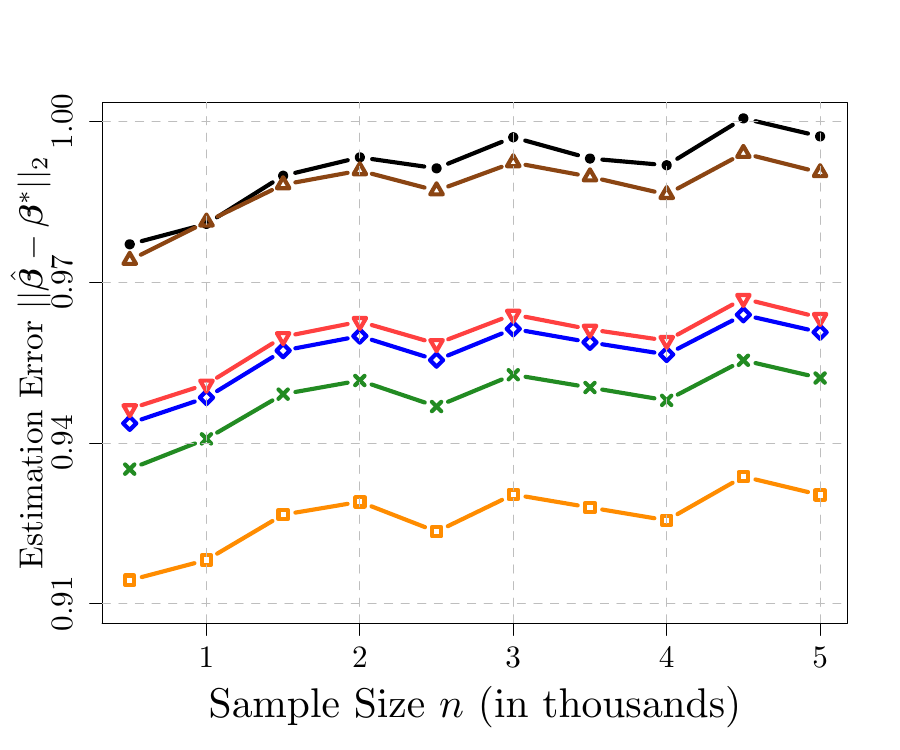}}
  \subfigure[Model \eqref{model.linear} under $\tau = 0.7$.]{\includegraphics[width=0.32\textwidth]{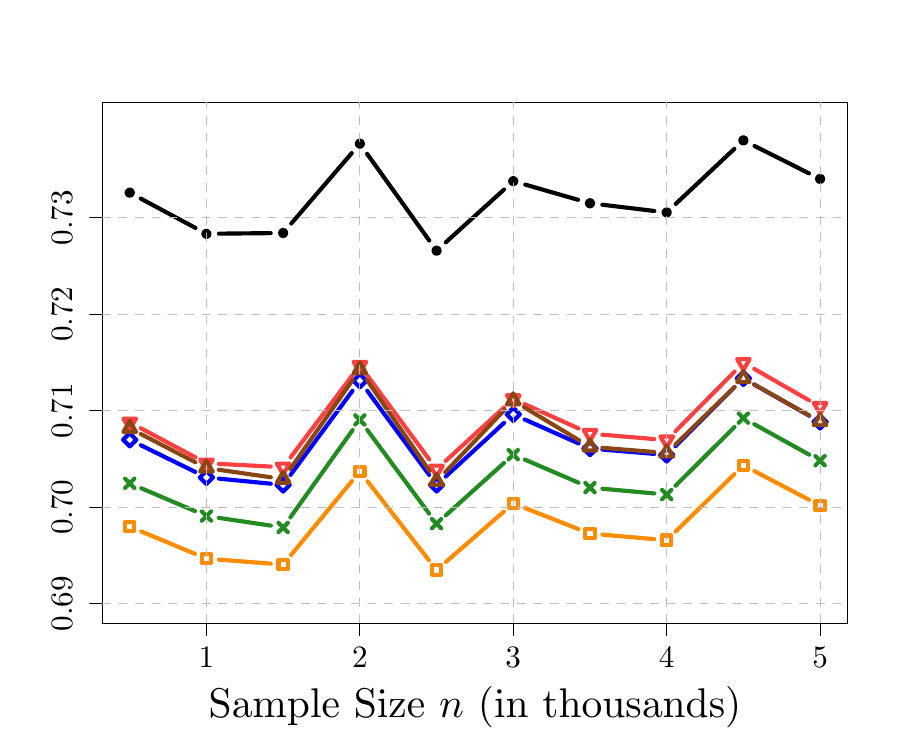}}
  \subfigure[Model \eqref{model.quad} under $\tau = 0.7$.]{\includegraphics[width=0.32\textwidth]{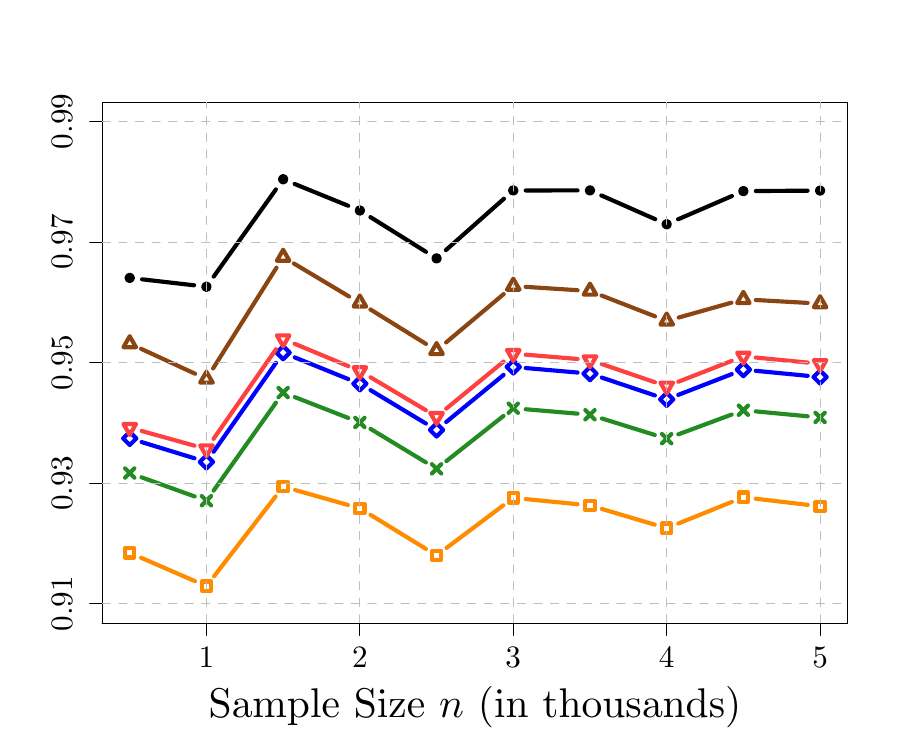}}
\caption{Results under models \eqref{model.homo}--\eqref{model.quad} in Section~\ref{sec:numerical} with $\tau \in \{0.1,0.3, 0.5, 0.7\}$ and $\mathcal{N}(0, 4)$ noise, averaged over 500 simulations. This figure extends the first row of Figure~\ref{est.9} to other quantile levels.}
  \label{est.1357.normal}
\end{figure}

\begin{figure}[!htp]
  \centering
  \subfigure[Model \eqref{model.homo} under $\tau = 0.1$.]{\includegraphics[width=0.32\textwidth]{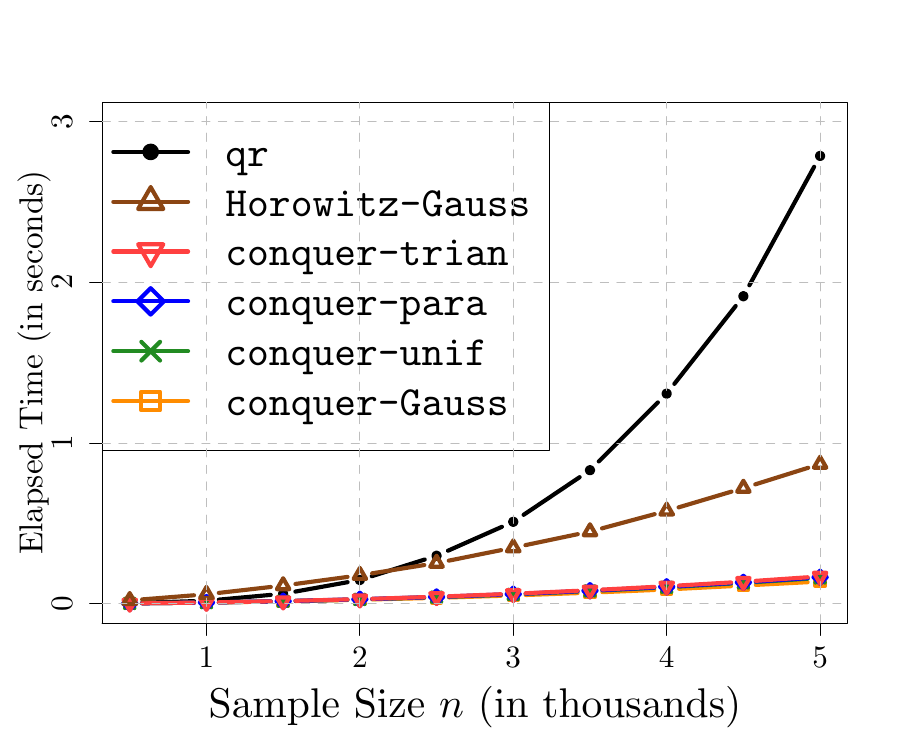}} 
  \subfigure[Model \eqref{model.linear} under $\tau = 0.1$.]{\includegraphics[width=0.32\textwidth]{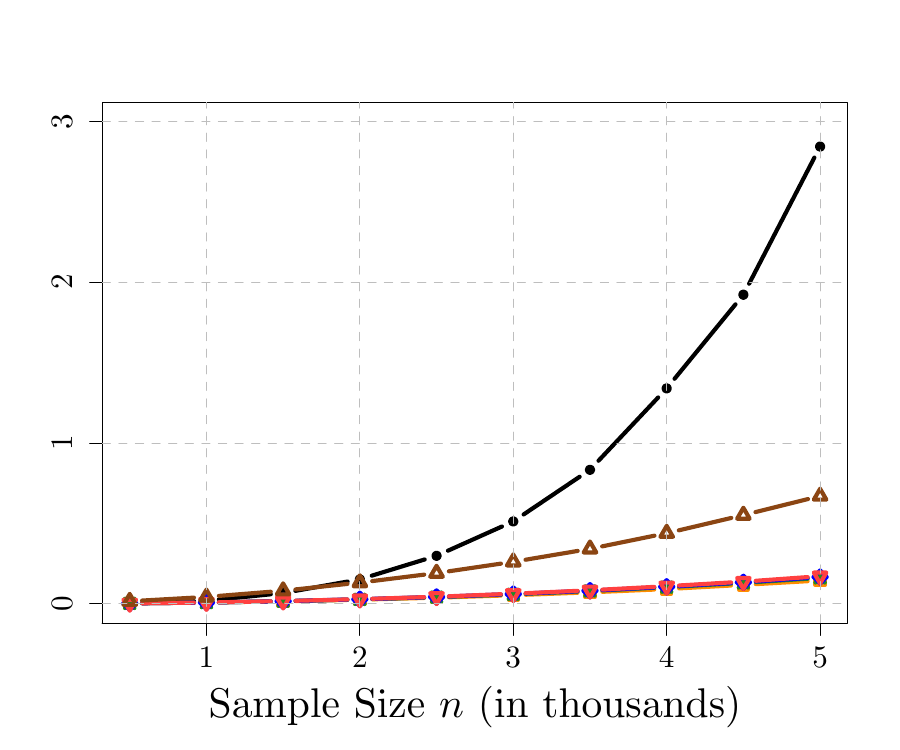}} 
  \subfigure[Model \eqref{model.quad} under $\tau = 0.1$.]{\includegraphics[width=0.32\textwidth]{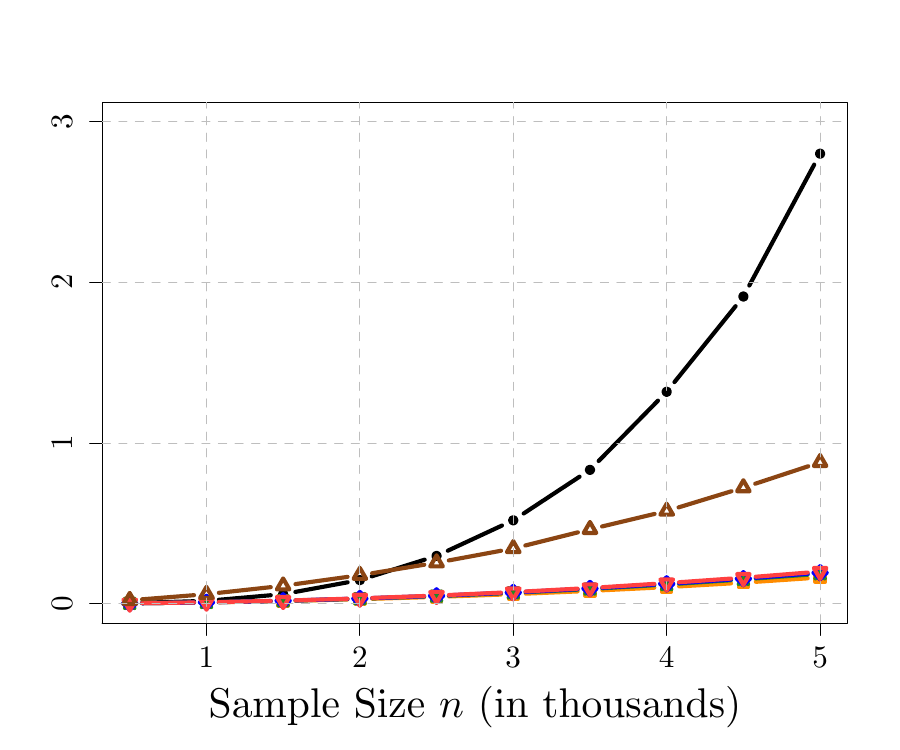}}
  \subfigure[Model \eqref{model.homo} under $\tau = 0.3$.]{\includegraphics[width=0.32\textwidth]{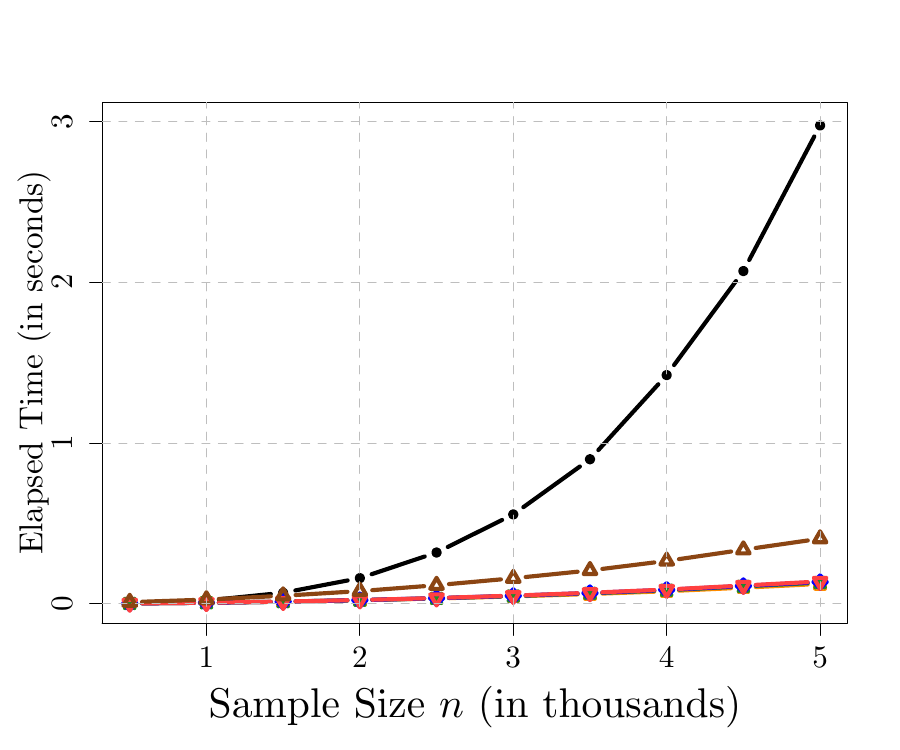}}
  \subfigure[Model \eqref{model.linear} under $\tau = 0.3$.]{\includegraphics[width=0.32\textwidth]{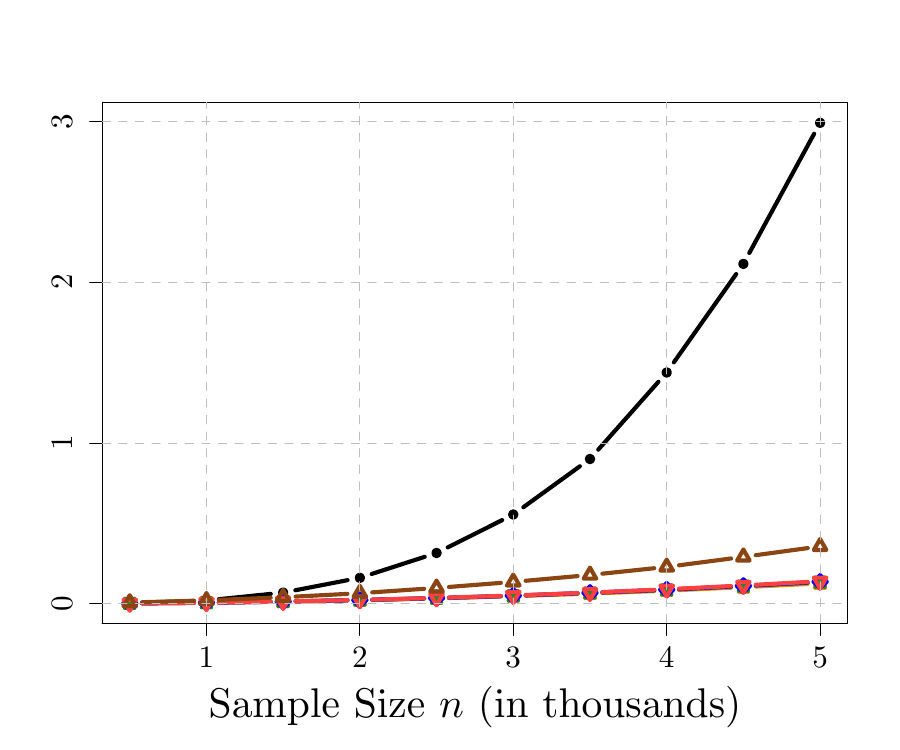}}
  \subfigure[Model \eqref{model.quad} under $\tau = 0.3$.]{\includegraphics[width=0.32\textwidth]{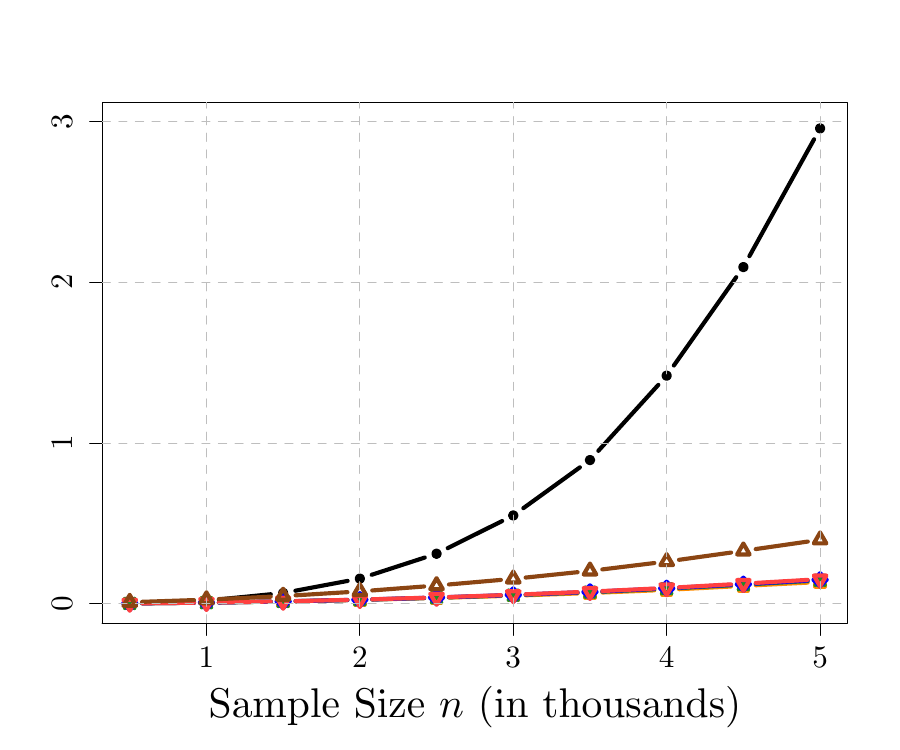}}
  \subfigure[Model \eqref{model.homo} under $\tau = 0.5$.]{\includegraphics[width=0.32\textwidth]{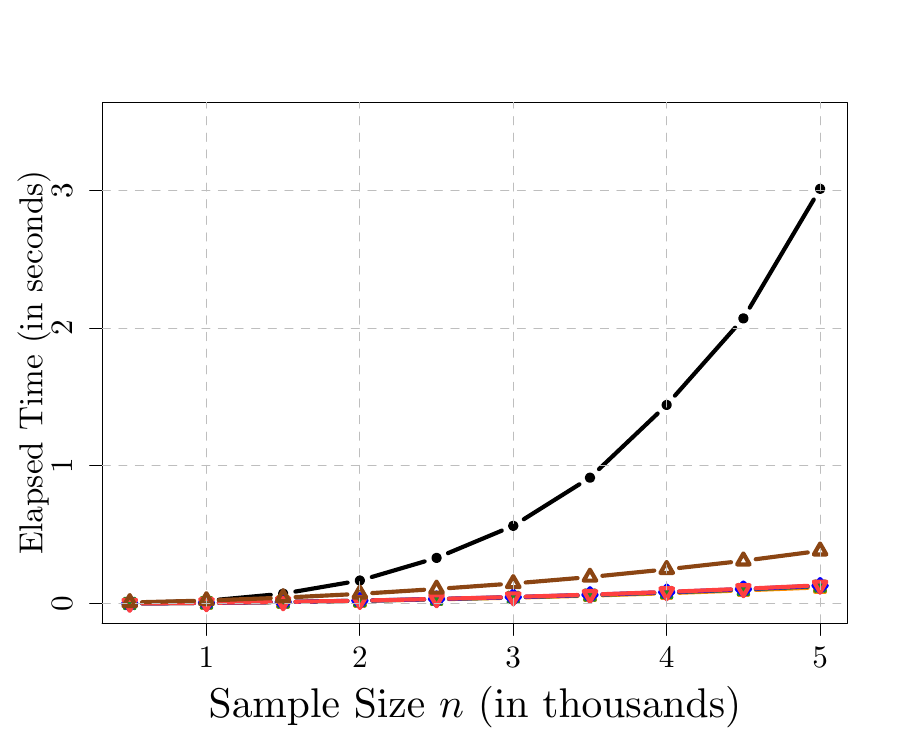}}
  \subfigure[Model \eqref{model.linear} under $\tau = 0.5$.]{\includegraphics[width=0.32\textwidth]{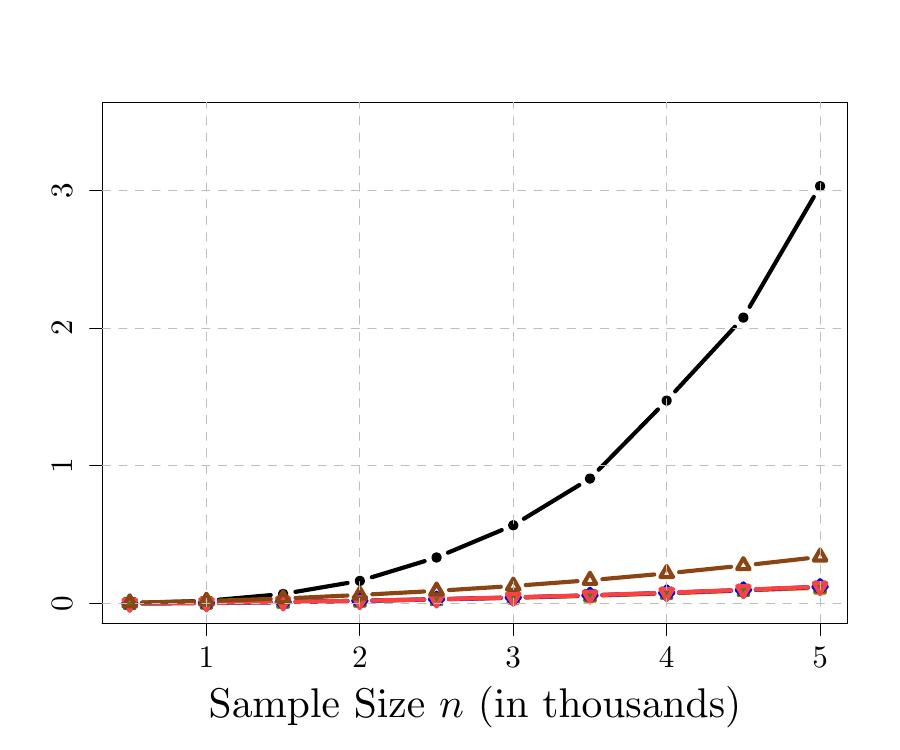}}
  \subfigure[Model \eqref{model.quad} under $\tau = 0.5$.]{\includegraphics[width=0.32\textwidth]{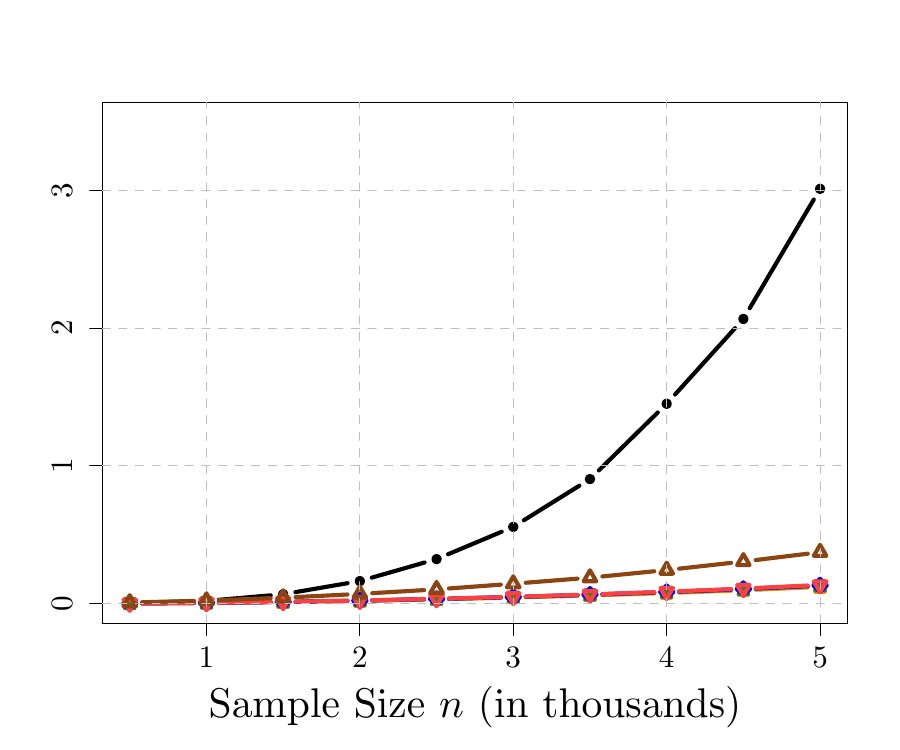}}
  \subfigure[Model \eqref{model.homo} under $\tau = 0.7$.]{\includegraphics[width=0.32\textwidth]{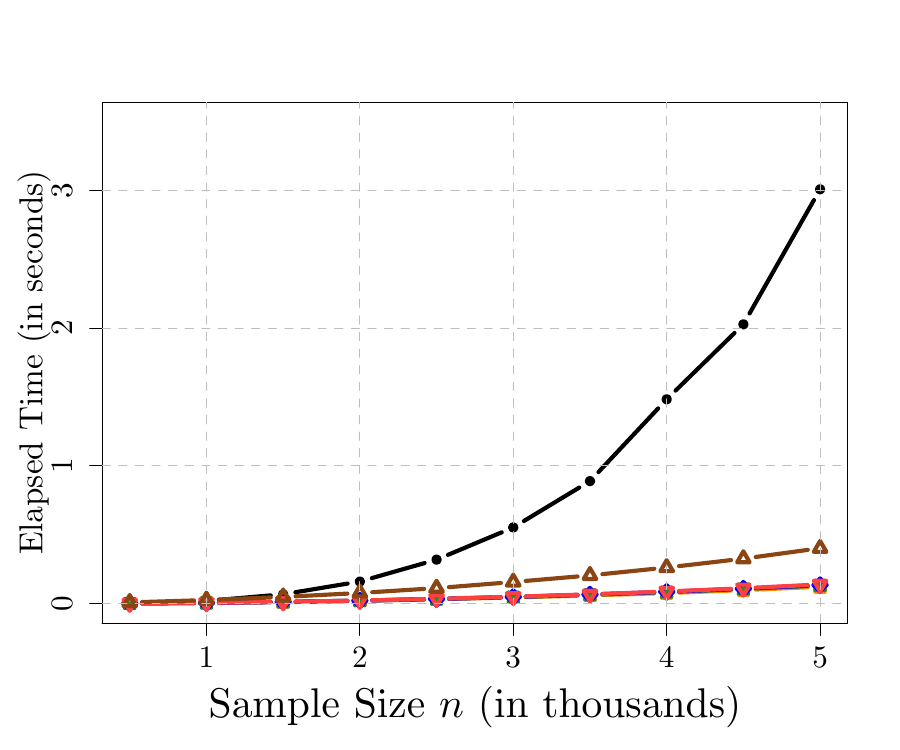}}
  \subfigure[Model \eqref{model.linear} under $\tau = 0.7$.]{\includegraphics[width=0.32\textwidth]{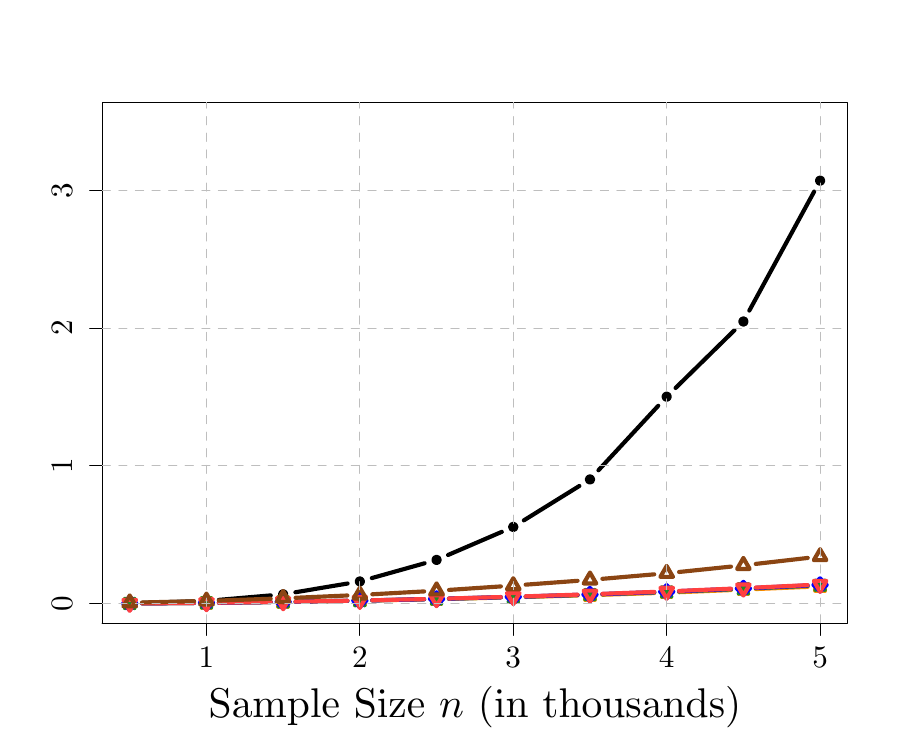}}
  \subfigure[Model \eqref{model.quad} under $\tau = 0.7$.]{\includegraphics[width=0.32\textwidth]{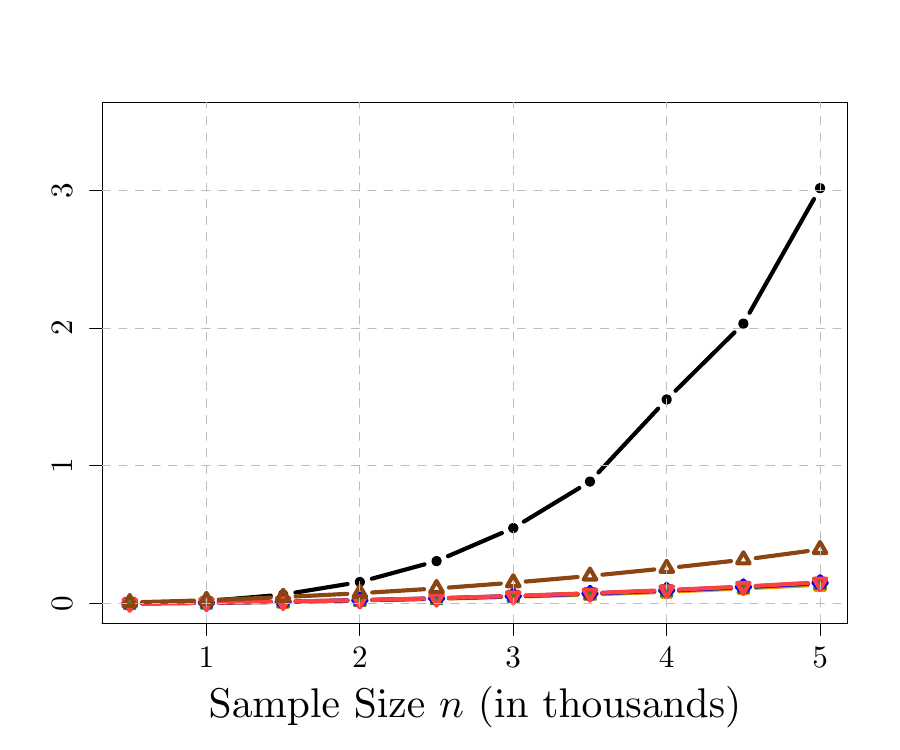}}
\caption{Elapsed time of standard QR, Horowitz's smoothing and conquer under $\mathcal{N}(0, 4)$ noise setting when $\tau \in \{0.1, 0.3, 0.5, 0.7\}$. This figure extends the first row of Figure~\ref{time.9} to other quantile levels.}
  \label{time.1357.normal}
\end{figure}

\begin{figure}[!htp]
  \centering
  \subfigure[Model \eqref{model.homo} under $\tau = 0.1$.]{\includegraphics[width=0.32\textwidth]{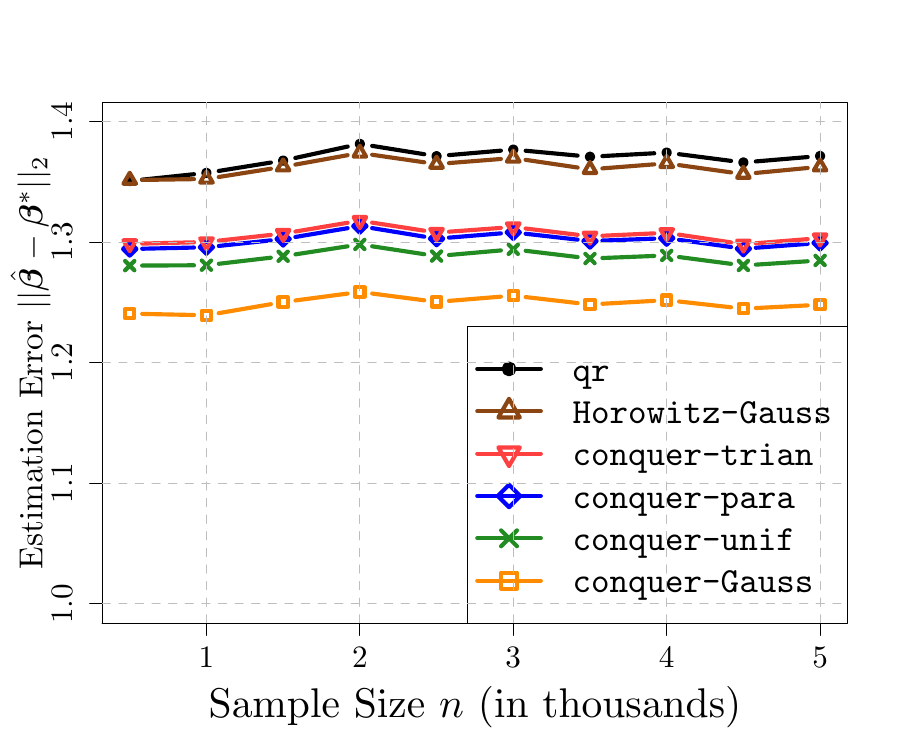}} 
  \subfigure[Model \eqref{model.linear} under $\tau = 0.1$.]{\includegraphics[width=0.32\textwidth]{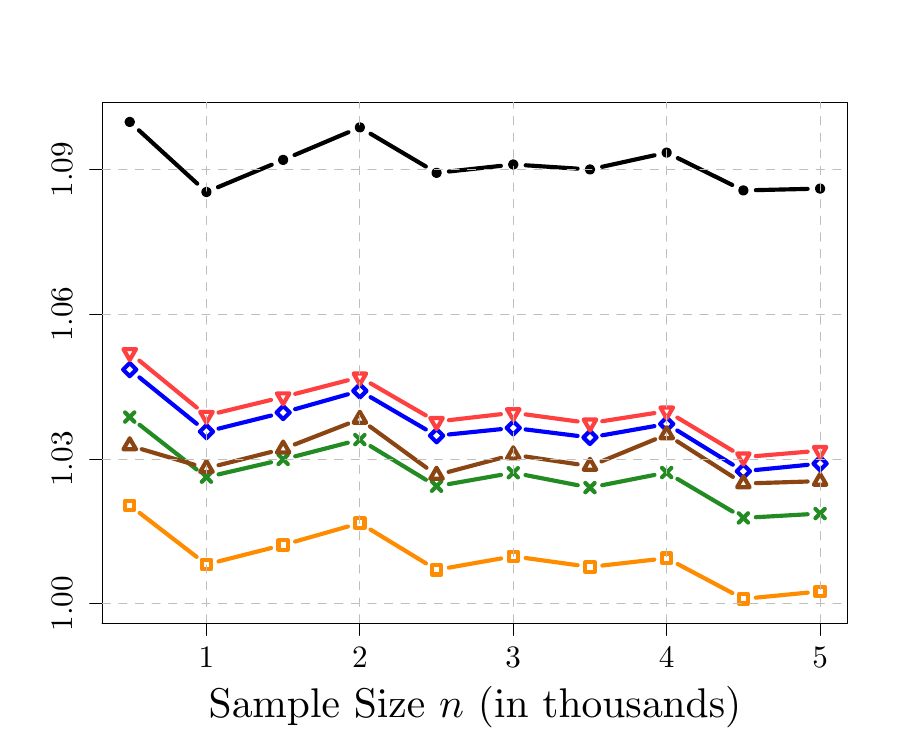}} 
  \subfigure[Model \eqref{model.quad} under $\tau = 0.1$.]{\includegraphics[width=0.32\textwidth]{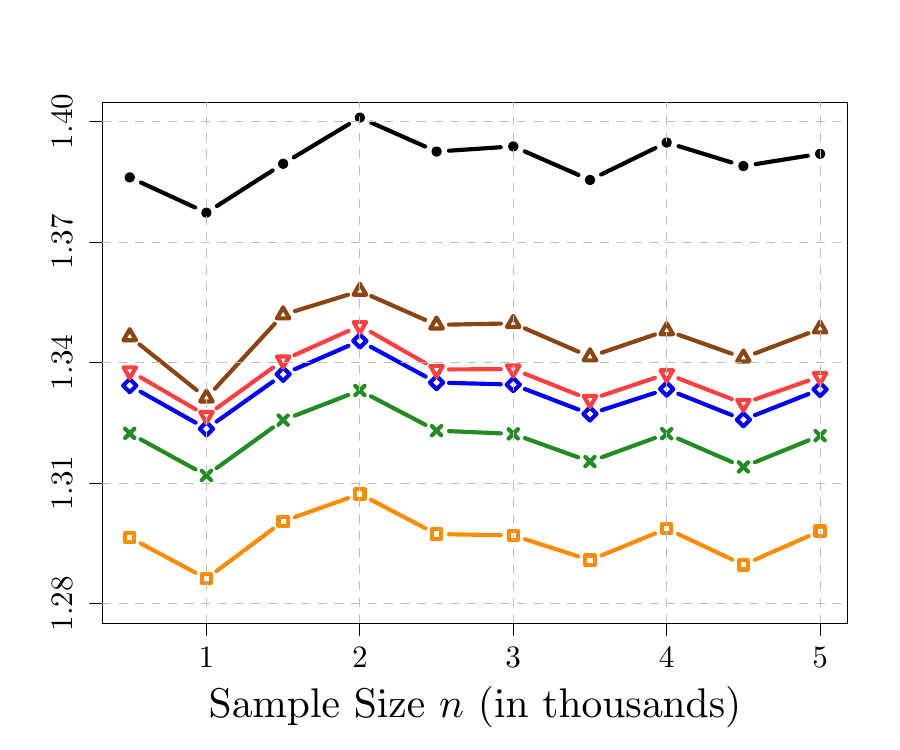}}
  \subfigure[Model \eqref{model.homo} under $\tau = 0.3$.]{\includegraphics[width=0.32\textwidth]{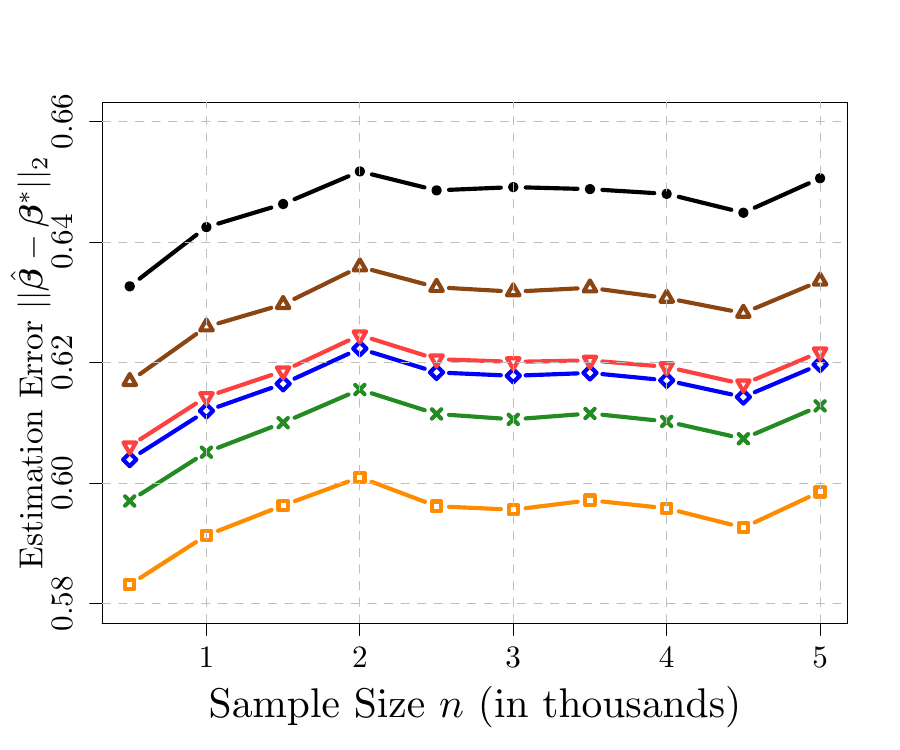}}
  \subfigure[Model \eqref{model.linear} under $\tau = 0.3$.]{\includegraphics[width=0.32\textwidth]{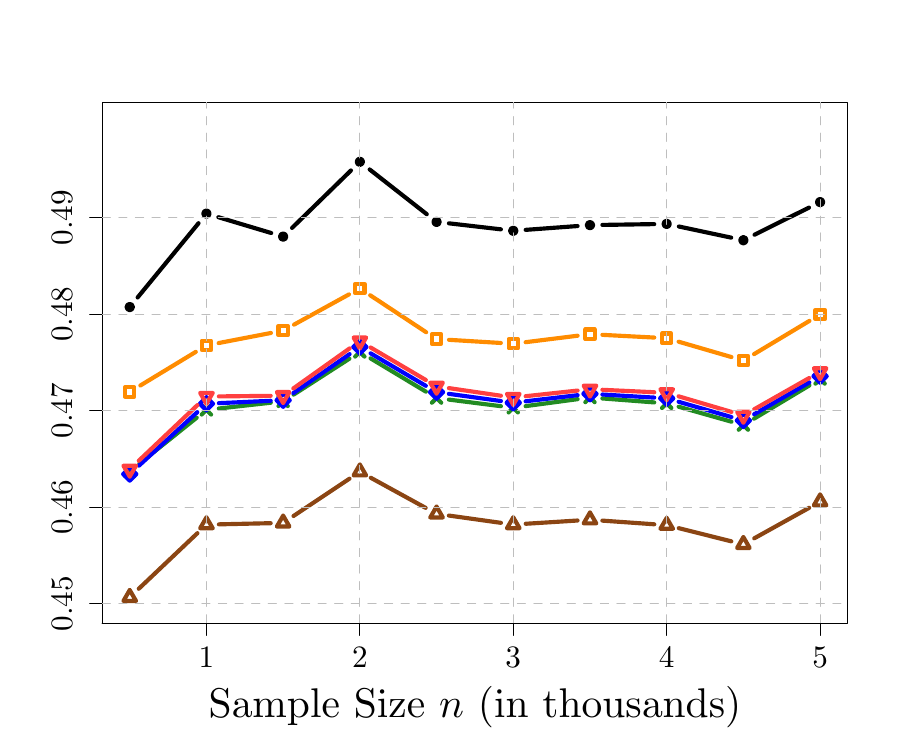}}
  \subfigure[Model \eqref{model.quad} under $\tau = 0.3$.]{\includegraphics[width=0.32\textwidth]{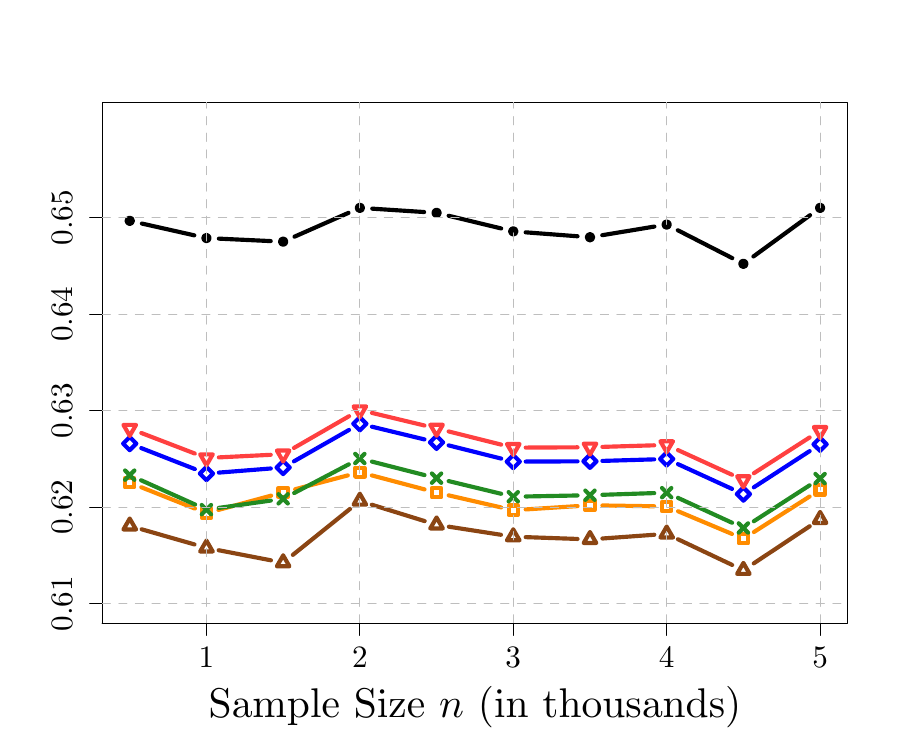}}
  \subfigure[Model \eqref{model.homo} under $\tau = 0.5$.]{\includegraphics[width=0.32\textwidth]{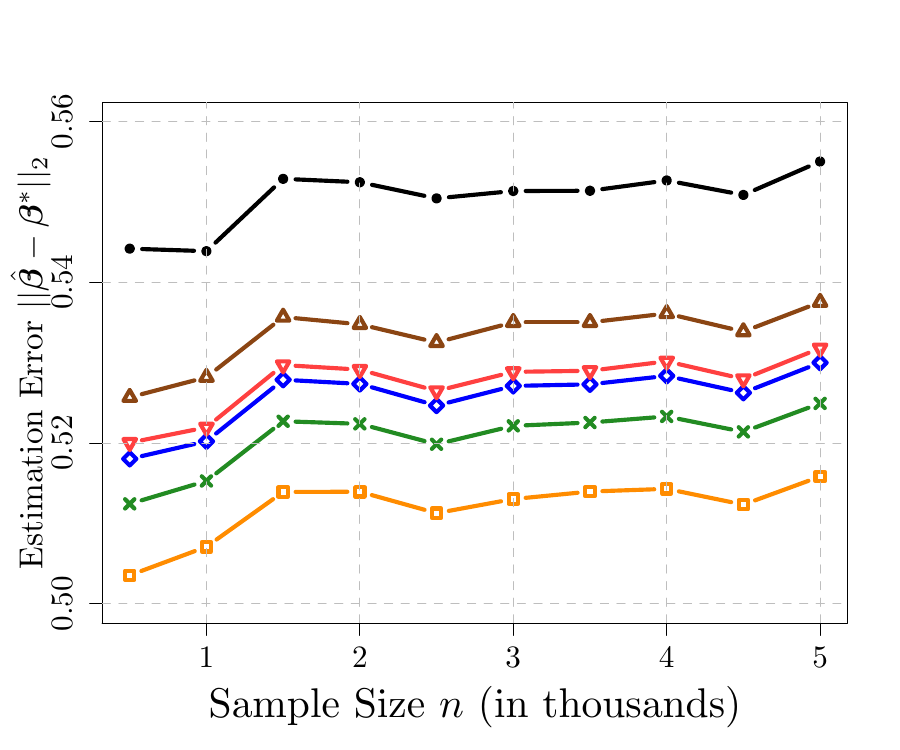}}
  \subfigure[Model \eqref{model.linear} under $\tau = 0.5$.]{\includegraphics[width=0.32\textwidth]{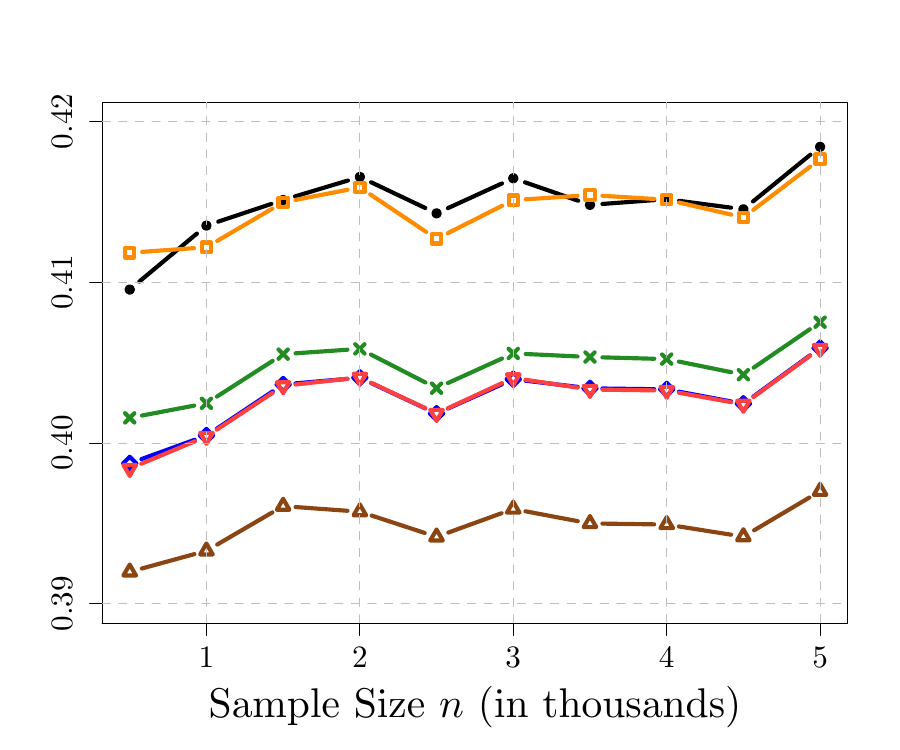}}
  \subfigure[Model \eqref{model.quad} under $\tau = 0.5$.]{\includegraphics[width=0.32\textwidth]{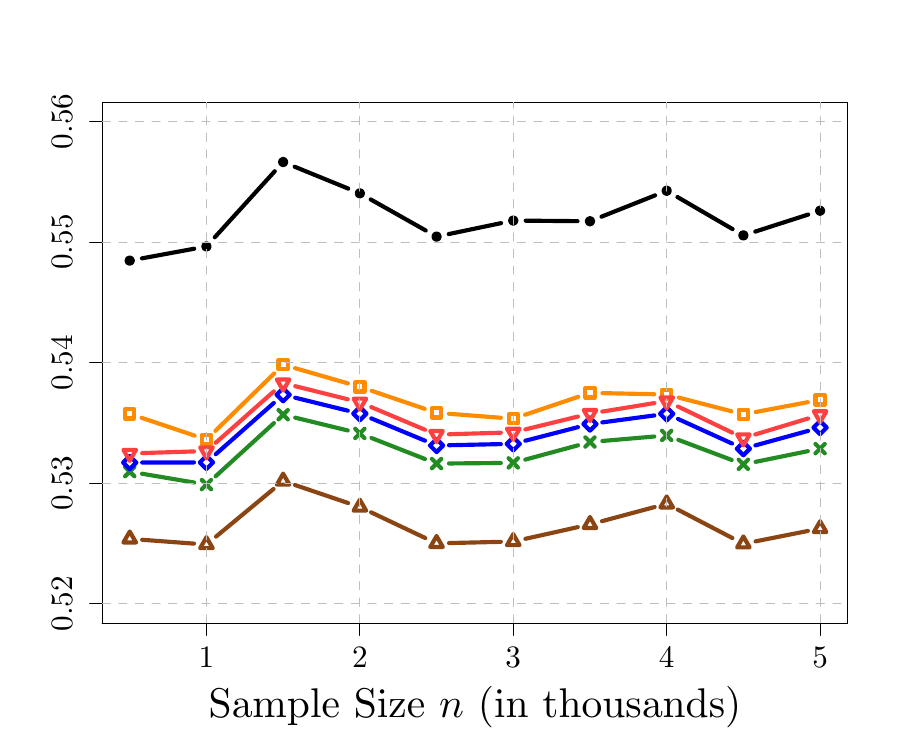}}
  \subfigure[Model \eqref{model.homo} under $\tau = 0.7$.]{\includegraphics[width=0.32\textwidth]{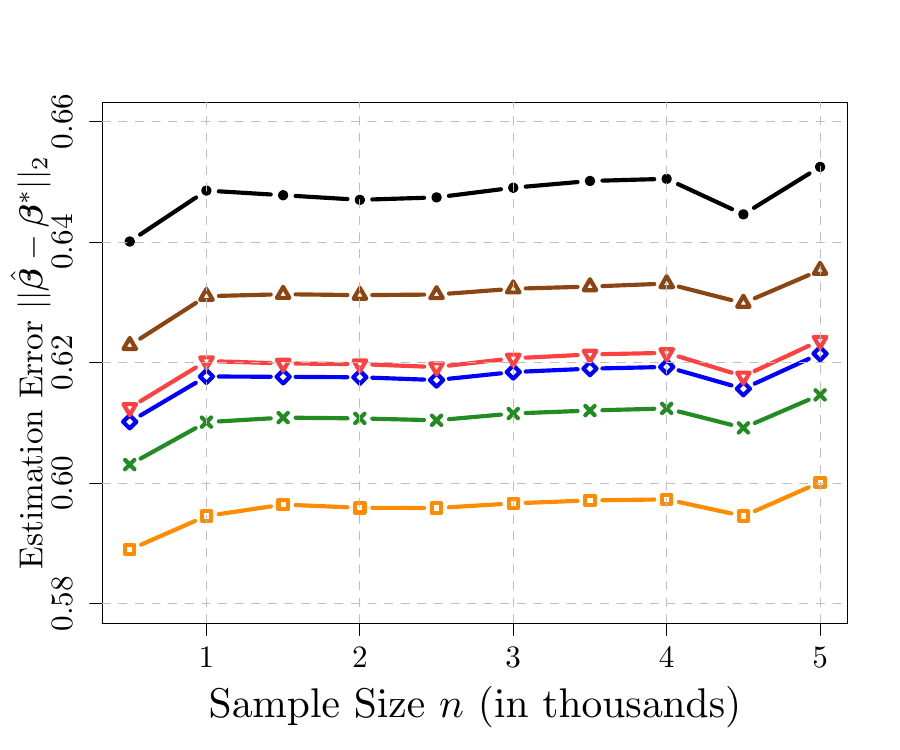}}
  \subfigure[Model \eqref{model.linear} under $\tau = 0.7$.]{\includegraphics[width=0.32\textwidth]{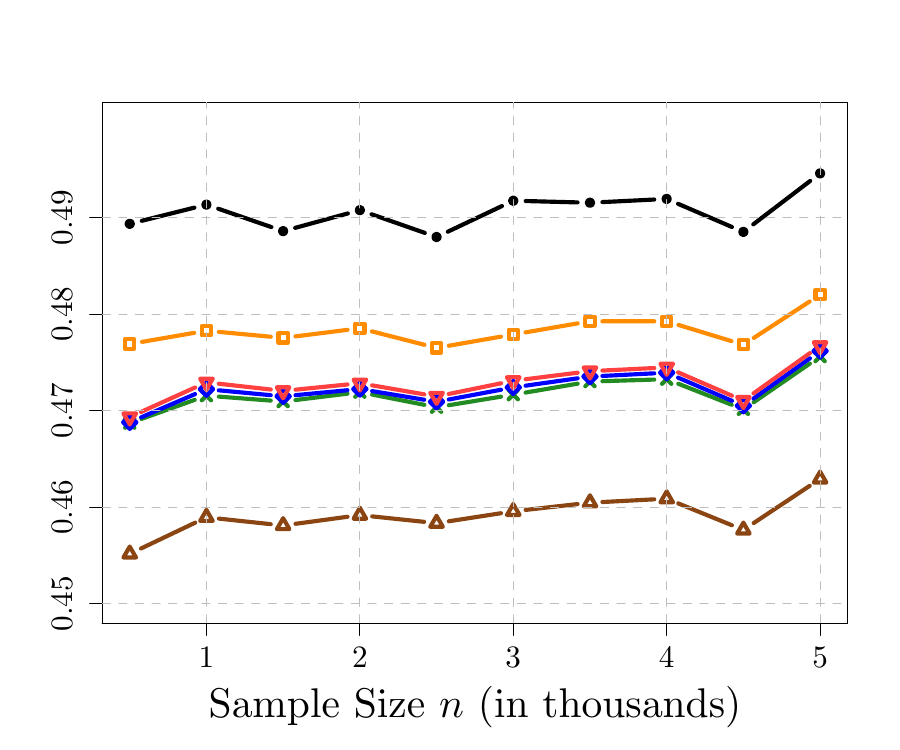}}
  \subfigure[Model \eqref{model.quad} under $\tau = 0.7$.]{\includegraphics[width=0.32\textwidth]{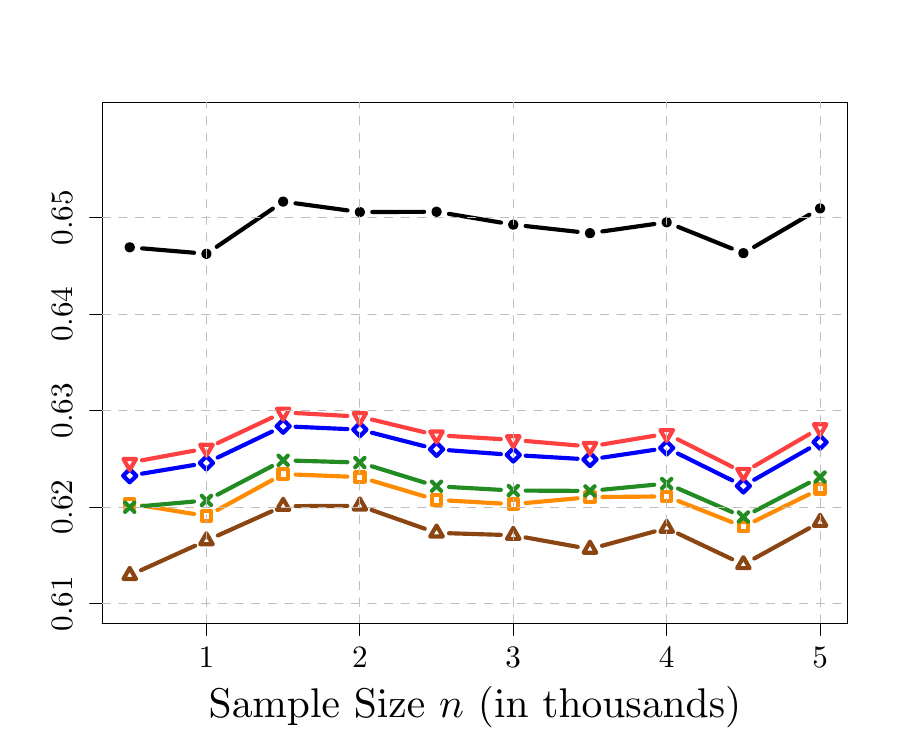}}
\caption{Results under models \eqref{model.homo}--\eqref{model.quad} in Section~\ref{sec:numerical} with $\tau \in \{0.1,0.3, 0.5, 0.7\}$ and $t_2$ noise, averaged over 500 simulations. This figure extends the last row of Figure~\ref{est.9} to other quantile levels.}
  \label{est.1357}
\end{figure}

\begin{figure}[!htp]
  \centering
  \subfigure[Model \eqref{model.homo} under $\tau = 0.1$.]{\includegraphics[width=0.32\textwidth]{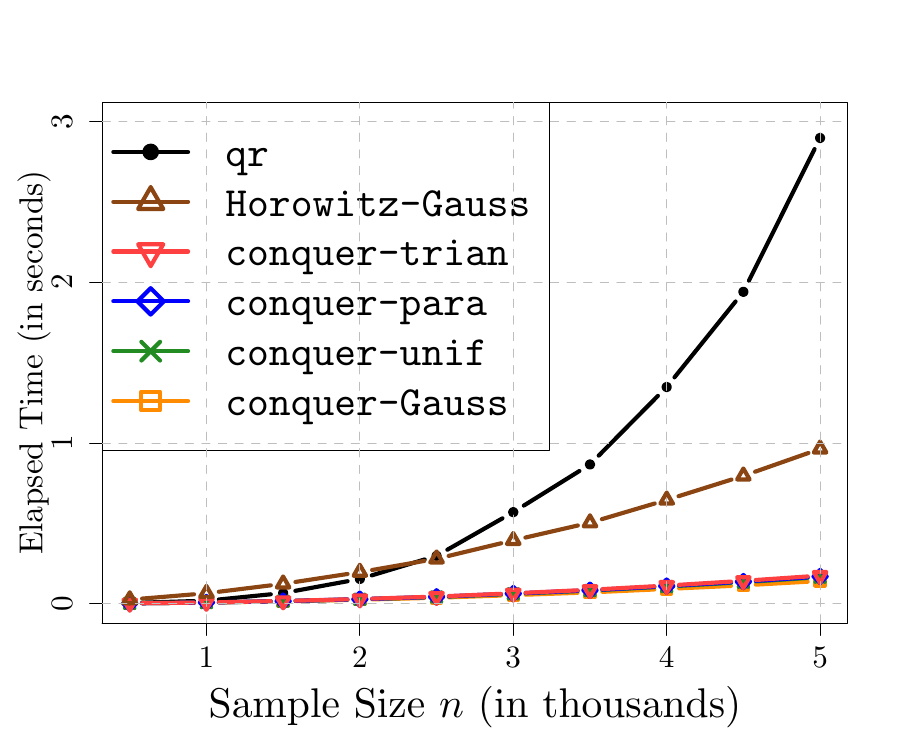}} 
  \subfigure[Model \eqref{model.linear} under $\tau = 0.1$.]{\includegraphics[width=0.32\textwidth]{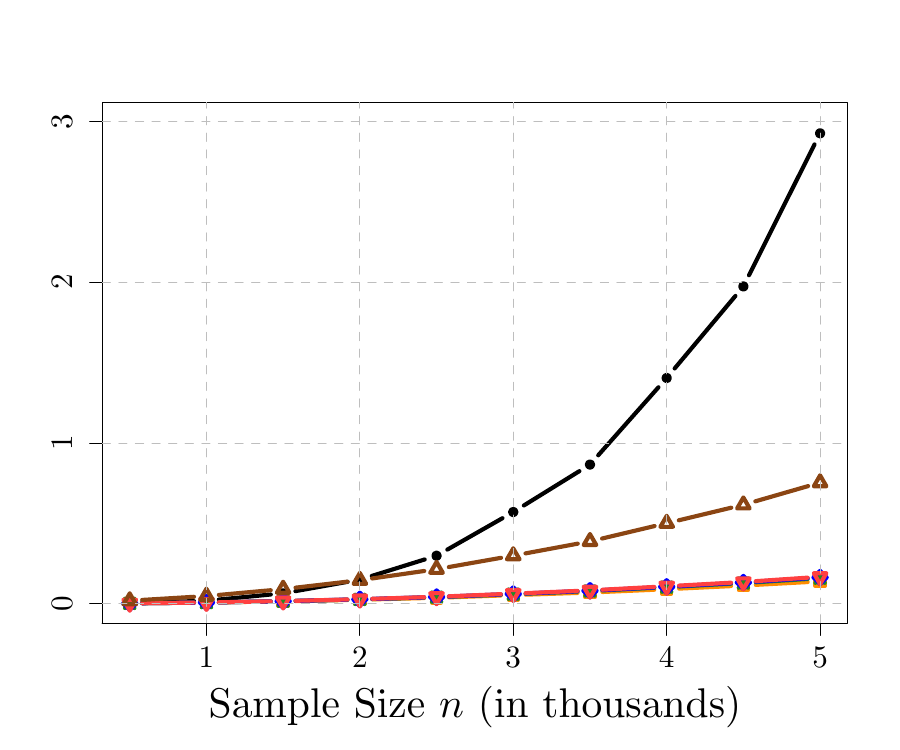}} 
  \subfigure[Model \eqref{model.quad} under $\tau = 0.1$.]{\includegraphics[width=0.32\textwidth]{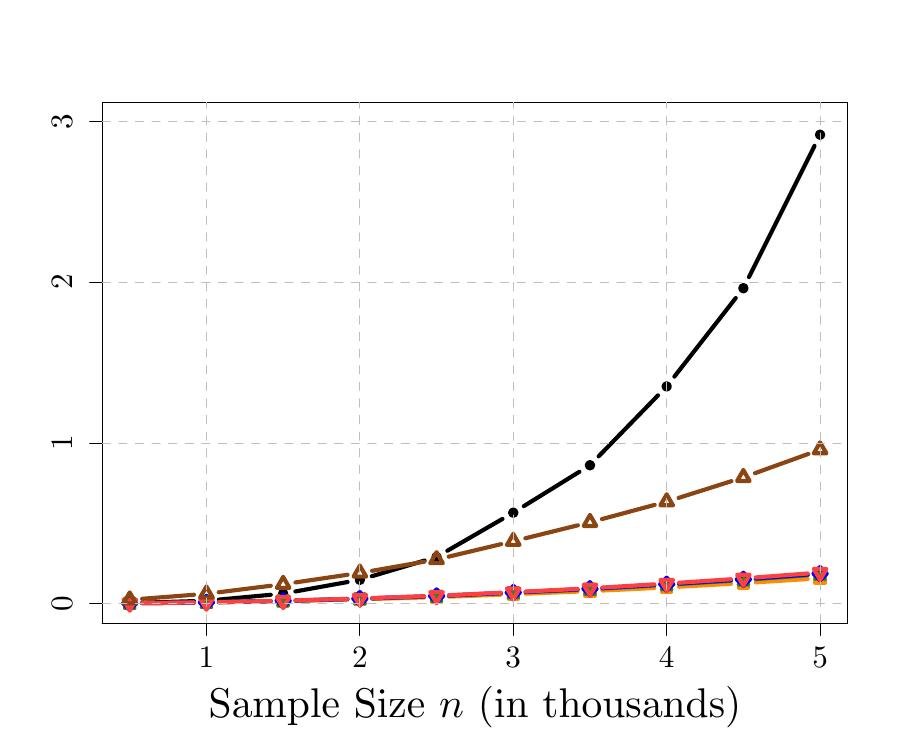}}
  \subfigure[Model \eqref{model.homo} under $\tau = 0.3$.]{\includegraphics[width=0.32\textwidth]{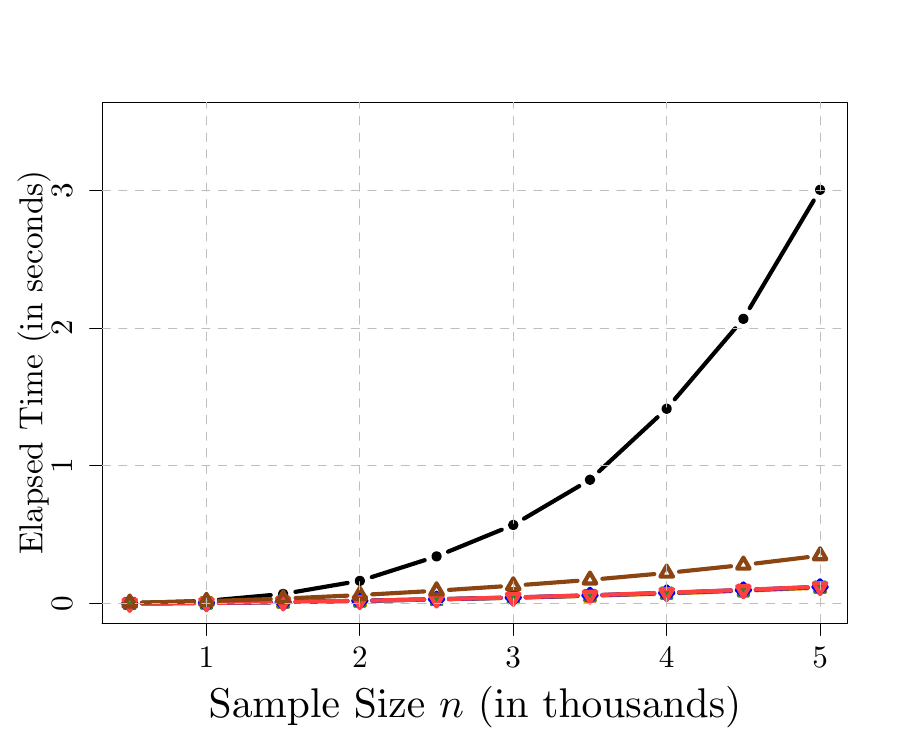}}
  \subfigure[Model \eqref{model.linear} under $\tau = 0.3$.]{\includegraphics[width=0.32\textwidth]{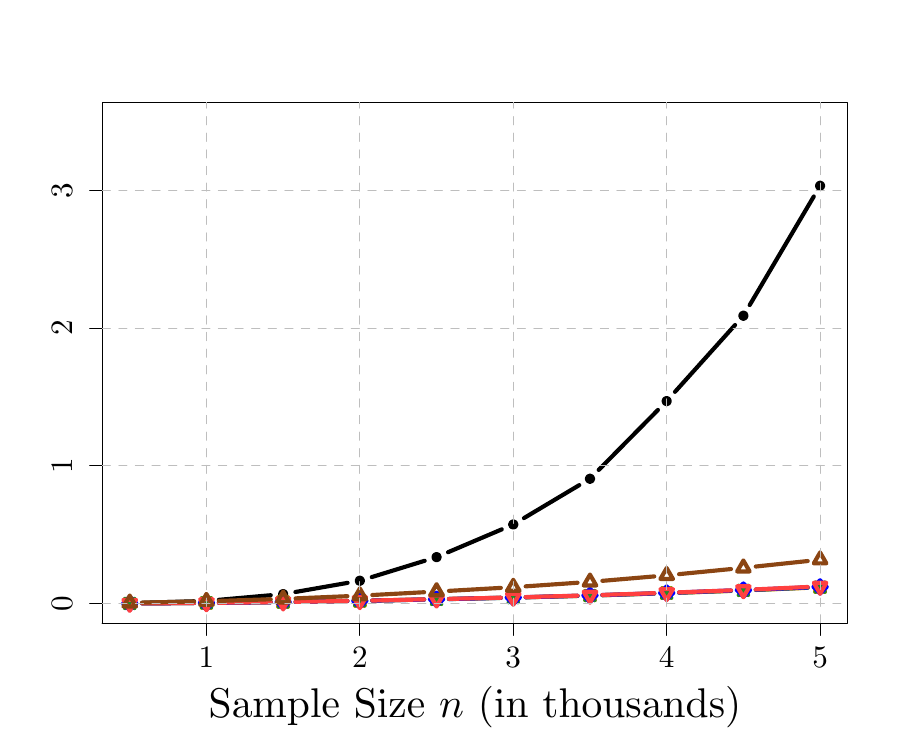}}
  \subfigure[Model \eqref{model.quad} under $\tau = 0.3$.]{\includegraphics[width=0.32\textwidth]{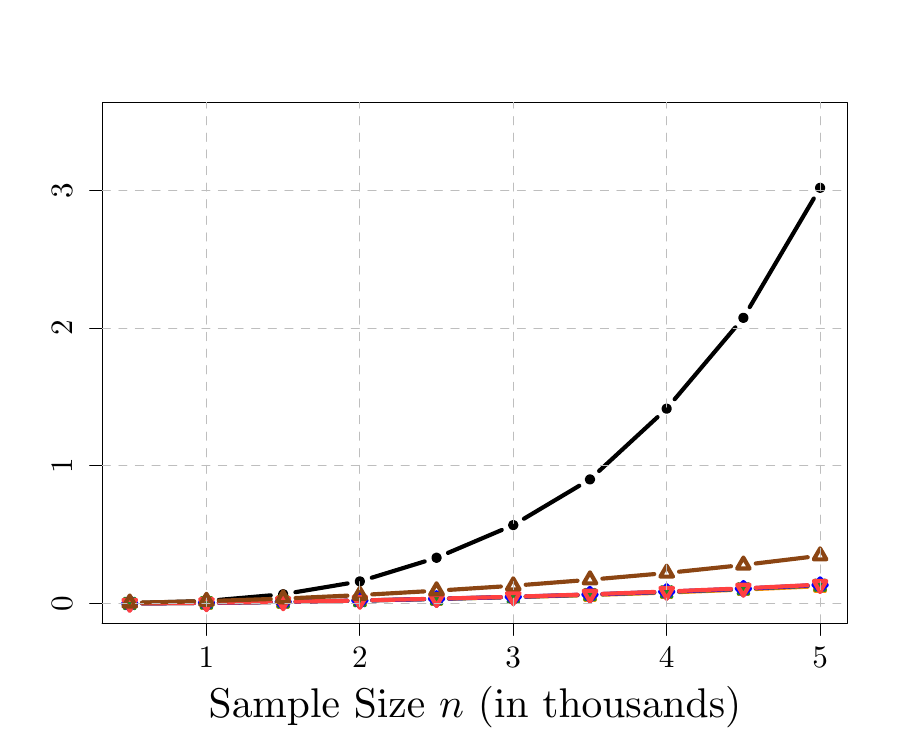}}
  \subfigure[Model \eqref{model.homo} under $\tau = 0.5$.]{\includegraphics[width=0.32\textwidth]{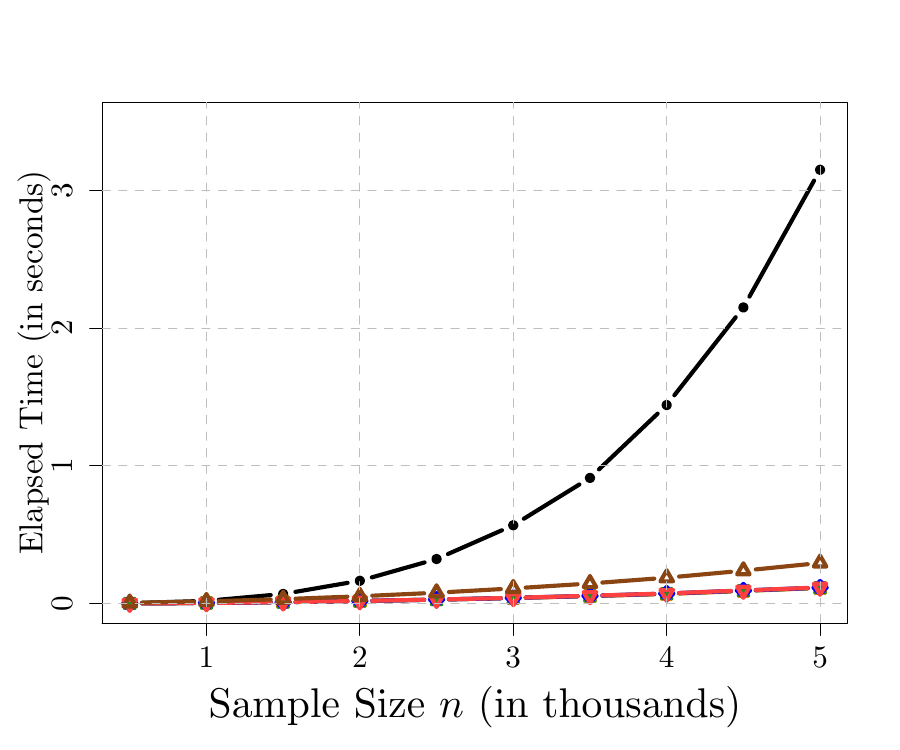}}
  \subfigure[Model \eqref{model.linear} under $\tau = 0.5$.]{\includegraphics[width=0.32\textwidth]{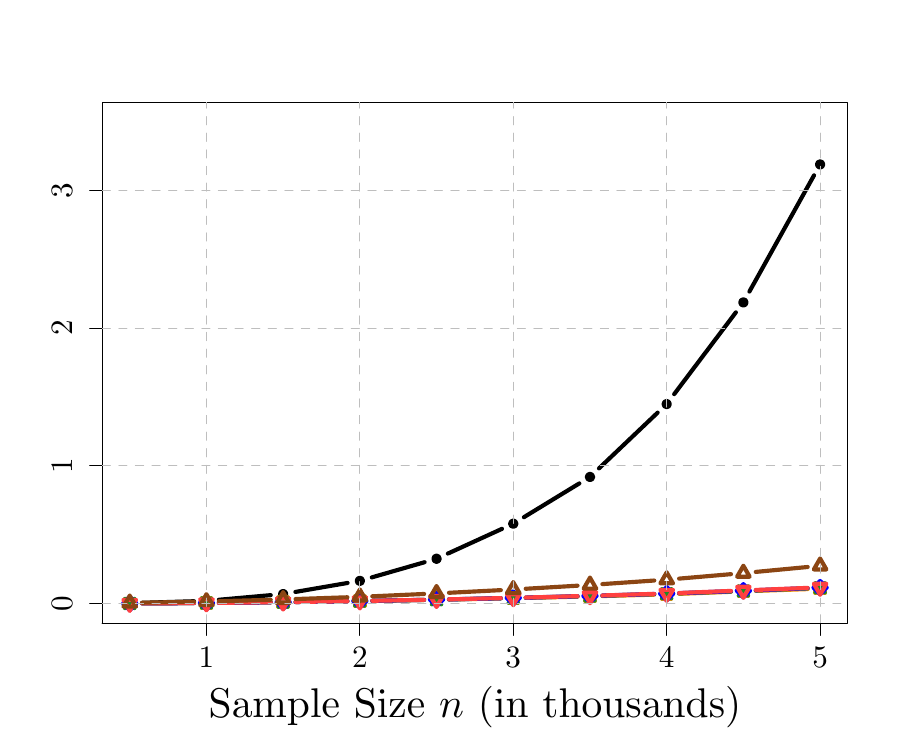}}
  \subfigure[Model \eqref{model.quad} under $\tau = 0.5$.]{\includegraphics[width=0.32\textwidth]{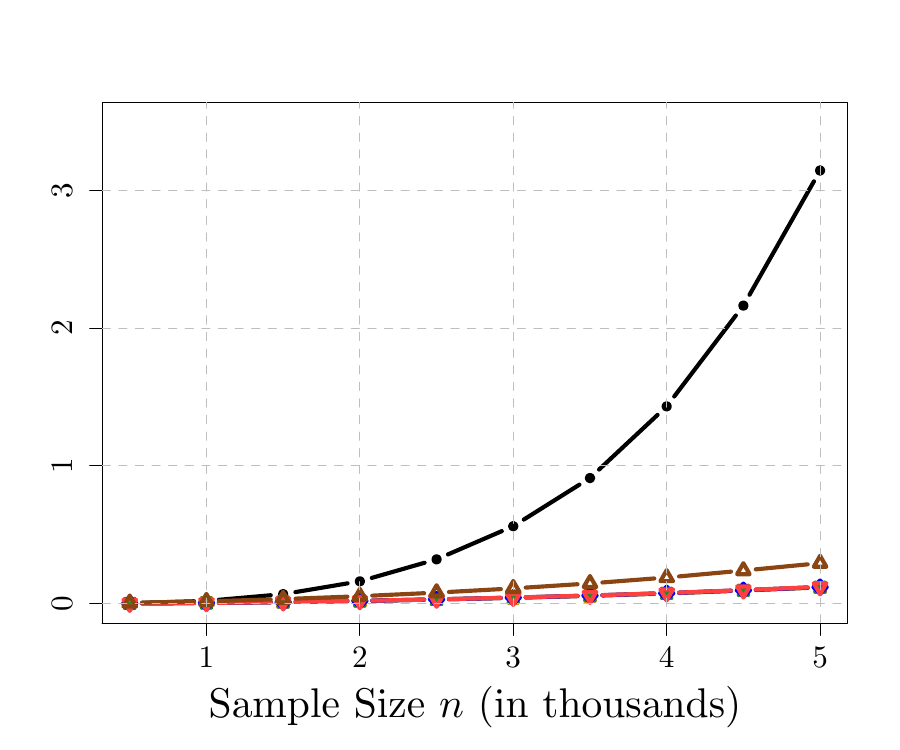}}
  \subfigure[Model \eqref{model.homo} under $\tau = 0.7$.]{\includegraphics[width=0.32\textwidth]{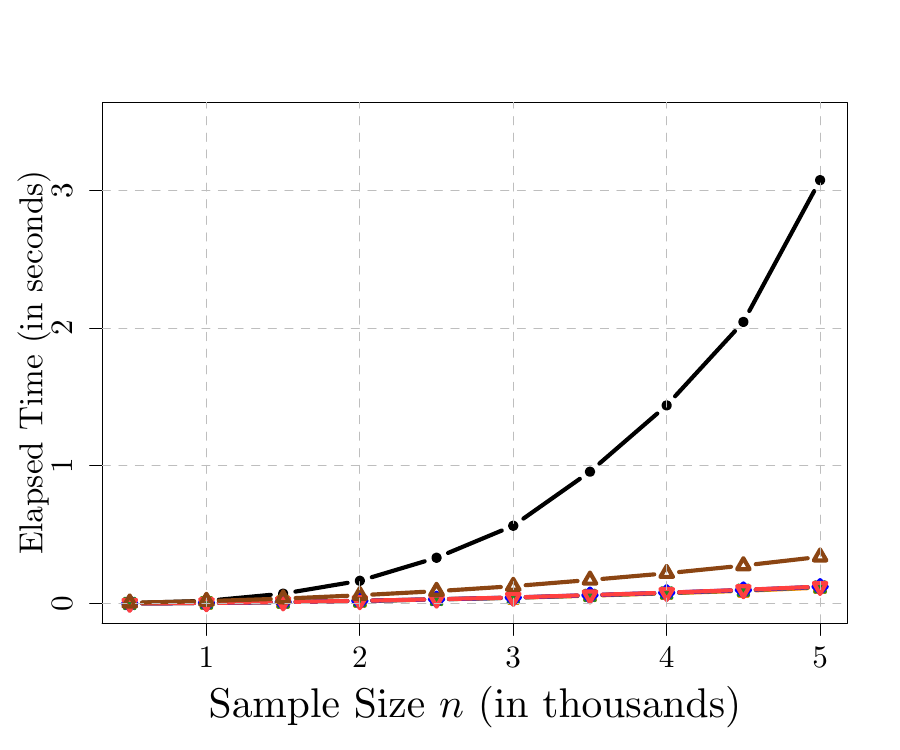}}
  \subfigure[Model \eqref{model.linear} under $\tau = 0.7$.]{\includegraphics[width=0.32\textwidth]{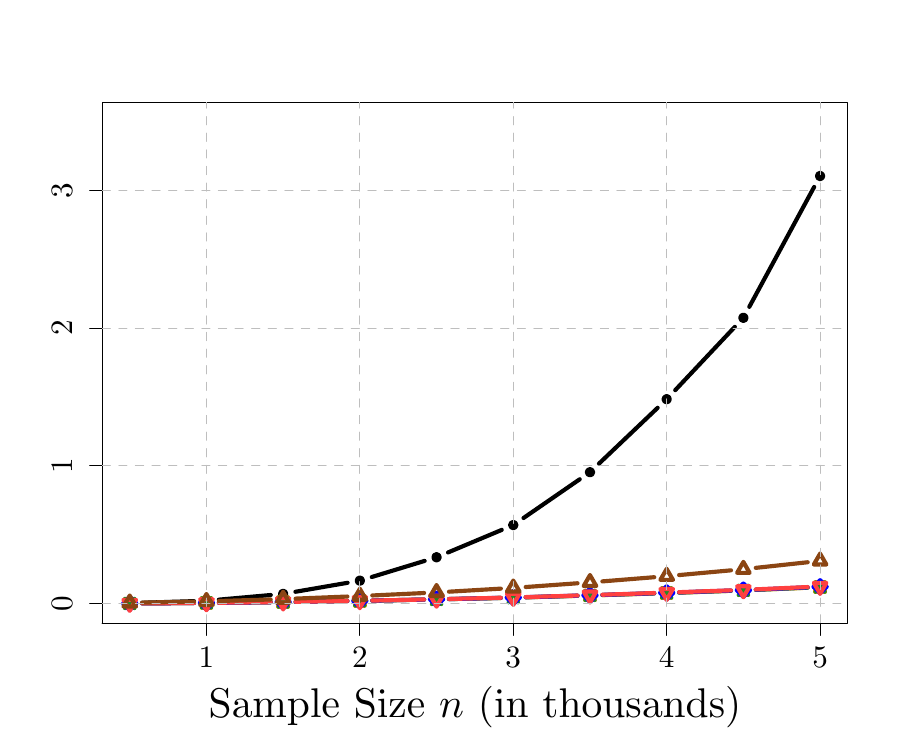}}
  \subfigure[Model \eqref{model.quad} under $\tau = 0.7$.]{\includegraphics[width=0.32\textwidth]{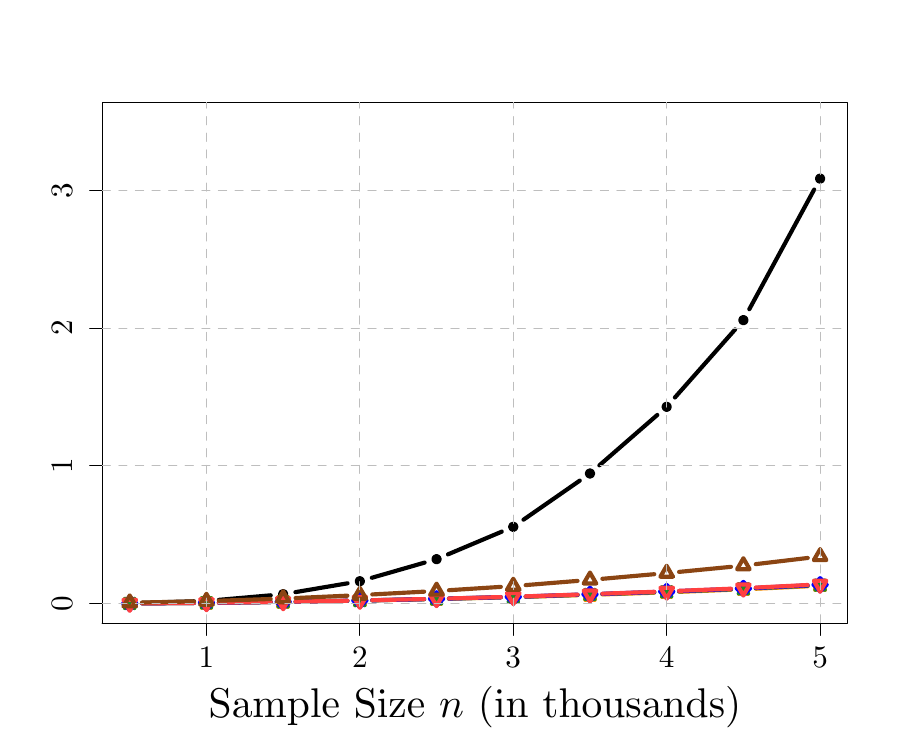}}
\caption{Elapsed time of standard QR, Horowitz's smoothing and conquer under $t_2$ noise setting when $\tau \in \{0.1, 0.3, 0.5, 0.7\}$. This figure extends the last row of Figure~\ref{time.9} to other quantile levels.}
  \label{time.1357}
\end{figure}

\begin{figure}[!htp]
  \centering
  \subfigure[Model \eqref{model.homo} under $\tau = 0.1$.]{\includegraphics[width=0.32\textwidth]{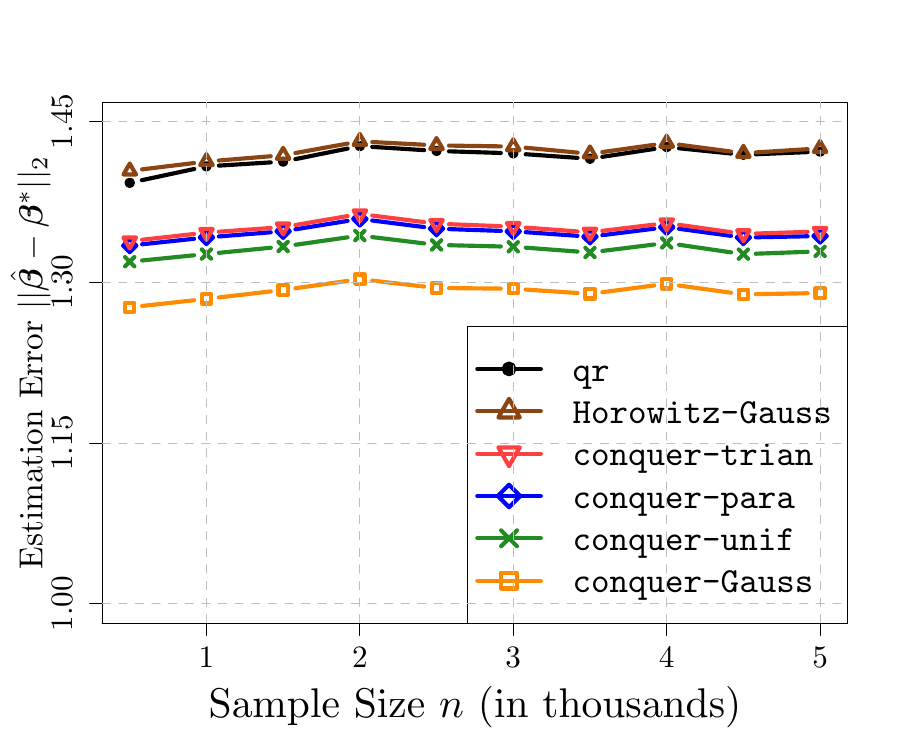}} 
  \subfigure[Model \eqref{model.homo} under $\tau = 0.3$.]{\includegraphics[width=0.32\textwidth]{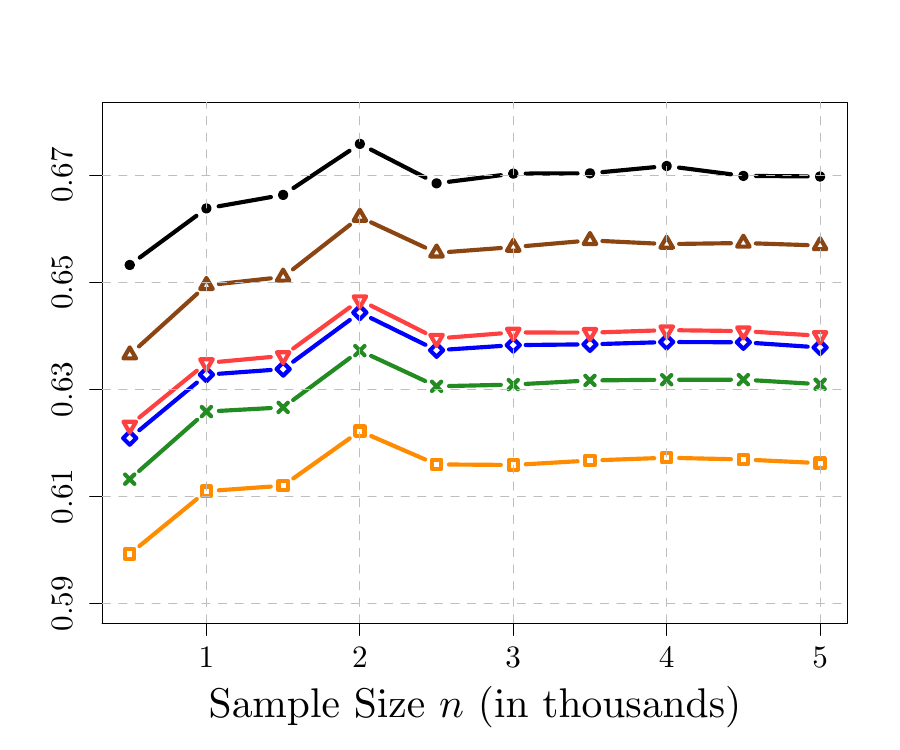}} 
  \subfigure[Model \eqref{model.homo} under $\tau = 0.5$.]{\includegraphics[width=0.32\textwidth]{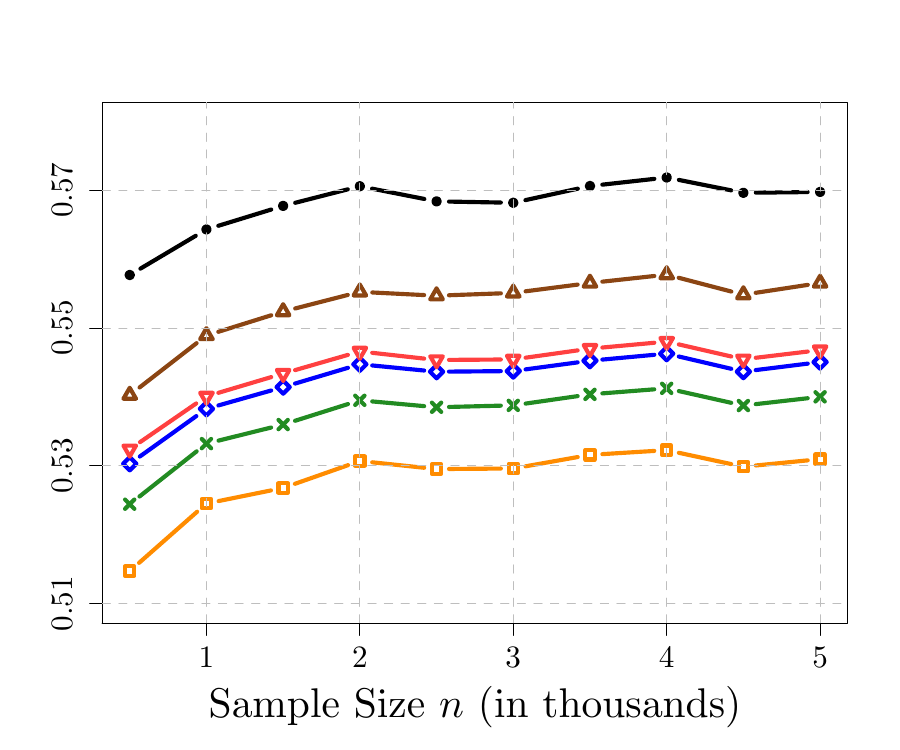}}
  \subfigure[Model \eqref{model.homo} under $\tau = 0.7$.]{\includegraphics[width=0.32\textwidth]{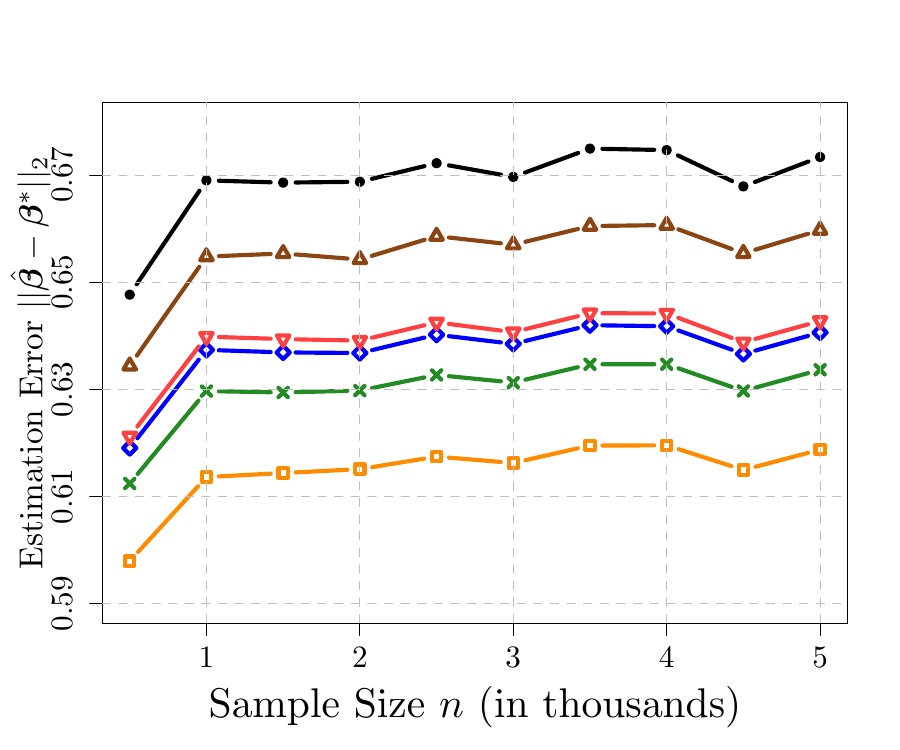}}
  \subfigure[Model \eqref{model.homo} under $\tau = 0.9$.]{\includegraphics[width=0.32\textwidth]{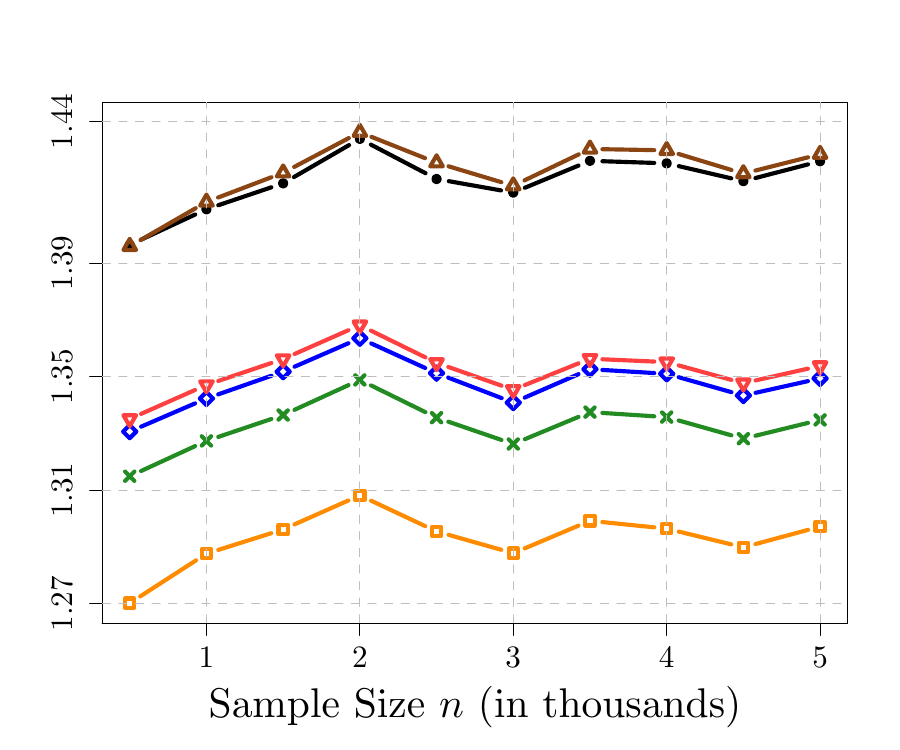}} \\
  \subfigure[Model \eqref{model.quad} under $\tau = 0.1$.]{\includegraphics[width=0.32\textwidth]{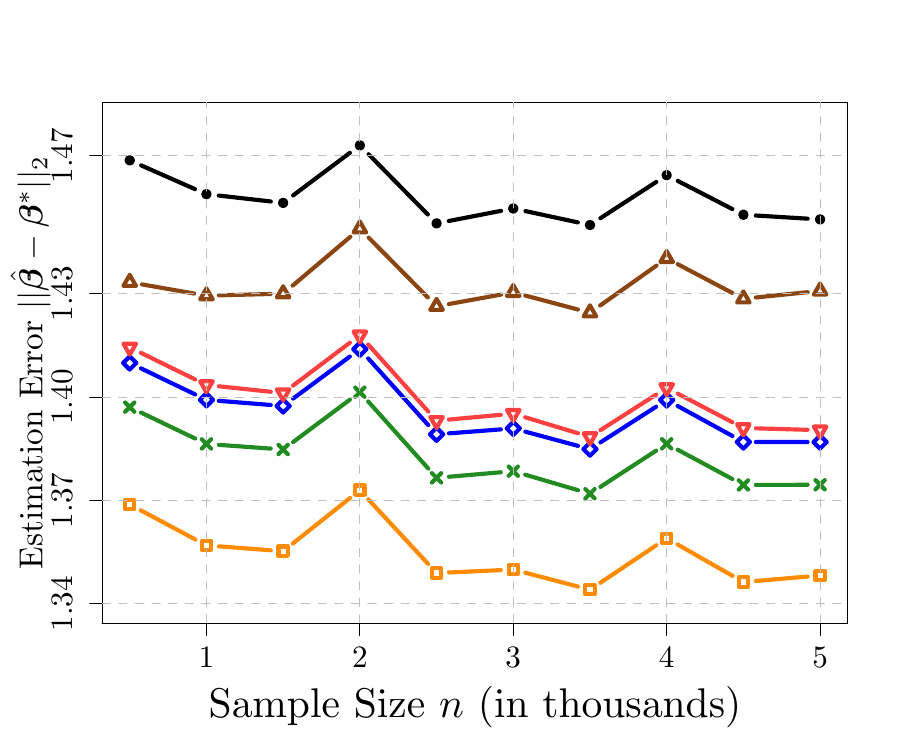}}
  \subfigure[Model \eqref{model.quad} under $\tau = 0.3$.]{\includegraphics[width=0.32\textwidth]{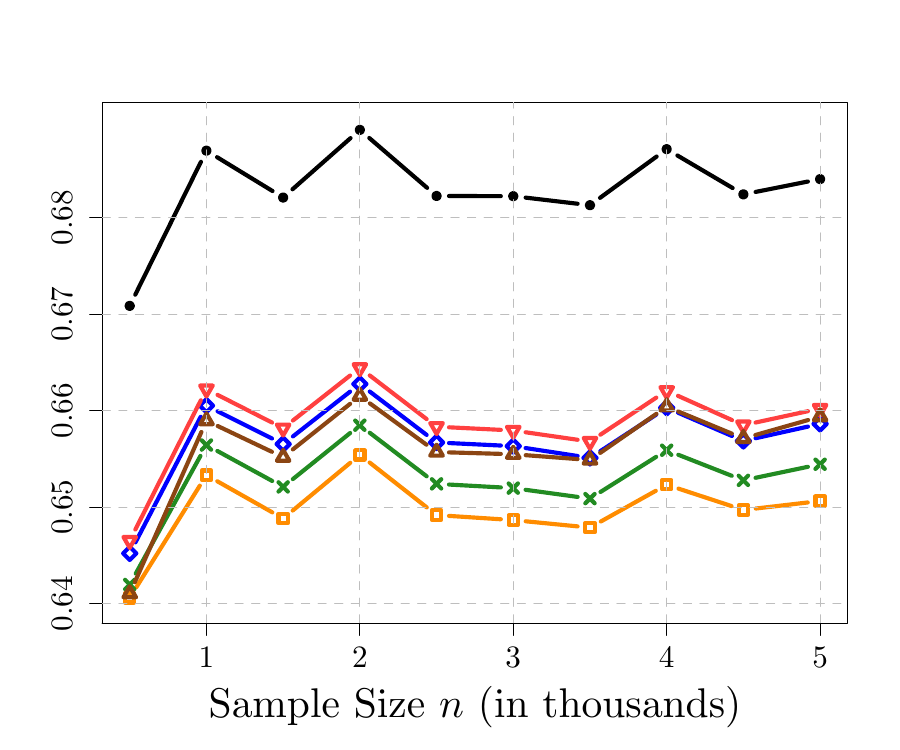}}
  \subfigure[Model \eqref{model.quad} under $\tau = 0.5$.]{\includegraphics[width=0.32\textwidth]{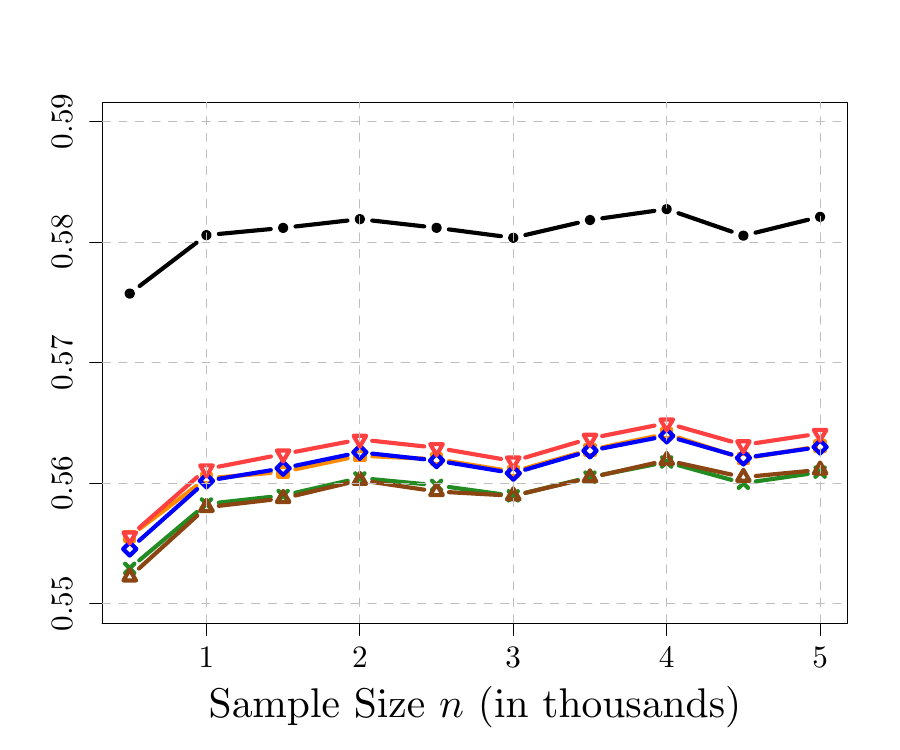}}
  \subfigure[Model \eqref{model.quad} under $\tau = 0.7$.]{\includegraphics[width=0.32\textwidth]{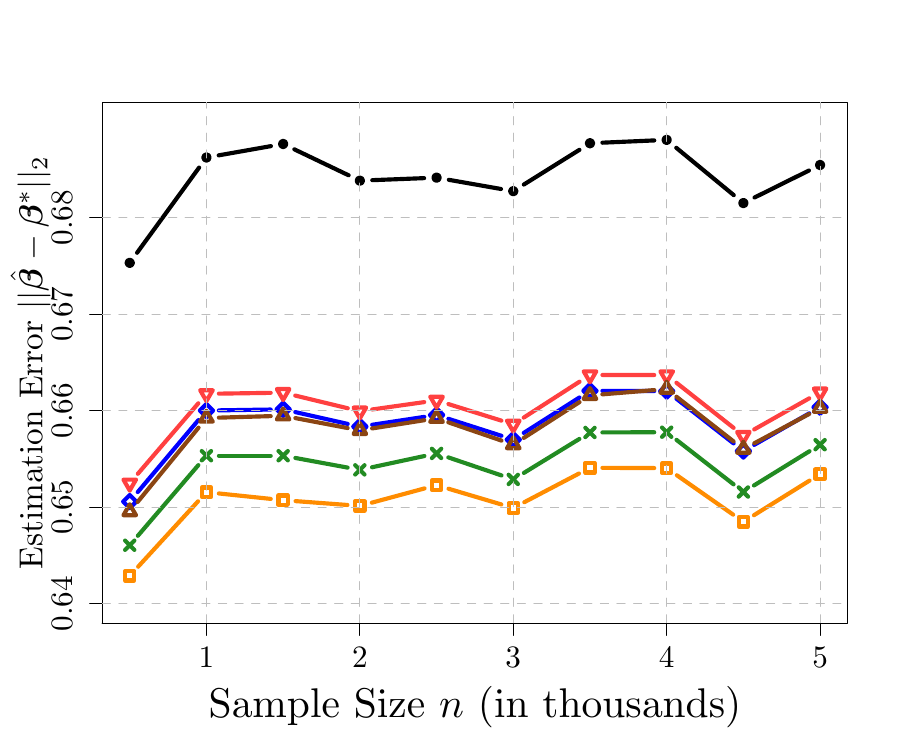}}
  \subfigure[Model \eqref{model.quad} under $\tau = 0.9$.]{\includegraphics[width=0.32\textwidth]{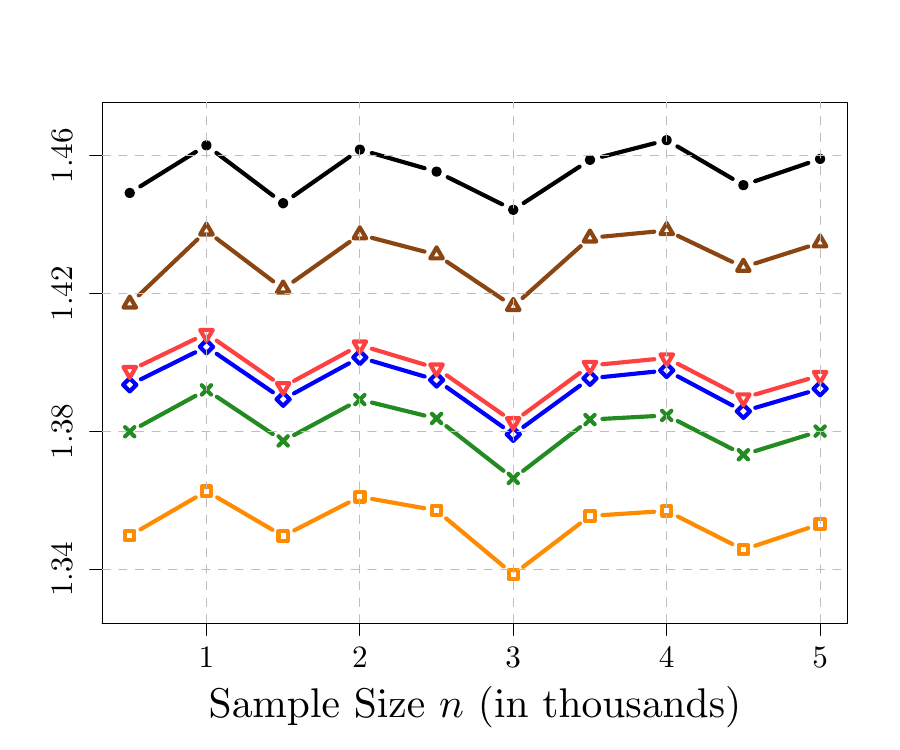}}
\caption{Estimation errors of standard QR, Horowitz's smoothing and conquer under multivariate normal random design and $t_2$ noise when $\tau \in \{0.1, 0.3, 0.5, 0.7, 0.9\}$.}
  \label{est.unbd}
\end{figure}

\begin{figure}[!htp]
  \centering
  \subfigure[Model \eqref{model.homo} under $\tau = 0.1$.]{\includegraphics[width=0.32\textwidth]{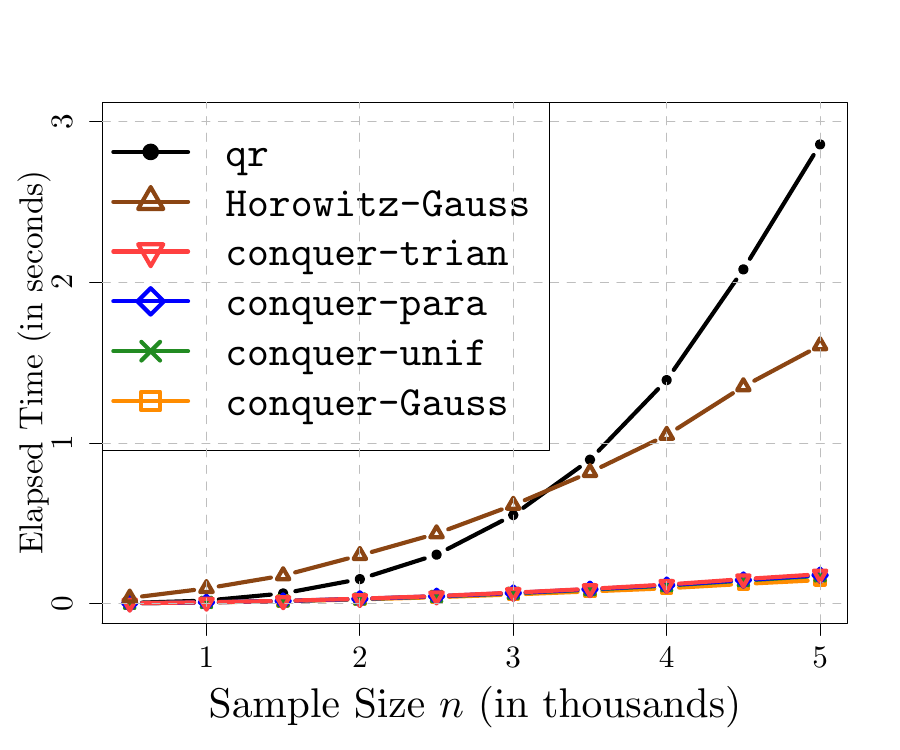}} 
  \subfigure[Model \eqref{model.homo} under $\tau = 0.3$.]{\includegraphics[width=0.32\textwidth]{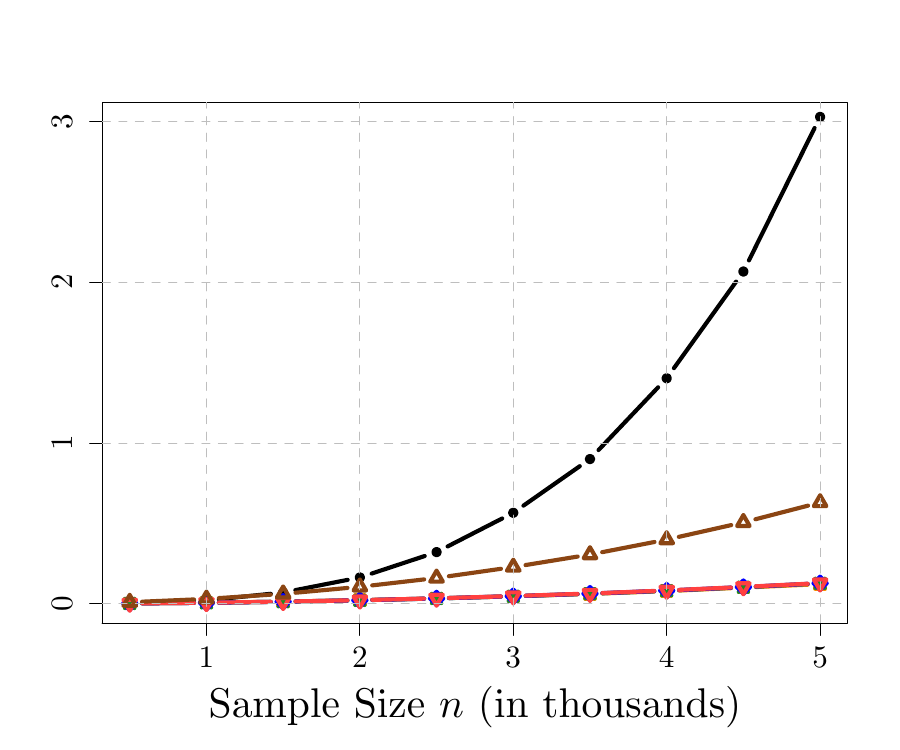}} 
  \subfigure[Model \eqref{model.homo} under $\tau = 0.5$.]{\includegraphics[width=0.32\textwidth]{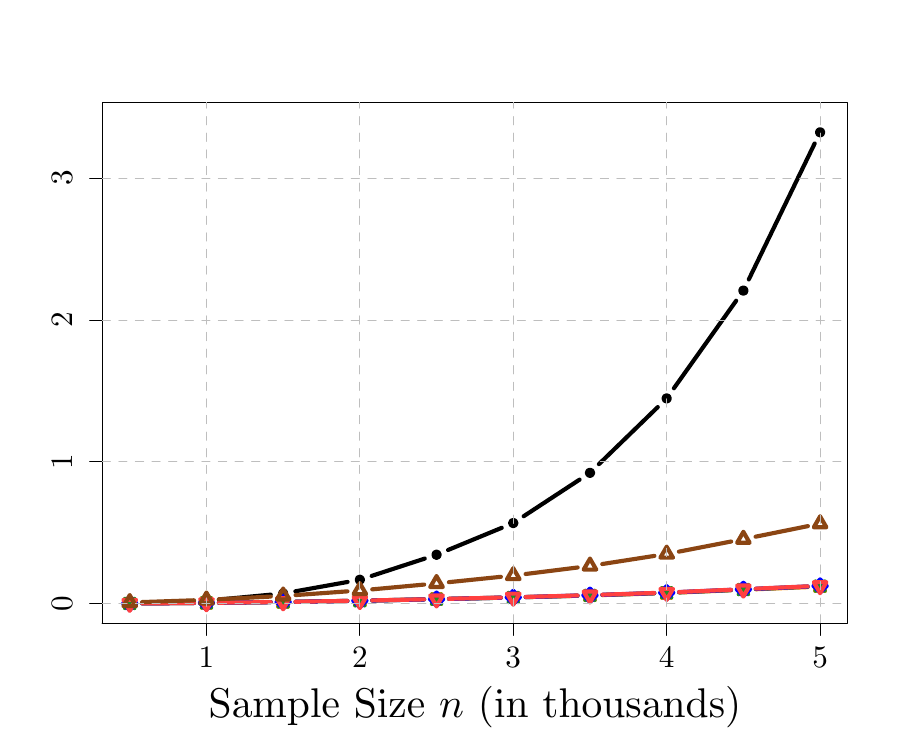}}
  \subfigure[Model \eqref{model.homo} under $\tau = 0.7$.]{\includegraphics[width=0.32\textwidth]{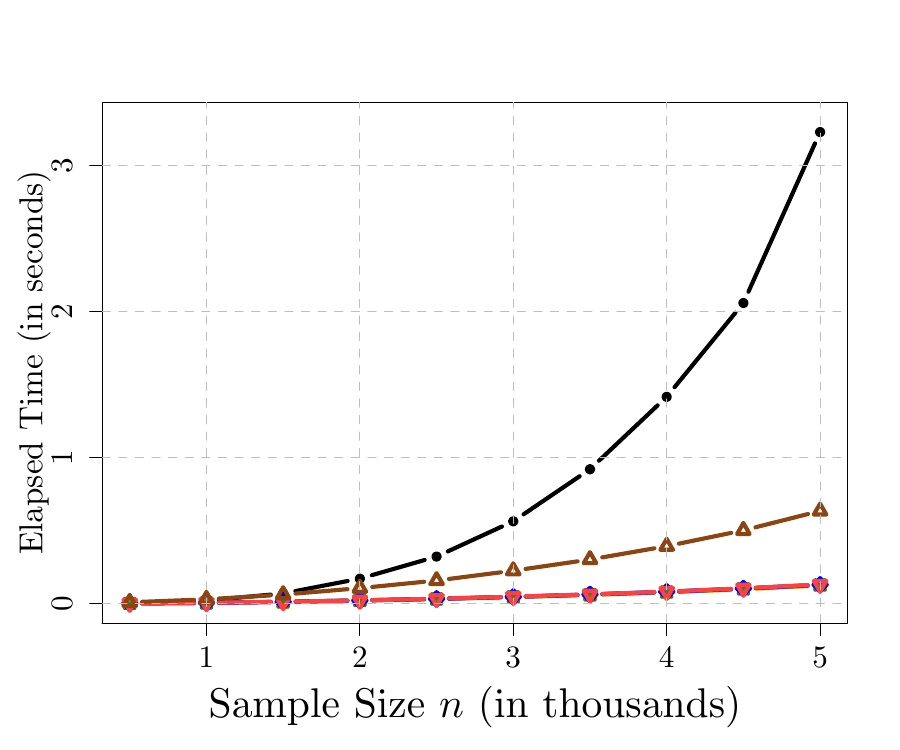}}
  \subfigure[Model \eqref{model.homo} under $\tau = 0.9$.]{\includegraphics[width=0.32\textwidth]{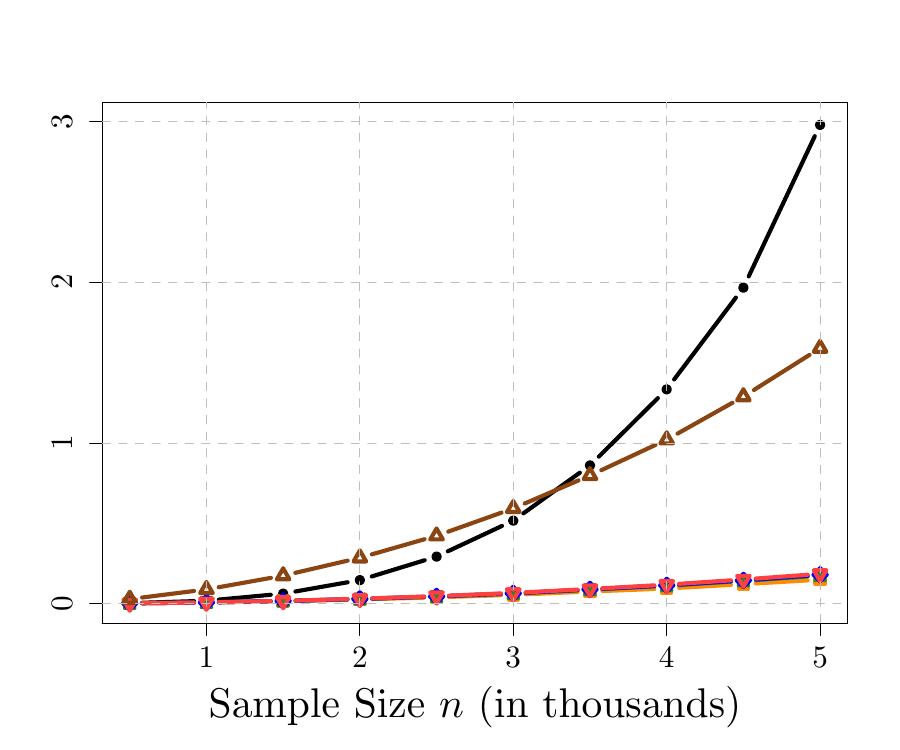}} \\
  \subfigure[Model \eqref{model.quad} under $\tau = 0.1$.]{\includegraphics[width=0.32\textwidth]{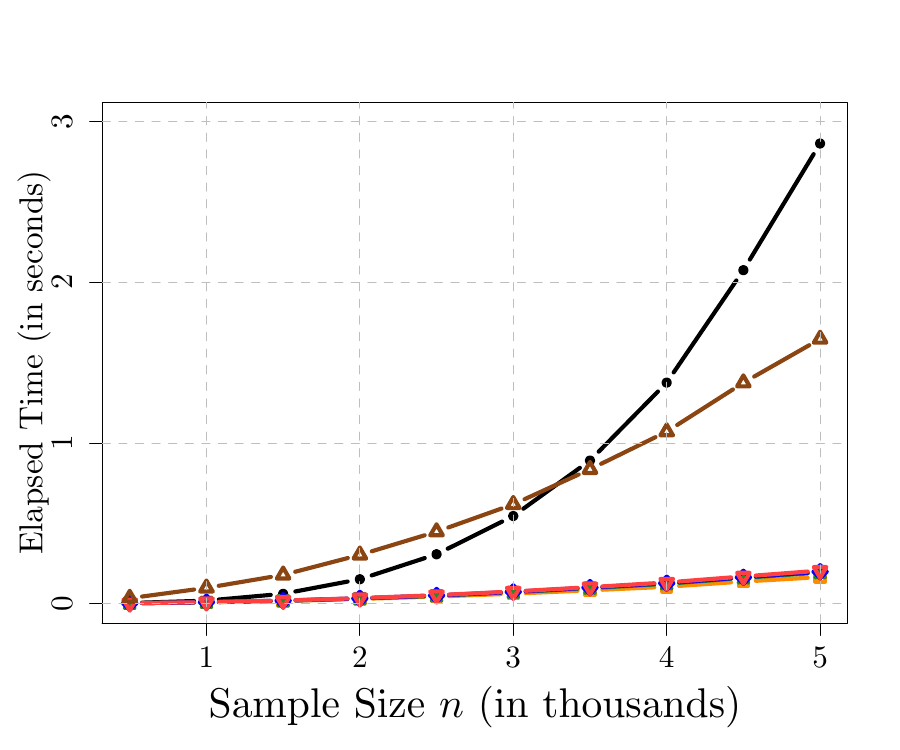}}
  \subfigure[Model \eqref{model.quad} under $\tau = 0.3$.]{\includegraphics[width=0.32\textwidth]{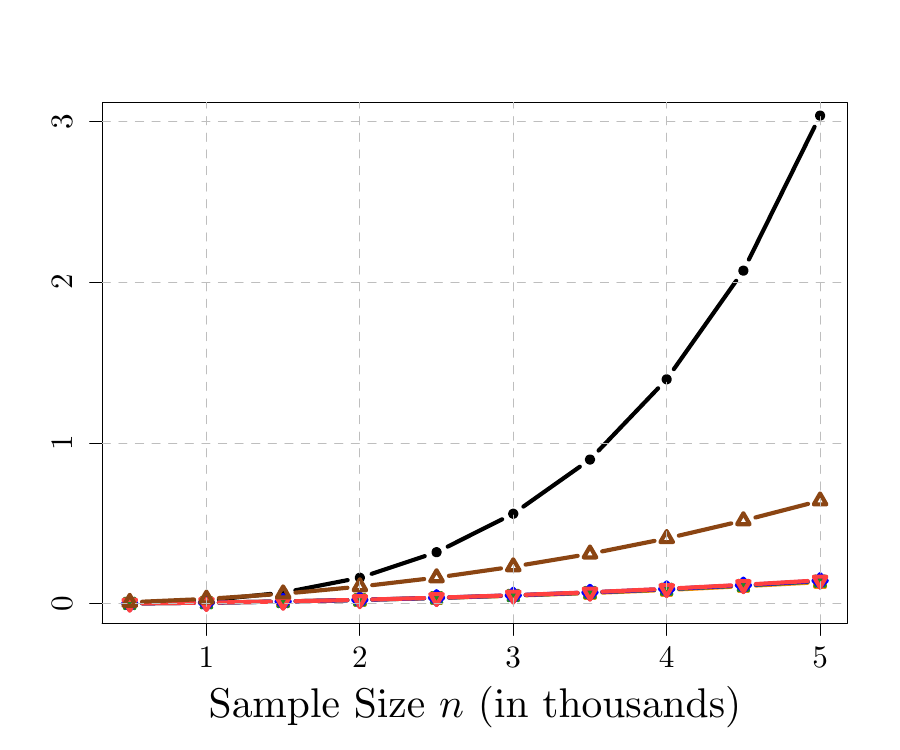}}
  \subfigure[Model \eqref{model.quad} under $\tau = 0.5$.]{\includegraphics[width=0.32\textwidth]{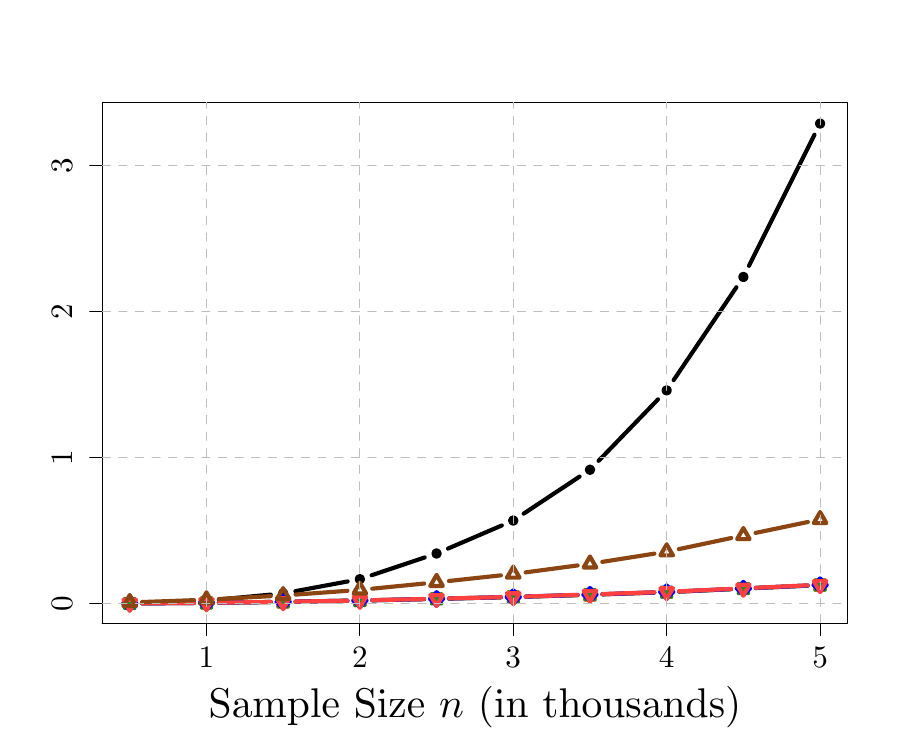}}
  \subfigure[Model \eqref{model.quad} under $\tau = 0.7$.]{\includegraphics[width=0.32\textwidth]{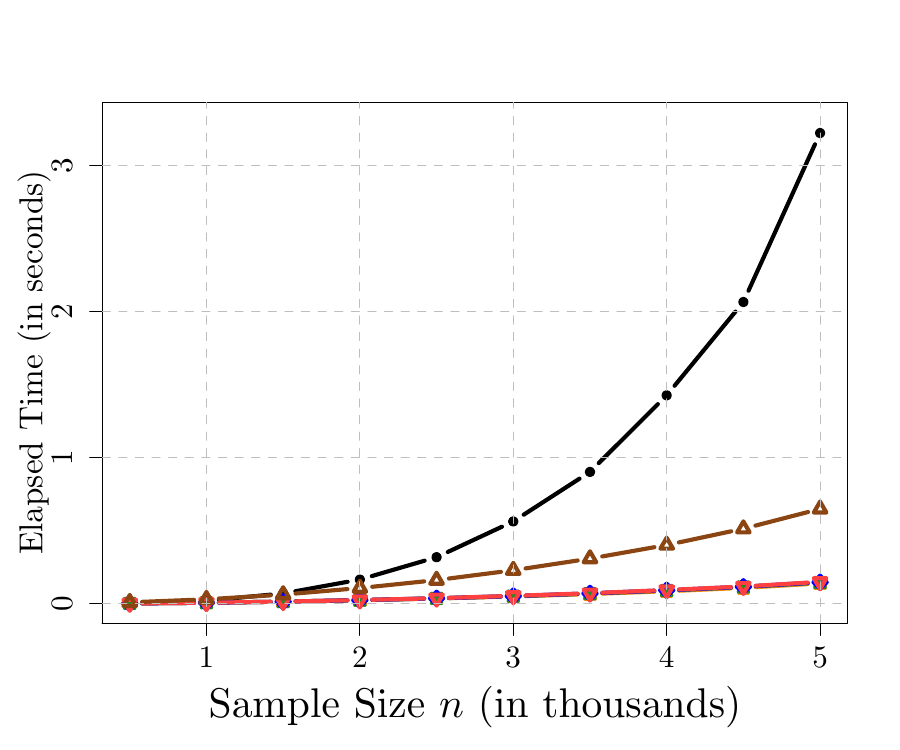}}
  \subfigure[Model \eqref{model.quad} under $\tau = 0.9$.]{\includegraphics[width=0.32\textwidth]{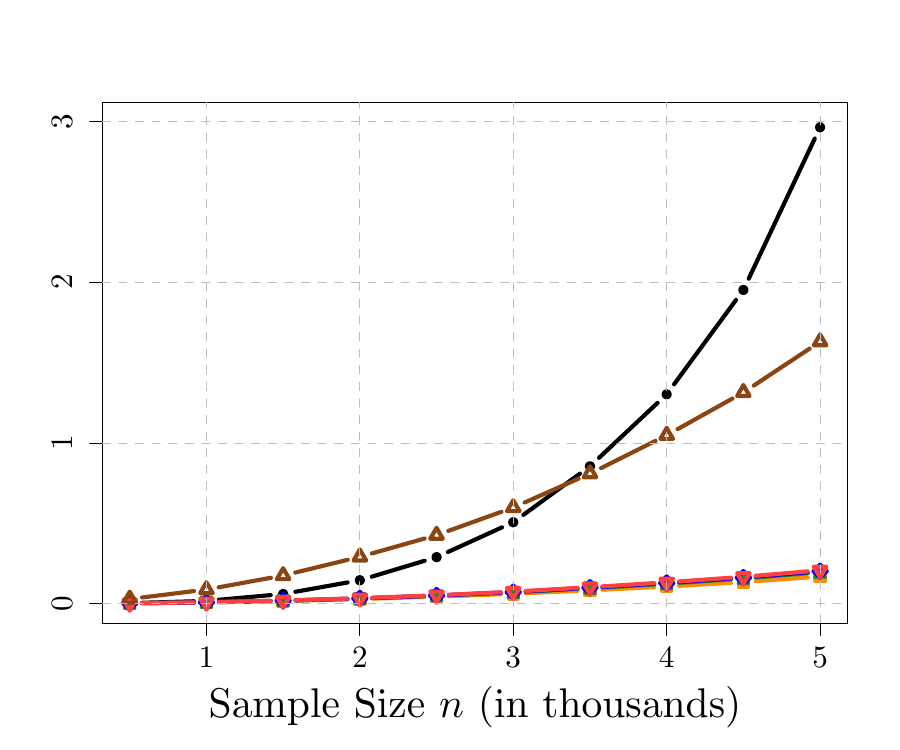}}
\caption{Elapsed time of standard QR, Horowitz's smoothing and conquer under multivariate normal random design and $t_2$ noise when $\tau \in \{0.1, 0.3, 0.5, 0.7, 0.9\}$.}
  \label{time.unbd}
\end{figure}

\begin{figure}[!htp]
  \centering
  \subfigure[$(n, p) = (2000, 25)$.]{\includegraphics[width=0.32\textwidth]{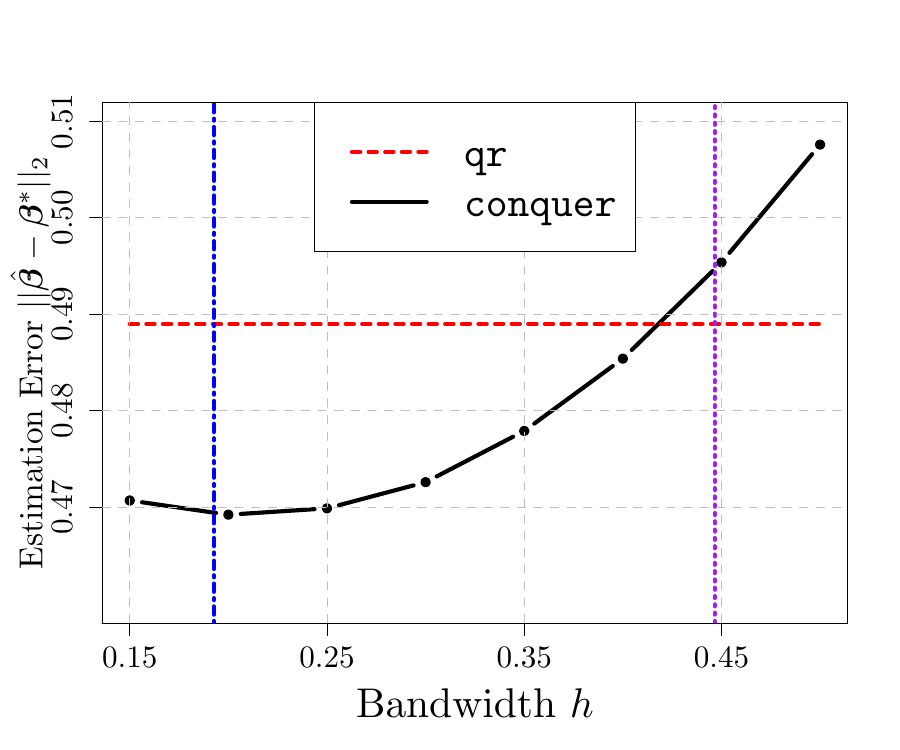}} 
  \subfigure[$(n, p) = (2000, 50)$.]{\includegraphics[width=0.32\textwidth]{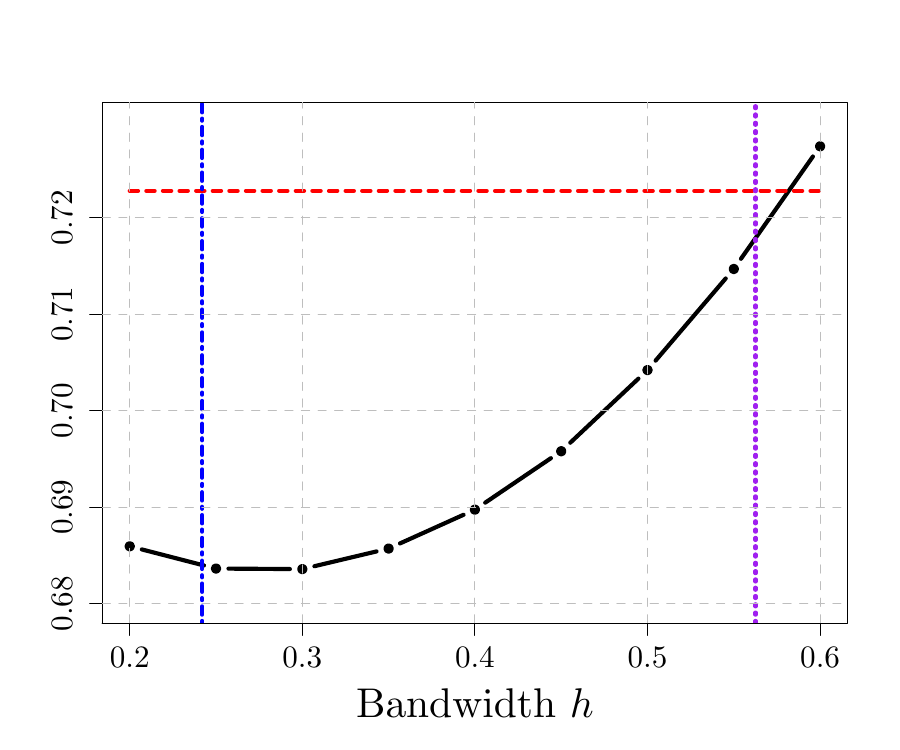}} 
  \subfigure[$(n, p) = (2000, 75)$.]{\includegraphics[width=0.32\textwidth]{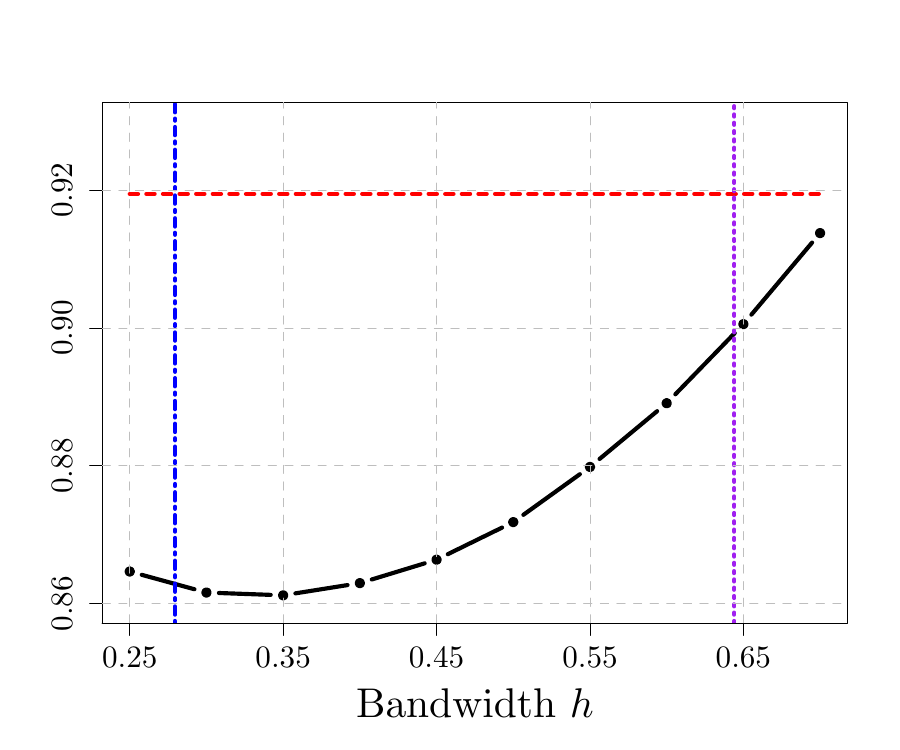}}
  \subfigure[$(n, p) = (2000, 100)$.]{\includegraphics[width=0.32\textwidth]{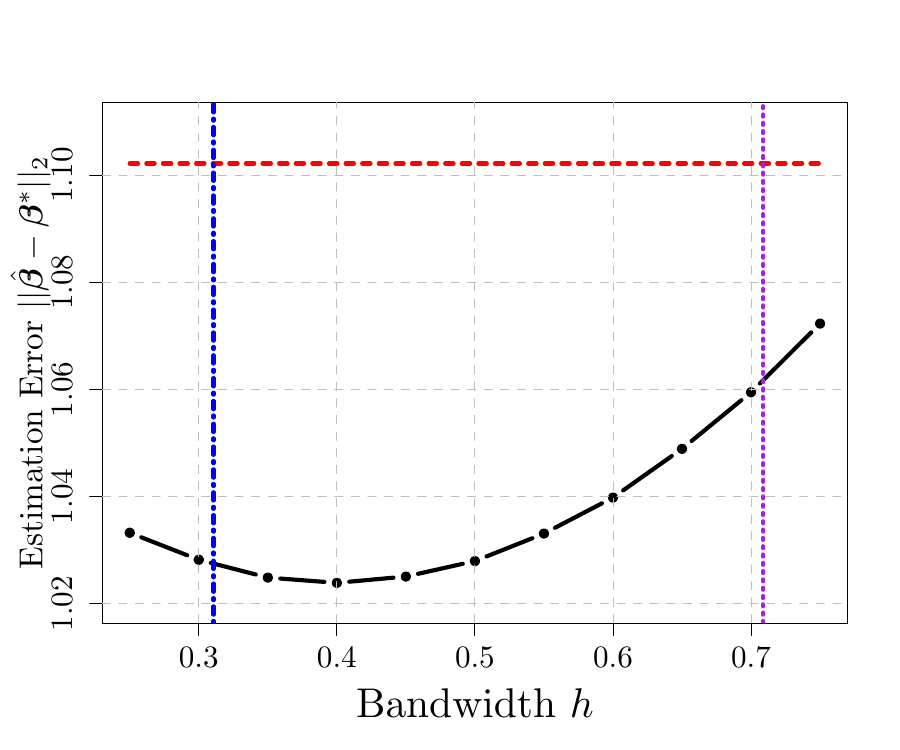}}
  \subfigure[$(n, p) = (2000, 125)$.]{\includegraphics[width=0.32\textwidth]{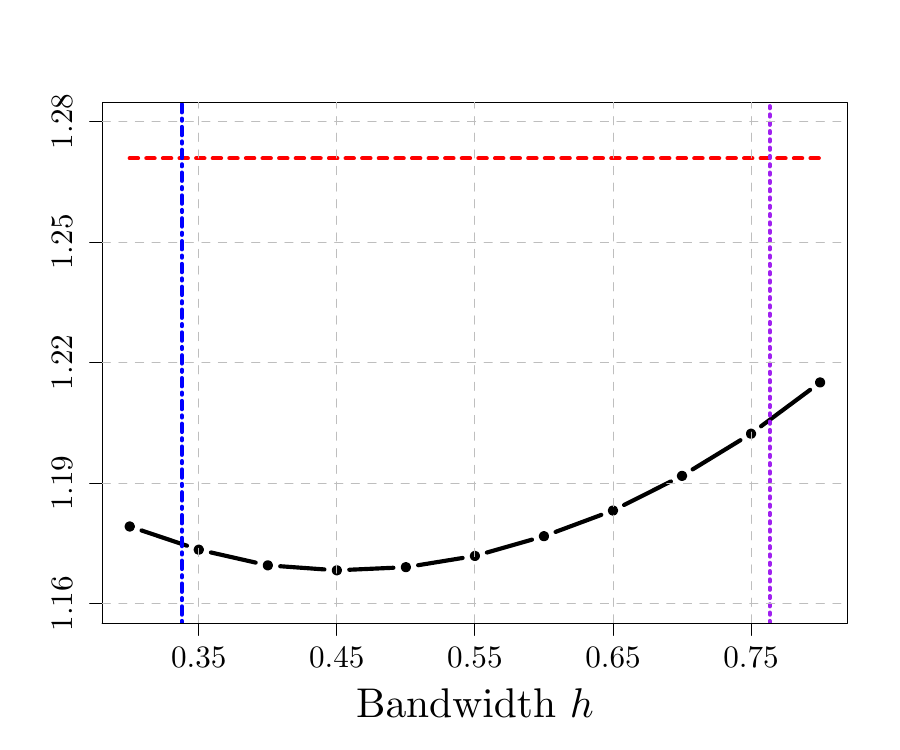}}
  \subfigure[$(n, p) = (2000, 150)$.]{\includegraphics[width=0.32\textwidth]{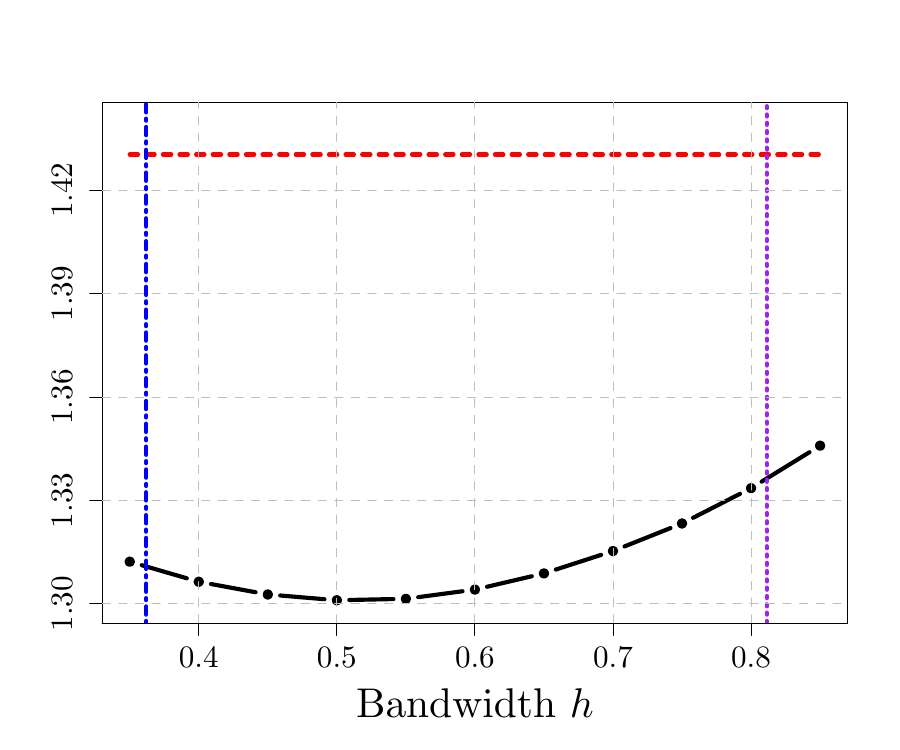}}
\caption{Sensitivity analysis of conquer using a Gaussian kernel under model \eqref{model.linear} and $t_2$ error with a range of bandwidth parameter $h$. Results for $n = 2000$ and $p \in \{25, 50, 75, 100, 125, 150\}$, averaged over 500 datasets. The blue vertical dash line represents our default choice $h_{{\rm de}} = \{ (p+\log n) /n\}^{2/5}$, the purple vertical dash line refers to the choice $h_{\mathrm{AMSE}}$ from \eqref{AMSE.bandwidth2} with some oracle knowledge, and the red horizontal dash line represents the estimation error of standard QR.}
  \label{sensitivity.grow}
\end{figure}

\begin{figure}[!htp]
  \centering
  \subfigure[Coverage under model \eqref{model.homo}.]{\includegraphics[width=0.32\textwidth]{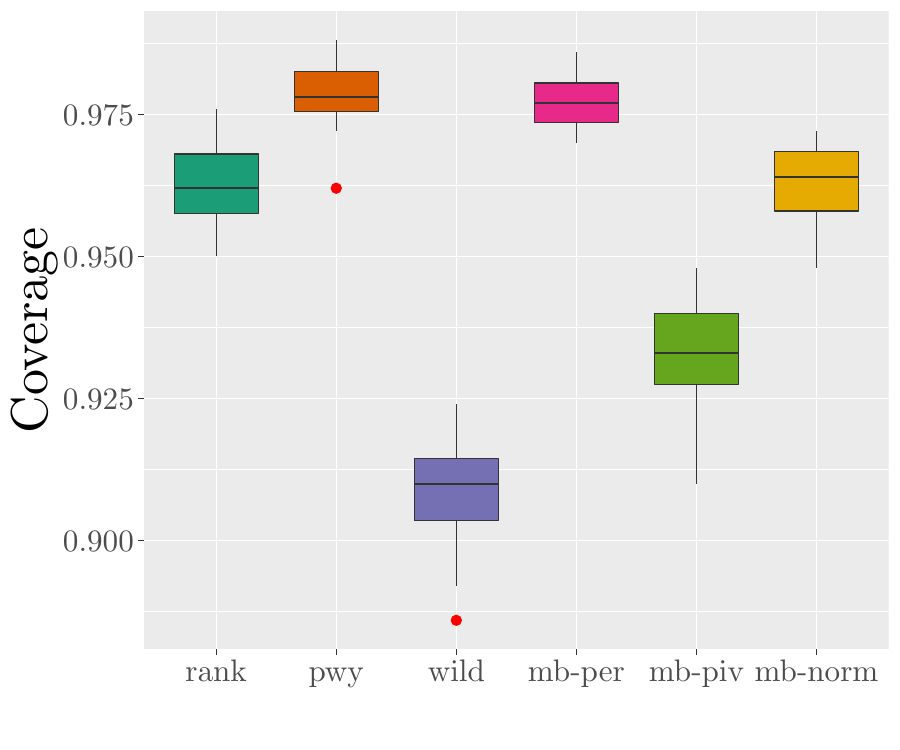}} 
  \subfigure[Coverage under model \eqref{model.linear}.]{\includegraphics[width=0.32\textwidth]{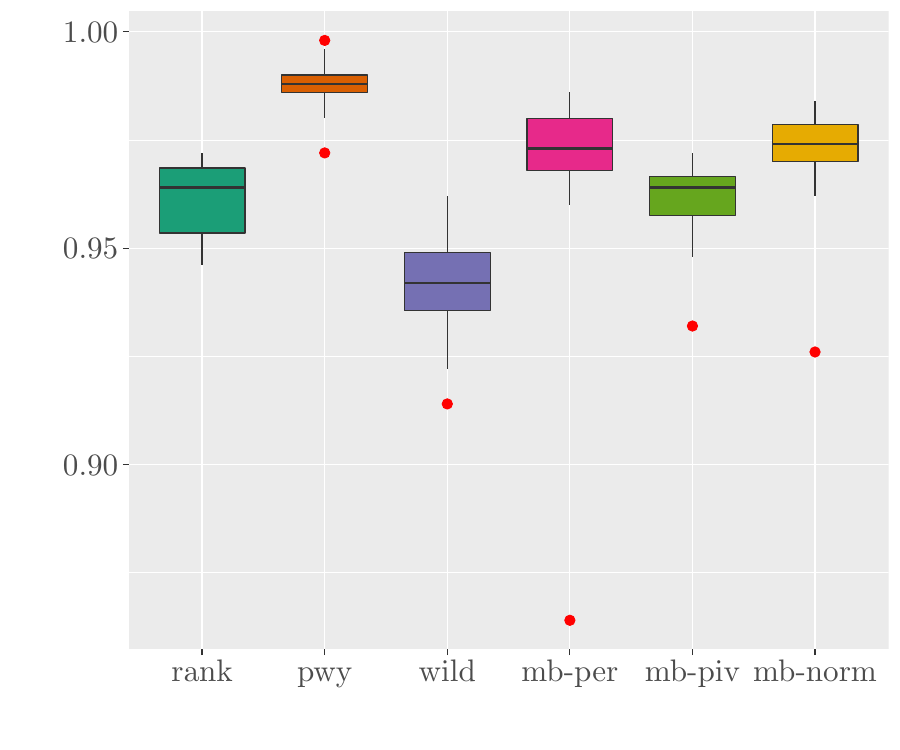}} 
  \subfigure[Coverage under model \eqref{model.quad}.]{\includegraphics[width=0.32\textwidth]{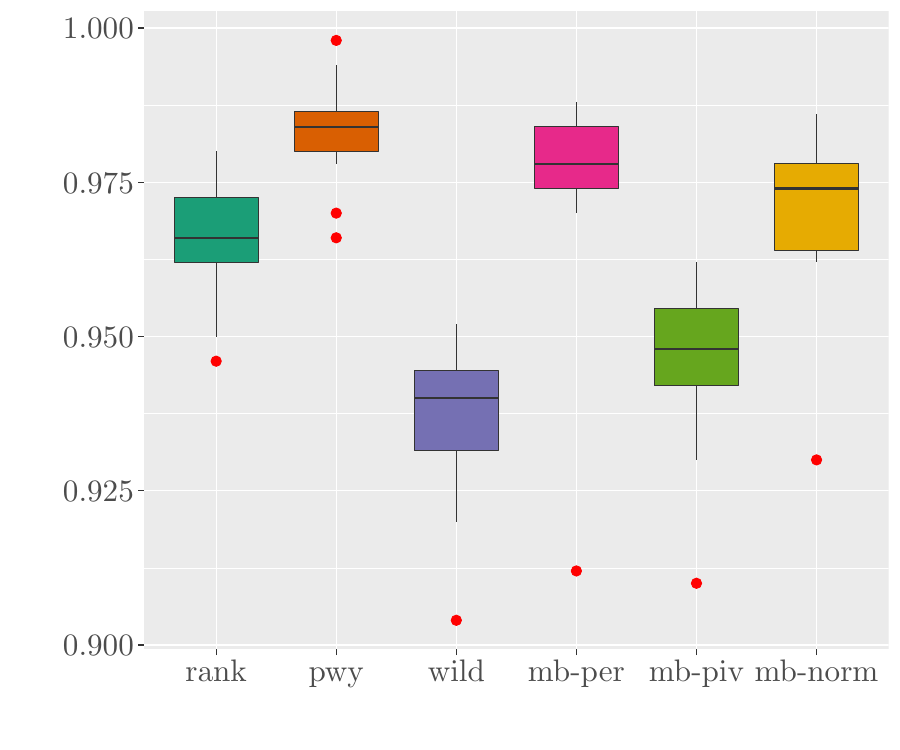}}
  \subfigure[CI width under model \eqref{model.homo}.]{\includegraphics[width=0.32\textwidth]{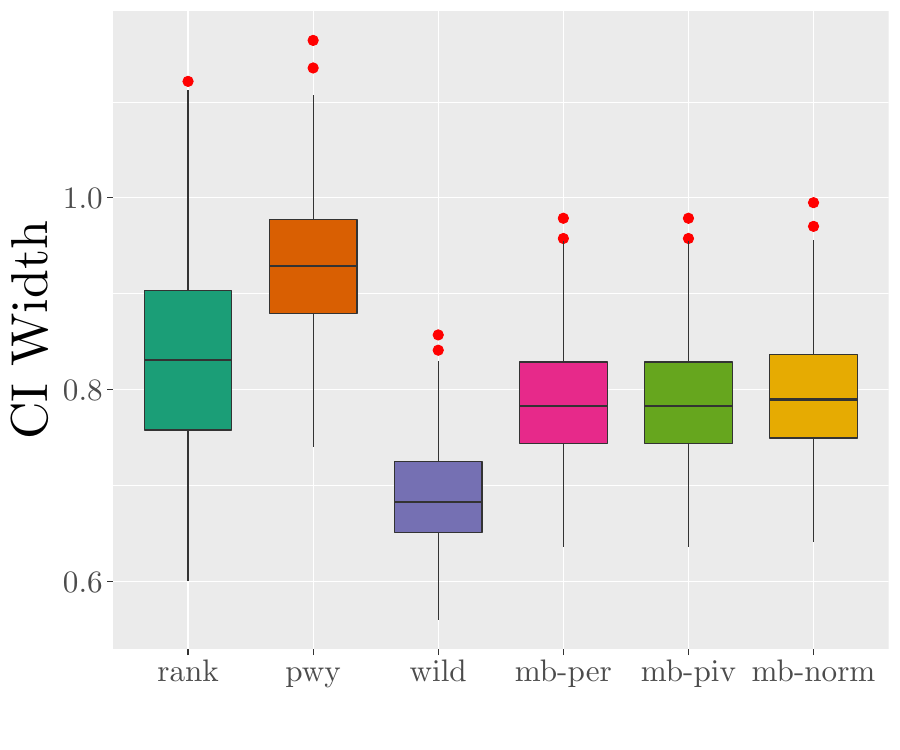}}
  \subfigure[CI width under model \eqref{model.linear}.]{\includegraphics[width=0.32\textwidth]{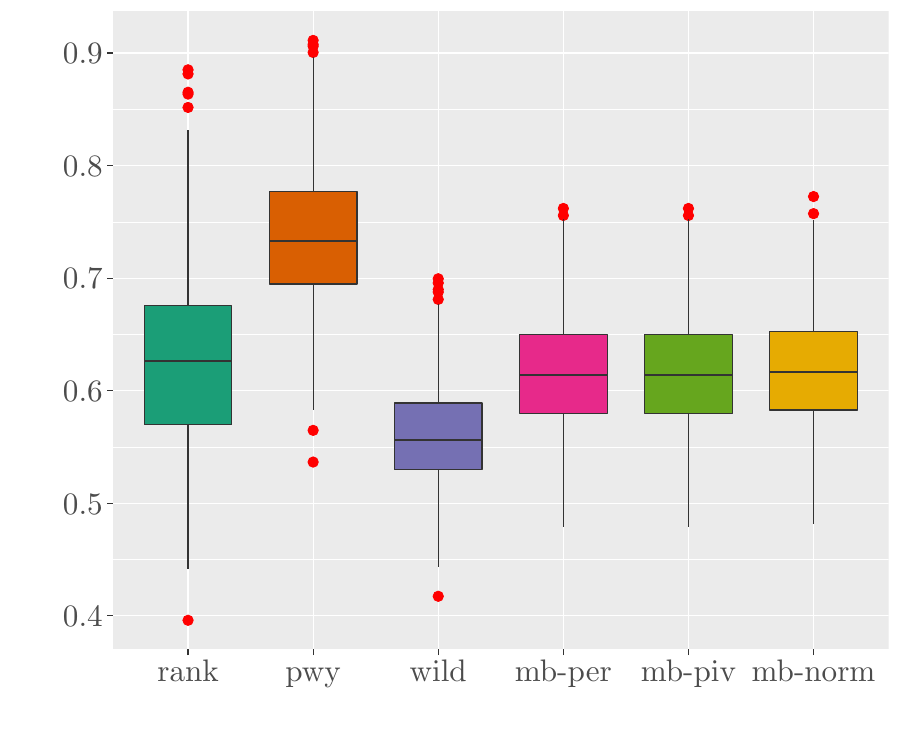}}
  \subfigure[CI width under model \eqref{model.quad}.]{\includegraphics[width=0.32\textwidth]{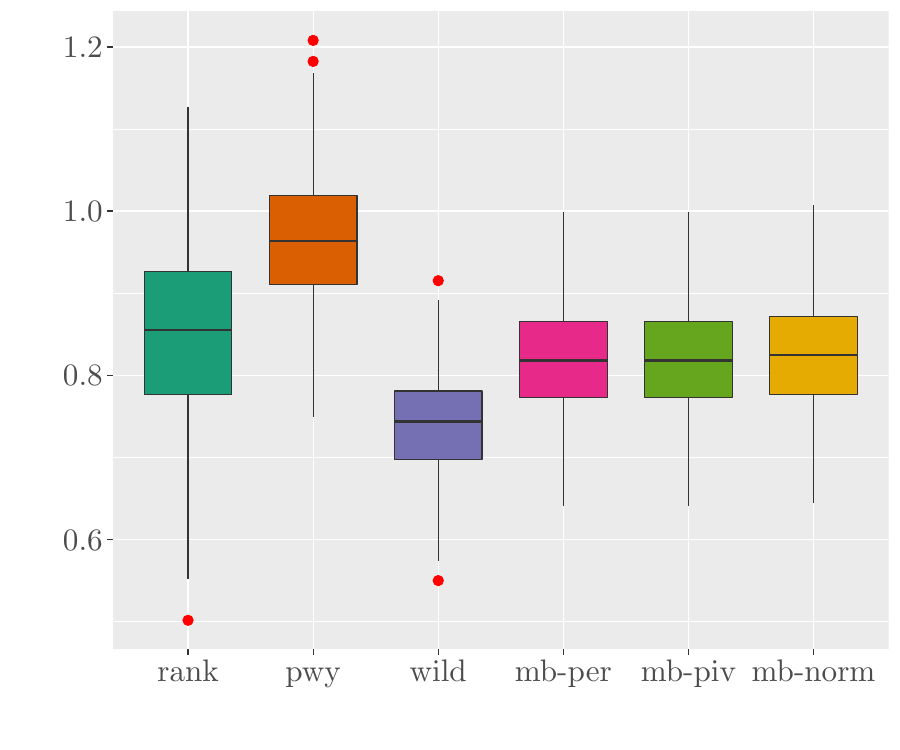}}
  \subfigure[Runtime under model \eqref{model.homo}.]{\includegraphics[width=0.32\textwidth]{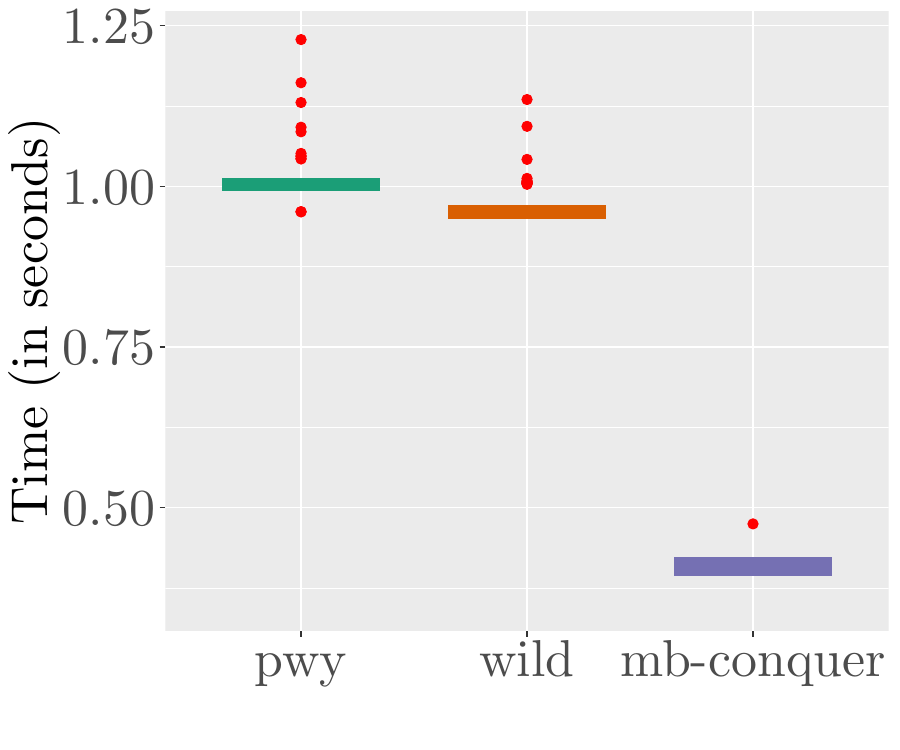}}
  \subfigure[Runtime under model \eqref{model.linear}.]{\includegraphics[width=0.32\textwidth]{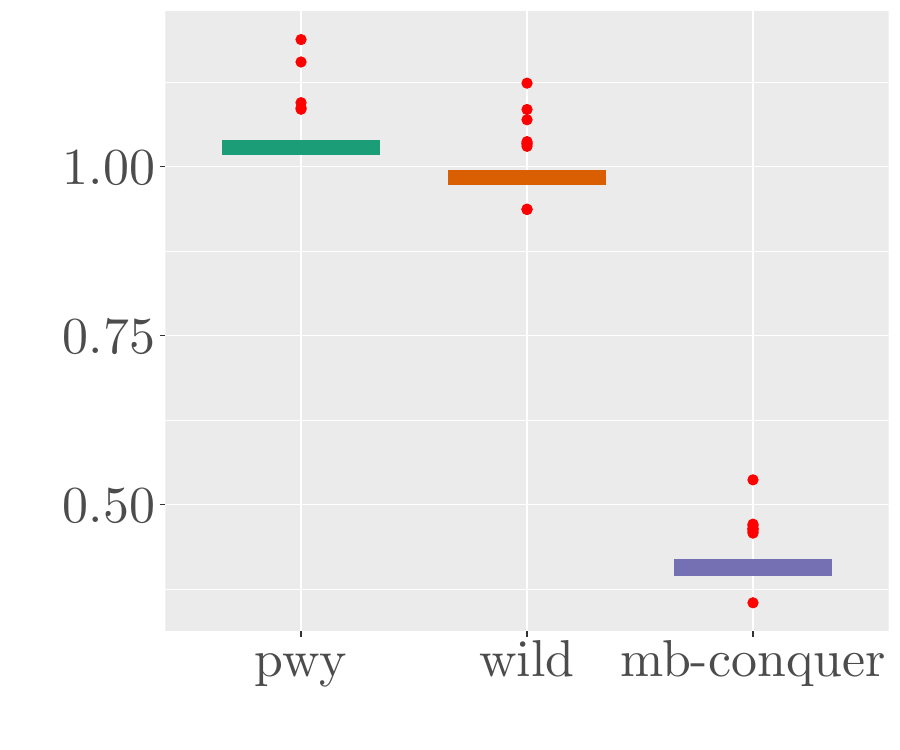}}
  \subfigure[Runtime under model \eqref{model.quad}.]{\includegraphics[width=0.32\textwidth]{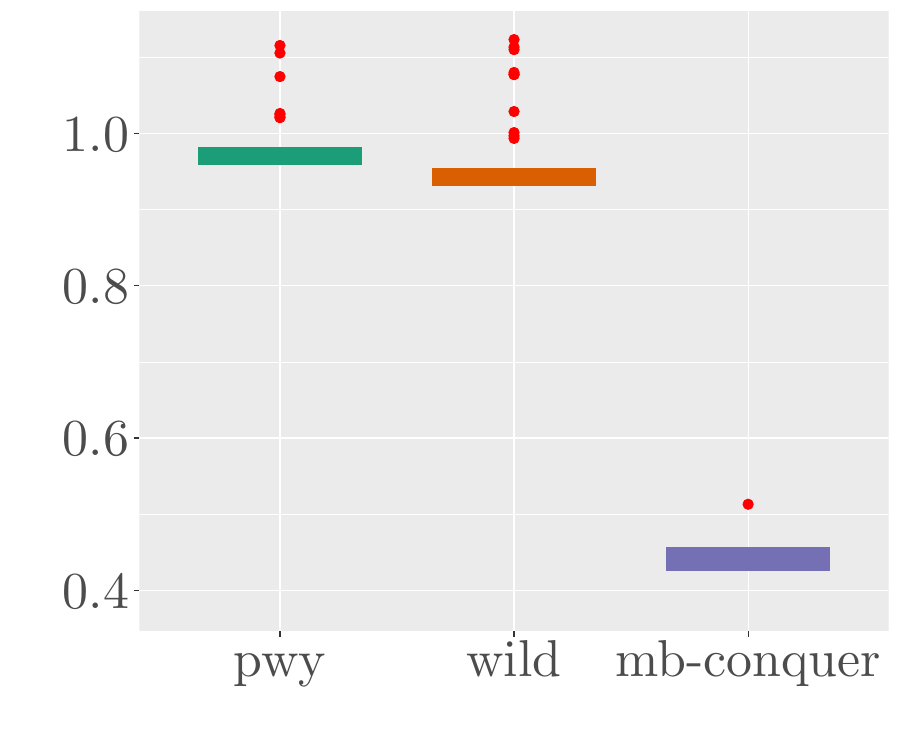}}
\caption{Empirical coverage, confidence interval width, and elapsed time of $six$ methods with $\mathcal{N}(0, 4)$ errors under $\tau = 0.9$. Other details are as in Figure~\ref{inf.t.9}.}
  \label{inf.normal.9}
\end{figure}

\begin{figure}[!htp]
  \centering
  \subfigure[Coverage under model \eqref{model.homo}.]{\includegraphics[width=0.32\textwidth]{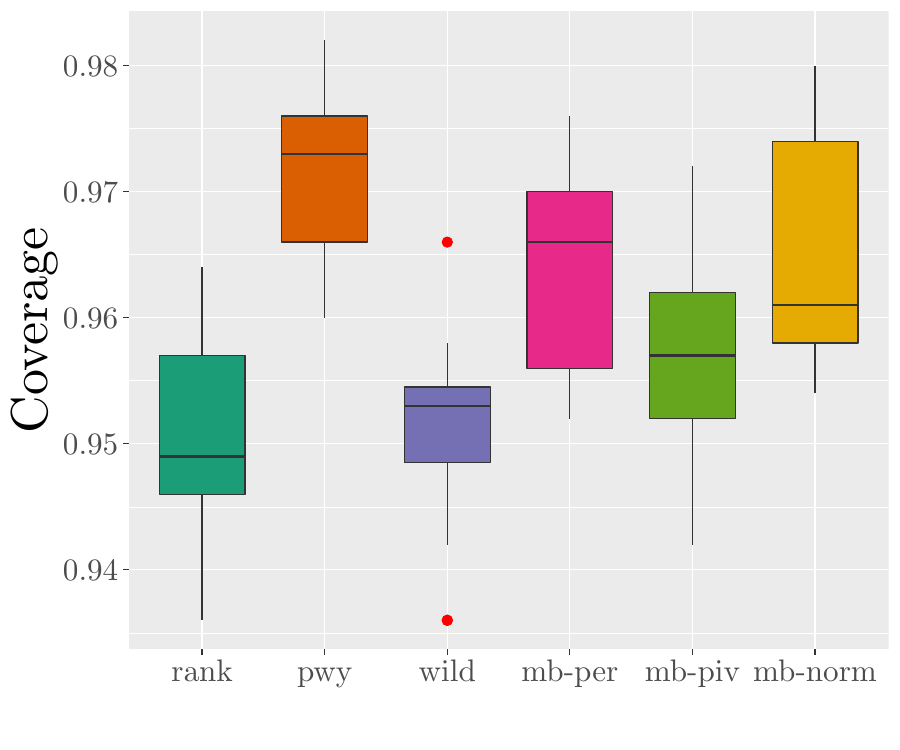}} 
  \subfigure[Coverage under model \eqref{model.linear}.]{\includegraphics[width=0.32\textwidth]{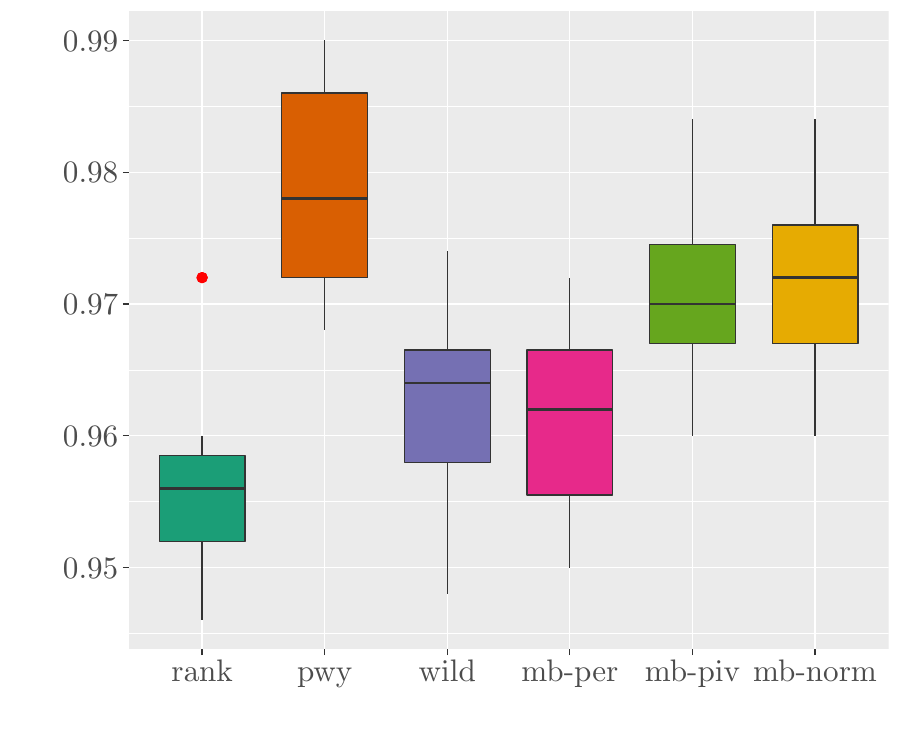}} 
  \subfigure[Coverage under model \eqref{model.quad}.]{\includegraphics[width=0.32\textwidth]{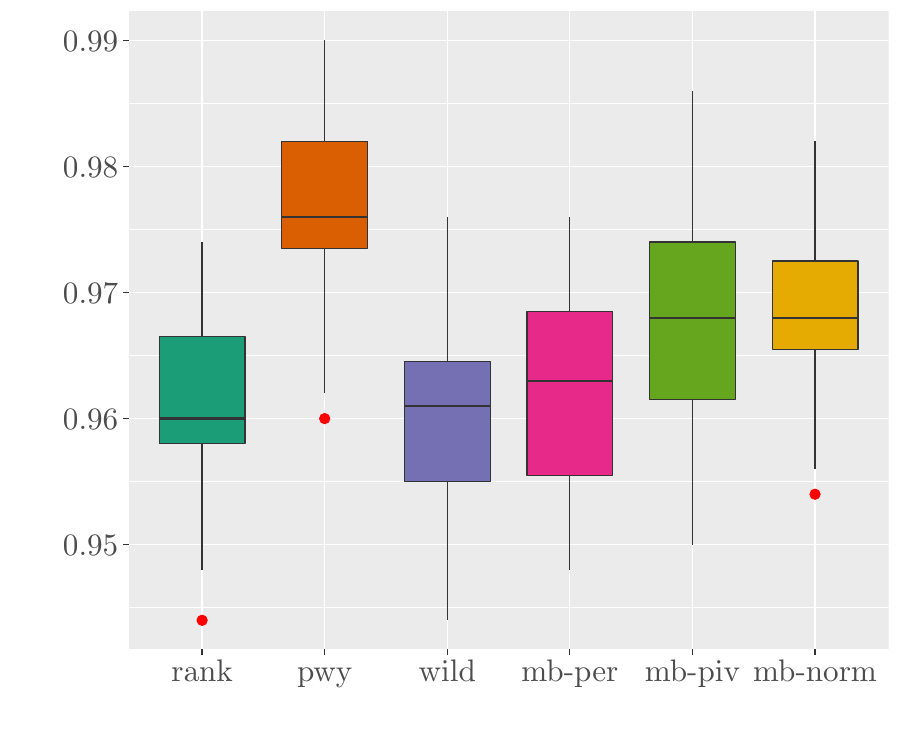}}
  \subfigure[CI width under model \eqref{model.homo}.]{\includegraphics[width=0.32\textwidth]{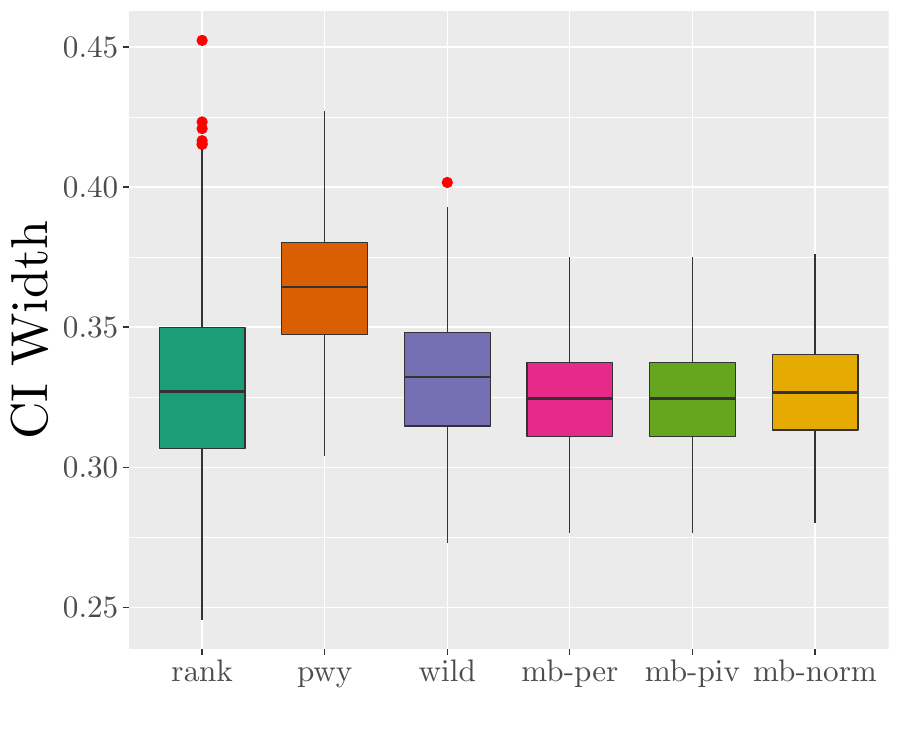}}
  \subfigure[CI width under model \eqref{model.linear}.]{\includegraphics[width=0.32\textwidth]{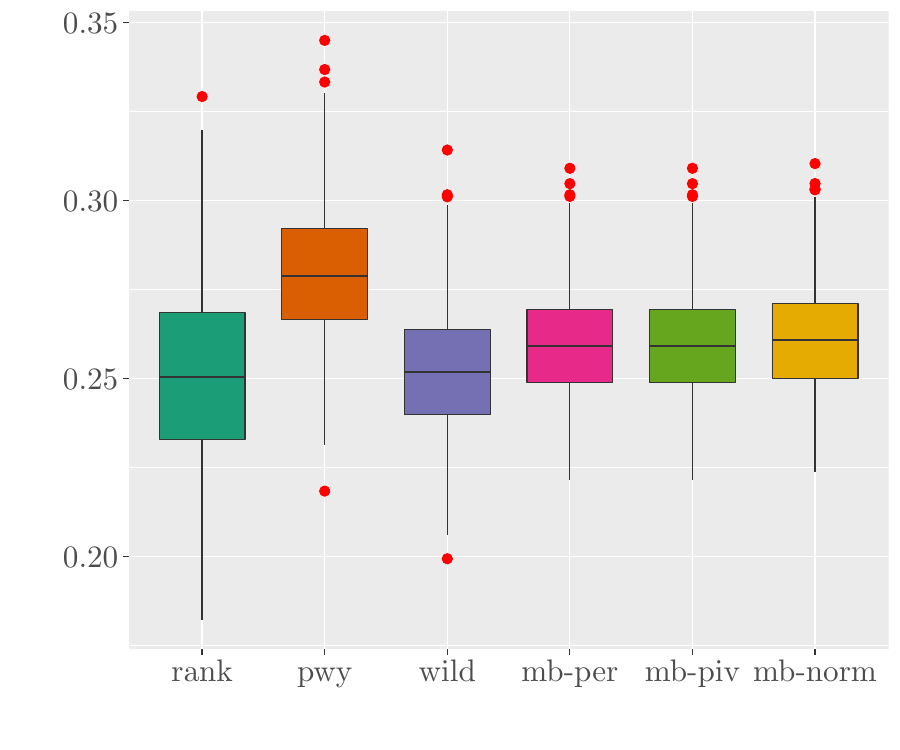}}
  \subfigure[CI width under model \eqref{model.quad}.]{\includegraphics[width=0.32\textwidth]{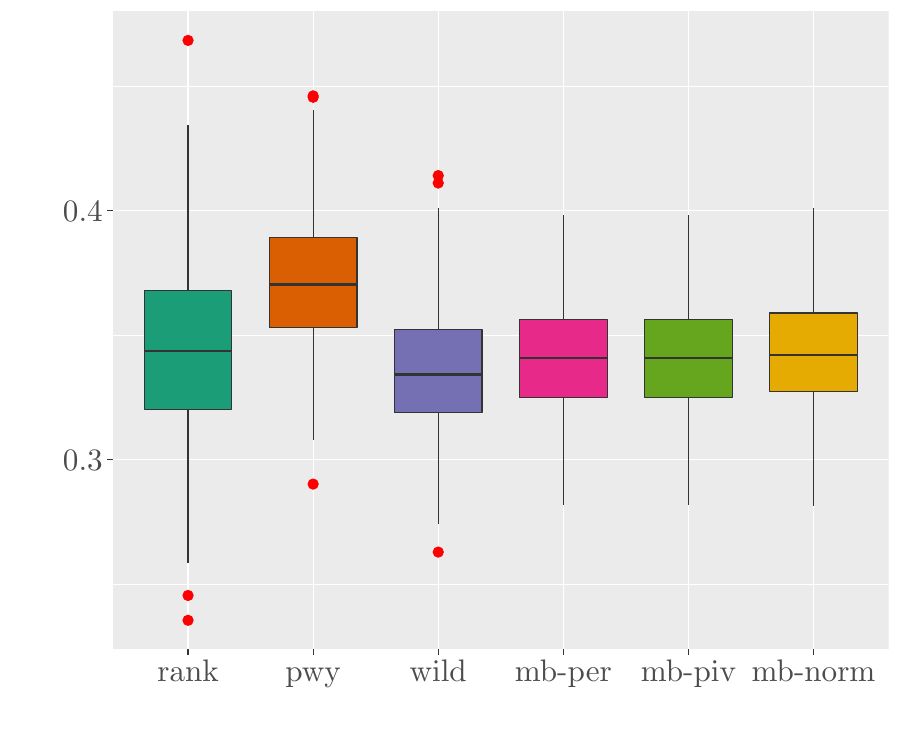}}
  \subfigure[Runtime under model \eqref{model.homo}.]{\includegraphics[width=0.32\textwidth]{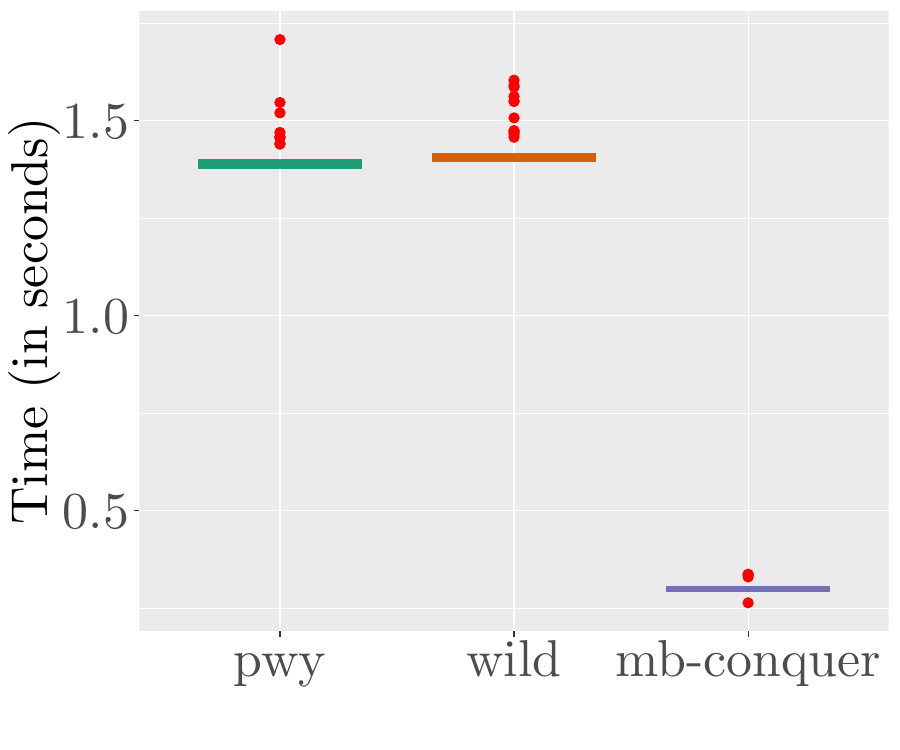}}
  \subfigure[Runtime under model \eqref{model.linear}.]{\includegraphics[width=0.32\textwidth]{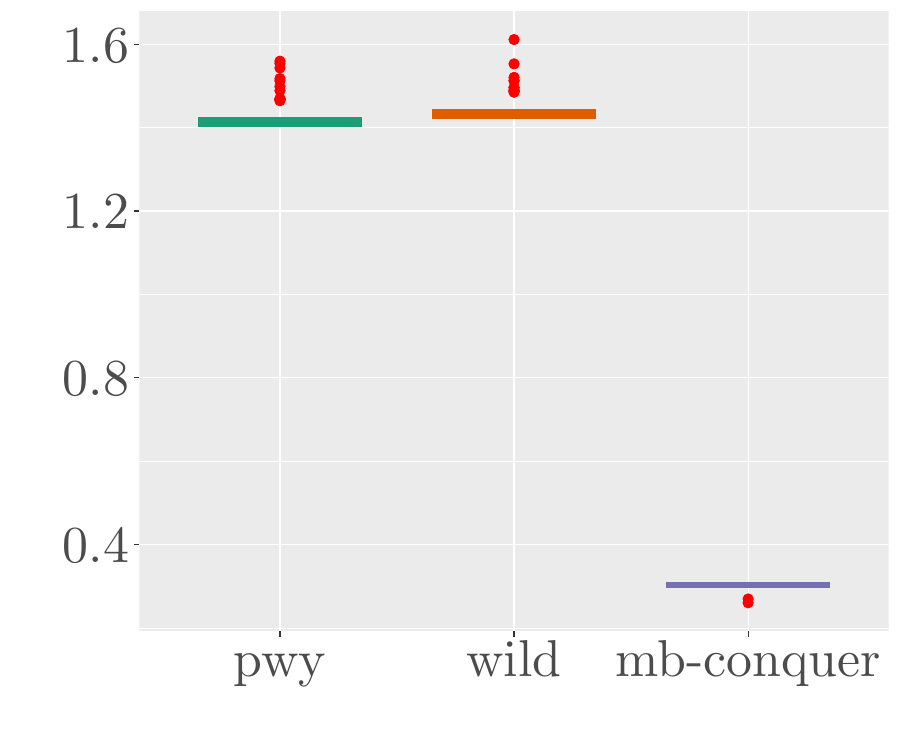}}
  \subfigure[Runtime under model \eqref{model.quad}.]{\includegraphics[width=0.32\textwidth]{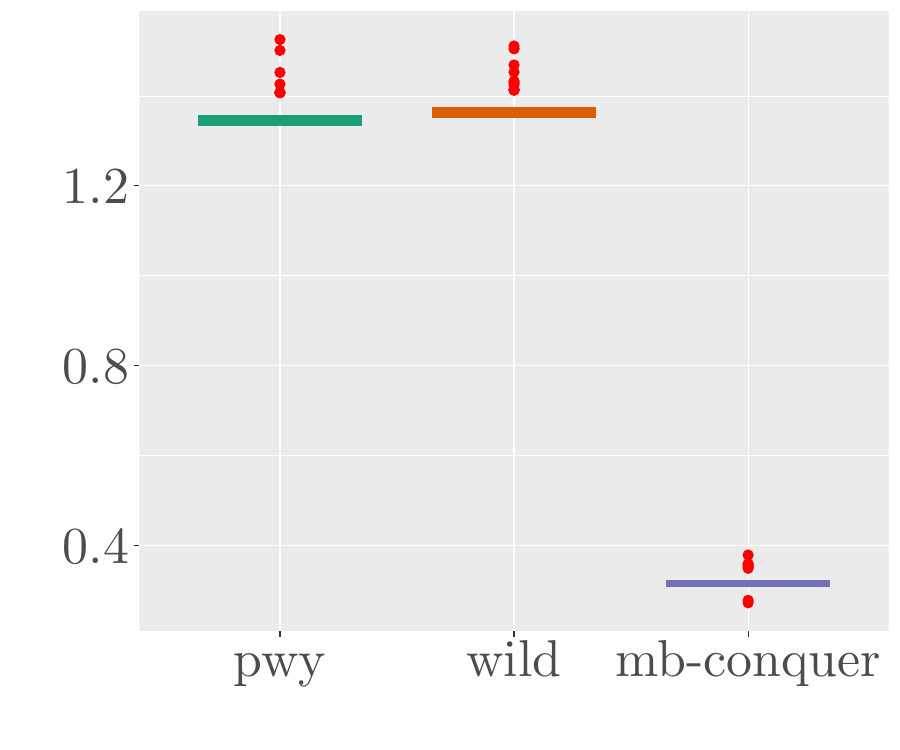}}
\caption{Empirical coverage, confidence interval width, and elapsed time of $six$ methods with $t_2$ errors under $\tau = 0.5$. Other details are as in Figure~\ref{inf.t.9}.}
  \label{inf.t.5}
\end{figure}

\begin{figure}[!htp]
  \centering
  \subfigure[Coverage under model \eqref{model.homo}.]{\includegraphics[width=0.32\textwidth]{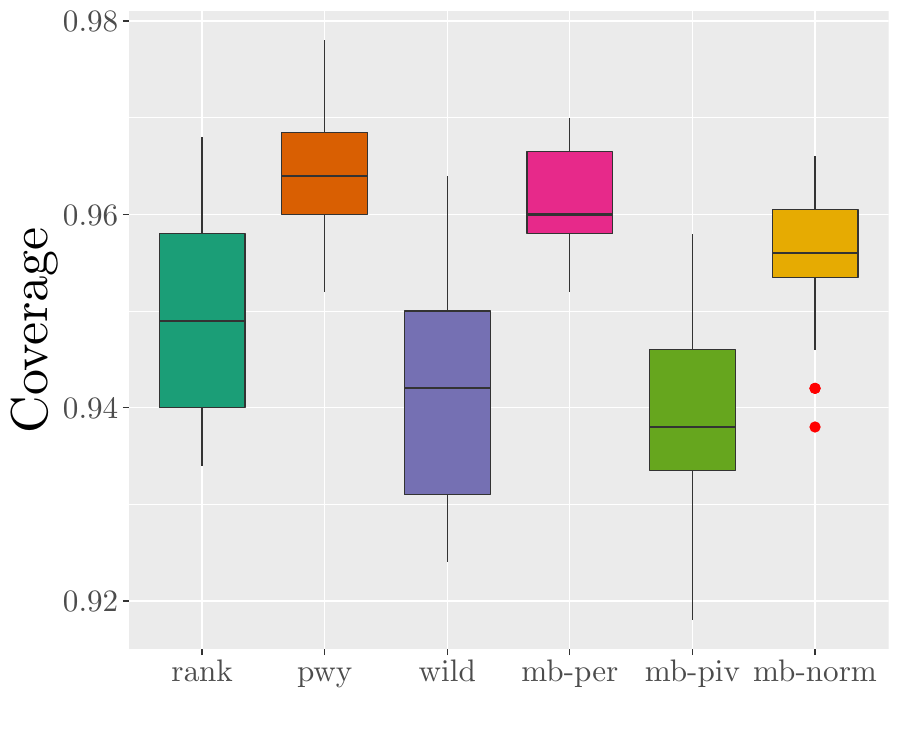}} 
  \subfigure[Coverage under model \eqref{model.linear}.]{\includegraphics[width=0.32\textwidth]{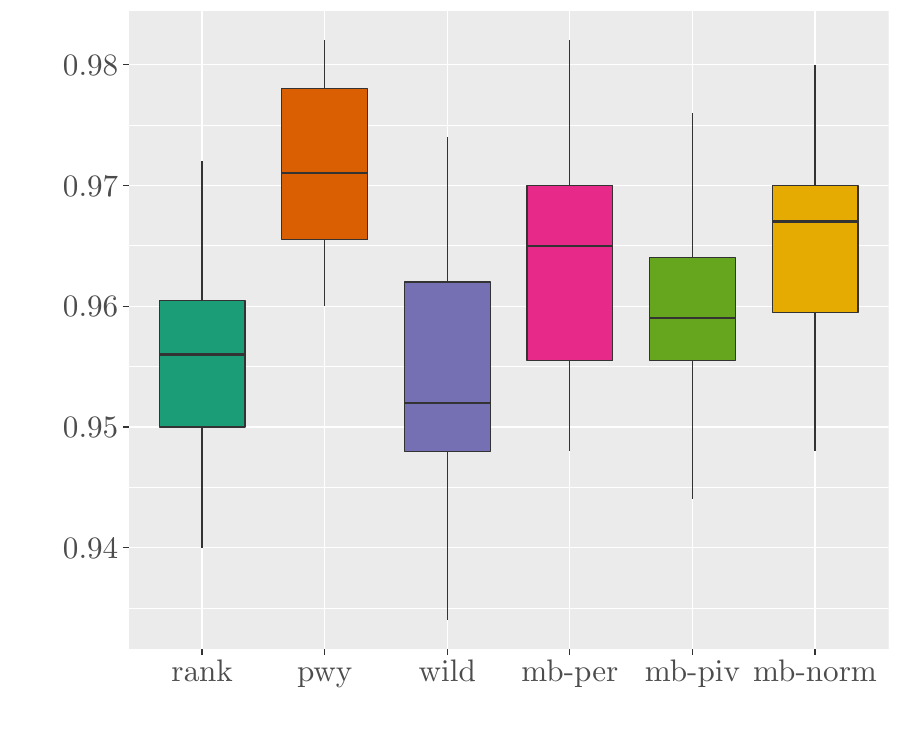}} 
  \subfigure[Coverage under model \eqref{model.quad}.]{\includegraphics[width=0.32\textwidth]{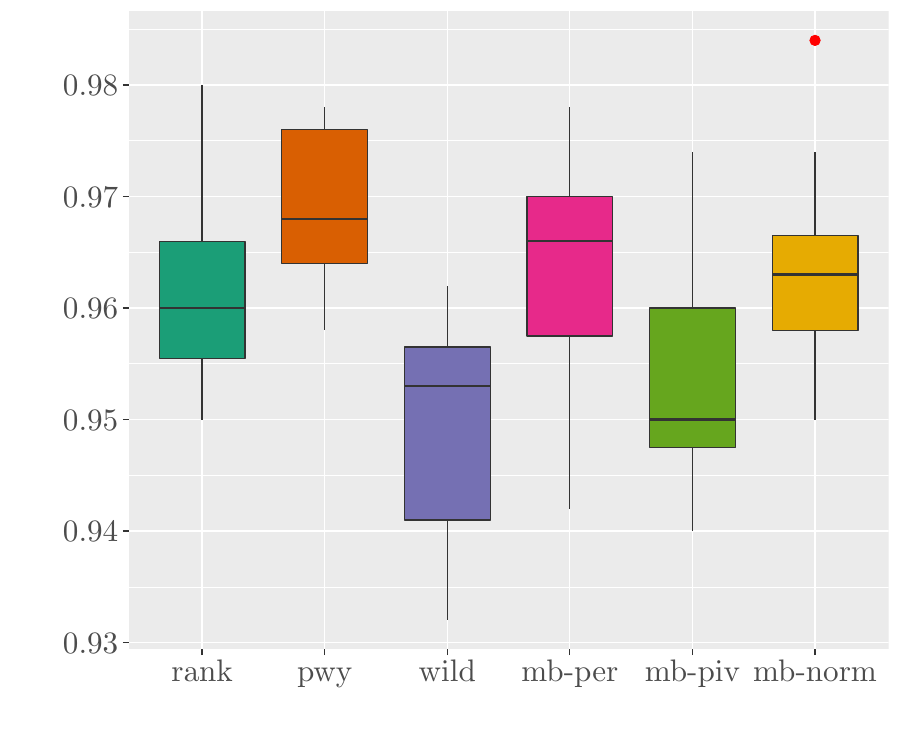}}
  \subfigure[CI width under model \eqref{model.homo}.]{\includegraphics[width=0.32\textwidth]{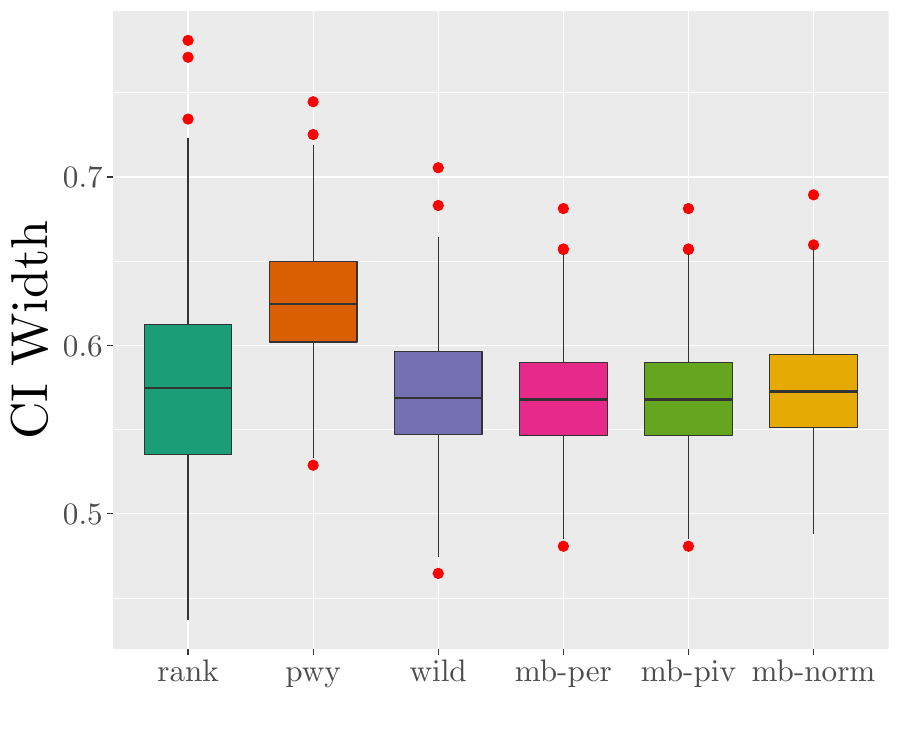}}
  \subfigure[CI width under model \eqref{model.linear}.]{\includegraphics[width=0.32\textwidth]{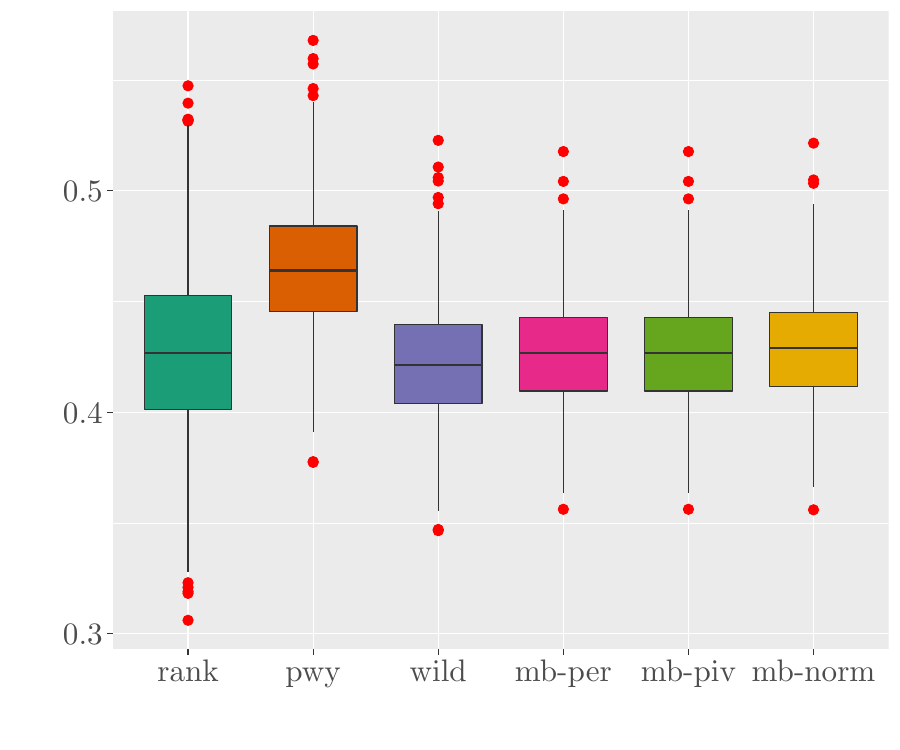}}
  \subfigure[CI width under model \eqref{model.quad}.]{\includegraphics[width=0.32\textwidth]{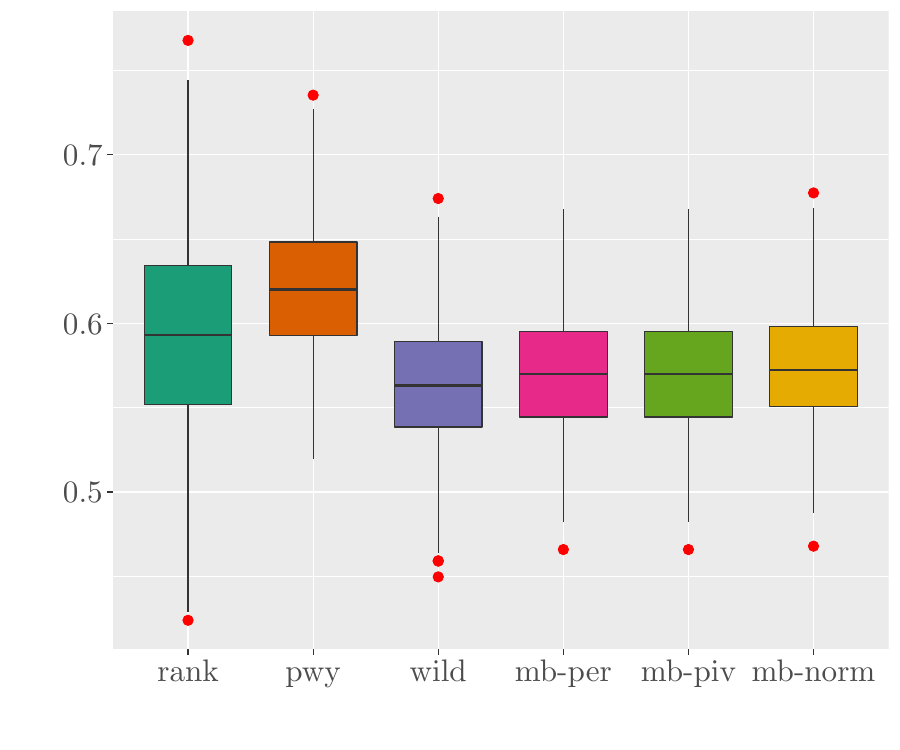}}
  \subfigure[Runtime under model \eqref{model.homo}.]{\includegraphics[width=0.32\textwidth]{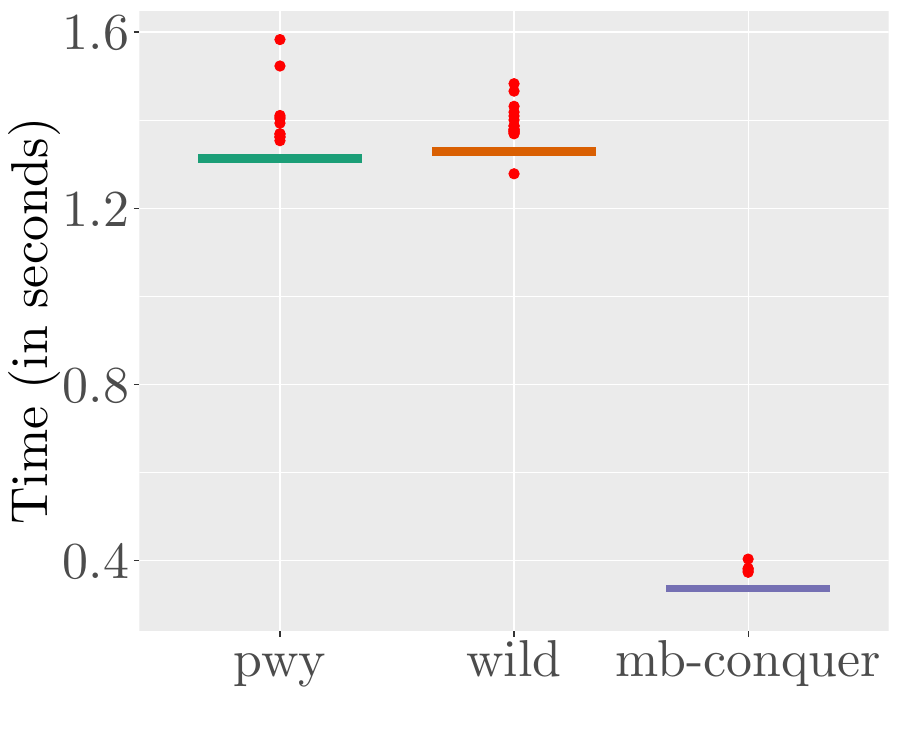}}
  \subfigure[Runtime under model \eqref{model.linear}.]{\includegraphics[width=0.32\textwidth]{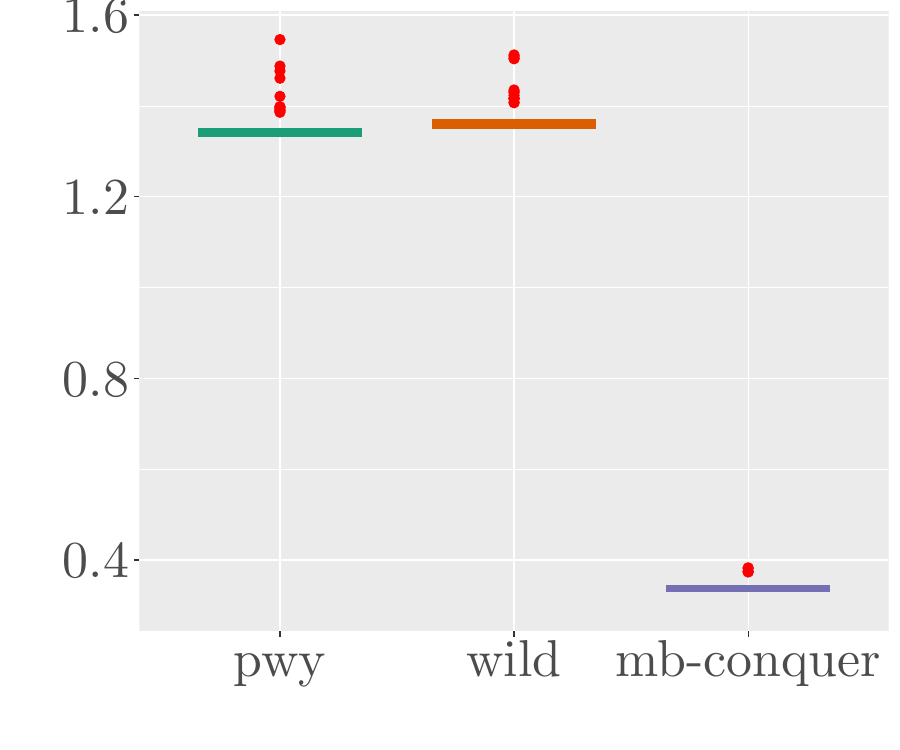}}
  \subfigure[Runtime under model \eqref{model.quad}.]{\includegraphics[width=0.32\textwidth]{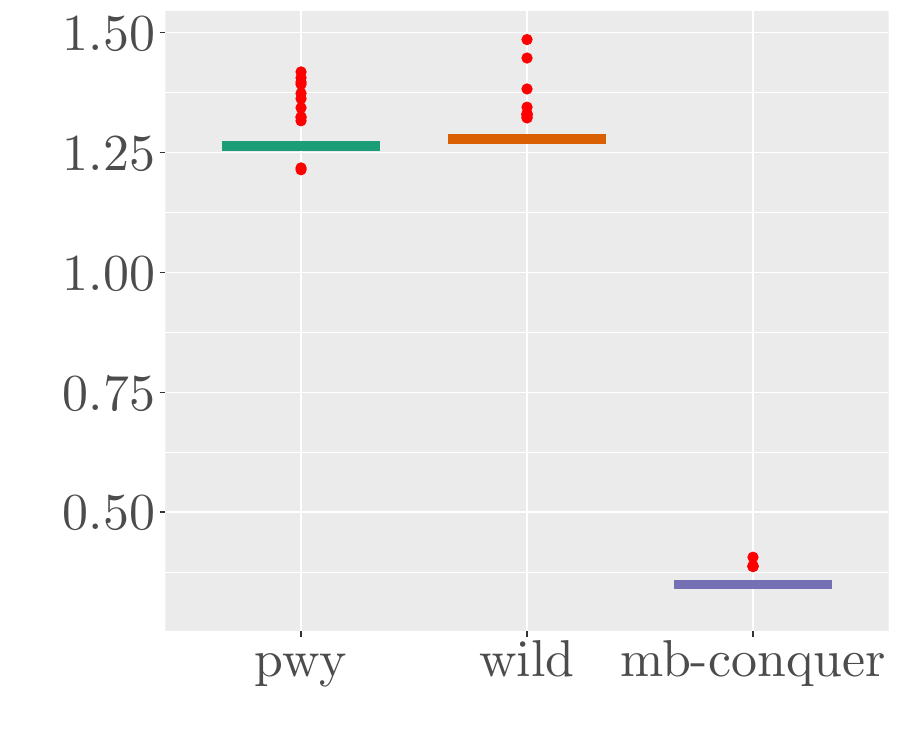}}
\caption{Empirical coverage, confidence interval width, and elapsed time of $six$ methods with $\mathcal{N}(0, 4)$ errors under $\tau = 0.5$. Other details are as in Figure~\ref{inf.t.9}.}
  \label{inf.normal.5}
\end{figure}

 In this section, we present additional results for the numerical studies as described in Section~\ref{sec:numerical} in Figures~\ref{est.1357.normal}--\ref{inf.normal.5}.

 %%%%%%%%%%%%%%%%%%%%%%%%%%%%%%%%%%%%%%%%%
%%%%%%%%%%%%%%%%%%%%%%%%%%%%%%%%%%%%%%%%%
% Bibliography
%%%%%%%%%%%%%%%%%%%%%%%%%1%%%%%%%%%%%%%%%%
%%%%%%%%%%%%%%%%%%%%%%%%%%%%%%%%%%%%%%%%%

\end{document}